%% file: upravl1.tex
\documentclass[russian,12pt]{report}
\usepackage[utf8]{inputenc}
\usepackage[russian]{babel}
\usepackage{graphicx}
\usepackage{amsmath,amssymb}
\begin{document}
\input{vved1}

\input{ch11}

\input{ch22}

\input{ch33}

\input{liter1}

\end{document}

%% file: vved1.tex
\begin{titlepage}
\vskip40pt 
\begin{center}
\Large 
КАЗАНСКИЙ ГОСУДАРСТВЕННЫЙ УНИВЕРСИТЕТ
\end{center}
\vskip80pt
\hfill {\large \bf На правах рукописи}
\vskip20pt
\hfill {\large \bf УДК  515.124.4}
\vskip50pt
\centerline{\large \bf  Сосов Евгений Николаевич}
\vskip30pt
\begin{center}
\Large \bf
ГЕОМЕТРИИ ВЫПУКЛЫХ  И КОНЕЧНЫХ МНОЖЕСТВ ГЕОДЕЗИЧЕСКОГО ПРОСТРАНСТВА
\end{center}
\vskip40pt
\begin{center}
\large \bf
01.01.04 --- геометрия и топология
\end{center}
\vskip20pt
\begin{center}
\large \bf
Д и с с е р т а ц и я 
\vskip10pt
на соискание ученой степени доктора 
\vskip10pt
физико-математических наук
\end{center}
\vskip130pt
\begin{center}
\bf
КАЗАНЬ --- 2010
\end{center}
\end{titlepage}
\newpage
\large
\setcounter{page}{2}
\tableofcontents
\pagebreak
\newpage
\newpage
\addtocontents{toc}{Введение \hfill {4} \par}

\section*{\bf Введение}

\textbf {Объектом исследования} настоящей работы являются
проблемы геометрии выпуклых (конечных) множеств геодезического пространства.
\\
\textbf{Актуальность.}  Метрическая геометрия возникла в 20--30-е годы двадцатого века в работах К. Менгера, П. С. Урысона, А. Вальда, С. Э. Кон-Фоссена, К. Куратовского,  Ф. Хаусдорфа, И. Шёнберга и других математиков. 
В этот начальный период метрическая геометрия еще не приобрела известность, а сам термин <<метрическая геометрия>> имел более узкий смысл. В 40--60-е годы были созданы основы метрической геометрии в фундаментальных работах Г. Буземана, А. Д. Александрова, В. А. Ефремовича, Л. М. Блюменталя, В. А. Залгаллера, Ю. Г. Решетняка, Ю. Д. Бураго и их учеников.
С 70х годов начался современный этап развития метрической геометрии, достижения которого отражены в монографиях Г. Буземана \cite {Bus1};  Г. Буземана, B. B. Phadke \cite {Bus2}; М. Л. Громова \cite {Grom}; А. В. Погорелова \cite {Pog1}; W. Ballmann \cite {Ballm}; A. Papadopoulos \cite {Papad}; M. Bridson, A. Haeffliger \cite {Brid}; С. В. Буяло, В. Шрёдер \cite {Buyal1}; М. М. Деза, М. Лоран \cite {Deza}; в первом учебнике на русском языке Ю. Д. Бураго, Д. Ю. Бураго, С. В. Иванова \cite {Bur}; в обзорах  Ю. Г. Решетняка \cite {Resh}; В. Н. Берестовского, И. Г. Николаева \cite {Beres} и некоторых других обзорах и монографиях. Кроме того, большое количество новых результатов пока не описано в обзорах, монографиях и учебниках, они содержатся лишь в научных статьях, число которых стабильно растет.
В настоящее время установилось много взаимосвязей метрической геометрии с комбинаторной геометрией, геометрией близости \cite {Efrem1}, римановой геометрией <<в целом>> \cite {Bur2}, теорией гиперболических групп, теорией фракталов, геометрической теорией меры \cite {Feder}, нелинейным функциональным анализом, субдифференциальным исчислением \cite {Ioffe}, выпуклым анализом \cite {Polov, Pog}, теорией некорректных задач, теорией вероятностей, теорией графов \cite {Deza}, теорией приближений \cite {Tikh, Send} и другими разделами математики \cite {Deza1}. Эти взаимосвязи  поддерживают актуальность метрической геометрии и постоянный приток в нее новых задач. Кроме того, развитие метрической геометрии связано с важностью и распространенностью метрических свойств объектов, исследуемых в различных разделах математики и прикладных науках, а также с тем, что по мере накопления геометрических фактов, полученных другими методами (например, методами математического анализа), проясняется метрическая природа многих из них. 

Г. Буземан метризовал группу всех движений метрического пространства и доказал, что в случае конечной компактности метрического пространства (сейчас чаще употребляются термины: собственное метрическое пространство \cite [с. 2]{Brid} или ограниченно-компактное метрическое пространство \cite [с. 17]{Bur}) эта группа является конечно-компактным метрическим пространством \cite [c. 30, 32]{Bus}. В. Н. Берестовским было доказано, что метрическая топология Буземана в группе всех движений конечно-компактного метрического пространства эквивалентна ее компактно-открытой топологии и группа всех движений однородного $G$-пространства Буземана является группой Ли \cite {Beres5, Beres2}. Для группы подобий аналогичные исследования проведены в [s2], \cite {Beres6} и \cite {Gun}.
   
Известно, что одним из условий в определениях $G$-пространства Буземана \cite [с. 54]{Bus} и хордового пространства \cite [c. 23]{Bus2} является условие конечной компактности метрического пространства. Актуальной задачей является ослабление условия конечной компактности метрического пространства, поскольку оно исключает из исследования многие геометрические объекты гильбертовых многообразий (в частности, гильбертовых пространств) и бесконечномерных банаховых многообразий, а также ограничивает общность исследования наилучших аппроксимирующих множеств (например, чебышевских центров или наилучших $N$-сетей) ограниченных множеств, возникающих при решении геометрических задач или задач теории приближений. На этом пути были исследованы некоторые свойства выпуклых множеств в обобщенно хордовом пространстве (обобщенном $G$-пространстве Буземана) [s4, s6], обобщающие соответствующие свойства выпуклых множеств в хордовом пространстве ($G$-пространстве Буземана) \cite [с. 65, 74, 75, 79, 80-82]{Bus2}, \cite [с. 154, 157, 160]{Bus}. В более общей ситуации (то есть при отсутствии собственности метрического пространства) были использованы условие неположительности кривизны по Буземану и понятие дифференцируемого пространства по Буземану в точке, что позволило обобщить и начать исследование касательного пространства по Буземану в точке геодезического пространства [s18]. Другие подходы проработаны более глубоко и основаны на понятиях конуса над пространством направлений \cite {Alexan}, касательного конуса по Громову--Хаусдорфу \cite [с. 328]{Bur} и их модификаций с использованием дифферецируемости в метрическом пространстве, ультрасходимости метрических пространств и отображений \cite {Lyt}. Известно, что некоторые метрические свойства (например, аппроксимативные свойства) множеств в равномерно выпуклом банаховом пространстве связаны со свойствами метрической или обобщенной метрической проекции на эти множества. Эти свойства исследовались Б. Секефальви--Надь \cite[теорема 3.35]{Brudn}, Ю. А. Брудным, Е. А. Гориным \cite {Brudn},
С. Б. Стечкиным, Н. В. Ефимовым, Л. П. Власовым (см. обзор \cite {Vlas}), А. В. Мариновым \cite {Marin, Marin1, Marin2}, I. Singer \cite {Singer} и другими математиками. Обобщенная метрическая проекция имеет также важное значение для исследования $\varepsilon$-квазирешений и квазирешений операторных уравнений первого рода \cite {Marin, Lisk, Lisk1}. Оказалось, что многие из таких свойств допускают обобщение на геодезические пространства, удовлетворяющие дополнительным условиям, обеспечивающим в совокупности аналог свойства равномерной выпуклости банахова пространства [s7, s12, s13]. Аналог равномерной выпуклости в геодезическом пространстве дает возможность получить достаточные условия существования и единственности чебышевского центра ограниченного множества и исcледовать геометрические свойства наилучших $N$-сетей ограниченных множеств. В банаховых и гильбертовых пространствах свойства чебышевских центров и наилучших $N$-сетей ограниченных множеств исследовали А. Л. Гаркави \cite {Gark,Gark1}, П. К. Белобров  \cite {Belob, Belob1}, D. Amir, J. Mach, K. Saatkamp \cite {Amir,Amir1}, L. Vesely \cite {Vesel}, A. Wisnicki и J. Wosko \cite {Wisnicki}, В. С. Балаганский \cite {Balag}. В метрическом пространстве свойства чебышевского центра ограниченного множества исследовались при более сильных ограничениях на пространство \cite {Nar, Buyal2}, \cite [c. 26]{Ballm}. 
Аналогично ситуации с чебышевским центром, некоторые результаты С. И. Дудова и И. В. Златорунской \cite {Dudov, Dudov1} о наилучшем приближении в метрике Хаусдорфа выпуклого компакта банахова пространства шаром допускают обобщение на случай специального геодезического пространства неположительной кривизны по Буземану [s17].

Таким образом, важность установления новых связей метрической геометрии с теорией приближений, выпуклым анализом и функциональным анализом делает тему диссертации актуальной. Кроме того, есть много внутренних нерешенных проблем метрической геометрии. Например, проблема Буземана о том, является ли $G$-пространство Буземана топологическим многообразием \cite [c. 69]{Bus}. Г. Буземан \cite {Bus3}, В. Krakus \cite {Krak} и P. Thurston \cite {Thurst} доказали, что $G$-пространство Буземана является топологическим многообразием в размерностях $2$, $3$ и $4$ соответственно. В. Н. Берестовский установил некоторые достаточные условия конечномерности $G$-пространства Буземана \cite {Beres4}. В общем случае проблема остается открытой. 

В данной работе при исследовании геометрических свойств прост\-ранств выпуклых (конечных) множеств метрического пространства используются в основном прямые, синтетические методы Буземана \cite {Bus, Bus1, Bus2, Konf, Papad, Beres1, Beres2, Beres3, And1, And2, And3, And}, стандартные методы теории метрических пространств \cite {Kurat, Nadl} и методы теории приближений \cite {Tikh, Brudn, Vlas, Gark, Dudov, Marin, Marin1, Marin2}. Треугольники сравнения и верхние углы по А. Д. Александрову \cite {Beres} используются для исследования одного метрического аналога слабой сходимости последовательности в вещественном гильбертовом пространстве.

\textbf{Целью} настоящей работы является исследование геометрии выпуклых (конечных) множеств геодезического пространства.

\textbf{Научная новизна.} Все основные результаты, представленные в настоящей работе и выносимые на защиту, являются новыми. Перечислим эти результаты.
\\
1. Найдены необходимые и достаточные условия, при которых пространства $(X_N,\alpha_p)$, $(\Sigma_2 (X),\alpha)$ являются пространствами с внутренней метрикой, а также метрически выпуклыми (выпуклыми по Менгеру, собственными, геодезическими) пространствами. Получены достаточные условия, при которых пространство $(X^*_N,\alpha_p)$ является геодезическим пространством (удовлетворяет локальному условию неположительности кривизны по Буземану). Найдены необходимые и достаточные условия, при которых пространство  $(\Sigma_N (X),\alpha_{p,R})$ является пространством с внутренней метрикой, а также собственным (собственным метрически выпуклым, собственным выпуклым по Менгеру, собственным геодезическим) пространством. 
\\
2. Установлено, что одулярные структуры прямого $G$-пространства Буземана и геометрии Гильберта являются топологическими одулярными структурами. Исследованы геометрические свойства выпуклых $U$-множеств обобщенного хордового пространства.
\\
3. Получены оценки изменения относительного чебышевского радиуса $R_W (M)$ при изменении непустых ограниченных множеств $M, \, W$ метрического пространства. Найдены замыкание и внутренность множества всех $N$-сетей, каждая из которых обладает принадлежащим ей единственным относительным чебышевским центром, в множестве всех $N$-сетей специального геодезического пространства относительно метрики Хаусдорфа. Получены достаточные условия существования и единственности чебышевского центра, а также принадлежности чебышевского центра замыканию выпуклой оболочки непустого ограниченного множества специального геодезического пространства.
\\
4. Теоремы Б. Секефальви - Надь, С. Б. Стечкина и Н. В. Ефимова об
аппроксимативных свойствах множеств, а также теоремы Л. П. Власова и А. В. Маринова о непрерывности и связности метрической $\delta$-проекции в равномерно выпуклом банаховом пространстве обобщены на случай специального геодезического пространства. В специальном метрическом пространстве получены обобщения теорем \mbox{П. К. Белоброва} и  А. Л. Гаркави о наилучших $N$-сетях непустых ограниченных замкнутых выпуклых множеств в гильбертовом и в специальном банаховом пространствах. Для каждого непустого ограниченного множества бесконечномерного пространства Лобачевского доказано существование наилучшей $N$-сети и наилучшего $N$-сечения, а также установлена сильная устойчивость чебышевского центра. 
\\
5. Получена оценка сверху для расстояния Хаусдорфа от непустого ограниченного множества до множества всех замкнутых шаров специального геодезического пространства $X$ неположительной кривизны по Буземану. Доказано, что множество всех центров $\chi (M)$ замкнутых шаров, наилучшим образом приближающих в метрике
Хаусдорфа выпуклый компакт $M \subset X$, непустое и принадлежит $M$.
\\
6. Установлено, что метрика на касательном пространстве в произвольной точке пространства неположительной кривизны по Буземану (дифференцируемого по \mbox{Буземану} метрического пространства) внутренняя. Доказано, что касательное пространство в произвольной точке локально полного дифференцируемого по Буземану метрического пространства является полным пространством, а также, что касательное пространство в произвольной точке локально компактного пространства неположительной  кривизны по Буземану
является собственным геодезическим пространством.
\\
7. Доказано, что пространство всех слабо ограниченных гомеоморфизмов с метрикой Куратовского, каждый из которых равномерно непрерывен на произвольном замкнутом шаре с центром в фиксированной точке метрического пространства вместе со своим обратным гомеоморфизмом, является паратопологической группой (топологической группой при связности произвольного замкнутого шара с центром в данной фиксированной точке), непрерывно действующей на метрическом пространстве $X$. Теорема Банаха об обратном операторе и принцип равностепенной непрерывности для  $F$-пространств обобщены на случай специального геодезического отображения специальных геодезических пространств. 
\\
8. Доказано, что пространство $(H_B (X,Y,\alpha),\delta_p)$ всех отображений из метрического пространства $X$ в метрическое пространство $Y$, удовлетворяющих равномерному условию Гельдера с фиксированными показателем и коэффициентом, является полным (собственным) метрическим пространством, если $Y$ --- полное метрическое пространство ($X$, $Y$ --- собственные метрические пространства). Установлено, что если $X$ --- собственное метрическое пространство, то топология пространства $(H_B (X,Y,\alpha),\delta_p)$ совпадает как с топологией поточечной сходимости, так и с компактно-открытой топологией. В специальном метрическом пространстве введены два аналога слабой сходимости последовательности в вещественном гильбертовом пространстве и исследованы их геометрические свойства.
\\
9. Доказано, что: 

-- если $X$, $Y$ --- полные (собственные) метрические пространства, то пространство $Sim (X,Y)\cup Const (X,Y)$, состоящее из всех подобий и всех постоянных отображений из $X$ в $Y$, с метрикой Буземана $\delta_p$ является полным (собственным); 

-- если $X$ --- собственное метрическое пространство, то топология пространства $(Sim (X,Y)\cup Const (X,Y),\delta_p)$ совпадает как с топологией поточечной сходимости, так и с компактно-открытой топологией;

-- $(Sim (X),\delta_p)$ --- топологическая группа, действующая непрерывно на пространстве $X$;

-- группы подобий $Sim (X)$ и изометрий $Iso (X)$ с метрикой Куратовского $\delta$ являются топологическими группами, непрерывно действующими на пространстве $X$. Найдено замыкание группы подобий полного метрического пространства в объемлющем пространстве отображений $\Phi (X,X)$ с метрикой Буземана $\delta_p$.

\textbf{Методы исследования.} Основными методами исследования, применяемыми в настоящей работе, являются:
\\
-- синтетические, прямые методы Буземана;
\\
-- стандартные методы из теории метрических пространств;
\\
-- методы теории приближений.

\textbf{Достоверность} полученных в диссертации результатов
обусловлена тем, что: 
\\
-- применяются проверенные, точные и строго обоснованные методы исследования; 
\\
-- многие результаты диссертации являются обобщением полученных ранее результатов и совпадают с этими результатами в частных случаях;
\\
-- все основные результаты диссертации доказаны и опубликованы.

\textbf{Апробация.} Основные результаты, изложенные в
диссертации, опубликованы в 21 публикации.

Результаты докладывались 
\\
-- ежегодно на научных семинарах кафедры геометрии Казанского государственного университета в 1993-2009 г.г.; 
\\
-- на итоговых научных конференциях Казанского педагогического университета в 1994-2000 г.г.;
\\
-- на итоговых научных конференциях Казанского государственного университета в 2001-2009 г.г.; 
\\
-- на научных семинарах НИИ математики и механики им. Н. Г. Чеботарева (Казань, 2001-2007 г.г.);
\\
-- на международном геометрическом семинаре имени Н. И. Лобачевского <<Современная геометрия и теория физических полей>> (Казань, 4-6 февраля 1997 г.);
\\
-- на международной научной конференции <<Актуальные проблемы математики и механики>> в НИИ математики и механики им. Н. Г. Чеботарева (Казань, 1-3 октября 2000 г.); 
\\
-- на международной научной конференции <<Topology, Analysis and Rela\-ted Topics>>, посвященной шестидесятилетию А. С. Мищенко (Московский гос. ун-т, 29-31 августа 2001 г.);  
\\
-- на международной научной конференции <<Геометрия <<в целом>>, топология и их приложения>>, посвященной девяностолетию со дня рождения А. В. Погорелова (Харьковский национальный ун-т, 22-27 июня 2009 г.);
\\
-- на Восьмой научной школе-конференции <<Лобачевские чтения 2009>> (Казань, 1-6 ноября 2009 г.);
\\
-- на научном семинаре кафедры дифференциальной геометрии и приложений Московского государственного университета (Москва, 14 декабря 2009 г.);
\\
-- на геометрическом семинаре им. А. Д. Александрова Санкт-Петер\-бургского отделения Математического института им. В. А. Стеклова РАН (Санкт-Петербург, 4 марта 2010 г.).

\textbf{Краткое описание содержания работы по главам}.
Диссертация состоит из введения, трех глав и списка использованной литературы.

В первой главе диссертации исследуются геометрические свойства выпуклых и конечных множеств в геодезическом пространстве.

В параграфе 1.1 рассматривается внутренняя метрика Хаусдорфа. Доказано, что метрика Хаусдорфа на множестве всех непустых замкнутых ограниченных подмножеств метрического пространства $(X,\rho)$ является внутренней метрикой тогда и только тогда, когда метрика $\rho$ --- внутренняя (теорема 1.1.1). Установлено, что пространство $(X,\rho)$ --- метрически выпукло тогда и только тогда, когда для любых двух ограниченных множеств существования пространства $X$ найдется середина этих множеств относительно метрики Хаусдорфа. Указано, как построить такую середину в метрически выпуклом пространстве (теорема 1.1.2). В метрическом пространстве с внутренней метрикой получена верхняя оценка для хаусдорфова расстояния между обобщенными шарами (лемма 1.1.2).

В параграфе 1.2 при заданном метрическом пространстве $(X,\rho)$ рассматриваются множество всех $N$-сетей $\Sigma_N (X)$ пространства $X$, его подмножество $\Sigma^*_N (X)$, элементами которого служат произвольные $N$-сети мощности $N$,  симметризованная степень $X_N$ порядка $N$ пространства $X$, отождествленная с множеством всех $N$-сетей с повторениями, и его подмножество $X^*_N$, равномощное множеству $\Sigma^*_N (X)$.
Множество $X_N$ наделим метрикой $\alpha_p:$
$$\alpha_p ([(x_1,\ldots ,x_N)],[(y_1,\ldots ,y_N)])=$$
$$\min \{\rho_{N,p} ((x_1,\ldots ,x_N),(y_{\sigma (1)},\ldots ,y_{\sigma (N)})) : \sigma \in S(N)\},$$
где $S(N)$ --- группа всех подстановок множества из $N \geq 1$ элементов,
$$\rho_{N,p} ((x_1,\ldots ,x_N),(y_1,\ldots ,y_N)) = (|x_1y_1|^p + \ldots + |x_Ny_N|^p)^{1/p}$$
при $p \in [1,\infty)$  и
$$\rho_{N,\infty} ((x_1,\ldots ,x_N),(y_1,\ldots ,y_N)) = \max \{|x_1y_1|,\ldots ,|x_Ny_N|\}$$
при $p =\infty$. Рассмотрим на множестве $X_N$ следующее отношение эквивалентности $R$:
$$[(x_1,\ldots ,x_N)]R[(y_1,\ldots ,y_N)], \, \mbox{если} \quad \{x_1,\ldots ,x_N\} = \{y_1,\ldots ,y_N\}.$$
Множество $\Sigma_N (X)$ отождествим с фактор-пространством $X_N / R$ и наделим его фактор-метрикой $\alpha_{p,R}$ или метрикой Хаусдорфа $\alpha$.

Доказано, что пространства $(X_N,\alpha_p)$, $(\Sigma_2 (X),\alpha)$ являются пространствами с внутренней метрикой или метрически выпуклыми (выпуклыми по Менгеру, собственными, геодезическими) пространствами тогда и только тогда, когда сответствующими свойствами обладает пространство $(X,\rho)$ (теорема 1.2.1, следствие 1.2.1). Получены достаточные условия, при которых пространство $(X^*_N,\alpha_p)$ является геодезическим пространством или удовлетворяет локальному условию неположительности кривизны по Буземану (теорема 1.2.2). Доказано, что пространство  $(\Sigma_N (X),\alpha_{p,R})$ является пространством с внутренней метрикой или собственным (собственным метрически выпуклым, собственным выпуклым по Менгеру, собственным геодезическим) пространством тогда и только тогда, когда соответствующими свойствами обладает пространство $(X,\rho)$ (теорема 1.2.3). Для собственного геодезического пространства $(X,\rho)$ получены достаточные условия, при которых произвольные две $N$-сети пространства $(\Sigma_N (X),\alpha)$ могут быть соединены сегментом (следствие 1.2.4). В метрически выпуклом пространстве, удовлетворяющем
глобальному условию неположительности кривизны по Буземану, найдены геометрическое свойство отображения, сопоставляющее произвольному сегменту его середину (лемма 1.2.1), и некоторые геометрические свойства пространства $(\Sigma_2 (X),\alpha)$, связанные в основном со свойством
метрической выпуклости этого пространства (лемма 1.2.2).

В параграфе 1.3 исследуются геометрические свойства выпуклых множеств в обобщенном хордовом пространстве. Введено понятие обобщенного хордового пространства (обобщенного $G$-пространства Буземана). Получены достаточные условия выпуклости замыкания (внутренности) выпуклого $U$-множества обобщенного хордового пространства, а также совпадения замыкания (внутренности) выпуклого 
\\
$U$-множества с замыканием внутренности (внутренностью замыкания внутренности) этого множества (теорема 1.3.1). Найдены условия, эквивалентные условию выпуклости всех замкнутых шаров в обобщенном хордовом пространстве, являющимся $U$-множеством, а также в обобщенном 
$G$-пространстве Буземана (теоремы 1.3.2, 1.3.3). Установлены достаточные условия, при которых обобщенное хордовое пространство является обобщенным 
$G$-пространством Буземана (теорема 1.3.4). Геометрически охарактеризовано множество всех точек на всех опорных (полукасательных) хордах в произвольной граничной точке выпуклого множества из открытого шара обобщенного хордового пространства (теорема 1.3.5). Получены геометрические характеристики для трансверсалей произвольной прямой обобщенного прямого хордового пространства, все открытые шары которого выпуклы (теорема 1.3.6).

В параграфе 1.4 исследуются одулярные структуры геометрии Гильберта и прямого $G$-пространства Буземана. 
Получены достаточные условия, при которых одулярная структура геодезического пространства, через каждые две
различные точки которого можно провести единственную прямую, является топологической одулярной структурой (лемма 1.4.1). Установлено, что одулярные структуры прямого $G$-пространства Буземана и геометрии Гильберта являются топологическими одулярными структурами (теоремы 1.4.$1'$, 1.4.3). Доказано, что метрика геометрии Гильберта в открытом шаре $B$ строго выпуклого банахова пространства топологически эквивалентна индуцированной метрике в $B$ банахова пространства, а также, что эти метрики липшицево эквивалентены на каждом замкнутом шаре, содержащимся в $B$ (теорема 1.4.2). Получен явный вид для основных операций одулярной структуры геометрии Гильберта (лемма 1.4.2). Вычислены некоторые пределы функций, связанных с основными операциями одулярной структуры геометрии Лобачевского (теорема 1.4.4).

Основные результаты по перечисленным темам опубликованы в статьях
автора [s3, s4, s6, s11, s14, s16, s20].

Во второй главе диссертации исследуются аппроксимативные свойства множеств в геодезическом пространстве. 

В параграфе 2.1 рассматриваются относительные чебышевский центр,  чебышевский радиус и множество всех диаметральных точек ограниченного множества метрического пространства. Получены оценки изменения относительного чебышевского радиуса $R_W (M)$ при изменении непустых ограниченных множеств $M, \, W$ метрического пространства (теорема 2.1.1). Доказано, что из всякой последовательности компактных множеств метрического пространства, сходящейся относительно метрики Хаусдорфа к некоторому компактному множеству $M$, можно выбрать подпоследовательность, для которой последовательность множеств всех относительных чебышевских центров (всех диаметральных точек)  ее элементов сходится относительно отклонения Хаусдорфа к множеству всех относительных чебышевских центров (всех диаметральных точек) множества $M$ (теорема 2.1.2).

В параграфе 2.2 найдены замыкание и внутренность множества всех $N$-сетей, каждая из которых обладает принадлежащим ей единственным относительным чебышевским центром, в множестве всех $N$-сетей специального геодезического пространства относительно метрики Хаусдорфа (теорема 2.2.1).
Доказано, что относительно метрики Хаусдорфа при $N > 2$ граница множества всех $N$-сетей, каждая из которых обладает не более, чем $N-2$ принадлежащими ей относительными чебышевскими центрами, совпадает с множеством всех диаметральных $N$-сетей в множестве всех $N$-сетей специального геодезического пространства (следствие 2.2.1). 

В параграфе 2.3 получены достаточные условия существования и единственности чебышевского центра непустого ограниченного множества геодезического пространства (теорема 2.3.1, следствие 2.3.1).

В параграфе 2.4 теоремы Б. Секефальви - Надь \cite [теорема 3.35]{Brudn},
С. Б. Стечкина и Н. В. Ефимова \cite [теоремы 1.1 и 1.2]{Vlas} об аппроксимативных свойствах множеств в равномерно выпуклом банаховом пространстве обобщены на случай специального геодезического пространства (теоремы 2.4.1, 2.4.2, 2.4.3).

В параграфе 2.5 теоремы Л. П. Власова \cite{Vlas1, Vlas} и А. В. Маринова \cite {Marin1} о непрерывности и связности метрической $\delta$-проекции в равномерно выпуклом банаховом пространстве обобщены на случай специального геодезического пространства (теоремы 2.5.1 -- 2.5.4). Одним из простых следствий такого обобщения является справедливость аналогичных результатов в пространствах Лобачевского (включая бесконечномерные).

В параграфе 2.6 доказано, что теоремы А. В. Маринова из \cite{Marin, Marin2} о непрерывности метрической $\delta $-проекции на выпуклое множество в линейном нормированном пространстве остаются верными в специальном метрическом пространстве (теоремы 2.6.1, 2.6.2).

В параграфе 2.7 в специальном метрическом пространстве получены обобщения некоторых результатов \mbox{П. К. Белоброва \cite{Belob, Belob1}} и  А. Л. Гаркави \cite{Gark1} о наилучших $N$-сетях непустых
ограниченных замкнутых выпуклых множеств в гильбертовом и в специальном банаховом пространствах (теоремы 2.7.1, следствие 2.7.2), а также о принадлежности чебышевского центра замыканию выпуклой оболочки данного множества (лемма 2.7.2).

В параграфе 2.8 доказано, что некоторые результаты А. Л. Гаркави \cite {Gark, Gark1} и
П. К. Белоброва \cite {Belob1, Belob} о наилучшей сети, наилучшем сечении и чебышевском
центре ограниченного множества в специальном банаховом пространстве верны и в бесконечномерном пространстве Лобачевского. А именно, для каждого непустого ограниченного множества бесконечномерного пространства Лобачевского доказано существование наилучшей $N$-сети (теорема 2.8.1) и наилучшего $N$-сечения (теорема 2.8.2), а также установлена сильная устойчивость чебышевского центра (теорема 2.8.3).

В параграфе 2.9 рассматривается наилучшее приближение выпуклого компакта геодезического пространства шаром. Получена оценка сверху для расстояния Хаусдорфа от непустого ограниченного множества до множества всех замкнутых шаров специального геодезического пространства $X$ неположительной кривизны по Буземану (теорема 2.9.1). Доказано, что множество всех центров $\chi (M)$ замкнутых шаров, наилучшим образом приближающих в метрике Хаусдорфа выпуклый компакт $M \subset X$, непустое и содержится в $M$ (теорема 2.9.2). Исследованы геометрические свойства множества $\chi (M)$. Таким образом, теоремы С. И. Дудова и И. В. Златорунской \cite {Dudov, Dudov1} обобщены на случай специального геодезического пространства неположительной кривизны по Буземану.

В параграфе 2.10 исследуются метрические свойства касательного пространства для метрического пространства 
более общего, чем дифференцируемое G-пространство Буземана. Установлено,
что метрика на касательном пространстве в произвольной точке пространства
неположительной кривизны по Буземану (дифференцируемого по \mbox{Буземану}
метрического пространства) внутренняя. Доказано, что касательное
пространство в произвольной точке локально полного дифференцируемого
по Буземану метрического пространства является полным пространством, а
также, что касательное пространство в произвольной точке локально
компактного пространства неположительной  кривизны по Буземану
является собственным геодезическим пространством (теоремы 2.10.1, 2.10.2).

Основные результаты по перечисленным темам этой главы опубликованы в статьях автора [s7, s8, s9, s12, s13, s14, s17, s18, s19, s21].

В третьей главе исследуются специальные отображения метрических пространств.

В параграфе 3.1 рассматривается метрическое пространство слабо ограниченных отображений метрических пространств с метрикой Куратовского $\delta$. Доказано, что пространство всех слабо ограниченных гомеоморфизмов $(HB(X),\delta)$, каждый из которых равномерно непрерывен на произвольном замкнутом шаре с центром в фиксированной точке метрического пространства $X$ вместе со своим обратным гомеоморфизмом, является паратопологической группой, непрерывно действующей на пространстве $X$. Установлено, что $(HB(X),\delta)$ является топологической группой при связности произвольного замкнутого шара с центром в фиксированной точке метрического пространства $X$ (теорема 3.1.2).

В параграфе 3.2 рассматриваются геодезические отображения специальных геодезических пространств. Исследованы некоторые геометрические свойства таких отображений. Теорема Банаха об обратном операторе и принцип равностепенной непрерывности для  $F$-пространств \cite [c. 99, 104]{Sadov} обобщены на случай специальных геодезических отображений специальных геодезических пространств (теоремы 3.2.1 -- 3.2.5). 

В параграфе 3.3 метрика Буземана $\delta_p$ распространяется на множество $\Phi (X, Y)$ всех непрерывных отображений 
$$f : (X,\rho) \rightarrow (Y, d),$$ 
удовлетворяющих следующему условию: для любых $x, \, y \in X$
$$d (f (x), f (y)) \leq B_f e^{|xy|},$$
где $B_f$ --- неотрицательная константа. Доказано, что пространство 
\\
$(H_B (X,Y,\alpha),\delta_p)$ всех отображений из метрического пространства $X$ в метрическое пространство $Y$, удовлетворяющих равномерному условию Гельдера с фиксированными показателем $\alpha \in (0,1]$ и коэффициентом $B \in \mathbb{R}_{+}$, является полным (собственным) метрическим пространством, если $Y$ --- полное метрическое пространство ($X$, $Y$ --- собственные метрические пространства) (теоремы 3.3.1, 3.3.3). Установлено, что если $X$ --- собственное метрическое пространство, то топология пространства
$(H_B (X,Y,\alpha),\delta_p)$ совпадает как с топологией поточечной сходимости, так и с компактно-открытой топологией (теорема 3.3.2).

В параграфе 3.4 рассматривается метрическое пространство всех подобий $(Sim (X, Y),\delta_p)$ из метрического пространства $X$ на метрическое пространство $Y$ с метрикой Буземана $\delta_p$. Доказано, что: 

-- если группа всех подобий действует транзитивно на полном метрическом пространстве, то и группа изометрий действует на нем транзитивно (лемма 3.4.1);

-- если $X$, $Y$ --- полные (собственные) метрические пространства, то пространство $(Sim (X,Y)\cup Const (X,Y),\delta_p)$ --- полное (собственное), где $Const (X,Y)$ --- множество всех постоянных отображений из $X$ в $Y$ (теоремы 3.4.1, 3.4.3); 

-- если $X$ --- собственное метрическое пространство, то топология пространства $(Sim (X,Y)\cup Const (X,Y),\delta_p)$ совпадает как с топологией поточечной сходимости, так и с компактно-открытой топологией (теорема 3.4.2);
 
-- $(Sim (X),\delta_p)$ --- топологическая группа, действующая непрерывно на пространстве $X$ (теорема 3.4.4);

-- группы подобий $Sim (X)$ и изометрий $Iso (X)$ с метрикой Куратовского $\delta$ являются топологическими группами, непрерывно действующими на пространстве $X$. 

Найдено замыкание группы подобий полного метрического пространства в пространстве $(\Phi (X,X),\delta_p)$ (теорема 3.4.5).

В параграфе 3.5 в специальном метрическом пространстве введены два аналога слабой сходимости последовательности в вещественном гильбертовом пространстве и исследованы их геометрические свойства (теоремы 3.5.1 -- 3.5.4).

Основные результаты по перечисленным темам этой главы опубликованы в статьях автора [s1, s2, s5, s15].

Всюду в диссертации используется обозначение $|xy|$ для расстояния $\rho (x,y)$ между точками метрического пространства $(X,\rho)$. Теорема 1.2.3 обозначает теорему 3 из второго параграфа главы 1.

\pagebreak

%% file: ch11.tex
\chapter{Выпуклые и конечные множества в геодезическом пространстве}
\vskip20pt

В параграфе 1.1 рассматривается внутренняя метрика Хаусдорфа. Доказано, что метрика Хаусдорфа на множестве всех непустых замкнутых ограниченных подмножеств
метрического пространства $(X,\rho)$ является внутренней метрикой тогда и
только тогда, когда метрика $\rho$ --- внутренняя. Установлено, что пространство $(X,\rho)$ --- метрически выпукло тогда и
только тогда, когда для любых двух ограниченных множеств существования пространства $X$ найдется середина этих множеств относительно метрики Хаусдорфа. Указано, как построить такую середину в метрически выпуклом пространстве. В метрическом пространстве с внутренней метрикой получена верхняя оценка для хаусдорфова расстояния между обобщенными шарами.

В параграфе 1.2 при заданном метрическом пространстве $(X,\rho)$ рассматриваются множество всех $N$-сетей $\Sigma_N (X)$ пространства $X$, его подмножество $\Sigma^*_N (X)$, элементами которого служат произвольные $N$-сети мощности $N$,  симметризованная степень порядка $N$ пространства $X$, отождествленная с 
множеством всех $N$-сетей с повторениями $X_N$, и его подмножество $X^*_N$, равномощное множеству $\Sigma^*_N (X)$.
Множество $X_N$ наделяется метрикой $\alpha_p$, где $p \in [1,\infty]$, а множество $\Sigma_N (X)$ наделяется фактор-метрикой $\alpha_{p,R}$ или метрикой Хаусдорфа $\alpha$.
Найдены необходимые и достаточные условия, при которых пространства $(X_N,\alpha_p)$, $(\Sigma_2 (X),\alpha)$ являются пространствами с внутренней метрикой, а также метрически выпуклыми (выпуклыми по Менгеру, собственными, геодезическими) пространствами. Получены достаточные условия, при которых пространство $(X^*_N,\alpha_p)$ является геодезическим пространством или удовлетворяет локальному условию неположительности кривизны по Буземану. Найдены необходимые и достаточные условия, при которых пространство  $(\Sigma_N (X),\alpha_{p,R})$ является пространством с внутренней метрикой, а также собственным (собственным метрически выпуклым, собственным выпуклым по Менгеру, собственным геодезическим) пространством. Для собственного геодезического пространства $(X,\rho)$ получены достаточные условия, при которых произвольные две $N$-сети пространства $(\Sigma_N (X),\alpha)$ могут быть соединены сегментом. В метрически выпуклом пространстве, удовлетворяющем
глобальному условию неположительности кривизны по Буземану, найдены геометрическое свойство отображения, сопоставляющее произвольному сегменту его середину, и некоторые геометрические свойства пространства $(\Sigma_2 (X),\alpha)$, связанные в основном со свойством
метрической выпуклости этого пространства.

В параграфе 1.3 исследуются геометрические свойства выпуклых множеств в обобщенном хордовом пространстве. Введено понятие обобщенного хордового пространства (обобщенного $G$-пространства Буземана). Получены достаточные условия выпуклости замыкания (внутренности) выпуклого $U$-множества обобщенного хордового пространства, а также совпадения замыкания (внутренности) выпуклого 
\\
$U$-множества с замыканием внутренности (внутренностью замыкания внутренности) этого множества. Найдены условия, эквивалентные условию выпуклости всех замкнутых шаров в обобщенном хордовом пространстве, являющимся $U$-множеством, а также в обобщенном 
\\
$G$-пространстве Буземана. Установлены достаточные условия, при которых обобщенное хордовое пространство является обобщенным 
\\
$G$-пространством Буземана. Геометрически охарактеризовано множество всех точек на всех опорных (полукасательных) хордах в произвольной граничной точке выпуклого множества из открытого шара обобщенного хордового пространства. Получены геометрические характеристики для трансверсалей произвольной прямой обобщенного прямого хордового пространства, все открытые шары которого выпуклы.

В параграфе 1.4 исследуются одулярные структуры геометрии Гильберта и прямого $G$-пространства Буземана.
Получены достаточные условия, при которых одулярная структура геодезического пространства, через каждые две
различные точки которого можно провести единственную прямую, является топологической одулярной структурой. Установлено, что одулярные структуры прямого $G$-пространства Буземана и геометрии Гильберта являются топологическими одулярными структурами. Доказано, что метрика геометрии Гильберта в открытом шаре $B$ строго выпуклого банахова пространства топологически эквивалентна индуцированной метрике в $B$ банахова пространства, а также, что эти метрики липшицево эквивалентены на каждом замкнутом шаре, содержащимся в $B$. Получен явный вид для основных операций одулярной структуры геометрии Гильберта. Вычислены некоторые пределы функций, связанных с основными операциями одулярной структуры геометрии Лобачевского.

Основные результаты по перечисленным темам опубликованы в статьях
автора [s3, s4, s6, s11, s14, s16, s20].

\section{Внутренняя метрика Хаусдорфа}

\vskip20pt

Напомним наиболее общее определение пространства с внутренней метрикой.
Метрическое пространство $(X,\rho)$ называется {\it пространством с внутренней метрикой}, если для любых $x, \, y \in X$, $\varepsilon > 0$ найдется конечная последовательность точек
$$z_0 =x,\, z_1, \ldots ,\, z_k=y \in X$$
такая, что для каждого $i \in \{0,\ldots, k-1\}$
$$|z_iz_{i+1}| < \varepsilon \, \, \mbox{и}\, \, |z_0z_1| +\ldots + |z_{k-1}z_k| < |xy| +\varepsilon,$$
где $|xy| =\rho (x,y)$ для $x, \, y \in X$ \cite{Bur1}.

Покажем, что это условие внутренности
метрики $\rho$ равносильно следующему условию $(A)$ \cite [c. 121]{Ioffe}.
\vskip7pt
$(A)\quad$ Для любых $x, \, y \in X$, $\varepsilon > 0$ множество
$$\omega (x,y,\varepsilon) = \{z \in X : 2 \max \{|xz|, |zy|\} < |xy| +
\varepsilon\}$$
непустое.
\vskip7pt
{\bf Лемма 1} [s11]. {\it Метрическое пространство $(X,\rho)$ является
пространством с внутренней метрикой тогда и только тогда, когда
выполняется условие $(A)$.}
\vskip7pt
{\bf Доказательство леммы 1}.
\vskip7pt
Необходимость. Пусть для любых
$$x, \, y \in X,\quad \varepsilon > 0$$
найдется конечная последовательность точек
$$z_0=x,\, z_1, \ldots ,\, z_k=y \in X$$
такая, что для каждого $i \in \{0,\ldots, k-1\}$
$$|z_iz_{i+1}| < \dfrac{1}{2}\varepsilon\, \, \mbox{и}\, \, |z_0z_1| +\ldots + |z_{k-1}z_k| < |xy| + \dfrac{1}{2}\varepsilon.$$
Если $x=y$, то доказательство очевидно. Пусть $x\neq y$. Кроме того,
можно считать, что точки
$$z_0=x,\, z_1, \ldots ,\, z_k = y$$
попарно различны. Тогда найдется такая точка $z_i\in X$, что
$$|z_0z_i| \leq \dfrac{1}{2}|xy|,\quad |z_0z_{i+1}| \geq \dfrac{1}{2}|xy|.$$
Из условия леммы 1 и неравенства треугольника получим неравенства
$$|z_0z_{i+1}| + |z_{i+1}z_k| \leq
|z_0z_1| +\ldots + |z_{k-1}z_k| < |xy| + \dfrac{1}{2}\varepsilon.$$
Тогда
$$|z_{i+1}z_k| < |xy| + \dfrac{1}{2}\varepsilon -
|z_0z_{i+1}| \leq \dfrac{1}{2}|xy| + \dfrac{1}{2}\varepsilon,$$
$$|z_0z_{i+1}| \leq |z_0z_i| +
|z_iz_{i+1}| < \dfrac{1}{2}|xy| + \dfrac{1}{2}\varepsilon.$$
Следовательно, $z_{i+1} \in \omega (x,y,\varepsilon)$.
\vskip7pt
Достаточность. Пусть
$$x, \, y \in X, \quad \varepsilon > 0.$$
Тогда найдется
натуральное число $n$ такое, что $|xy| < \varepsilon(2^n - 1)$. Положим
$$\delta = \dfrac{\varepsilon}{2^n - 1},\quad z_{00} = x,\quad z_{01} = y.$$
Используя условие $(A)$, на первом шаге найдем точки:
$$z_{10}=z_{00},\quad
z_{11} \in \omega (z_{00},z_{01},\delta),\quad z_{12}=z_{01},$$
где первый номер индекса --- номер шага, а второй --- номер точки на этом шаге.
Продолжая применять условие $(A)$, на шаге с номером $n$ получим $2^n +
1$ точек:
$$z_{n0}=z_{n-1,0}\,,\ldots,\, z_{n,2j}= z_{n-1,j},$$
$$z_{n,2j+1} \in
\omega (z_{n-1,j},z_{n-1,j+1},\delta),\ldots,\, z_{n,2m}=z_{n-1,m},$$
где
$m=2^{n-1}$, $0\leq j < m$. Тогда
$$|z_{n,2j}z_{n,2j \pm 1}| <
\dfrac{1}{2}|z_{n-1,j}z_{n-1,j \pm 1}| + \dfrac{1}{2}\delta < \ldots < |xy|2^{-n} + \delta(1 - 2^{-n}),$$
где $0\leq j < m$, в случае верхнего знака в индексе, и $1\leq j
\leq m$, в случае нижнего знака в индексе. Следовательно,
$$|z_{n0}z_{n1}|
+\ldots + |z_{n,2m-1}z_{n,2m}| < 2^n (|xy|2^{-n} + \delta(1- 2^{-n}))=
|xy| + \varepsilon.$$ 
Таким образом, лемма 1 доказана.
\vskip7pt
{\bf Замечание 1.} {\it В метрической геометрии часто используется следующее менее общее определение
пространства с внутренней метрикой {\rm (}см., например, {\rm \cite [с. 6]{Grom})}. Метрическое пространство $(X,\rho)$ называется пространством с внутренней метрикой,
если для любых  $x$, $y \in  X$ расстояние $|xy|$ равно инфимуму длин кривых, соединяющих точки $x$, $y$. Другое доказательство аналога леммы 1 для такого определения пространства с внутренней метрикой можно найти в} \cite [theorem 1.8, с. 6]{Grom}.
\vskip7pt
Примем следующие обозначения.
\vskip7pt
$\mathbb{R}_{+}$ --- множество всех неотрицательных вещественных чисел.
\vskip7pt
$B[X]$ ($B(X)$, $\Sigma (X)$) --- множество всех непустых ограниченных и
замкнутых (непустых ограниченных, непустых) подмножеств пространства $(X,\rho)$.
\vskip7pt
$B[x,r] \quad (B(x,r))$ --- замкнутый (открытый) шар с центром в точке
$x \in X$, радиуса $r > 0$.
$$|MW| = \inf \{|xy| : x \in M, \, y \in W\},$$ 
где $M$, $W \in \Sigma (X)$.
$$\beta (M,W) = \sup \{|xW| : x \in M\}$$
--- {\it отклонение} множества $M\in
B(X)$ от множества $W \subset X$ \cite [с. 223]{Kurat}.
$$\alpha : B(X)\times B(X)\rightarrow \mathbb{R}_{+},\quad \alpha(M,W) = \max\{\beta (M,W),\beta (W,M)\}$$
--- {\it псевдометрика Хаусдорфа} на множестве $B(X)$ (ограничение которой на
множество $B[X]$ является {\it метрикой Хаусдорфа}) \cite [с. 223]{Kurat}.
$$P :  X\times \Sigma (X) \times \mathbb{R}_{+} \rightarrow B(X),\quad P(x,M,\varepsilon) =  M \cap B[x,|xM| + \varepsilon]$$
--- {\it оператор метрического $\varepsilon$-проектирования} ({\it оператор
метрического проектирования} при
$$\varepsilon = 0, \quad M \cap B[x, |xM|] \neq \emptyset)$$
\cite{Marin1} или \cite{Vlas}. В дальнейшем вместо $$\bigcup \limits_{x \in W}{P(x,M)}$$
будем писать $P(W, M)$, где
$$W \in \Sigma (X),\quad P(x,M) = P(x,M,0).$$
$P(x,M)$ --- множество всех {\it оснований} точки $x$ на множестве $M$ \cite [с. 10] {Bus2}.
Множество $M\subset X$ называется {\it множеством существования}, если для каждого $x \in X$ $P(x,M) \neq \emptyset$ \cite{Vlas}.
\vskip7pt
Используем лемму 1 для доказательства следующей теоремы, содержащей
необходимое и достаточное условие внутренности метрики Хаусдорфа на $B[X]$.
\vskip7pt
{\bf Теорема 1} [s11]. {\it Метрическое пространство $(B[X],\alpha)$ является
пространством с внутренней метрикой тогда и только тогда, когда $(X,\rho)$
является пространством с внутренней метрикой.}
\vskip7pt
{\bf Доказательство теоремы 1}.
\vskip7pt
Необходимость теоремы непосредственно
следует из леммы 1 и определения метрики Хаусдорфа.
\vskip7pt
Достаточность. Рассмотрим множество
$$\Omega (M, W, \varepsilon) = \bigcup \{\omega (x, v, \varepsilon) \cup
\omega (y, u, \varepsilon) :$$
$$x \in M,\, y \in W,\, v \in P (x, W, \varepsilon ),\,
u \in P(y, M, \varepsilon ) \},$$
где $\varepsilon > 0$.
\vskip7pt
$1.\, \,$ Докажем, что
$$2\beta (M,\Omega (M, W, \varepsilon)) \leq
\beta (M, W) + 2 \varepsilon.$$
Для каждого $\delta > 0$ найдется точка $x \in M$ такая, что
$$\beta (M,\Omega (M ,W, \varepsilon)) < |x\Omega (M, W, \varepsilon)| + \delta.$$
Следовательно, для каждого
 $\delta > 0$ найдется точка $x \in M$ такая, что для каждого
$v \in P (x, W, \varepsilon)$ верны неравенства
$$\beta (M,\Omega (M,W,\varepsilon)) < \dfrac{1}{2}|xv|
+ \dfrac{1}{2}\varepsilon + \delta \leq  \dfrac{1}{2}|xW| + \varepsilon + \delta.$$
Тогда для каждого $\delta  > 0$
$$\beta (M,\Omega (M, W, \varepsilon)) < \dfrac{1}{2}\beta (M, W)
+ \varepsilon + \delta.$$
В силу произвольности $\delta  > 0$ 
$$\beta (M, \Omega (M, W, \varepsilon)) \leq \dfrac{1}{2}\beta (M, W) + \varepsilon.$$

$2.\, \,$ Докажем теперь неравенство
$$2\beta (\Omega (M, W, \varepsilon), M) \leq
\alpha (M, W) + 3\varepsilon$$
методом от противного.
Пусть найдутся множества $M,\, W \in B[X]$ такие, что
$$2\beta (\Omega (M, W, \varepsilon), M) > \alpha (M, W) + 3 \varepsilon.$$
Тогда найдутся точки
$$z \in \Omega (M, W, \varepsilon),\quad u \in P (z, M, \varepsilon/2)$$
такие, что $$2 |zu| > \alpha (M, W) + 3 \varepsilon.$$
Теперь возможны лишь два случая.
\vskip7pt
Случай 1. Найдутся точки
$$v \in W,\quad g \in P (v, M, \varepsilon)$$
такие, что $z \in \omega (v, g, \varepsilon)$. Тогда
$$|vM| \geq  |vg| - \varepsilon > 2|zg| - 2\varepsilon \geq$$
$$2|zM| - 2\varepsilon \geq 2|zu| - 3\varepsilon > \alpha (M,W).$$
Получили противоречие.
\vskip7pt
Случай 2. Найдутся точки
$$t \in M,\quad q \in P (t,W, \varepsilon)$$
такие, что $z \in \omega (t,q,\varepsilon)$. Тогда
$$|tW| \geq  |tq| -
\varepsilon > 2|zt| - 2\varepsilon  \geq$$
$$2|zM| - 2\varepsilon \geq 2|zu| - 3\varepsilon > \alpha (M,W).$$
Снова получили противоречие. Таким образом,
неравенство
$$2\beta (\Omega (M,W,\varepsilon),M) \leq \alpha (M,W) + 3 \varepsilon$$
установлено. Из полученных неравенств следует неравенство
$$2\alpha (\Omega
(M,W,\varepsilon),M) < \alpha (M,W)+ 4\varepsilon.$$ Аналогично
доказывается, что
$$2\alpha (\Omega (M,W,\varepsilon),W) <
\alpha (M,W) + 4\varepsilon.$$
Из полученных неравенств и леммы 1 следует
справедливость теоремы 1.
\vskip7pt
Рассмотрим теперь более сильное условие $(A_0)$ на метрическое пространство
$(X,\rho)$.
\vskip7pt
$(A_0)\quad$ Для любых $x, \, y \in X$ множество
$$\omega (x,y) = \{z \in X : 2\max \{|xz|, |zy|\} \leq |xy|\}$$ непустое.
\vskip7pt
Нетрудно заметить, что это условие можно представить в следующей форме.

Для любых $x, \, y \in X$ найдется такой элемент $z \in X$, что
$$|xz| = |zy|= \dfrac{1}{2}|xy|.$$
Метрическое пространство, удовлетворяющее условию $(A_0)$, называется
{\it метрически выпуклым пространством} \cite {Bing}, \cite [с. 38]{Nadl}. 
Напомним, что метрическое пространство $X$ называется {\it выпуклым по Менгеру} \cite [с.43]{Bus}, если для любых различных точек $x$, $y \in X$ найдется отличная от них такая точка $z \in X$, что
$$|xz| + |zy| = |xy|.$$
Очевидно, что метрически выпуклое пространство является  метрическим пространством, выпуклым по Менгеру. 
В следующей теореме, в частности, указано, как построить середину для двух ограниченных множеств существования в метрически выпуклом пространстве.
\vskip7pt
{\bf Теорема 2} [s11]. {\it Пространство  $(X,\rho)$ метрически выпукло тогда и
только тогда, когда для каждых ограниченных множеств существования $M, \,W \subset X$ найдется множество $\Omega \in B(X)$ такое, что
$$\alpha(M,\Omega) = \alpha(\Omega,W) = \dfrac{1}{2}\alpha(M,W).$$
В качестве множества $\Omega$ достаточно взять следующее множество
$$\Omega =\bigcup \{\omega (x,v)\cup \omega (y,u) :$$
$$x \in M,\, y \in W,\, v \in P (x, W),\, u
\in P (y, M) \}.$$}

{\bf Доказательство теоремы 2}.
\vskip7pt
Достаточность теоремы непосредственно
следует из определения метрики Хаусдорфа. 
\vskip7pt
Необходимость. 
Докажем, что для выбранного множества $\Omega$ имеет место неравество
$$2\beta (M,\Omega) \leq \beta (M,W).$$
Для каждого $\delta  > 0$ найдется точка $x \in M$ такая, что
$$\beta (M,\Omega) < |x\Omega| +\delta.$$
Следовательно, для каждого
$\delta > 0$ найдется точка $x \in M$ такая, что для каждого
$v \in P (x, W)$ верны неравенства
$$\beta (M,\Omega) < \dfrac{1}{2}|xv| + \delta = \dfrac{|xW|}{2} +
\delta.$$
Тогда для каждого $\delta > 0$ верно неравенство
$$\beta (M,\Omega) < \dfrac{1}{2}\beta (M,W) + \delta.$$
В силу произвольности
$\delta > 0$,
$$\beta (M,\Omega) \leq \dfrac{1}{2} \beta (M,W).$$
Докажем теперь неравенство
$$2\beta (\Omega,M) \leq \alpha (M,W)$$
методом от противного. Пусть найдутся ограниченные множества существования $M, \, W \subset X$ такие, что
$$2\beta (\Omega,M) > \alpha (M,W).$$
Тогда найдутся точки
$$z \in \Omega,\quad u \in P (z, M)$$
такие, что верно неравенство
$$2|zu| > \alpha (M,W).$$
Теперь возможны лишь два случая.
\vskip7pt
$1.\, \,$ Найдутся точки
$$v \in W,\quad g \in P (v, M)$$
такие, что $z \in \omega (v,g)$.
Тогда
$$|vM| = |vg| \geq 2 |zu| > \alpha (M,W).$$
Получили противоречие.
\vskip7pt
$2.\, \,$ Найдутся точки
$$t \in M,\quad q \in P (t, W)$$
такие, что $z \in \omega (t,q)$.
Тогда
$$|tW| = |tq| \geq 2|zu| > \alpha (M,W).$$ 
Снова получили противоречие. Таким образом, неравенство
$$2\beta (\Omega,M) \leq \alpha (M,W)$$
установлено. Из полученных неравенств следует неравенство
$$2\alpha (\Omega,M) \leq \alpha (M,W).$$
Аналогично доказывается, что
$$2\alpha (\Omega,W) \leq \alpha (M,W).$$ 
Справедливость утверждения следует теперь из неравенства треугольника.
\vskip7pt
{\bf Следствие 1} [s11]. {\it Пусть пространство $(X,\rho)$ удовлетворяет условию $(A_0)$.
Если для каждых $x,\, y \in X$ множество $\omega (x,y)$ конечное, то метрическое пространство
$(\Sigma_{\mathbb{N}} (X),\alpha )$ удовлетворяет условию $(A_0)$, где $\Sigma_{\mathbb{N}} (X)$ --- множество всех непустых конечных подмножеств в $X$.}
\vskip7pt
Пусть $r \geq 0$. Обозначим через
$$B[M,r] = \{ x \in X : |xM| \leq r \}$$
{\it обобщенный замкнутый шар} для непустого подмножества $M \subset X$ \cite [с. 219]{Kurat}. В пространстве с внутренней метрикой имеет место простая оценка сверху расстояния в метрике Хаусдорфа между произвольными обобщенными замкнутыми шарами непустых ограниченных множеств пространства $X$.
\vskip7pt
{\bf Лемма 2} [s11]. {\it Пусть $(X,\rho)$ --- пространство с внутренней метрикой. Тогда имеет место неравенство:
$$\alpha (B[M,r],B[W,R]) \leq \alpha (M,W) + |R - r|$$
для всех
$$M,\, W \in B(X), \quad R >0,\quad r > 0.$$}

{\bf Доказательство леммы 2}.
\vskip7pt
$1.\, \,$ Докажем сначала, что
$$|xB[W,R]| = \max \{|xW| - R, 0\}$$
для каждого $x \in X$. Из неравенства треугольника для отклонения $\beta$ \cite{Marin1} следует:
$$|xW| - R \leq |xB[W,R]| + \beta (B[W,R],W) - R \leq |xB[W,R]|$$
для каждого $x \in X$ и значит, для каждого $x \in X$ верно неравенство
$$|xB[W,R]| \geq  \max \{|xW| - R, 0\}.$$
Используем метод доказательства от противного. Допустим, что равенства нет. Тогда найдется
точка $x \in X \backslash B[W,R]$ такая, что
$$|xW| - R < |xB[W,R]|.$$
Следовательно, найдется такое
$\delta \in (0,R)$, что
$$|xW| + \delta  - R < |xB[W,R]|.$$
Тогда найдется точка $y \in W$ такая, что
$$|xy| + \delta  - R < |xB[W,R]|.$$
Из условия внутренности метрики следует, что для каждого
$\varepsilon \in (0,\delta/2)$ найдется конечная последовательность точек
$$z_0 =x,\, z_1, \ldots ,\, z_k =y$$
такая, что для каждого $i \in \{0, 1, \ldots, k-1\}$ 
$$|z_iz_{i+1}| < \varepsilon \, \, \mbox{и}\, \, |z_0z_1| +\ldots + |z_{k-1}z_k| < |xy| +\varepsilon.$$
Тогда найдется точка
$$z_i \in \{z_1,\ldots, z_{k-1}\}$$ 
такая, что
$$z_i \in B[W,R],\quad |z_iy| \geq R - \varepsilon.$$
Cледовательно,
$$|xz_i| < |xy| - |z_iy| + \varepsilon  \leq$$
$$|xy| - R + 2\varepsilon <
|xB[W,R]| + 2\varepsilon - \delta < |xB[W,R]|.$$
Получили противоречие. 
\vskip7pt
$2.\,\,$ Используя установленное равенство, получим
$$\alpha  (B[M,r],B[W,R]) =
\max\{\sup\limits_{z \in B[M,r]}{|zB[W,R]|}; \sup\limits_{w \in B[W,R]}{|wB[M,r]|}\} = $$
 $$\max\{\sup\limits_{z \in B[M,r]}{\max\{|zW| - R,0\}}; \sup\limits_{w \in B[W,R]}{\max\{|wM|
 - r,0\}}\} \leq $$
$$\max\{\max\{\beta (B[M,r],W) - R,0\}; \max\{\beta (B[W, R],M) - r,0\}\} \leq $$
$$\max\{\max\{\beta (B[M,r],M) +\beta (M,W) - R, 0\};$$
$$\max\{\beta (B[W, R],W) + \beta (W,M) - r, 0\}\} \leq $$
$$\max\{\max\{\beta (M,W) + r - R, 0\};\max\{\beta (W,M) + R - r, 0\}\} \leq$$
$$\alpha (M,W) + |R - r|.$$
Таким образом, лемма 2 доказана.
\vskip7pt
Из теоремы 1 и леммы 2 следует
\vskip7pt
{\bf Следствие 2} [s11]. {\it Пусть $(X,\rho)$ --- пространство с
внутренней метрикой. Тогда  имеет место неравенство:
$$\alpha_B(B_B [M,r],B_B [W,R]) \leq \alpha (M,W) + |R - r|$$
для всех 
$$M, \, W \in B(X), \quad R >0,\quad r > 0,$$
где $\alpha_B : B[B(X)]\times B[B(X)]\rightarrow \mathbb{R}_{+}$
--- метрика Хаусдорфа, $B_B [M,r]$ --- замкнутый шар с центром в точке $M \in B(X)$ пространства $(B[B(X)], \alpha_B)$.}
\vskip7pt
{\bf Лемма 3} [s11]. {Пусть  $F$ --- подмножество пространства $(B[X],\alpha)$ с индуцированной
метрикой. Тогда имеет место неравенство:
$$\alpha (H,Q) \leq \alpha_F (H,Q)$$
для всех $H,\, Q \in F$, где 
$$\alpha_F : B(F)\times B(F)\rightarrow \mathbb{R}_{+}$$
--- псевдометрика Хаусдорфа.}
\vskip7pt
{\bf Доказательство леммы 3}.
\vskip7pt
Пусть $x \in H$ и $\underline{\in }$ --- знак принадлежности элемента в пространстве  $(F,\alpha )$.
Тогда найдется элемент $Z \underline{\in } H$ такой, что $x \in Z$. Кроме того,
$$|xQ| \leq \inf\limits_{Y \underline{\in} Q}{|xY|} \leq \inf\limits_{Y \underline{\in} Q}{\alpha  (Z,Y)} \leq$$
$$\sup\limits_{Z \underline{\in } H} \inf\limits_{Y \underline{\in} Q}{\alpha  (Z,Y)} = \beta _F (H,Q).$$
Следовательно,
$$\beta (H,Q) = \sup\limits_{x \in H}{|xQ|} \leq \beta_F (H,Q).$$
Неравенство $$\beta (Q,H) \leq \beta_F (Q,H)$$
устанавливается аналогично.

\section{Пространство всех $N$-сетей и симметризованная степень порядка $N$ метрического пространства}

\vskip20pt

Пусть $p \in [1,\infty]$ и $S(N)$ --- группа всех подстановок множества из $N \geq 1$ элементов. Рассмотрим на декартовом произведении $X^N$ из $N$ экземпляров метрического пространства
$(X,\rho)$ метрику
$$\rho_{N,p}: X^N \times  X^N \, \rightarrow \mathbb{R}_{+},$$
$$\rho_{N,p} ((x_1,\ldots ,x_N),(y_1,\ldots ,y_N)) = (|x_1y_1|^p + \ldots + |x_Ny_N|^p)^{1/p}$$
при $p \in [1,\infty)$  и
$$\rho_{N,\infty} ((x_1,\ldots ,x_N),(y_1,\ldots ,y_N)) = \max \{|x_1y_1|,\ldots ,|x_Ny_N|\}$$
при $p =\infty.$
Зададим на $X^N$ следующее отношение эквивалентности $\sim$:
$$(x_1,\ldots ,x_N) \sim (y_1,\ldots ,y_N),$$
если найдется такое $\sigma \in S(N)$, что
$$y_1=x_{\sigma (1)},\ldots ,y_N=x_{\sigma (N)}.$$
Полученное факторпространство $X_N = X^N /\!\!\sim$ множества $X^N$ по этому отношению эквивалентности есть симметризованная степень порядка $N$ пространства $X$, а элементы этого множества можно отождествить с $N$-сетями с повторениями пространства $X$. Рассмотрим на этом множестве метрику
$$\alpha_p  : X_N \times X_N \, \rightarrow \mathbb{R}_{+},$$
$$\alpha_p ([(x_1,\ldots ,x_N)],[(y_1,\ldots ,y_N)])=$$
$$\min \{\rho_{N,p} ((x_1,\ldots ,x_N),(y_{\sigma (1)},\ldots ,y_{\sigma (N)})) : \sigma \in S(N)\},$$
где $p \in [1,\infty]$ (при $p =\infty $, см. \cite{Fedor}). Ограничение этой метрики на множество $$X^*_N = \{[(x_1,\ldots,x_N)] \in X_N : \, \mbox{card\,} \{x_1,\ldots,x_N\} = N\},$$
где {\rm card\,}$\{x_1,\ldots,x_N\}$ --- мощность множества $\{x_1,\ldots,x_N\}$,
будем обозначать тем же символом $\alpha_p$.  Отметим также, что псевдометрика
$$\alpha_{*}  : X_N \times X_N \, \rightarrow \mathbb{R}_{+},$$
$$\alpha_{*} ([(x_1,\ldots ,x_N)],[(y_1,\ldots ,y_N)])=
\alpha (\{x_1,\ldots ,x_N\},\{y_1,\ldots ,y_N\}),$$
индуцирует на множестве $X^*_N$ метрику.
\vskip7pt
Прежде, чем исследовать некоторые свойства метрики $\alpha_p$, напомним следующие определения.
\vskip7pt
Кривая, соединяющая точки $x$, $y \in X$, длина которой равна расстоянию между этими точками, называется {\it сегментом} $[x,y]$  с концами $x,\, y \in  X$ \cite [с.42]{Bus}.
\vskip7pt
Метрическое пространство называется {\it геодезическим пространством}, если любые две его точки можно соединить сегментом \cite [с.4]{Brid}. Напомним, что полное метрически выпуклое пространство является геодезическим пространством \cite{Efrem}.
\vskip7pt
Метрическое пространство $(X,\rho)$ называется {\it собственным}, если любой его замкнутый шар компактен \cite [с. 2]{Brid}. 
\vskip7pt

Для $M$, $W \in B(X)$ обозначим
$$D(M,W) = \sup \{|xy| : x \in M, \, y \in W \}.\quad D(M) = D(M,M)$$
--- диаметр множества $M$.
\vskip7pt
{\bf Теорема 1} [s20]. {\it Пусть $p \in [1,\infty]$. Тогда верны следующие утверждения.
\vskip7pt
$(i)\, \,$ $X^*_N$ --- открытое множество пространства $(X_N,\alpha_p)$. На множестве $X_N$ имеет место неравенство $\alpha_{*} \leq \alpha_p$. Кроме того, при $N > 2$ для каждого
$$S = [(x_1,\ldots ,x_N)] \in X^*_N$$
имеет место равенство
$\alpha_{*}(S_1,S_2) = \alpha_{\infty}(S_1,S_2)$,
где $S_1$, $S_2$ --- произвольные элементы из открытого шара $B(S,\varepsilon) \subset (X_N,\alpha_{*})$ радиуса 
$$\varepsilon =
\dfrac{1}{4} \min \{|x_ix_j| : i \neq j; \, i, \, j \in \{1,\ldots,N\}\},$$ 
а при $N \leq 2$ имеет место равенство $\alpha_{*} = \alpha_{\infty}$.
\vskip7pt
$(ii)\, \,$ $(X,\rho)$ является собственным пространством тогда и только тогда, когда $(X_N,\alpha_p)$~--- собственное пространство. 
\vskip7pt
$(iii)\, \,$ $(X,\rho)$  является пространством с внутренней метрикой тогда и только тогда, когда
$(X_N,\alpha_p)$, $(X^*_N,\alpha_p)$ --- пространства с внутренней метрикой. Кроме того, $X^*_N$ --- всюду плотное подмножество пространства $(X_N,\alpha_p)$ при {\rm card\,}$(X) > 1$.
\vskip7pt
$(iv)\, \,$ $(X,\rho)$ является геодезическим {\rm {(}}метрически выпуклым, выпуклым по Менгеру{\rm {)}} пространством тогда и только тогда, когда $(X_N,\alpha_p)$ --- геодезическое {\rm {(}}метрически выпуклое, выпуклое по Менгеру{\rm {)}} пространство.}
\vskip7pt
{\bf Доказательство теоремы 1.}
\vskip7pt
$(i)\, \,$ Нетрудно проверить, что для каждого
$$S = [(x_1,\ldots ,x_N)] \in X^*_N$$
открытый шар $B(S,\delta)$ пространства $(X_N,\alpha_p)$, где
$$\delta = \dfrac{1}{2} \min \{|x_ix_j| : i \neq j; \, i, \, j \in \{1,\ldots,N\}\},$$
принадлежит множеству $X^*_N$. Следовательно, $X^*_N$ --- открытое множество пространства $(X_N,\alpha_p)$. Докажем, что $\alpha_{*} \leq \alpha_p$. Для каждого $\sigma \in S(N)$ имеем, что
$$\alpha_{*} ([(x_1,\ldots ,x_N)],[(y_1,\ldots ,y_N)])= \max \{\max \{|x_i\{y_1,\ldots ,y_N\}| :$$
$$1 \leq i \leq N\}, \max \{|y_j\{x_1,\ldots ,x_N\}| : 1 \leq j \leq N\}\} \leq$$
$$\max \{\rho_{N,p}((x_1,\ldots ,x_N),(y_{\sigma (1)},\ldots ,y_{\sigma (N)})),
\rho_{N,p}((y_1,\ldots ,y_N),(x_{\sigma (1)},\ldots ,x_{\sigma (N)}))\}.$$
Осталось взять минимум по всем $\sigma \in S(N)$ в правой части неравенства и использовать определение метрики $\alpha_p$.

Пусть
$$N > 2, \quad S_1,\, S_2 \in B(S,\varepsilon) \subset (X_N,\alpha_{*}).$$
Тогда
$$\alpha_{*}(S_1,S_2) < 2 \varepsilon.$$
Учитывая определения псевдометрики $\alpha_{*}$ и $\varepsilon$, получим, что
$$S_1,\, S_2 \in X^*_N \quad \alpha_{*}(S_1,S_2) = \alpha_{\infty}(S_1,S_2).$$
Пусть теперь $N = 2$ (случай, когда $N = 1$, очевиден) и для определенности положим
$$S = [(x, y)],\quad S_1 = [(u, v)],\quad |xu| =  D(S,S_1).$$
Тогда
$$\alpha_{*} (S,S_1) = \max \{\max \{|xv|, |y\{u,v\}|\};\max \{|yu|, |v\{x,y\}|\}\}.$$

Рассмотрим два случая.
\vskip7pt
$1.\, \,$ Если $|yu| > |yv|$, то
$$\alpha_{*} (S,S_1) = \max \{\max \{|xv|, |yv|\};|yu|\} = $$
$$\max \{|xv|, |yu|\} = \min \{\max \{|xu|, |yv|\};\max \{|xv|, |yu|\}\} = \alpha_{\infty}(S,S_1).$$

$2.\, \,$ Если $|yu| \leq |yv|$, то
$$\alpha_{*} (S,S_1) = \max \{\max \{|xv|, |yu|\};\max \{|yu|, |v\{x,y\}|\}\} =$$
$$\max \{|xv|, |yu|\} = \min \{\max \{|xu|, |yv|\};\max \{|xv|, |yu|\}\} = \alpha_{\infty} (S,S_1).$$

$(ii)\, \,$ Доказательство достаточности очевидно. Докажем необходимость.
Пусть
$$S = [(x_1,\ldots ,x_N)] \in X_N \, \, \mbox {и} \, \, (S_n = [(y^n_1,\ldots ,y^n_N)])$$
--- ограниченная последовательность пространства $(X_N,\alpha_p)$, то есть найдется такая вещественная константа $c > 0$, что для каждого $n \in \mathbb{N}$ справедливо неравенство:
$$\alpha_p (S_n,S) \leq c.$$
Тогда для каждого $n \in \mathbb{N}$ найдется такое $\sigma_n \in S(N)$, что
$$\rho_{N,p} ((x_1,\ldots ,x_N),(y^{n}_{\sigma_n (1)},\ldots ,y^{n}_{\sigma_n (N)})) \leq c.$$
Следовательно, для каждого $k \in \{1,\ldots,N\}$ найдется подпоследовательность
$(y^m_{\sigma_m (k)})$ последовательности $(y^n_{\sigma_n (k)}),$
сходящаяся к некоторому
$u_k \in X$ при $m \rightarrow \infty$, поскольку пространство $X$ --- собственное. Положим
$$S_0 = [(u_1,\ldots ,u_N)].$$
Тогда 
$$\alpha_p (S_m,S_0) \leq
\rho_{N,p} ((u_1,\ldots ,u_N),(y^{m}_{\sigma_m (1)},\ldots ,y^{m}_{\sigma_m (N)}))
\rightarrow 0$$
при $m \rightarrow \infty$.  Таким образом, $(X_N,\alpha_p)$ --- собственное пространство.
\vskip7pt
$(iii)\, \,$ Доказательство достаточности для пространства $(X_N,\alpha_p)$ очевидно. Докажем необходимость. Пусть произвольно выбраны
$$S = [(x_1,\ldots ,x_N)],\quad S_1 = [(y_1,\ldots ,y_N)] \in X_N.$$
Используя лемму 1.1.1, для каждого $\varepsilon > 0$ выберем
$$S_2 = [(z_1,\ldots ,z_N)],$$
где $z_i \in \omega (x_i,y_{\sigma(i)},\varepsilon_1)$ для каждого $i \in \{1,\ldots,N\}$, 
$\varepsilon_1 = \varepsilon N^{-1/p}$
при $p \in [1,\infty)$, $\varepsilon_1 = \varepsilon$ при $p = \infty$ и $\sigma \in S(N)$ такое, что  
$$\rho_{N,p} ((x_1,\ldots ,x_N),(y_{\sigma (1)},\ldots ,y_{\sigma (N)})) = \alpha_p (S,S_1).$$
Тогда
$$2 \max \{\alpha_p (S,S_2), \alpha_p (S_2,S_1)\} \leq$$
$$2 \max \{\rho_{N,p}((x_1,\ldots ,x_N),(z_1,\ldots ,z_N)),
\rho_{N,p}((z_1,\ldots ,z_N),(y_{\sigma (1)},\ldots ,y_{\sigma (N)}))\} <$$
$$\rho_{N,p}((x_1,\ldots ,x_N),(y_{\sigma (1)},\ldots ,y_{\sigma (N)})) + \varepsilon = \alpha_p (S,S_1) + \varepsilon.$$
В силу леммы 1.1.1, $(X_N,\alpha_p)$ --- пространство с внутренней метрикой.

С помощью леммы 1.1.1 нетрудно доказать, что $X^*_N$ --- всюду плотное подмножество пространства $(X_N,\alpha_p)$ при {\rm card}$X > 1$. Из этих двух доказанных утверждений и леммы 1.1.1 следует теперь, что $(X^*_N,\alpha_p)$ --- пространство с внутренней метрикой.
\vskip7pt
$(iv)\, \,$ Доказательство достаточности очевидно. Докажем необходимость. Пусть произвольно выбраны
$$S = [(x_1,\ldots ,x_N)],\, \, S_1 = [(y_1,\ldots ,y_N)] \in X_N$$
и пусть $\sigma \in S(N)$ такое, что
$$\rho_{N,p} ((x_1,\ldots ,x_N),(y_{\sigma (1)},\ldots ,y_{\sigma (N)})) = \alpha_p (S,S_1).$$
Для каждого
$i \in \{1,\ldots,N\}$ в силу геодезичности пространства $X$ найдется некоторый сегмент $[x_i,y_{\sigma (i)}]$. Выберем для каждого $i \in \{1,\ldots,N\}$ и для каждого $\lambda \in [0,1]$ такую точку $z_i(\lambda) \in [x_i,y_{\sigma (i)}]$, что
$$|x_iz_i(\lambda)| = \lambda |x_iy_{\sigma (i)}|.$$
Для каждого $\lambda \in [0,1]$ положим
$$S(\lambda) = [(z_1(\lambda),\ldots ,z_N(\lambda))].$$
Тогда для каждого $\lambda \in [0,1]$,
используя неравенство треугольника, нетрудно получить равенство
$$\alpha_p (S,S(\lambda)) + \alpha_p (S(\lambda),S_1) = \alpha_p (S,S_1).$$
Следовательно, $S(\lambda)$ при изменении параметра $\lambda$ на отрезке $[0,1]$ есть параметризация некоторого сегмента с концами $S$, $S_1$ в пространстве $(X_N,\alpha_p)$, и это пространство геодезическое. Утверждения этого пункта теоремы 1, приведенные в скобках, теперь нетрудно доказать, используя аналогию с доказанным случаем. Таким образом, теорема 1 доказана.
\vskip7pt
Обозначим через $\Sigma^{*}_N(X)$ (через $\Sigma_N(X)$) множество всех (непустых) подмножеств в $(X,\rho )$, состоящих (не более чем) из $N$ точек, с индуцированной метрикой Хаусдорфа $\alpha$. Причем обозначение $\Sigma_1(X)$ будем заменять на $X$. Элементы множества $\Sigma_N(X)$ называются $N$-сетями \cite{Gark}.

Рассмотрим сюръекцию
$$f : X_N \, \rightarrow \Sigma_N(X), \quad f([(x_1,\ldots ,x_N)]) =
\{x_1,\ldots ,x_N\}.$$
Очевидно, что при $N \leq 2$ сюръекция $f$ является изометрией пространства $(X_2,\alpha_{\infty})$ на пространство $(\Sigma_2(X),\alpha)$.
Используя эту изометрию, из теоремы 1 получаем
\vskip7pt
{\bf Следствие 1} [s20]. {\it $(X,\rho)$ --- собственное пространство {\rm {(}}пространство с внутренней метрикой, геодезическое пространство, метрически выпуклое пространство, выпуклое по Менгеру пространство{\rm {)}} тогда и только тогда, когда
$(\Sigma_2(X),\alpha)$ --- собственное пространство {\rm {(}}пространство с внутренней метрикой, геодезическое пространство, метрически выпуклое пространство, выпуклое по Менгеру пространство{\rm{)}}.}
\vskip7pt
Для того, чтобы при $N > 2$ пространство $(X^*_N,\alpha_p)$ было геодезическим пространством, на пространство $(X,\rho)$ приходится налагать более жесткие условия.
\vskip7pt
{\bf Теорема 2} [s20]. {\it Пусть $p \in [1,\infty]$ и $(X,\rho)$ --- геодезическое пространство, удовлетворяющее следующим двум условиям $(A_1)$, $(A_2)$.
\vskip7pt
$(A_1)\, \,$ Каждые две различные точки пространства $(X,\rho)$ можно соединить единственным сегментом.
\vskip7pt
$(A_2)\, \,$ Если два сегмента имеют общий конец и общую внутреннюю точку, то один из этих сегментов есть подмножество другого сегмента.
\vskip7pt
Тогда $(X^*_N,\alpha_p)$ --- геодезическое пространство.
\\
Если кроме условий $(A_1)$, $(A_2)$ выполнено следующее глобальное условие $(A_3)$ неположительности кривизны по Буземану {\rm {(см. \cite [с. 304]{Bus})}}
$$2|\omega (z,x)\omega (z,y)| \leq |xy| \leqno (A_3)$$
для любых $x,\, y,\, z \in  X$, то при $N>1$ для каждого
$$S = [(x_1,\ldots ,x_N)] \in (X^*_N,\alpha_p)$$
и для любых
$W$, $T$, $D \in B (S,\varepsilon)$, где $B (S,\varepsilon)$ --- открытый шар пространства
$(X^*_N,\alpha_p)$ радиуса
$$\varepsilon =
\dfrac{1}{4} \min \{|x_ix_j| : i \neq j; \, i, \, j \in \{1,\ldots,N\}\},$$
верно неравенство
$$2\alpha_p (\omega (W,D),\omega (T,D)) \leq \alpha_p (W,T),$$
то есть при $N>1$ пространство
$(X^*_N,\alpha_p)$ удовлетворяет локальному условию
неположительности кривизны по Буземану.}
\vskip7pt
{\bf Доказательство теоремы 2.}
\vskip7pt
При {\rm card\,}$(X) = 1$ или $N =1$ теорема 2, очевидно, верна. Предположим, что {\rm card\,}$(X) > 1$ и $N > 1$.
Пусть произвольно выбраны различные
$$S = [(x_1,\ldots ,x_N)], \quad S_1 = [(y_1,\ldots ,y_N)] \in X^*_N.$$
Рассмотрим сначала частный случай, когда все точки
$$\{x_1,\ldots ,x_N,y_1,\ldots ,y_N\}$$
принадлежат одному сегменту и $p \in [1,\infty)$. Без потери общности можно считать, что точки упорядочены следующим образом:
$$x_1 > \ldots > x_N,\quad y_1 > \ldots > y_N,\quad x_1 \geq y_1.$$
Докажем, что
$$\rho_{N,p} ((x_1,\ldots ,x_N),(y_1,\ldots ,y_N)) = \alpha_p (S,S_1),$$
индукцией по $N$.
Пусть $N = 2$. Если $x_2 \leq y_1$, то нетрудно получить неравенство
$$\rho_{2,p} ((x_1,x_2),(y_1,y_2)) \leq \rho_{2,p} ((x_1,x_2),(y_2,y_1)),$$
из которого следует требуемое равенство. Пусть $x_2 > y_1$. Заметим, что функция
$$f : \mathbb{R}_{+} \rightarrow \mathbb{R}_{+},\quad f(t) = (a + t)^p - (b + t)^p,$$
где $b \leq a$, $a$, $b \in \mathbb{R}_{+}$, 
неубывающая.  С учетом этого, требуемое равенство следует теперь из неравенств
$$|x_1y_1|^p - |x_2y_1|^p \leq (|x_1y_1| + |y_1y_2|)^p - (|x_2y_1| + |y_1y_2|)^p = |x_1y_2|^p - |x_2y_2|^p.$$
Итак, при $N = 2$ рассматриваемое равенство установлено. Предположим, что это равенство верно для всех $N \leq n -1$ и докажем его для $N = n$. Пусть $\sigma \in S (N)$. Рассмотрим два случая.
\vskip7pt
$(i)\, \,$ Если $\sigma (1) = 1$, то по предположению индукции
$$\rho_{N-1,p} ((x_2,\ldots ,x_N),(y_2,\ldots ,y_N)) \leq \rho_{N-1,p} ((x_2,\ldots ,x_N),(y_{\sigma (2)},\ldots ,y_{\sigma (N)}))$$ 
или
$$\rho_{N-1,p} ((x_2,\ldots ,x_N),(y_2,\ldots ,y_N)) \leq \rho_{N-1,p} ((y_2,\ldots ,y_N),(x_{\sigma (2)},\ldots ,x_{\sigma (N)})).$$ 
Следовательно, 
$$\rho_{N,p} ((x_1,\ldots ,x_N),(y_1,\ldots ,y_N)) \leq \rho_{N,p} ((x_1,\ldots ,x_N),(y_{\sigma (1)},\ldots ,y_{\sigma (N)})).$$ 
или
$$\rho_{N,p} ((x_1,\ldots ,x_N),(y_1,\ldots ,y_N)) \leq \rho_{N,p} ((y_1,\ldots ,y_N),(x_{\sigma (1)},\ldots ,x_{\sigma (N)})).$$ 

$(ii)\, \,$ Если $\sigma (1) > 1$ и $\sigma (k) = 1$, где $k \in \{2,\ldots ,N\}$, то, используя установленные неравенства пункта $(i)$ и предположение индукции, нетрудно получить неравенства:
$$\rho_{N,p} ((x_1,\ldots ,x_N),(y_1,\ldots ,y_N)) \leq$$ 
$$\rho_{N,p} ((x_1,\ldots ,x_N),(y_{\sigma (k)},\ldots ,y_{\sigma (k-1)},y_{\sigma (1)},y_{\sigma (k+1)},\ldots ,y_{\sigma (N)})) \leq$$
$$\rho_{N,p} ((x_1,\ldots ,x_N),(y_{\sigma (1)},\ldots ,y_{\sigma (N)})).$$ 
Следовательно, требуемое равенство доказано.  

Выберем для каждого $\lambda \in [0,1]$ и для каждого $k \in \{1,\ldots,N\}$ такую точку $z_k(\lambda) \in [x_k,y_k]$, что
$$|x_kz_k(\lambda)| = \lambda |x_ky_k|.$$
Тогда для каждого $\lambda \in [0,1]$ справедливы неравенства:
$$z_1(\lambda) > \ldots > z_N(\lambda).$$ 
Положим для каждого $\lambda \in [0,1]$ 
$$S(\lambda) = [(z_1(\lambda),\ldots ,z_N(\lambda))].$$
Так же как и при доказательстве утверждения $(iv)$ теоремы 1, используя доказанные выше неравенства, получаем, что $S(\lambda)$ при изменении параметра $\lambda$ на отрезке $[0,1]$ есть параметризация некоторого сегмента с концами $S$, $S_1$ в пространстве $(X^*_N,\alpha_p)$. Таким образом, в частном случае теорема 2 верна. 

Рассмотрим общий случай. Пусть $p \in [1,\infty]$ и
$$\Xi = \{\sigma \in S(N) : \rho_{N,p} ((x_1,\ldots ,x_N),(y_{\sigma (1)},\ldots ,y_{\sigma (N)})) = \alpha_p (S, S_1)\}.$$
Индукцией по $N$ докажем вспомогательное утверждение $(H)$.
\vskip7pt
$(H)\quad$ Найдется такое $\sigma \in \Xi$,
что для любых различных $i$, $j \in \{1,\ldots,N\}$, для каждого 
$\lambda \in [0,1]$ справедливо неравенство $z_i(\lambda) \neq z_j(\lambda),$ 
где точка 
$$z_k(\lambda) \in [x_k,y_{\sigma (k)}]$$ 
при $k \in \{i,j\}$ такая, что  
$$|x_kz_k(\lambda)| = \lambda |x_ky_{\sigma (k)}|,$$
а сегмент $[x_k,y_{\sigma (k)}]$ существует и является единственным в силу условия $(A_1)$. 
\vskip7pt
Пусть $N = 2$. Предположим противное. Тогда для каждого $\sigma \in \Xi$ найдется такое $\lambda \in (0,1)$, что  $z_1(\lambda) = z_2(\lambda)$. Используя определение точек
$z_1(\lambda)$, $z_2(\lambda)$ и неравенство треугольника, получим
$$\rho_{2,p} ((x_1,x_2),(y_{\sigma (2)},y_{\sigma (1)})) \leq \rho_{2,p} ((x_1,x_2),(y_{\sigma (1)},y_{\sigma (2)})).$$
Если полученное неравенство строгое, то в силу определения множества
$\Xi$ и метрики $\alpha_p$ получили противоречие. Если имеет место равенство, то в силу условий
$(A_1)$, $(A_2)$ точки
$$x_1,\quad x_2,\quad y_{\sigma (1)},\quad y_{\sigma (2)}$$
принадлежат одному сегменту и $\Xi = S(2)$. Тогда в силу определений точек $z_1(\lambda)$, $z_2(\lambda)$ и метрики $\alpha_p$ снова получили противоречие.

Допустим, что утверждение $(H)$ верно для всех $N \leq n -1$ и докажем его для $N = n$.
Предположим противное. Тогда для каждого $\sigma \in \Xi$ найдутся такие различные $i$, $j \in \{1,\ldots,N\}$ и такое $\lambda \in (0,1)$, что  $z_i(\lambda) = z_j(\lambda)$. В этом случае
пару сегментов 
$$\{[x_i,y_{\sigma (i)}],[x_j,y_{\sigma (j)}]\}$$ 
назовем отмеченной.

Пусть $\sigma \in \Xi$ и $p \in [1,\infty)$. В семействе сегментов
$$D = ([x_k,y_{\sigma (k)}])_{k \in \{1,\ldots,N\}}$$
заменим отмеченную пару 
$$\{[x_i,y_{\sigma (i)}],[x_j,y_{\sigma (j)}]\}$$
на пару 
$$\{[x_i,y_{\sigma (j)}],[x_j,y_{\sigma (i)}]\}.$$ 
Это приводит к новой подстановке $\pi$, полученной перестановкой $\sigma (i)$ и $\sigma (j)$ в подстановке $\sigma$, в частности, $\pi (i) = \sigma (j)$, $\pi (j) = \sigma (i)$. Тогда, учитывая установленное неравенство при $N = 2$, получим
$$\rho_{N,p}((x_1,\ldots ,x_N),(y_{\pi (1)},\ldots ,y_{\pi (N)})) \leq \rho_{N,p}((x_1,\ldots ,x_N),(y_{\sigma (1)},\ldots ,y_{\sigma (N)})).$$
Если полученное неравенство строгое, то в силу определения множества $\Xi$ и метрики $\alpha_p$ получили противоречие. Если имеет место равенство, то в силу условий $(A_1)$, $(A_2)$ точки
$$x_i,\quad y_{\pi (i)},\quad x_j,\quad y_{\pi (j)}$$
принадлежат одному сегменту. Кроме того, $\pi \in \Xi$ и пара 
$$\{[x_i,y_{\pi (i)}], [x_j,y_{\pi (j)}]\}$$ 
--- неотмеченная. Учитывая предположение индукции, можно считать, что, если  удалить из семейства сегментов 
$$\Delta_1 = ([x_k,y_{\pi (k)}])_{k \in \{1,\ldots,N\}}$$ 
сегменты 
$$[x_i,y_{\pi (i)}],\quad [x_j,y_{\pi (j)}],$$ 
то в оставшемся семействе не останется отмеченных пар сегментов. 

Если в семействе $\Delta_1$ нет отмеченных пар сегментов, то в силу определения множества $\Xi$ и метрики $\alpha_p$ получаем противоречие. Предположим для определенности, что пара 
$$\{[x_i,y_{\pi (i)}], [x_l,y_{\pi (l)}]\}$$ 
--- отмеченная. Заменив ее в семействе $\Delta$ на пару 
$$\{[x_i,y_{\pi (l)}],[x_l,y_{\pi (i)}]\},$$ 
получим новое семейство сегментов $\Delta_2$. Аналогично предыдущему, вводя новую подстановку и повторяя рассуждения, мы либо сразу получаем противоречие, либо получаем, что точки
$$x_i,\quad y_{\pi (i)},\quad x_j, \quad y_{\pi (j)}, \quad x_l, \quad y_{\pi (l)}$$
принадлежат одному сегменту $L$. Учитывая доказанный частный случай, можно считать, что на сегменте $L$ нет отмеченных пар сегментов. А по предположению индукции можно считать, что нет отмеченных пар сегментов и среди сегментов семейства $\Delta_1$, расположенных вне сегмента $L$.  Продолжая такую процедуру замены отмеченных пар сегментов, мы через конечное число шагов получим противоречие, поскольку для точек, расположенных на одном сегменте, теорема 2 доказана и группа подстановок --- конечная. 

Пусть $\sigma \in \Xi$ и $p = \infty$. В семействе сегментов $D$ выберем сегмент $[x_i,y_{\sigma (i)}]$ длины $\alpha_p (S,S_1)$. Учитывая предположение индукции, можно считать, что, если удалить этот сегмент из семейства $D$, то в оставшемся семействе не останется отмеченных пар сегментов. Пусть пара 
$$\{[x_i,y_{\sigma (i)}],[x_j,y_{\sigma (j)}]\}$$ 
отмеченная. Если в $D$ такого сегмента $[x_j,y_{\sigma (j)}]$ не существует, то в силу определения множества $\Xi$ и метрики $\alpha_p$, сразу получаем противоречие. Заменив в семействе $D$ отмеченную пару 
$$\{[x_i,y_{\sigma (i)}],[x_j,y_{\sigma (j)}]\}$$
на пару 
$$\{[x_i,y_{\sigma (j)}],[x_j,y_{\sigma (i)}]\},$$ 
получим семейство сегментов $D(1)$. Это приводит к новой подстановке $\pi$, полученной перестановкой $\sigma (i)$ и $\sigma (j)$ в подстановке $\sigma$. Длины сегментов 
$$[x_i,y_{\sigma (j)}],\quad [x_j,y_{\sigma (i)}]$$ 
строго меньше, чем $\alpha_p (S,S_1)$. В семействе $D(1)$ выберем новый сегмент длины $\alpha_p (S,S_1)$ и повторим рассуждения. Через конечное число шагов получим противоречие, поскольку группа подстановок --- конечная и на каждом шаге сегменты из отмеченных пар заменяются на сегменты меньшей длины, чем $\alpha_p (S,S_1)$. Таким образом, утверждение $(H)$ доказано.

Пусть $\sigma \in \Xi$ такое, как описано в $(H)$.
Для каждого $\lambda \in [0,1]$ положим
$$S(\lambda) = [(z_1(\lambda),\ldots ,z_N(\lambda))].$$
Так же, как и в доказательстве утверждения $(iv)$ теоремы 1, получаем, что $S(\lambda)$ при изменении параметра $\lambda$ на отрезке $[0,1]$ есть параметризация некоторого сегмента с концами $S$, $S_1$ в пространстве $(X^*_N,\alpha_p)$. 

Докажем теперь последнее утверждение теоремы. Пусть $S$,
$W$, $T$, $D$ выбраны в соответствии с условиями теоремы. Тогда из
неравенства
$$\alpha_p (S,T) < \varepsilon,$$
неравенства треугольника и определения $\varepsilon$  следует, что для каждого $t \in f(T)$ найдется
единственный элемент $x(t) \in f(S)$ такой, что
$$|tx(t)| < \varepsilon,$$
а также для каждого $x \in f(S)$ найдется такой элемент $u \in f(T)$, что
$$|ux| < \varepsilon.$$
Учтем также, что $S \in X^*_N$. Тогда получим,
что для каждого $x \in f(S)$ в шаре $B(x,\varepsilon)$ содержится точно один
элемент $t(x) \in f(T)$. Аналогичное рассуждение справедливо и для $W$, $D$. Соответствующие элементы обозначим $w(x) \in f(W)$, $d(x) \in f(D)$.
Тогда для каждого $x \in f(S)$ имеем:
$$2|\omega (w(x),d(x))\omega (t(x),d(x))| \leq |w(x)t(x)|.$$
Теперь нетрудно проверить, что
$$2\alpha_p (\omega (W,D),\omega (T,D)) \leq \alpha_p (W,T).$$
Таким образом, теорема 2 доказана.
\vskip7pt
Используя сюръекцию $f$, для каждого $p \in [1,\infty]$ построим симметрику
$$\hat \alpha_p  : \Sigma_N(X) \times  \Sigma_N(X) \, \rightarrow \mathbb{R}_{+}, \qquad \hat \alpha_p (S,S_1)= \alpha_p (f^{-1}(S),f^{-1}(S_1)).$$
Ограничение $f_1$ сюръекции $f$ на подмножество $\Delta_N \cup X^*_N \subset X_N$, где
$$\Delta_N = \{[(x_1,\ldots ,x_N)] \in X_N : x_1 = \ldots = x_N\},$$
является биекцией на свой образ $X\cup \Sigma^{*}_N (X) \subset \Sigma_N(X)$. Следовательно, ограничение симметрики $\hat \alpha_p$ на множество
$$(X\cup \Sigma^{*}_N (X))\times (X \cup \Sigma^{*}_N (X))$$
является метрикой и
$$f_1 : (\Delta_N \cup X^*_N,\alpha_p) \rightarrow (X\cup \Sigma^{*}_N (X),\hat \alpha_p)$$
--- изометрия. Отметим также, что $\Sigma^{*}_N(X)$ является открытым подмножеством пространства $(\Sigma_N(X),\alpha)$.
Действительно, пусть $N > 1$ и $S \in \Sigma^{*}_N(X)$,
тогда открытый шар $B(S,\varepsilon)$ пространства $(\Sigma_N(X),\alpha)$, где
$$\varepsilon =\dfrac{1}{2} \min \{|xy| : x, \, y \in S, \, x \neq y \},$$
принадлежит подмножеству $\Sigma^{*}_N(X)$.
Используя изометрию $f_1$, получаем следствие 2 теоремы 1 и следствие 3 теоремы 2.
\vskip7pt
{\bf Следствие 2} [s20]. {\it Пусть $p \in [1,\infty]$. Тогда верны следующие утверждения
\vskip7pt
$(i)\, \,$ На множестве $X\cup \Sigma^{*}_N (X)$ имеет место неравенство $\alpha \leq \hat \alpha_p$. Кроме того, при $N > 2$ для каждого
$$S = \{x_1,\ldots ,x_N\} \in \Sigma^{*}_N (X)$$
имеет место равенство
$$\alpha(S_1,S_2) = \hat \alpha_{\infty}(S_1,S_2),$$
где $S_1$, $S_2$ --- произвольные элементы из открытого шара $B(S,\varepsilon) \subset (\Sigma_N (X),\alpha)$ радиуса 
$$\varepsilon =
\dfrac{1}{4} \min \{|x_ix_j| :\, i \neq j; \, i, \, j \in \{1,\ldots,N\}\},$$ 
а при $N \leq 2$ имеет место равенство $\alpha = \hat \alpha_{\infty}$. 
\vskip7pt
$(ii)\, \,$ $(X,\rho)$ является пространством с внутренней метрикой тогда и только тогда, когда
$(\Sigma^{*}_N (X),\hat \alpha_p)$ --- пространство с внутренней метрикой.}
\vskip7pt
{\bf Следствие 3} [s20]. {\it Пусть $p \in [1,\infty]$ и $(X,\rho)$ --- геодезическое пространство, удовлетворяющее условиям $(A_1)$, $(A_2)$.
Тогда $(\Sigma^{*}_N (X),\hat \alpha_p)$ --- геодезическое пространство. Если, кроме условий $(A_1)$, $(A_2)$, выполнено условие $(A_3)$,
то при $N>1$ для каждого
$$S = \{x_1,\ldots ,x_N\} \in (\Sigma^{*}_N (X),\hat \alpha_p)$$
и для любых
$W$, $T$, $D \in B (S,\varepsilon)$, где $B (S,\varepsilon)$ --- открытый шар пространства
$(\Sigma^{*}_N (X),\hat \alpha_p)$ радиуса
$$\varepsilon =
\dfrac{1}{4} \min \{|x_ix_j| : i \neq j; \, i, \, j \in \{1,\ldots,N\}\},$$
имеет место неравенство
$$2\hat \alpha_p (\omega (W,D),\omega (T,D)) \leq \hat \alpha_p (W,T),$$
то есть при $N>1$ пространство $(\Sigma^{*}_N (X),\hat \alpha_p)$ удовлетворяет локальному условию неположительности кривизны по Буземану.}
\vskip7pt
Введем на множестве $X_N$ следующее отношение эквивалентности $R$:
$$[(x_1,\ldots ,x_N)]R[(y_1,\ldots ,y_N)],$$
если
$$\{x_1,\ldots ,x_N\} = \{y_1,\ldots ,y_N\}.$$
Рассмотрим для $p \in [1,\infty]$ псевдометрику (см.\cite [с. 73]{Bur})
$$\alpha_{p,R} : X_N \times X_N \rightarrow \mathbb{R}_{+},$$
$$\alpha_{p,R} (S,\hat S) = \inf \{\alpha_p (S_1,\tilde S_1) + \ldots + \alpha_p (S_k,\tilde S_k)\},$$ 
где инфимум берется по всем $k \in \mathbb{N}$ и таким наборам $\{S_i\}$, $\{\tilde S_i\},$ $1 \leq i \leq k$, что $S_1 = S$, $\tilde S_k = \hat S$ и точка $\tilde S_i$ $R$-эквивалентна точке $S_{i+1}$ при всех $i= 1,\ldots,k-1$. Стандартным образом сопоставим псевдометрическому пространству $(X_N,\alpha_{p,R})$ метрическое пространство $(X_N/ R,\alpha_{p,R})$, отождествляя точки, находящиеся на нулевом расстоянии и сохраняя обозначение $\alpha_{p,R}$ для полученной фактор-метрики.
Отображение
$$g : X_N / R \rightarrow \Sigma_N (X), \quad g([[(x_1,\ldots ,x_N)]]_R) = \{x_1,\ldots ,x_N\}$$ является биекцией, с помощью которой мы отождествим фактор-множество $X_N / R$ с множеством всех $N$-сетей $\Sigma_N (X)$ пространства $X$ и наделим множество $\Sigma_N (X)$ фактор-метрикой $\alpha_{p,R}$.
\vskip7pt
{\bf Теорема 3} [s20]. {\it Пусть $p \in [1,\infty]$. Тогда верны следующие утверждения
\vskip7pt
$(i)\, \,$  На множестве $\Sigma_N(X)$ справедливо неравенство $\alpha \leq \alpha_{p,R}$. Кроме того, при $N > 2$ для каждого 
$$S = \{x_1,\ldots ,x_N\} \in \Sigma^{*}_N (X)$$ 
имеет место равенство
$$\alpha(S_1,S_2) = \alpha_{\infty,R}(S_1,S_2),$$ 
где $S_1$, $S_2$ --- произвольные элементы из открытого шара $B(S,\varepsilon) \subset (\Sigma_N (X),\alpha)$ радиуса
$$\varepsilon = \dfrac{1}{4} \min \{|x_ix_j| : i \neq j; \, i, \, j \in \{1,\ldots,N\}\},$$ 
а при  $N \leq 2$ верно равенство $\alpha =  \alpha_{\infty,R}$.
\vskip7pt
$(ii)\, \,$ $(X,\rho)$ является собственным пространством тогда и только тогда, когда $(\Sigma_N(X),\alpha_{p,R})$ --- собственное пространство.
\vskip7pt
$(iii)\, \,$ $(X,\rho)$ является пространством с внутренней метрикой тогда и только тогда, когда
$(\Sigma_N(X),\alpha_{p,R})$, $(\Sigma^*_N(X),\alpha_{p,R})$ --- пространства с внутренней метрикой. Кроме того, $\Sigma^*_N(X)$ --- всюду плотное подмножество пространства
$(\Sigma_N(X),\alpha_{p,R})$ при {\rm card\,}$X > 1$.
\vskip7pt
$(iv)\, \,$ $(X,\rho)$ является геодезическим {\rm {(}}метрически выпуклым, выпуклым по Менгеру{\rm {)}} собственным пространством тогда и только тогда, когда
$(\Sigma_N (X),\alpha_{p,R})$ --- геодезическое {\rm {(}}метрически выпуклое, выпуклое по Менгеру{\rm {)}} собственное пространство.}
\vskip7pt
{\bf Доказательство теоремы 3.}
\vskip7pt
$(i)\, \,$ Нетрудно понять, что для доказательства  неравенства $\alpha \leq \alpha_{p,R}$ достаточно доказать неравенство $\alpha_{*} \leq \alpha_{p,R}$ для псевдометрик на множестве $X_N$. Выберем произвольно $\varepsilon > 0$. Пусть $S$, $\hat S \in X_N$. В силу определения псевдометрики $\alpha_{p,R}$, найдутся такие
$$S_1, \ldots,S_k,\tilde S_1,\ldots, \tilde S_k \in X_N,$$
что  $S_1 = S$, $\tilde S_k = \hat S$, $\tilde S_iRS_{i+1}$ при всех $i= 1,\ldots,k-1$ и
$$\alpha_p (S_1,\tilde S_1) + \ldots + \alpha_p (S_k,\tilde S_k) \leq \alpha_{p,R}(S,\hat S) + \varepsilon.$$
Тогда, используя неравенство треугольника и утверждение $(i)$ теоремы 1, получим
$$\alpha_{*} (S,\hat S) \leq \alpha_{*}(S_1,\tilde S_1) + \alpha_{*}(\tilde S_1,S_2) + \alpha_{*} (S_2,\tilde S_2)+ \ldots + \alpha_{*} (\tilde S_{k-1},S_k) +$$
$$\alpha_{*} (S_k,\tilde S_k) = \alpha_{*}(S_1,\tilde S_1) + \alpha_{*} (S_2,\tilde S_2) + \ldots  +
\alpha_{*} (S_k,\tilde S_k) \leq$$
$$\alpha_p (S_1,\tilde S_1) + \ldots + \alpha_p (S_k,\tilde S_k) \leq \alpha_{p,R}(S,\hat S) + \varepsilon.$$
Используя произвол в выборе $\varepsilon > 0$,
получим требуемое неравенство. Оставшиеся равенства следуют из соответствующих равенств
утверждения $(i)$ следствия 2, если учесть доказанное неравенство и следующее простое неравенство для псевдометрики и метрики на $X_N$: $\alpha_{p,R} \leq \alpha_p$ для $p \in [1,\infty]$.
\vskip7pt
$(ii)\, \,$ Доказательство достаточности очевидно. Для доказательства необходимости достаточно доказать, что произвольный замкнутый шар
$B[[S]_R, r]$ в пространстве $(\Sigma_N(X),\alpha_{p,R})$ компактен. Рассмотрим полный прообраз
$$M = \pi^{-1}(B[[S]_R, r])$$
этого шара относительно канонической проекции
$$\pi : X_N \rightarrow  \Sigma_N(X),\quad \pi (S) = [S]_R.$$
Из непрерывности канонической проекции следует, что множество $M$ замкнуто. Кроме того,
$$M \subset B[S,r_1] \subset (X_N,\alpha_p),$$
где
$$S \in [S]_R,\quad r_1 = (N)^{1/p} (r + D(f(S))).$$
Действительно, используя определения метрик
$\alpha$, $\alpha_p$, сюръекции $f$, доказанное утверждение $(i)$ и неравенство треугольника, для каждого $\hat S \in M$ получим неравенства
$$\alpha_p(S,\hat S) \leq (N)^{1/p} \alpha_{\infty}(S,\hat S) \leq
(N)^{1/p} (\alpha (f(S),f(\hat S)) + D(f(S))) \leq$$
$$(N)^{1/p} (\alpha_{p,R} (f(S)),f(\hat S)) + D(f(S))) \leq r_1.$$
Следовательно, множество $M$ ограничено и замкнуто в пространстве $(X_N,\alpha_p)$. В силу утверждения $(ii)$ теоремы 1, $M$ компактно в пространстве $(X_N,\alpha_p)$. Тогда его образ $\pi (M) = B[[S]_R,r]$ относительно непрерывного отображения $\pi$ является компактным множеством в пространстве $(\Sigma_N(X),\alpha_{p,R})$.
\vskip7pt
$(iii)\, \,$ Доказательство достаточности для пространства $(\Sigma_N(X),\alpha_{p,R})$ очевидно. Докажем необходимость. Пусть произвольно выбраны
$$S = [(x_1,\ldots ,x_N)],\quad \hat S = [(y_1,\ldots ,y_N)]$$
в псевдометрическом пространстве $(X_N,\alpha_{p,R})$. Для произвольно выбранного $\varepsilon > 0$ найдутся такие
$$S_1, \ldots,S_k,\tilde S_1,\ldots, \tilde S_k \in X_N,$$
что
$$S_1 = S,\quad \tilde S_k = \hat S,\quad \tilde S_iRS_{i+1}$$
при всех $i= 1,\ldots,k-1$ и
$$C = \alpha_p (S_1,\tilde S_1) + \ldots + \alpha_p (S_k,\tilde S_k) < \alpha_{p,R}(S,\hat S) + \varepsilon.$$
Найдем такое $j \in \{1,\ldots,k-1\}$, что
$$\alpha_p (S_1,\tilde S_1) + \ldots + \alpha_p (S_j,\tilde S_j) \leq C/2$$
и
$$\alpha_p (S_1,\tilde S_1) + \ldots + \alpha_p (S_{j+1},\tilde S_{j+1}) > C/2.$$
Пусть
$$\varepsilon_1 = \alpha_{p,R}(S,\hat S) + \varepsilon - C > 0.$$
Используя утверждение $(iii)$ теоремы 1 и лемму 1.1.1, выберем такое
$$S^{*} = [(z_1,\ldots ,z_N)] \in X_N,$$
что
$$2 \max \{\alpha_p (S_{j+1},S^{*}),\alpha_p (S^{*},\tilde S_{j+1})\} < \alpha_p (S_{j+1},\tilde S_{j+1}) + \varepsilon_1.$$
Тогда
$$\alpha_p (S_1,\tilde S_1) + \ldots + \alpha_p (S_j,\tilde S_j) + \alpha_p (S_{j+1},S^{*}) +
\alpha_p (S^{*},\tilde S_{j+1}) +$$
$$\alpha_p (S_{j+2},\tilde S_{j+2}) + \ldots + \alpha_p (S_k,\tilde S_k) < C + \varepsilon_1 = \alpha_{p,R}(S,\hat S) + \varepsilon.$$
Теперь нетрудно проверить, что
$$2 \max \{\alpha_{p,R} (S,S^{*}),\alpha_{p,R} (S^{*},\hat S)\} < \alpha_{p,R}(S,\hat S) + \varepsilon.$$
Следовательно, в силу леммы 1.1.1, $(\Sigma_N(X),\alpha_{p,R})$ --- пространство с внутренней метрикой.

С помощью леммы 1.1.1 нетрудно доказать, что $\Sigma^*_N(X)$ --- всюду плотное подмножество пространства $(\Sigma_N(X),\alpha_{p,R})$ при {\rm card}$X > 1$. Из этих двух доказанных утверждений и леммы 1.1.1 следует  теперь верность утвержения для пространства $(\Sigma_N(X),\alpha_{p,R})$.
\vskip7pt
$(iv)\, \,$ Доказательство достаточности очевидно, а необходимость является следствием доказанного утверждения $(iii)$ и известной теоремы Хопфа-Ринова и Кон-Фоссена (см. \cite {Bur}, с. 60). Таким образом, теорема 3 доказана.
\vskip7pt
Отметим для сравнения, что из предложения 1
\cite{Foert} следует, что в геодезическом пространстве, удовлетворяющим условиям $(A_1)$, $(A_2)$, $(A_3)$, множество всех непустых ограниченных
замкнутых выпуклых множеств с метрикой Хаусдорфа является геодезическим пространством, удовлетворяющем глобальному условию неположительности кривизны по Буземану.
В следующем cледствии 4 теорем 1, 3 установлены полезные достаточные условия, при которых две $N$-сети могут быть соединены сегментом в пространстве $(\Sigma_N (X),\alpha)$. 
\vskip7pt
{\bf Следствие 4} [s20]. {\it Пусть $(X,\rho)$ --- геодезическое пространство.
\vskip7pt
$(i)\, \,$ Если $S_1$, $S_2 \in \Sigma^{*}_N(X)$ и
$$\alpha(S_1,S_2) = \hat \alpha_{\infty}(S_1,S_2),$$
то $N$-сети $S_1$, $S_2$ могут быть соединены сегментом в пространстве $(\Sigma_N(X),\alpha)$.
\vskip7pt
$(ii)\, \,$ Если $(X,\rho)$ --- собственное пространство, $S_1$, $S_2 \in \Sigma_N(X)$ и
$$\alpha (S_1,S_2) = \alpha_{\infty,R} (S_1,S_2),$$
то $N$-сети $S_1$, $S_2$ могут быть соединены сегментом в пространстве $(\Sigma_N(X),\alpha)$.}
\vskip7pt
{\bf Примеры.}
\vskip7pt
$I.\, \,$ Рассмотрим в евклидовой плоскости $\mathbb{R}^2$ две $3$-сети
$$M= \{(0;0),(-a;-a),(-a;a)\}, \quad W=\{(0;0),(a;a),(a;-a)\},$$
где $a \in \mathbb{R}_{+}$.
Тогда нетрудно найти, что
$$\alpha (M,W)=\sqrt{2} a,\quad \hat \alpha_{\infty} (M,W)= \alpha_{\infty,R}(M,W) = 2a$$
и множество
$$\Omega =\{(0;0),(-a/2;-a/2),(-a/2;a/2),(a/2;a/2),(a/2;-a/2)\} \in \Sigma_5(\mathbb{R}^2)$$
из теоремы 1.1.2.
Кроме того,
$$2\alpha(M,T)=2\alpha(W,T)=\alpha(M,W),$$
$$2\hat \alpha_{\infty}(M,D)=2\hat \alpha_{\infty}(W,D)= \hat \alpha_{\infty}(M,W),$$
где
$$T=\{(-a/2;-a/2),(-a/2;a/2),(a/2;a/2),(a/2;-a/2)\} \in \Sigma_4(\mathbb{R}^2),$$
$$D=\{(0;0),(0;a),(0;-a)\} \in \Sigma_3(\mathbb{R}^2).$$
Теперь нетрудно установить что, в любой окрестности $3$-сети
$$O=\{(0;0)\} \in (\Sigma_3(\mathbb{R}^2),\alpha)$$
существуют две $3$-сети, несоединимые сегментом в пространстве
$$(\Sigma_3(\mathbb{R}^2),\alpha)\, \, \mbox{и} \, \, \alpha \neq \hat \alpha.$$

$II.\, \,$  Рассмотрим на прямой $\mathbb{R}$ три $3$-сети
$$M= \{0,2a,3a + b\}, \quad W=\{a,2a + b,4a + b\}, \quad T = \{a,3a + b\},$$
где $a$, $b \in \mathbb{R}_{+}$. В следующих трех случаях нетрудно подсчитать расстояния.
\vskip7pt
$1.\, \,$ Если $b \leq a$, то
$$\alpha (M,W)= \alpha_{\infty,R}(M,W) = \hat \alpha_{\infty} (M,W)= a.$$

$2.\, \,$ Если $a < b \leq 2a$, то
$$\alpha (M,W)= a < b = \hat \alpha_{\infty} (M,W)= \alpha_{\infty,R}(M,W).$$

$3.\, \,$ Пусть $2a < b$. Если $a = 0$, то $M = W = T$. Если $a > 0$, то
$$\alpha (M,W)= a < 2a = \alpha (M,T) + \alpha (T,W) =$$
$$\alpha_{\infty,R}(M,W) < b = \hat \alpha_{\infty} (M,W),$$
и не существует сегмента с концами $M$, $W$ в пространстве $(\Sigma_3(\mathbb{R}),\alpha)$.
\vskip7pt
Сравнивая полученные результаты, нетрудно заметить, что пространство $(\Sigma_N(X),\alpha)$ лишь при $N = 2$ обладает достаточно хорошими геометрическими свойствами, зависящими от аналогичных свойств пространства $(X,\rho)$. Поэтому исследуем этот случай немного подробнее.

Рассмотрим метрическое пространство $(X,\rho)$ с выделенным семейством
сегментов $S$, удовлетворяющее следующим условиям.
\vskip7pt
$(B_1)\, \,$  Каждый подсегмент произвольного сегмента из $S$ принадлежит $S$.
\vskip7pt
$(B_2)\, \,$  Для любых  $x$, $y \in X$ существует единственный сегмент $[x, y] \in  S$.
\vskip7pt
Для $\lambda \in  [0,1]$, $x, y \in X$ обозначим через $\omega_{\lambda}(x,y)$
точку сегмента $[x,y] \in  S$ такую, что
$$|x\omega_{\lambda}(x,y)| = \lambda |xy|.$$

Кроме условий $(B_1)$, $(B_2)$ будем использовать следующее глобальное условие неположительности кривизны по Буземану пространства $(X,\rho)$ \cite [c. 63]{Bus2}.
\vskip7pt
$(B_3)\, \,$ Для любых $$x, \, y,\, z \in  X,\quad 2|\omega_{1/2}(z,x)\omega_{1/2}(z,y)| \leq |xy|.$$
\vskip7pt
{\bf Лемма 1} [s16]. {\it Пусть метрическое пространство $X$
удовлетворяет условиям $(B_1)$, $(B_2)$, $(B_3)$. Тогда
\\
отображение
$$\pi : (\Sigma_2 (X),\alpha) \rightarrow (X,\rho),\quad
\pi (\{x,y\}) = \omega_{1/2} (x,y)$$
обладает следующими свойствами.
\vskip7pt
$(i)\, \,$ Для любых $S,\, S_1 \in \Sigma_2 (X)$ справедливы неравенства:
$$|\pi (S)\pi (S_1)| \leq \alpha (S,S_1) \leq |\pi (S)\pi (S_1)| +
\dfrac{1}{2}(D(S) + D(S_1)).$$
$(ii)\, \,$ Для любых $x,\, y \in X$ имеет место равенство:
$$\inf \{\alpha (S,S_1) : S \in \pi^{-1}(x),\, S_1 \in \pi^{-1}(y)\} = |xy|.$$}

{\bf Доказательство леммы 1}.
\vskip7pt
$(i)\,\,$ Пусть
$$|xu| =  D(S,S_1),\quad S = \{x, y\},\quad S_1 = \{u, v\}.$$
Тогда в силу условий $(B_1)$, $(B_2)$, $(B_3)$ и изометричности пространств $(X_2,\alpha_{\infty})$, $(\Sigma_2(X),\alpha)$ получим:
$$|\pi (S)\pi (S_1)| = |\omega_{1/2}(x,y)\omega_{1/2}(u,v)| \leq$$
$$|\omega_{1/2}(x,y)\omega_{1/2}(y,v)| + |\omega_{1/2}(y,v)\omega_{1/2}(u,v)| \leq \dfrac{1}{2}|xv| + \dfrac{1}{2}|yu| \leq$$
$$\max \{|xv|, |yu|\} = \min \{\max \{|xu|, |yv|\};\max \{|xv|, |yu|\}\} = \alpha (S,S_1).$$
С другой стороны
$$\alpha (S,S_1) = \min \{\max \{|xu|, |yv|\};\max \{|xv|, |yu|\}\} = \max \{|xv|, |yu|\} \leq$$
$$\max \{|x\pi (S)| + |\pi (S)\pi (S_1)| + |\pi (S_1)v|; |u\pi (S_1)| + |\pi (S_1)\pi (S)| + |\pi (S)y|\} =$$
$$|\pi (S)\pi (S_1)| + \dfrac{1}{2} (D(S) + D(S_1)).$$
Доказательство $(ii)$ сразу следует из $(i)$. Таким образом, лемма 1 доказана.
\vskip7pt
Введем следующие обозначения.
$$\Delta = \{\{x,y\} \in \Sigma_2(X) :  x = y \},$$
$$\tau(S) =\{S_1 \in \Sigma_2(X) :  \alpha (S,S_1) = D(S,S_1) \},$$
где $S \in \Sigma_2 (X)$.
Для $S = \{x, y\} \in \Sigma_2(X)$ положим $S(x) = x$ при $x =y$ и $S(x)= y$ при $x \neq y$.

Пусть пространство $X$ удовлетворяет условиям $(B_1)$, $(B_2)$, $(B_3)$ и $S, \, S_1 \in \Sigma_2(X)$. Введем множество
$$\Omega (S,S_1) = \{\{\omega_{1/2} (x,S_1 (u)),\omega_{1/2} (S(x),u)\} \in \Sigma_2(X) :$$
$$x \in S, \, u \in S_1, \, |xu| = D(S,S_1)\}$$
и рассмотрим cледующее условие на метрическое пространство $X$.
\vskip7pt
$(B_4)\, \,$ Для любых различных точек $x, \, y \in X$ найдется единственная точка $z \in X$ такая, что
$y = \omega_{1/2} (x,z)$.
\vskip7pt
{\bf Лемма 2} [s16]. {\it Пусть метрическое пространство $X$
удовлетворяет условиям $(B_1)$, $(B_2)$, $(B_3)$. Тогда верны следующие свойства.
\vskip7pt
$(i)\, \,$ Если $S_1 \not \in \tau(S)\backslash \Delta $, то существует
единственная $2$-сеть
$$\Omega (S,S_1)  \in \Sigma_2(X).$$

$$S \in \Sigma_2(X) \backslash \Delta,\quad S_1 \in \tau(S)\backslash \Delta \leqno (ii)$$
тогда и только тогда, когда множество
$$\Omega (S,S_1) \subset \Sigma_2 (X)$$
состоит из двух различных $2$-сетей.
\vskip7pt
$(iii)\, \,$ Если пространство $X$ удовлетворяет также условию $(B_4)$, то
для каждого $x \in \Delta $ и для каждого
$$S \in \Sigma_2 (X)\backslash \{x\}$$
существует единственная $2$-сеть $S_1 \in \Sigma_2 (X)$ такая, что
$$S = \Omega (x,S_1).$$}

{\bf Доказательство леммы 2}.
\vskip7pt
$(i)\, \,$ Пусть
$$S = \{x, y\},\quad S_1 = \{u, v\},\quad |xu| = D(S,S_1) > \alpha (S,S_1).$$
Тогда, в силу изометричности пространств $(X_2,\alpha_{\infty})$, $(\Sigma_2 (X),\alpha)$,
$$\alpha (S,S_1) = \max \{|xv|, |yu|\}.$$
Следовательно,
$$\Omega (S,S_1) = \{\omega_{1/2}(x,v),\,\omega_{1/2}(y,u)\}.$$

$(ii)\, \,$ Пусть $$|xu| = \alpha (S,S_1) = D(S,S_1).$$
Тогда
$$\alpha (S,S_1) = \max \{|xv|, |yu|\}.$$
Следовательно,
$$\Omega (S,S_1) = \{\{\omega_{1/2}(x,v),\omega_{1/2}(y,u)\},\{\omega_{1/2}(x,u),\omega_{1/2}(y,v)\}\}.$$
$(iii)\, \,$ Пусть $S = \{u, v\}$. В силу условий $(B_1)$, $(B_2)$, $(B_3)$, $(B_4)$, найдутся единственная точка
$z \in X$ и единственная точка $w \in X$ такие,
что
$$u = \omega_{1/2}(x,z),\quad v = \omega_{1/2}(x,w).$$
Следовательно, $S = \Omega (x,S_1)$, где $S_1 = \{z, w\}$. Таким образом, лемма 2 доказана.
\vskip7pt
{\bf Примеры.}
\vskip7pt
$III.\, \,$ Пусть полуплоскость
$$P = \{(x;y) \in \mathbb{R}^{2} : x \geq y \}$$
наделена метрикой:
$$d : P\times P \rightarrow \mathbb{R},\quad d((x; y),(u; v))
= \max \{|x - u|, |y - v|\}.$$
Тогда отображение
$$f : (P,d) \rightarrow (\Sigma_2(\mathbb{R}),\alpha),\quad f((x;y)) = \{x, y\}$$
является изометрией, поскольку для всех
$(x; y), \, (u; v) \in P$
$$\min \{\max \{|x - u|, |y - v|\}, \max \{|x - v|, |y - u|\}\} =
\max \{|x - u|, |y - v|\}.$$
При этом $f ^{-1} (\Omega (S,S_1))$ изображается серединой отрезка
$$[f ^{-1}(S),f ^{-1}(S_{1})],$$
где $S,\, S_1 \in \Sigma_2 (\mathbb{R})$.
\vskip7pt
$IV.\, \,$ Пусть $(X,\rho)$ --- евклидово пространство и $S \in \pi^{-1}(x)$. Тогда
$$\alpha (S,\pi^{-1}(y)) = |xy|$$
для всех   $x,\, y \in X$.
\vskip7pt
$V.\, \,$ Пусть в открытом шаре $B(0,1)$ вещественного гильбертового
пространства задана метрика $\rho$ следующим образом:
$$|xy|= Arch \left(\dfrac{1- (x,y)}{((1 -x^2)(1 - y^2))^{
1/2}}\right),$$
где $(x,y)$ --- скалярное произведение векторов $x, \, y $ из
$B(0,1)$. Таким образом, получили бесконечномерный вариант известной
интерпретации Бельтрами--Клейна геометрии Лобачевского \cite [с. 48]{Nut}. В пространстве $(B(0,1),\rho)$ рассмотрим такие невырожденные
перпендикулярные отрезки $[x,y]$, $[a,b]$, что $x = \omega_{1/2}(a,b)$.
Тогда, используя соображения симметрии и тригонометрию плоскости
Лобачевского, нетрудно получить равенство
$$\sh (\alpha (\{a,b\},\pi^{-1}(y))) =
 \ch \left(\dfrac{|ab|}{2}\right)\sh (|xy|).$$
Следовательно, в этом случае имеет место неравенство:
$$\alpha (\{a,b\},\pi^{-1}(y)) > |xy|.$$

\section{Выпуклые множества в обобщенном хордовом пространстве}

\vskip20pt

{\it Обобщенным хордовым пространством} [s6] назовем метрическое пространство $(X,\rho)$ с выделенным семейством $S$ сегментов, называемых хордами, удовлетворяющее следующим условиям  $(C_1-C_8)$.
\vskip7pt
$(C_1)\, \, $  Каждый подсегмент произвольной хорды из $S$ является хордой.
\vskip7pt
$(C_2)\, \, $  Для любых  $x$, $y \in X$ существует хорда $T(x, y) \in  S$ с концами $x$, $y$.
\vskip7pt
$(C_3)\, \, $ Если в параметризации пропорциональной длине дуги с областью определения равной отрезку $[0,1]$ последовательность хорд поточечно сходится к некоторому сегменту, то этот сегмент является хордой.
\vskip7pt
$(C_4)\, \, $ Для каждой точки $p \in X$ существует открытый шар $B(p,r_p)$ $(r_p > 0)$ такой, что всякие две различные точки из этого шара можно соединить единственной хордой.
\vskip7pt
$(C_5)\, \, $ Для каждой точки $p \in X$ всякую хорду $T(x, y)$, соединяющую различные точки из шара $B(p,r_p)$, можно единственным образом продолжить в обе стороны до большей хорды, т. е. найдется такая хорда $T(u,v) \subset B(p,r_p)$, что
$$T(x, y) \subset T(u,v) \backslash \{u,v\}$$
и если произвольные две хорды из шара $B(p,r_p)$ имеют общий конец и общую внутреннюю точку, то одна из этих хорд есть подмножество другой хорды.
\vskip7pt
Часто вместо условия $(C_5)$ будем использовать более сильное условие $(PG)$.
\vskip7pt
$(PG)\, \, $ Всякую хорду можно единственным образом продолжить до такой геодезической линии, у которой каждый принадлежащий ей сегмент пространства $X$ является хордой. Здесь и далее область параметризации геодезической линии совпадает с $\mathbb{R}$ \cite [c. 48]{Bus},  \cite [c. 3]{Bus2}.
\vskip7pt
Для компактной записи следующих условий введем несколько обозначений. Если выполнено условие $(PG)$, то $K = \mathbb{R}$, в ином случае $K = [0,1]$.
Для произвольного $x \in X$ $W_x$ --- множество всех $y \in X$ таких, что существует единственная хорда $T(x,y)$. $T_{\lambda}(x,y)$ --- точка, которая находится из следующих условий:
\vskip7pt
$(i)\, \,$ если $\lambda \in [0,1]$, то
$$T_{\lambda}(x,y) \in T(x,y),\quad |xT_{\lambda}(x,y)| = \lambda |xy|;$$
\vskip7pt
$(ii)\, \,$ если выполнено условие $(PG)$ и $\lambda > 1$ $(\lambda < 0)$, то
$$|xT_{\lambda}(x,y)| = |\lambda| |xy|$$
и точка $T_{\lambda}(x,y)$ принадлежит геодезической, содержащей хорду $T(x,y)$, таким образом, чтобы точка $y$  $(x)$ лежала между точкой $x$  $(y)$ и точкой $T_{\lambda}(x,y)$.
\vskip7pt
Для каждого $\lambda \in K$ положим
$$\omega_{\lambda} (x,y) = \{T_{\lambda}(x,y) : T(x,y) \in S\} \eqno (1).$$

$(C_6)\, \, $ Из условий
$$x = \lim \limits_{n \rightarrow \infty} x_n,\quad y = \lim \limits_{n \rightarrow \infty} y_n$$ следует, что для каждого
$\lambda \in K$ верно равенство
$$\lim \limits_{n \rightarrow \infty} |\omega_{\lambda}(x_n,y_n)\omega_{\lambda}(x,y)| = 0.$$

$(C_7)\, \, $ Для каждого $x \in X$ и для каждого $\lambda \in (0,1)$ отображение
$$\lambda_x : W_x \rightarrow W_x,\quad \lambda_x (y) = \omega_{\lambda} (x,y),$$
отображает открытые множества в $X$, принадлежащие множеству $W_x$, в открытые множества в $X$.
\vskip7pt
$(C_8) \, \, $ Для каждого $p \in X$ существует непрерывное отображение
$$\tau_p : B(p,r_p) \rightarrow B(p,r_p),$$
удовлетворяющее условию
$$p = \omega_{1/2}(x,\tau_p (x))$$
для каждого $x \in B(p,r_p)$.
\vskip7pt
Нетрудно проверить, что условия $(C_1-C_8)$ выполняются в хордовом пространстве \cite [c. 23]{Bus2}, а условия $(C_7)$, $(C_8)$ следуют из условий $(C_1-C_4)$, $(C_6)$, $(PG)$.
\vskip7pt
Обобщенное хордовое пространство назовем {\it обобщенным \\
$G$-пространством Буземана}, если оно удовлетворяет условию $(G)$.
\vskip7pt
$(G) \quad$ Cемейство $S$ совпадает с семейством всех сегментов пространства $X$.
\vskip7pt
Нетрудно проверить, что $G$-пространство Буземана \cite [c. 54]{Bus} является обобщенным $G$-пространством Буземана.
\vskip7pt
В дальнейшем нам понадобятся следующие определения.
\vskip7pt
Множество $M$ пространства $X$, удовлетворяющего условиям $(C_1)$, $(C_2)$, назовем {\it $U$-множеством}, если для каждых двух различных точек из множества $M$ существует единственная хорда с концами в этих точках. Это аналог $U$-множества в хордовом пространстве \cite [c. 73]{Bus2}.
\vskip7pt
{\it Обобщенным прямым хордовым пространством} ({\it обобщенным прямым
$G$-пространством Буземана}) назовем обобщенное хордовое пространство (обобщенное $G$-пространство Буземана), которое является $U$-множеством и удовлетворяет условию $(PG)$. Это аналоги прямого хордового пространства и прямого $G$-пространства Буземана \cite [c. 26]{Bus2}. Банахово пространство, в котором в качестве хорд выбраны прямолинейные сегменты, является примером обобщенного прямого хордового пространства.
\vskip7pt
{\it Прямой} в пространстве $X$, удовлетворяющим условиям $(C_1)$, $(C_2)$, назовем геодезическую линию, у которой каждый принадлежащий ей сегмент пространства $X$ является хордой.
\vskip7pt
Множество $M$ пространства $X$, удовлетворяющего условиям  $(C_1)$, $(C_2)$, назовем {\it выпуклым}, если для каждых двух различных точек из $M$ все хорды с концами в этих точках принадлежат множеству $M$  \cite [с. 72]{Bus2}.
\vskip7pt
Множество $M$ с непустой внутренностью обобщенного хордового пространства назовем {\it строго выпуклым}, если множество $\bar {M}$ является $U$-множеством и для каждых двух различных точек $x$, $y \in \bar {M}$ каждая точка $z$, удовлетворяющая условию $[xzy]$, является внутренней точкой множества $M$ (ср. \cite [с. 153]{Bus}). Здесь $[xzy]$ означает, что точки $x$, $y$, $z$ попарно различные и точка $z$ принадлежит единственной хорде с концами $x$, $y$.
\vskip7pt
Пусть $E$ --- выпуклое множество вещественной прямой $\mathbb{R}$. Функция
$$f : E \rightarrow \mathbb{R}$$
называется {\it безвершинной}, если она непрерывна и
$$f(y) \leq \max \{f(x), f(z)\}$$
для всех $x$, $y$, $z \in E$ таких, что $x < y < z$ \cite [с. 144]{Bus}.
\vskip7pt
Функция
$$f : X \rightarrow \mathbb{R}$$
называется {\it слабо безвершинной}, если она непрерывна, пространство $X$ удовлетворяет условиям $(C_1)$, $(C_2)$ и
$$f(y) \leq \max \{f(x), f(z)\}$$
для всех $x$, $y$, $z$ таких, что $[xyz]^W$, где $[xyz]^W$ означает, что точки $x$, $y$, $z$ попарно различные и точка $y$ принадлежит некоторой хорде с концами $x$, $z$ \cite [с. 73, с. 149]{Bus2}.
\vskip7pt
Прямую $\Lambda_1$ в обобщенном прямом хордовом пространстве назовем {\it перпендикуляром} к прямой $\Lambda$ в точке $x$, если $\{x\}= \Lambda\cap \Lambda_1$ и каждая точка прямой $\Lambda_1$
имеет своим основанием на прямой $\Lambda$ точку $x$. При этом прямая $\Lambda$ называется {\it трансверсалью} к прямой $\Lambda_1$ в точке $x$ (ср. \cite [с. 10]{Bus2}, \cite [с. 158]{Bus}).
\vskip7pt
Пусть $M$ --- выпуклое множество с непустой внутренностью обобщенного хордового пространства
такое, что $\bar M \subset B(p,r_p)$. Фиксируем точку $q \in \partial M$, где $\partial M$ --- граница множества $M$. {\it Опорной хордой} множества $M$ в точке $q$ назовем такую хорду $T(x,y)$, что 
$$x,\, y \in B(p,r_p), \, [xqy]\, \, \mbox{и}\, \, T(x,y)\cap \stackrel{\circ}{M} = \emptyset.$$

Обозначим множество всех точек на всех опорных хордах в точке $q$ через $V_q$ (ср. \cite [с. 79]{Bus2}). Введем множества (ср. \cite [с. 79]{Bus2})
$$A_q = \{x \in B(p,r_p) : \exists y \, (y \in  \stackrel{\circ}{M}, \, [xyq]) \},$$
$$C_q = \{x \in B(p,r_p) : \exists y \, (y \in  \stackrel{\circ}{M}, \, [xqy]) \},$$
где $\stackrel{\circ}{M}$ --- внутренность множества $M$.
\vskip7pt
{\it Полукасательной хордой} множества $M \subset X$ в точке $q$ назовем хорду $T(q,x)$, удовлетворяющую следующим условиям:
\vskip7pt
$(i)\, \,$ хорда $T(q,x)$ принадлежит некоторой опорной хорде множества $M$ в точке $q$;
\vskip7pt
$(ii)\, \,$ существует последовательность точек $(y_n)_{n \in \mathbb{N}}$ из $A_q$ такая, что
 $$\lim \limits_{n \rightarrow \infty} y_n = z, \quad z \in T(q,x)\backslash \{q\}.$$
Обозначим множество всех точек на всех полукасательных хордах множества $M$ в точке $q$ через $\tilde V_q$ (ср. \cite [с. 79]{Bus2}).
\vskip7pt
{\bf Замечание 1.} {\it В обобщенном прямом хордовом пространстве имеют смысл понятия опорных и полукасательных прямых, если в определениях множеств $V_q$, $A_q$, $C_q$, $\tilde V_q$ заменить открытый шар $B(p,r_p)$ на пространство $X$.}
\vskip7pt
В следующей лемме 1 приведены некоторые элементарные свойства рассмотренных понятий. Доказательства этих свойств те же, что и доказательства аналогичных свойств, приведенные в \cite [с. 10, 73]{Bus2}, \cite [с. 156]{Bus}.
\vskip7pt
{\bf Лемма 1} [s6]. {\it
\vskip7pt
$(i)\, \,$ Непрерывная функция
$$f : X \rightarrow \mathbb{R}$$
слабо безвершинная тогда и только тогда, когда для каждого $r \in \mathbb{R}$ множество
$$\{x \in X : f(x) \leq r\}$$
выпуклое.
\vskip7pt
$(ii)\, \,$ Все замкнутые шары пространства $X$, удовлетворяющего условиям $(C_1)$, $(C_2)$, выпуклые тогда и только тогда, когда для каждого $p \in X$ функция
$$f : X \rightarrow \mathbb{R},\quad f(x) = |px|$$
--- слабо безвершинная.
\vskip7pt
$(iii)\, \,$ Если пространство $X$, удовлетворяет условиям $(C_1)$, $(C_2)$, точка $y$ является основанием точки $p$ на множестве $M$ и выполнется условие $[pxy]^W$, то точка $y$ является основанием точки $x$ на множестве $M$. Если, кроме того, пространство $X$ удовлетворяет условиям $(C_5)$, $(G)$, то точка $y$ является единственным основанием точки $x$ на множестве $M$.
\vskip7pt
$(iv)\, \,$ Пусть $M$ является сегментом, лучом или геодезической линией пространства $X$ с представлением $x=x(t)$. Если для $u \in [u_1,u_2]$ $z=z(u)$
--- непрерывная параметризованная кривая такая, что множество  $P (z(u), M)$ связное для каждого $u \in
[u_1,u_2]$ и
$$x(t_1) \in P (z(u_1), M),\quad x(t_2) \in P (z(u_2), M),\quad t_3 \in [t_1,t_2],$$
то найдется число $u_3 \in [u_1,u_2]$ такое, что $x(t_3) \in P (z(u_3), M)$.}
\vskip7pt
{\bf Теорема 1} [s6], [s4]. {\it Пусть $M$ --- выпуклое $U$-множество пространства $X$, удовлетворяющего условиям $(C_1-C_5)$. Тогда верны следующие утверждения.
\vskip7pt
$(i)\, \,$ Если выполняется условие $(C_7)$ и
$$\stackrel{\circ}{M} \neq \emptyset,$$
то $\stackrel{\circ}{M}$ --- выпуклое множество и $\bar M = \bar {\stackrel{\circ}{M}}$.
\vskip7pt
$(ii)\, \,$ Если выполняется условие $(C_6)$ и $\bar M$ является $U$-множеством, то $\bar M$ --- выпуклое множество {\rm {(}}аналогичное утверждение для хордового пространства см. {\rm {\cite [с. 74]{Bus2})}}.
\vskip7pt
$(iii)\, \,$ Если $\bar M$ является $U$-множеством обобщенного хордового пространства $X$,
$$x_0 \in \stackrel{\circ}{M},\quad x \in \bar M,$$
то
$$(T(x_0,x)\backslash \{x\}) \subset \stackrel{\circ}{M},\quad \stackrel{\circ}{M} = \stackrel{\circ}{\bar {\stackrel{\circ}{M}}}.$$}

{\bf Доказательство теоремы 1}.
\vskip7pt
$(i)\, \,$ Пусть $x$, $y \in \stackrel{\circ}{M}$. Тогда найдется окрестность
$$U_y \subset \stackrel{\circ}{M}$$
точки $y$.
Возьмем произвольную точку $z$, удовлетворяющую условию $[xzy]$. Тогда $z \in \lambda_x(U_y)$, где 
$$\lambda = \dfrac{|xz|}{|xy|},$$
поскольку $\stackrel{\circ}{M} \subset W_x$. Из условия $(C_7)$ следует, что
$$\lambda_x (U_y) \subset \stackrel{\circ}{M}.$$ Поэтому $\stackrel{\circ}{M}$ --- выпуклое множество.
Пусть
$$x_0 \in \stackrel{\circ}{M},\quad x \in \bar M.$$
Возьмем произвольную окрестность $U_x$ точки $x$. Тогда существует точка $z \in U_x\cap M$. Кроме того, из условия $(C_7)$ следует, что
$$(T(x_0,z)\backslash \{z\}) \subset \stackrel{\circ}{M},$$
поскольку
$$x_0 \in \stackrel{\circ}{M}$$
и $z$ --- точка выпуклого множества $M \subset W_z$. Теперь очевидно, что
$$(T(x_0,z)\backslash \{z\})\cap U_x \neq \emptyset.$$
Поэтому
$$U_x\cap \stackrel{\circ}{M} \neq \emptyset.$$ 
Таким образом,
$$x \in \bar {\stackrel{\circ}{M}}, \quad \bar M \subset \bar {\stackrel{\circ}{M}}.$$
Обратное включение очевидно.
\vskip7pt
$(ii)\, \,$ Пусть $x$, $y \in \bar M$, $z \in T(x,y)$. Тогда эти точки можно представить в виде
$$x = \lim \limits_{n \rightarrow \infty} x_n,\quad y = \lim \limits_{n \rightarrow \infty} y_n,\quad
z = \omega_{\lambda} (x,y),$$ 
где
$$\lambda \in [0,1],\quad (x_n)_{n \in \mathbb{N}},\, \, (y_n)_{n \in \mathbb{N}}$$
--- последовательности из $M$.
Кроме того, для каждого натурального числа $n$
$$\{\omega_{\lambda} (x_n,y_n)\} \subset M,$$
поскольку $M$ --- выпуклое множество. Из условия $(C_6)$ следует теперь, что $z \in \bar M$. Таким образом, множество $\bar M$ выпуклое.
\vskip7pt
$(iii)\, \,$ Используем метод доказательства от противного. Пусть существует точка $y$ отличная от точки $x$ такая, что
$$y \in T(x_0,x)\backslash \stackrel{\circ}{M}.$$
Выберем $\lambda \in (0,1)$ такое, что
$$\lambda_y (x_0) \in B(y,r_y), \quad |y\lambda_y (x_0)| < |yx|.$$
Из условия непрерывности отображения $\lambda_y$ следует, что существует окрестность
$$U_{x_0} \subset \stackrel{\circ}{M}$$
точки $x_0$, обладающая следующим свойством.
$$\lambda_y (U_{x_0}) \subset B(y,r_y).$$
Из условий $(C_7)$, $(C_8)$ следует, что множество
$$\tau_y(\lambda_y (U_{x_0}))$$
открытое и
$$\tau_y(\lambda_y (U_{x_0}))\cap T(y,x) \neq \emptyset.$$
Но из утверждения $(i)$ следует, что
$$\bar M = \bar {\stackrel{\circ}{M}},$$ а из утверждения $(ii)$ следует, что
$$T(y,x) \in \bar M.$$
Поэтому существует точка
$$z \in \tau_y(\lambda_y (U_{x_0}))\cap \stackrel{\circ}{M}.$$
Пусть
$$z_1 = \lambda_y^{-1}(\tau_y^{-1}(z)).$$
Тогда
$$z_1 \in U_{x_0}, \quad y \in \mu_z(U_{x_0}),$$
где
$$\mu = \dfrac{|yz|}{|zz_1|}.$$
Из условия $(C_7)$ и утверждения $(i)$ следует теперь, что $y \in \stackrel{\circ}{M}$. Получили противоречие.
Пусть
$$x \in \stackrel{\circ}{\bar {\stackrel{\circ}{M}}}.$$ Тогда существуют окрестность $B(x,r)$, где $0 < r < r_x$, и точка $z$ такие, что
$$B(x,r) \subset \bar {\stackrel{\circ}{M}}, \quad z \in \stackrel{\circ}{M}\cap B(x,r).$$
Используя условие $(C_5)$, продолжим хорду $T(z,x)$ до большей хорды
$$T(z,y) \subset B(x,r).$$ Тогда из первой части доказательства следует, что
$$x \in \stackrel{\circ}{M}.$$
Следовательно,
$$\stackrel{\circ}{\bar {\stackrel{\circ}{M}}} \subset \stackrel{\circ}{M}.$$
Обратное включение очевидно. Таким образом, теорема 1 доказана.
\vskip7pt
В следующей теореме 2 содержатся утверждения, аналогичные утверждениям для хордового пространства (см. \cite [с. 74, 75]{Bus2}). Отметим также, что, к сожалению, в \cite [с. 74, (4)]{Bus2} имеются неточности в формулировке и в доказательстве.
\vskip7pt
{\bf Теорема 2} [s6]. {\it В обобщенном хордовом пространстве $X$, являющемся $U$-множеством, следующие утверждения эквивалентны.
\vskip7pt
$(i)\, \,$ Все замкнутые шары выпуклые.
\vskip7pt
$(ii)\, \,$ Для каждой точки $p \in X$ и каждой хорды $T \in S$ функция
$$f(t)= |px(t)|,$$
где $x= x(t)$ --- параметрическое представление хорды $T$, является безвершинной.
\vskip7pt
$(iii)\, \,$ Для каждой точки $p \in X$ и каждой хорды $T \in S$ множество $P(p, T)$ связное.
\vskip7pt
Кроме того, из этих утверждений следует, что пересечение произвольной сферы и произвольной хорды является либо пустым множеством, либо хордой, либо двумя точками.}
\vskip7pt
{\bf Доказательство теоремы 2}.
\vskip7pt
$(i) \Rightarrow (ii).\, \,$ Если функция
$$f : \mathbb{R} \rightarrow \mathbb{R},\quad f(t)= |px(t)|$$ не является безвершинной, то на хорде $T$ существуют такие точки $x$, $y$, $z$, что $[xyz]$ и либо  $f(y) > \max \{f(x), f(z)\}$, либо
$$f(y) = r = \max \{f(x), f(z)\},\quad r = f(x) > f(z).$$
Первый случай приводит к противоречию, поскольку замкнутый шар $B[p,r]$ выпуклый. Во втором случае из утверждения $(iii)$ теоремы 1 следует, что $y$ --- внутренняя точка шара $B[p,r]$. Получили противоречие, поскольку $f(y) = r$.
\vskip7pt
$(ii) \Rightarrow (iii). \, \,$ Это очевидно.
\vskip7pt
$(iii) \Rightarrow (i). \, \,$ Допустим, что замкнутый шар $B[p,r]$ невыпуклый. Тогда, используя условие
$(C_1)$ и непрерывность параметрического представления хорды, можно найти хорду $T(a,b)$, удовлетворяющую условиям
$$|ap| = |bp| = r,\quad |xp| > r$$
для каждой точки $x$ такой, что $[axb]$. Получили противоречие с условием $(iii)$.
\vskip7pt
Докажем последнее утверждение. Если $r < |pT|$, то
$$B[p,r]\cap T = \emptyset.$$
Если $r = |pT|$, то множество $B[p,r]\cap T$ является хордой, поскольку замкнутый шар $B[p,r]$ выпуклый. Если $r > |pT|$, то из утверждения $(iii)$ теоремы 1 следует, что множество $B[p,r]\cap T$ состоит не более чем из двух точек. Таким образом, теорема 2 доказана.
\vskip7pt
В следующей теореме 3 содержатся утверждения, аналогичные утверждениям для $G$-пространства (см. \cite [с. 157, (20.9)]{Bus}).
\vskip7pt
{\bf Теорема 3} [s6], [s4]. {\it Пусть $B(p,r_p)$ --- фиксированный открытый шар из условия $(C_4)$ обобщенного $G$-пространства Буземана $X$. Тогда следующие утверждения эквивалентны.
\vskip7pt
$(i)\, \,$ Все замкнутые шары, принадлежащие шару $B(p,r_p)$, выпуклые.
\vskip7pt
$(ii)\, \,$ Все замкнутые шары, принадлежащие шару $B(p,r_p)$, строго выпуклые.
\vskip7pt
$(iii)\, \,$ Каждая точка $x \in B(p,r_p)$ имеет единственное основание на каждом сегменте $T$, удовлетворяющем условию
$$|xT| < r_p - |px|.$$
Кроме того, из $(i)$, $(ii)$, $(iii)$ следует утверждение
\vskip7pt
$(iv)\, \,$ Произвольный сегмент и произвольная сфера, принадлежащая шару $B(p,r_p)$, пересекаются не более, чем в двух точках.
\vskip7pt
В свою очередь, из $(iv)$ следует
\vskip7pt
$(v)\, \,$ Каждая точка $x \in B(p,r_p)$ имеет единственное основание на каждом сегменте $T$, удовлетворяющем условию
$$\max \{|xy| : y \in T\} < r_p - |px|.$$}

{\bf Доказательство теоремы 3}.
\vskip7pt
$(i) \Rightarrow (ii). \, \,$ Допустим, что замкнутый шар $B[x,r_1] \subset B(p,r_p)$ не является строго выпуклым. Тогда существуют точки $y$, $z$, $q$, удовлетворяющие условиям
$$|xy|= |xz|= |xq| = r_1,\quad [yzq].$$
Используя условие $(C_5)$, продолжим сегмент $[z,x]$ до такого сегмента $[z,u]$, что
$$r_1 + |ux| < r_p - |up|.$$
Тогда
$$r_2 = \max \{|uy|,|uq|\} < |ux| + r_1 = |uz|$$
и точка $z$ не принадлежит замкнутому шару $B[u,r_2]$.
В силу условия $(i)$ получаем противоречие.
\vskip7pt
$(ii) \Rightarrow (iii), (iv). \, \,$ Это очевидно.
\vskip7pt
$(iii) \Rightarrow (i). \, \,$ Допустим, что найдется невыпуклый замкнутый шар
$$B[x,r_1] \subset B(p,r_p).$$
Тогда существуют точки $a$, $b$, удовлетворяющие условиям
$$|xa|= |xb| = r_1,\quad[a,b] \cap B[x,r_1] = \{a,b\}.$$
Но отсюда следует, что точка $x$ имеет два основания на сегменте $[a,b]$. Получили противоречие.
\vskip7pt
$(iv) \Rightarrow (v). \, \,$ Допустим, что точка $x \in B(p,r_p)$ имеет два основания на сегменте $[a,b]$. Предположение о том, что таких оснований больше двух сразу приводит к противоречию в силу условия $(iv)$. Пусть точка $z$ удовлетворяет условию
$$|xz| = \max \{|xq| : q \in [a,b]\}.$$
Непрерывно уменьшая $\lambda$ $(\lambda \leq 1)$, найдем первую точку $z_1$, отличную от точки $z$, которая принадлежит пересечению сферы
$$S(\lambda_z(x),\lambda |zx|)$$
с сегментом $[a,b]$. Если таких точек больше одной, то получаем противоречие в силу условия $(iv)$.
Пусть для определенности $z_1 \in [a,z]$. Очевидно, что найдется сфера $S(\lambda_z(x),r_1)$, радиус
которой удовлетворяет неравенствам
$$\max \{\min \{|\lambda_z(x)q| : q \in [z,b]\},\min \{|\lambda_z(x)q| : q \in [z_1,z]\}\} < r_1 < \lambda |xz|.$$
Эта сфера пересекает сегмент $[a,b]$, по меньшей мере, в трех точках. В силу условия $(iv)$ получили противоречие. Таким образом, теорема 3 доказана.
\vskip7pt
Следующая теорема 4 аналогична теореме 4 для хордового пространства из \cite [с. 65]{Bus2} и имеет такое же доказательство, если учесть теорему 1.
\vskip7pt
{\bf Теорема 4} [s6]. {\it Пусть для каждой точки $p$ обобщенного хордового пространства $X$ существует открытый шар $B(p,r_1(p))$ такой, что все открытые шары из этого шара строго выпуклые и $r_1(p) < r_p$. Тогда пространство $X$ является обобщенным $G$-пространством Буземана.}
\vskip7pt
Следующая теорема 5 аналогична предложениям (20.3), (20.4) из \cite [с. 154]{Bus} и предложению (1) из \cite [с. 79]{Bus2}.
\vskip7pt
{\bf Теорема 5} [s6], [s4]. {\it Пусть обобщенное хордовое пространство $X$ не имеет размерности $1$, $q$ --- граничная точка выпуклого множества $M \subset X$ такого, что
$$\stackrel{\circ}{M} \neq \emptyset, \quad \bar M \subset B(p,r_p).$$
Тогда
$$B(p,r_p) = V_q \cup A_q \cup C_q$$
и множество $V_q$ $(\tilde V_q)$ отделяет в шаре $B(p,r_p)$ множество $A_q$ от множества $C_q$.}
\vskip7pt
{\bf Доказательство теоремы 5}.
\vskip7pt
Из определения множеств $A_q$, $C_q$, $V_q$ и теоремы 1 следует, что
$$A_q \cap C_q = \emptyset, \quad B(p,r_p) = V_q \cup A_q \cup C_q.$$
Пусть
$$x = x(t),\quad t \in [0,1],$$ --- непрерывный путь такой, что
$x(0) \in A_q$, $x(1) \in C_q$. Обозначим через
$$t_0 = \sup \{t : t \in [0,1], \, \exists y_t \, (y_t \in \stackrel{\circ}{M}, \, [x(t)y_tq])\}.$$
Очевидно, что $t_0 > 0$. Допустим, что $x(t_0) \in A_q$. Из условия $(C_6)$ следует, что отображение $\lambda_q$ непрерывно в шаре $B(p,r_p)$ для каждого $\lambda \in (0,1)$. Следовательно, для малого
$t - t_0 > 0$ существует точка
$$y_t \in  \stackrel{\circ}{M},$$
удовлетворяющая условию $[x(t)y_tq]$.
В силу определения $t_0$, получили противоречие. Допустим, что $x(t_0) \in C_q$. Отображения
$$\tau_q : B(q,r_q)\rightarrow B(q,r_q),\quad \lambda_q : B(q,r_q)\rightarrow \lambda_q(B(q,r_q))$$ непрерывны для каждого $\lambda \in (0,1)$. Отсюда нетрудно получить, что для малого $t_0 - t > 0$
существует точка
$$\bar y_{t_1} \in  \stackrel{\circ}{M},$$
удовлетворяющая условию $[x(t_1)y_{t_1}q]$.
Тогда при
$$0 < t_0 - t_1 < t_0 - t$$
выполнено условие $[\bar y_{t_1}qy_{t_1}]$. Но по теореме 1 $\stackrel{\circ}{M}$ --- выпуклое множество. Значит, $q \in \stackrel{\circ}{M}$. Снова получили противоречие. Следовательно, $x(t_0) \in V_q$. Теперь нетрудно заметить, что $x(t_0) \in \tilde V_q$.
Таким образом, теорема 5 доказана.
\vskip7pt
Следующая теорема 6 аналогична предложениям (20.11), (20.12) из \cite [с. 160]{Bus} и предложениям (2), (3) из \cite [с. 80-82]{Bus2}. Она имеет такое же доказательство, если учесть теоремы 2, 5.
\vskip7pt
{\bf Теорема 6} [s6], [s4]. {\it Пусть обобщенное прямое хордовое пространство $X$ не имеет размерности $1$ и все открытые шары в нем выпуклые. Тогда верны следующие утверждения.
\vskip7pt
$(i)\, \,$ Для каждой точки $x \in X$ и для каждой прямой $\Lambda \subset X$ существует перпендикуляр к прямой $\Lambda$, содержащий точку $x$. Объединение всех перпендикуляров к прямой $\Lambda$, содержащих точку $x$, пересекает прямую $\Lambda$ по хорде. Причем, если основание точки $x$ на прямой $\Lambda$ единственно, то через точку $x$ проходит единственный перпендикуляр к прямой $\Lambda$.
\vskip7pt
$(ii)\, \,$ Трансверсали к данной прямой $\Lambda \subset X$ в данной точке $q \in \Lambda$ совпадают с опорными прямыми в граничной точке $q$ для шара с центром в произвольной точке $x \in \Lambda$, отличной от точки $q$, радиуса $|xq|$.
\vskip7pt
$(iii)\, \,$ Множество точек, расположенных на всех трансверсалях к данной прямой $\Lambda \subset X$ в данной точке $q \in \Lambda$, разбивает пространство $X$ на два линейно связных множества.}

\section{Одулярные структуры геометрии Гильберта и прямого $G$-пространства Буземана}

Рассмотрим геодезическое пространство $(X,\rho)$, через каждые две
различные точки которого можно провести единственную прямую, т. е. кривую, изометричную числовой прямой со стандартной метрикой. 
С помощью формулы (1.3.1) определено однопараметрическое семейство $(\omega_\lambda)_{\lambda \in \mathbb{R}}$ бинарных операций
$$\omega_\lambda : X\times X \rightarrow X \quad (\lambda \in \mathbb{R}).$$
Кроме того, для фиксированных $\lambda \in \mathbb{R}$, $p \in X$ определена унарная операция
$$\lambda_p : X \rightarrow X, \quad \lambda_p (x) = \omega_\lambda (p,x).$$
Нетрудно проверить, что при $\lambda \neq 1$ точка $\lambda_x (y)$ делит ориентированный сегмент $[x,y]$ с началом в точке $x$ и концом в точке $y$
в отношении $\lambda : (1-\lambda)$.
Отметим некоторые элементарные свойства введенных операций.
\vskip7pt
$1.\, \,$ Для всех $p$, $x$, $y \in X$, $\lambda$, $\mu \in \mathbb{R}$ имеют место равенства
$$\lambda_x (y) = (1-\lambda)_y (x), \quad |\lambda_p (x)\mu_p (x)| = |\lambda - \mu| |px|.$$

$2.\, \,$ Для всех $x$, $y \in X$, $\lambda \in \mathbb{R}\backslash \{-1\}$, $\mu   \in \mathbb{R}$ имеет место равенство
$$\mu_x (y)= (\mu + \lambda (\mu -1))_z (y),$$
где точка
$$z = \left(\dfrac{\lambda}{\lambda + 1}\right)_x (y)$$
делит ориентированный сегмент $[x,y]$ в отношении $\lambda$.
\vskip7pt
Фиксируем точку $p\in X$ и определим еще одну бинарную операцию
$$+_p : X\times X \rightarrow X, \quad x +_p\, y = 2_p (\left(\dfrac{1}{2}\right)_x (y)).$$
Нетрудно проверить, что алгебра
$$\langle X,\, +_p,\, p,\, \{\lambda_p\} _{\lambda \in \mathbb{R}}\rangle$$
с нейтральным элементом $p \in X$ является одулем (определения см. \cite [с. 801]{Sabin}). Обозначим этот одуль через $X_p$.

Семейство одулей $\{X_p\}_{p \in X}$ определяет на геодезическом пространстве $(X,\rho)$ одулярную структуру $W$. Из определения операций следует, что каждая изометрия пространства $X$ на себя является автоморфизмом одулярного пространства $(X,W)$. Найдем достаточные условия, при которых одулярная сруктура $W$ является топологической, т.е. для каждой точки $p \in X$ отображения
$$\omega^p : \mathbb{R}\times X\rightarrow X, \quad \omega^p (\lambda, x) = \lambda_p (x);$$
$$L^p : X\times X \rightarrow X, \quad L^p (x,y) = x +_p\,  y,$$
являются непрерывными.
\vskip7pt
{\bf Лемма 1} [s3]. {\it Одулярная структура $W$ геодезического пространства $X$, через каждые две
различные точки которого можно провести единственную прямую, является топологической одулярной структурой, если дополнительно выполнено,
по крайней мере, одно из следующих условий $(i)$, $(ii)$.
\vskip7pt
$(i)\, \,$ Для каждой точки $p \in X$ и для каждого $\lambda \in \mathbb{R}$ отображения
$$\omega_{1/2}: X\times X \rightarrow X, \quad \omega^p_{\lambda} : X\rightarrow X, \quad \omega^p_{\lambda} (x) = \lambda_p (x)$$
непрерывны.
\vskip7pt
$(ii)\, \,$ Для каждого $\lambda \in \mathbb{R}$ отображение
$$\omega_{\lambda} : X\times X \rightarrow X$$
непрерывно.}
\vskip7pt
{\bf Доказательство леммы 1}.
\vskip7pt
Доказательство сразу следует из неравенств
$$|\omega^p (\lambda, x) \omega^p (\mu, y)| \leq |\mu_p (x)\mu_p (y)| + |\lambda_p (x)\mu_p (x)| =$$
$$|\omega_{\mu} (p, x)\omega_{\mu} (p, y)| + |\lambda - \mu| |px|$$
и тождества
$$L^p (x,y) = \omega^p (2, \omega_{1/2} (x, y)).$$

{\bf Замечание 1.} {\it Нетрудно проверить, что геодезическое пространство $X$, через каждые две
различные точки которого можно провести единственную прямую, удовлетворяющее условию $(ii)$ леммы 1, является обобщенным прямым $G$-пространством Буземана.}
\vskip7pt
Следующая теорема 1 проясняет ситуацию в одном важном частном случае, когда геодезическое пространство
$X$ собственное.
\vskip7pt
{\bf Теорема 1} [s3]. {\it Одулярная структура $W$ собственного геодезического пространства $X$, через каждые две различные точки которого можно провести единственную прямую, является топологической одулярной структурой.}
\vskip7pt
{\bf Доказательство теоремы 1}.
\vskip7pt
Проверим для каждого $\lambda \in \mathbb{R}$ непрерывность отображения $\omega_{\lambda}$. Пусть
$$x = \lim \limits_{n \rightarrow \infty} x_n,\quad y = \lim \limits_{n \rightarrow \infty} y_n$$
Последовательность $(\omega_{\lambda} (x_n, y_n))_{n \in \mathbb{N}}$ ограничена, поскольку ограничена правая часть неравенства:
$$|x\omega_{\lambda} (x_n, y_n)| \leq |xx_n| + |x_n\omega_{\lambda} (x_n, y_n)| = |xx_n| + |\lambda| |x_ny_n|.$$
Пространство $X$ собственное, поэтому последовательность $(\omega_{\lambda} (x_n, y_n))_{n \in \mathbb{N}}$ имеет сходящуюся подпоследовательность, предел которой обозначим через $z$. Тогда точка $z$ удовлетворяет одному из следующих соотношений.
\vskip7pt
$1.\, \,$ $|xz| + |zy| = |xy|$, если $0 \leq \lambda \leq 1$.
\vskip7pt
$2.\, \,$ $|xy| + |yz| = |xz|$, если $\lambda > 1$.
\vskip7pt
$3.\, \,$ $|xz| + |xy| = |yz|$, если $\lambda < 0$.
\vskip7pt
Но точка $\omega_{\lambda} (x, y)$ удовлетворяет тому же соотношению, что и точка $z$, поэтому, в силу условий теоремы 1, эти точки совпадают. Теперь ясно также, что предел всякой сходящейся подпоследовательности последовательности  $(\omega_{\lambda} (x_n, y_n))_{n \in \mathbb{N}}$ совпадает с точкой $\omega_{\lambda} (x, y)$. Следовательно, для каждого $\lambda \in \mathbb{R}$ отображение $\omega_{\lambda}$ непрерывно, и теорема 1 верна в силу леммы 1.
\vskip7pt
{\bf Замечание 2.} {\it Нетрудно проверить, что геодезическое пространство, через каждые две различные точки которого можно провести единственную прямую, является собственным тогда и только тогда, когда оно является прямым $G$-пространством Буземана.}
\vskip7pt
В силу замечания 2 теорему 1 можно сформулировать следующим образом.
\vskip7pt
{\bf Теорема $1'$.} {\it Одулярная структура $W$ прямого $G$-пространства Буземана $X$ является топологической одулярной структурой.}
\vskip7pt
Рассмотрим открытый шар $X = B(0,r)$ с центром в нуле радиуса $r > 0$ в строго выпуклом банаховом пространстве с нормой $||.||$. Напомним (см. \cite [гл. 5, упр. 7]{Danf}), что банахово пространство называется {\it строго выпуклым}, если из равенства $||x + y|| = ||x|| + ||y||$ следует, что векторы $x$, $y$ банахова пространства линейно зависимы. В строго выпуклом банаховом пространстве сферы не содержат невырожденных отрезков. Введем обозначение
$$R(x, y, y_x, x_y)= \dfrac{||x - y_x|| ||y - x_y||}{||x - x_y|| ||y - y_x||},$$
где
$$x, \, y \in B(0,r), \quad x_y \quad (y_x)$$
--- точка пересечения луча с началом в точке $y\,$  $(x)$, содержащего точку $x\,$ $(y)$, со сферой
$S(0,r)$. Зададим в шаре $B(0,r)$ метрику геометрии Гильберта (см. \cite [с. 140]{Bus})
$$|xy| = \rho (x, y)  = \dfrac{k}{2} \ln R(x, y, y_x, x_y)$$
для всех $x, \, y \in B(0,r)$, где $k > 0$ --- константа. Выясним, как связана метрика геометрии Гильберта $\rho$ с метрикой банахова пространства $d$, где
$d (x, y) = ||x - y||$ для всех $x, \, y \in B(0,r)$.
\vskip7pt
{\bf Теорема 2} [s3]. {\it Метрики $\rho$ и $d$ топологически эквивалентны в открытом шаре $B(0,r)$
и липшицево эквивалентны на каждом замкнутом шаре $B[0,r_1]$ при $r_1 < r$. При этом верны
следующие неравенства
$$\dfrac{k (r - r_1) ||x - y||}{(r + r_1)^2} \leq |xy| \leq \dfrac{k r ||x - y||}{(r - r_1)^2}$$
для всех $x, \, y \in B[0,r_1]$ и равенство
$$\lim \limits_{r \rightarrow \infty} \dfrac{r}{k} |xy| = ||x - y||.$$}
\vskip7pt
{\bf Доказательство теоремы 2.}
\vskip7pt
Сначала получим оценки сверху и снизу для следующего выражения
$$R(x, y, y_x, x_y)= 1 + ||x - y||\left(\dfrac{1}{||x - x_y||} + \dfrac{1}{||y - y_x||} + \dfrac{||x - y||}{||x - x_y|| ||y - y_x||}\right).$$
С помощью неравенств
$$r - ||x|| \leq ||x - x_y||, \quad r - ||y|| \leq ||y - y_x||$$
для всех $x, \, y \in B(0,r)$ получим оценку сверху
$$R(x, y, y_x, x_y) \leq 1 + ||x - y||\left(\dfrac{1}{r - ||x||} + \dfrac{1}{r - ||y||} + \dfrac{||x - y||}{(r - ||x||) (r - ||y||)}\right).$$
А с помощью неравенств
$$||x - x_y|| \leq r + ||x||, \quad ||y - y_x|| \leq r + ||y||$$
для всех $x, \, y \in B(0,r)$ получим оценку снизу
$$R(x, y, y_x, x_y) \geq 1 + ||x - y||\left(\dfrac{1}{r + ||x||} + \dfrac{1}{r + ||y||} + \dfrac{||x - y||}{(r + ||x||) (r + ||y||)}\right).$$
Из свойств логарифмической функции и полученных оценок следует топологическая эквивалентность метрик
$\rho$ и $d$. Допустим теперь, что $x, \, y \in B[0,r_1]$. Тогда, используя полученную оценку сверху, получим оценку сверху для расстояния $|xy|$.
$$|xy|  = \dfrac{k}{2} \ln R(x, y, y_x, x_y) \leq$$ 
$$\dfrac{k}{2}||x - y||\left(\dfrac{1}{r - ||x||} + \dfrac{1}{r - ||y||} +
\dfrac{||x - y||}{(r - ||x||) (r - ||y||)}\right) \leq$$ 
$$\dfrac{k}{2}||x - y||\left(\dfrac{2}{r - r_1}  + \dfrac{2 r_1}{(r - r_1)^2}\right) =
 \dfrac{k r ||x - y||}{(r - r_1)^2}.$$
Используя полученную оценку снизу, получим оценку снизу для для расстояния $|xy|$.
$$|xy| \geq \dfrac{k}{2} \ln \left( 1 + \dfrac{2 ||x - y||}{r + r_1} \right).$$
Теперь используем следующее неравенство для  логарифмической функции
$$\ln (1 + t) \geq t - \dfrac{t^2}{2} = t\left(1 - \dfrac{t}{2}\right)$$
для $t \in \mathbb{R}_{+}$ и учтем, что $||x - y|| \leq 2 r_1$. 
$$|xy| \geq \dfrac{k ||x - y||}{r + r_1} \left( 1 - \dfrac{||x - y||}{r + r_1}\right) \geq
\dfrac{k (r - r_1) ||x - y||}{(r + r_1)^2}.$$
Из полученных оценок следует липшицева эквивалентность метрик
$\rho$ и $d$ на каждом замкнутом шаре $B[0,r_1]$ при $r_1 < r$.
Доказательство последнего равенства теперь нетрудно получить следующим образом. Используя полученные оценки для выражения $R(x, y, y_x, x_y)$, оценить сверху и снизу выражение  $\dfrac{r}{k} |xy|$  и вычислить пределы в правой и левой частях полученных неравенств по правилу Лопиталя.
Таким образом, теорема 2 доказана.
\vskip7pt
{\bf Лемма 2} [s3]. {\it 
\vskip3pt
$(i) \, \,$ Явная формула для операции 
$$\lambda_p : B(0,r) \rightarrow B(0,r)$$
для $\lambda \in \mathbb{R}$, $p \in B(0, r)$ имеет вид
$$\lambda_p (x) = p + \left(\dfrac{||p - x_p||^{\lambda} ||x - p_x||^{\lambda} - ||p - p_x||^{\lambda}||x - x_p||^{\lambda}}{||p - x_p||^{\lambda -1} ||x - p_x||^{\lambda} + ||p - p_x||^{\lambda -1}||x - x_p||^{\lambda}}\right)\dfrac{(x - p)}{||x - p||}.$$
$(ii) \, \,$ Отображение 
$$\lambda_p : B(0,r) \rightarrow B(0,r)$$
при $\lambda \in (0,1)$, $p \in B(0, r)$ уменьшает расстояния, то есть
для любых $x, \, y \in B(0,r)$
$$|\lambda_p (x) \lambda_p (y)| < |xy|.$$}

{\bf Доказательство леммы 2.}
\vskip7pt
$(i) \, \,$ Из определения $\lambda_p$ следует равенство
$$R(p, \lambda_p (x), x_p, p_x)= R^{\lambda}(p, x, x_p, p_x).$$
Тогда
$$||\lambda_p (x) - p|| = sign (\lambda) \left(\dfrac{||p - x_p||^{\lambda} ||x - p_x||^{\lambda} - ||p - p_x||^{\lambda}||x - x_p||^{\lambda}}{||p - x_p||^{\lambda -1} ||x - p_x||^{\lambda} + ||p - p_x||^{\lambda -1}||x - x_p||^{\lambda}}\right).$$
Отсюда нетрудно получить искомую формулу.
\vskip7pt
$(ii) \, \,$ Пусть 
$$z = \lambda_p (x), \quad q = \lambda_p (x),\quad u = [p, x_y] \cap [z_q, q_z],\quad v = [p, y_x] \cap [z_q, q_z].$$
Тогда
$$R(z, q, q_z, z_q) = \dfrac{||z - q_z|| ||q - z_q||}{||z - z_q|| ||q - q_z||} =$$
$$\dfrac{(||z - q|| + ||q -  q_z|| ) (||q - z|| + ||z - z_q||)}{||q - q_z|| ||z - z_q||}  <$$
$$\dfrac{(||z - q|| + ||q -  v|| ) (||q - z|| + ||z - u||)}{||q - v|| ||z - u||} =
R(z, q, v, u) = R(x, y, y_x, x_y).$$
Осталось применить натуральный логарифм, умноженный на $k/2$, к левой и правой частям полученного неравенства. Лемма 2 доказана.
\vskip7pt
Учитывая формулы
$$\rho (0,x) = \dfrac{k}{2} \ln \dfrac{r + ||x||}{r - ||x||} = k {\rm Arth} \left( \dfrac{||x||}{r}\right), \quad 
||x|| = r {\rm th} \left(\dfrac{\rho (0,x)}{k}\right),$$ 
из леммы 2 нетрудно получить
\vskip7pt
{\bf Следствие 1} [s3]. {\it Середина сегмента $[x,y] \subset B(0,r)$ может быть
найдена по формуле
$$\left(\dfrac{1}{2}\right)_x (y) = \dfrac{x + a y}{1 + a}$$
где 
$$a = \sqrt{\dfrac{||x - y_x|| ||x - x_y||}{||y - y_x|| ||y - x_y||}}.$$
Кроме того, при $p = 0$ верны равенства
$$||\lambda_0 (x)|| = r {\rm sign} (\lambda) \left(\dfrac{(r + ||x||)^{\lambda} - (r - ||x||)^{\lambda}}{(r + ||x||)^{\lambda} + (r - ||x||)^{\lambda}}\right) = r  {\rm th} \left( \dfrac{|\lambda| \rho (0,x)}{k}\right),$$
$$\lambda_0 (x) = x {\rm th} \left( \dfrac{|\lambda| \rho (0,x)}{k}\right) {\rm cth} \left( \dfrac{\rho (0,x)}{k}\right).$$}
\vskip7pt
Для исследования одулярной структуры геометрии Гильберта нам понадобится еще одна лемма 3.
\vskip7pt
{\bf Лемма 3} [s3]. {\it 
Отображение 
$$\varphi : (B(0,r)\times B(0,r)) \backslash D \rightarrow S(0,r), \quad \varphi (x, y) = y_x,$$
где $D$ --- диагональ в декартовом произведении $B(0,r)\times B(0,r)$, непрерывно относительно топологий, индуцированных топологией произведения банаховых пространств и топологией банахова пространства.}
\vskip7pt
{\bf Доказательство леммы 3.}
\vskip7pt
Пусть
$$x = \lim \limits_{n \rightarrow \infty} x_n,\quad y = \lim \limits_{n \rightarrow \infty} y_n,\quad \varphi (x_n, y_n) = x_n + t_n (y_n - x_n),$$
где 
$$n \in \mathbb{N}, \quad (x_n; y_n), \, (x; y) \in (B(0,r)\times B(0,r)) \backslash D, \quad t_n \in \mathbb{R}, \quad t_n > 1.$$
Последовательность $(t_n)_{n \in \mathbb{N}}$ ограничена. Это следует из того, что последовательности $(x_n)_{n \in \mathbb{N}}$, $(y_n)_{n \in \mathbb{N}}$
сходятся к разным пределам, и из неравенства
$$|t_n| = \dfrac{||\varphi (x_n, y_n) - x_n||}{||y_n - x_n||} \leq \dfrac{r + ||x_n||}{||y_n - x_n||}.$$
Следовательно, последовательность $(t_n)_{n \in \mathbb{N}}$ обладает сходящейся подпоследовательностью, предел которой обозначим через $t \in \mathbb{R}$. Предположим, что существует еще одна подпоследовательность, сходящаяся к другому пределу $t_1$. Тогда  из непрерывности операций в банаховом пространстве получим, что луч 
$\{x + \lambda (y - x) : \lambda \geq 0\}$ с началом во внутренней точке $x \in B(0,r)$ пересечет сферу
$S(0,r)$ в двух точках при $\lambda = t$ и $\lambda = t_1$. Получили противоречие. Следовательно,
$$\lim \limits_{n \rightarrow \infty} \varphi (x_n, y_n) = x + t (y - x) = \varphi (x,y).$$
Лемма 3 доказана.
\vskip7pt
{\bf Теорема 3} [s3]. {\it Одулярная структура $W$ геометрии Гильберта является топологической.}
\vskip7pt
{\bf Доказательство теоремы 3.}
\vskip7pt
Для каждой точки $p \in X$ и для каждого $\lambda \in \mathbb{R}$ отображение
$$\omega^{p}_{\lambda} : X\rightarrow X$$
непрерывно в точке $p$, поскольку 
$$|\omega^{p}_{\lambda} (p)\omega^{p}_{\lambda} (x)| = |\lambda| |px|$$
для каждого $x \in X$. Непрерывность отображения $\omega^{p}_{\lambda}$ на множестве $X\backslash \{p\}$ следует из теоремы 2 и лемм 2, 3. Непрерывность отображения $\omega_{1/2}$ в точках диагонали
$D$ множества $B(0,r)\times B(0,r)$ следует из следствия 1, теоремы 2 и неравенства
$$||\omega_{1/2} (x, y) - \omega_{1/2} (z, z)|| = ||x - z + \dfrac{a}{1 + a} (y - x)|| \leq ||x - z|| + ||y - z||,$$
где  
$$(x; y) \in B(0,r)\times B(0,r), \quad (z;z) \in D.$$
Непрерывность отображения $\omega_{1/2}$ на множестве $(B(0,r)\times B(0,r)) \backslash D$ следует из следствия 1, теоремы 2 и леммы 3. Теорема 3 следует теперь из леммы 1.
\vskip7pt
Пусть открытый шар $B(0,r)$ расположен в гильбертовом пространстве. В этом случае получаем известную модель Бельтрами-Клейна геометрии Лобачевского в шаре $B(0,r)$. Метрику $\rho$ нетрудно представить в следующей форме
$$|xy|  = \dfrac{k}{2} \ln \dfrac{A + \sqrt{A^2 - B}}{A - \sqrt{A^2 - B}},$$
где
$$A = r^2 - (x, y),\quad B = (r^2 - (x, x)) (r^2 - (y, y)), \quad (x, y)$$
--- скалярное произведение векторов $x, \, y \in B(0,r)$, или в эквивалентной форме
$$|xy|  = k {\rm Arch} \left( \dfrac{A}{\sqrt{B}}\right).$$
Заметим, что
$$B = r^4 {\rm ch}^{-2} \left(\dfrac{\rho (0,x)}{k}\right) {\rm ch}^{-2} \left(\dfrac{\rho (0,y)}{k}\right),$$
а формулы из следствия 1 преобразуются к виду
$$a = \sqrt{\dfrac{r^2 - (x, x)}{r^2 - (y, y)}} = \dfrac{{\rm ch} \left(\dfrac{\rho (0,y)}{k}\right)}{{\rm ch} \left(\dfrac{\rho (0,x)}{k}\right)},$$
$$\left(\dfrac{1}{2}\right)_x (y) = \dfrac{x {\rm ch} \left(\dfrac{\rho (0,x)}{k}\right) + y {\rm ch} \left( \dfrac{\rho (0,y)}{k}\right)}{{\rm ch} \left(\dfrac{\rho (0,x)}{k}\right) + {\rm ch} \left(\dfrac{\rho (0,y)}{k}\right)}.$$
В следующей теореме вычислены некоторые замечательные пределы для одуля $B_0(0,r)$ модели Бельтрами-Клейна геометрии Лобачевского.
\vskip7pt
{\bf Теорема 4} [s3]. {\it Для одуля $B_0(0,r)$ модели Бельтрами-Клейна геометрии Лобачевского имеют место следующие формулы.
$$\lim \limits_{\lambda \rightarrow \infty} \dfrac{|\lambda_0 (x)\lambda_0 (y)|}{\lambda} = \rho (0,x) + \rho (0,y) \leqno (i)$$
для всех $x,\, y \in B(0,r).$
$$\lim \limits_{\lambda \rightarrow 0+} \dfrac{|\lambda_0 (x)\lambda_0 (y)|}{\lambda} = ||y \dfrac{\rho (0,y)}{||y||} - x \dfrac{\rho (0,x)}{||x||}|| = \leqno (ii)$$
$$\sqrt{\rho^2 (0,x) + \rho^2 (0,y) - 2 \rho (0,x) \rho (0,y) \cos \alpha},$$
где $x,\, y \in B(0,r)\backslash \{0\}$, $\alpha$ --- величина угла между векторами $x,\, y$.
$$\lim \limits_{\lambda \rightarrow 0+} \dfrac{\rho (0,\omega_{1/2} (\lambda_0 (x),\lambda_0 (y)))}{\lambda} = \dfrac{1}{2}||y \dfrac{\rho (0,y)}{||y||} + x \dfrac{\rho (0,x)}{||x||}|| = \leqno (iii)$$ 
$$\dfrac{1}{2} \sqrt{\rho^2 (0,x) + \rho^2 (0,y) + 2 \rho (0,x) \rho (0,y) \cos \alpha},$$
где $x,\, y \in B(0,r)\backslash \{0\}$, $\alpha$ --- величина угла между векторами $x,\, y$.
$$\lim \limits_{\lambda \rightarrow 0+} \dfrac{|\lambda_0 (x +_0\,  y)(\lambda_0 (x) +_0\,  \lambda_0 (y))|}{\lambda} = \leqno (iv) $$
$$||x \dfrac{\rho (0,x)}{||x||} + y \dfrac{\rho (0,y)}{||y||} - 2 z \dfrac{\rho (0,z)}{||z||}||,$$
где 
$$x,\, y \in B(0,r)\backslash \{0\}, \quad x \neq - y, \quad z = \left(\dfrac{1}{2}\right)_x (y).$$
}
\vskip7pt
{\bf Доказательство теоремы 4.}
\vskip7pt
Доказательство состоит в применении правила Лопиталя и элементарных преобразований с использованием формул из следствия 1 и их преобразованных форм для модели Бельтрами-Клейна геометрии Лобачевского.
\vskip7pt

%% file: ch22.tex
\chapter{Аппроксимативные свойства множеств в геодезическом пространстве}
\vskip20pt

В параграфе 2.1 рассматриваются относительные чебышевский центр,  чебышевский радиус и множество всех диаметральных точек ограниченного множества метрического пространства. Получены оценки изменения относительного чебышевского радиуса $R_W (M)$ при изменении непустых ограниченных множеств $M, \, W$ метрического пространства. Доказано, что из всякой последовательности компактных множеств метрического пространства, сходящейся относительно метрики Хаусдорфа к некоторому компактному множеству $M$, можно выбрать подпоследовательность, для которой последовательность множеств всех относительных чебышевских центров (всех диаметральных точек)  ее элементов сходится относительно отклонения Хаусдорфа к множеству всех относительных чебышевских центров (всех диаметральных точек) множества $M$.

В параграфе 2.2 найдены замыкание и внутренность множества всех $N$-сетей, каждая из которых обладает принадлежащим ей единственным относительным чебышевским центром, в множестве всех $N$-сетей специального геодезического пространства относительно метрики Хаусдорфа.
При $N > 2$ найдена граница множества всех $N$-сетей, каждая из которых обладает не более, чем $N-2$ принадлежащими ей относительными чебышевскими центрами, в множестве всех $N$-сетей специального геодезического пространства относительно метрики Хаусдорфа.
Исследованы геометрические свойства относительных чебышевских центров конечного множества специального геодезического пространства, принадлежащих этому множеству.

В параграфе 2.3 получены достаточные условия существования и единственности чебышевского центра непустого ограниченного множества специального геодезического пространства.

В параграфе 2.4 теоремы Б. Секефальви - Надь \cite[теорема 3.35]{Brudn},
С. Б. Стечкина и Н. В. Ефимова \cite [теоремы 1.1 и 1.2]{Vlas} об
аппроксимативных свойствах множеств в равномерно выпуклом банаховом
пространстве обобщены на случай специального геодезического пространства.

В параграфе 2.5 теоремы Л. П. Власова \cite{Vlas1, Vlas} и А. В. Маринова \cite{Marin1} о непрерывности и связности метрической $\delta $-проекции в равномерно выпуклом банаховом пространстве обобщены на случай специального геодезического пространства. Одним из простых следствий такого обобщения является справедливость аналогичных результатов в пространствах Лобачевского (включая бесконечномерные).

В параграфе 2.6 доказано, что теоремы А. В. Маринова из \cite{Marin, Marin2} о непрерывности метрической
 $\delta $-проекции на выпуклое множество в линейном нормированном пространстве остаются верными в специальном метрическом пространстве.

В параграфе 2.7 в специальном метрическом пространстве получены обобщения некоторых результатов \mbox{П. К. Белоброва \cite{Belob, Belob1}} и  А. Л. Гаркави \cite{Gark1} о наилучших $N$-сетях непустых
ограниченных замкнутых выпуклых множеств в гильбертовом и в специальном банаховом пространствах, а также о принадлежности чебышевского центра замыканию выпуклой оболочки данного множества.

В параграфе 2.8 доказано, что некоторые результаты А. Л. Гаркави \cite {Gark, Gark1} и П. К. Белоброва \cite {Belob1, Belob} о наилучшей сети, наилучшем сечении и чебышевском
центре ограниченного множества в специальном банаховом пространстве верны и в бесконечномерном пространстве Лобачевского. А именно, для каждого непустого ограниченного множества бесконечномерного пространства Лобачевского доказано существование наилучшей $N$-сети и наилучшего $N$-сечения, а также установлена сильная устойчивость чебышевского центра.

В параграфе 2.9 рассматривается наилучшее приближение выпуклого компакта геодезического пространства шаром. Получена оценка сверху для расстояния Хаусдорфа от непустого ограниченного множества до множества всех замкнутых шаров специального геодезического пространства $X$ неположительной кривизны по Буземану. Доказано, что множество всех центров $\chi (M)$
замкнутых шаров, наилучшим образом приближающих в метрике
Хаусдорфа выпуклый компакт $M \subset X$, непустое и принадлежит
$M$. Исследованы геометрические свойства множества $\chi (M)$.
 Таким образом, теоремы С. И. Дудова и И. В. Златорунской \cite {Dudov, Dudov1} обобщены на случай специального геодезического пространства неположительной кривизны по Буземану.

В параграфе 2.10 исследуются метрические свойства
касательного пространства для метрического пространства более
общего, чем дифференцируемое $G$-пространство Буземана. Установлено,
что метрика на касательном пространстве в произвольной точке пространства
неположительной кривизны по Буземану (дифференцируемого по
\mbox{Буземану}
метрического пространства) внутренняя. Доказано, что касательное
пространство в произвольной точке локально полного дифференцируемого
по Буземану метрического пространства является полным пространством, а
также, что касательное пространство в произвольной точке локально
компактного пространства неположительной  кривизны по Буземану
является собственным геодезическим пространством.

Основные результаты по перечисленным темам этой главы опубликованы в статьях автора [s7, s8, s9, s12, s13, s14, s17, s18, s19, s21].

\section{Относительные чебышевский центр и чебышевский радиус ограниченного множества метрического пространства}

\vskip20pt
Напомним несколько определений.

Пусть в метрическом пространстве $(X,\rho)$ фиксировано некоторое семейство множеств $\Sigma \subset \Sigma (X)$. Множество $S^{*} \in \Sigma$ называется {\it наилучшим аппроксимирующим множеством} из семейства $\Sigma$ для множества $M \in B(X)$ \cite [c. 15]{Tikh}, если
$$\beta (M,S^{*}) = R_{\Sigma}(M),\, \, \mbox {где} \, \, R_{\Sigma }(M) = \inf\{\beta (M,S) :\, S \in \Sigma \}.$$

{\it Oтносительным чебышевским радиусом множества} $M \in B(X)$ по отношению к множеству $W \in \Sigma (X)$ называется число \cite {Wisnicki}
$$R_W (M) = \inf\{\beta (M,x) : x \in W \},$$
а множество
$$Z_W (M) = \{x \in W : \beta (M,x) = R_W (M)\}$$
называется {\it множеством всех относительных чебышевских центров множества} $M$ по отношению к множеству $W \in \Sigma (X)$ \cite {Amir}.

$R (M) = R_X (M)$ ($Z (M) = Z_X (M)$) --- {\it чебышевский радиус} ({\it множество всех чебышевских центров}) множества $M\in B(X)$.
$$H (M) = \{x \in M : \beta(M,x) = D(M)\}$$
--- {\it множество всех диаметральных точек множества} $M \in B(X)$ \cite {Khamsi}.
\vskip7pt
Кроме того, используем обозначения
$$R_0 (M) = R_M (M),\quad Z_0 (M) = Z_M (M),$$
$$Q_0 (M) = \{x \in M : \beta(Z_0 (M),x) = R_0 (M)\}.$$

$K (X)$ --- множество всех компактных подмножеств пространства $X$.
\vskip7pt
Различные свойства относительных чебышевских центров и относительного чебышевского радиуса, связанные, в основном, с задачами теории приближений исследовались многими авторами для множеств банахова пространства (см., например, статьи \cite {Wisnicki}, \cite {Amir} и литературу, указанную в них). В метрическом пространстве мы получим оценки изменения относительного чебышевского радиуса $R_W (M)$ при изменении множеств $M, \, W \in (B(X),\alpha)$. Кроме того, выясним поведение последовательности множеств относительных чебышевских центров $(Z_0 (M_n))_{n \in \mathbb{N}}$ (множеств диаметральных точек $(H(M_n))_{n \in \mathbb{N}}$, множеств $(Q_0(M_n))_{n \in \mathbb{N}}$), когда последовательность компактов $(M_n)_{n \in \mathbb{N}}$ сходится к компакту $M$ при $n \rightarrow \infty$ относительно метрики Хаусдорфа.

Сначала перечислим некоторые простые свойства и взаимосвязи рассмотренных понятий в метрическом пространстве $X$, доказательства первых четырех из которых непосредственно следуют из их определений.
\vskip7pt
$1.\, \,$ Если $M \in B(X)$, то условие $i)$ $H(M) = M$ равносильно условию
$ii)$ $Z_0 (M) = M$. Из этих условий следует, что $Q_0 (M) = M$ и $iii)$ $R_0 (M) = D(M)$. Если $M \in K(X)$,
то условия $i)$, $ii)$, $iii)$ равносильны.
\vskip7pt
$2.\, \,$ Пусть $M \in B(X)$. Если $M\cap Z(M) \neq \emptyset$, то
$$Z_0 (M) = M\cap Z(M),\quad R_0 (M) = R(M).$$

$3.\, \,$ Пусть $M$, $W$, $T \in B(X)$. Тогда $R_W (M) \leq R_T (M) + R_W (T)$ и
$$\theta : B(X)\times B(X) \rightarrow \mathbb{R}_{+},\quad \theta (M,W) = \max \{R_W (M), R_M (W)\}$$
--- функция, удовлетворяющая неравенству треугольника, а также неравенствам
$$\alpha (M,W) \leq \theta (M,W) \leq D(M,W).$$

$4.\, \,$ Для любых $M, \, W \in B[X]$ множество $Z_W (M)$ ($H(M)$, $Q_0(M)$) замкнуто в множестве
$W$ ($M$) с индуцированной метрикой.
\vskip7pt
$5.\, \,$ Если $M, \, W \in B(X)$, то
$$|D(M) - D(W)| \leq 2 \alpha (M,W).$$

Приведем доказательство этого утверждения (поскольку найти опубликованное доказательство непросто).
Для любого $\varepsilon > 0$ выберем такие элементы $x, \, y \in M$, $u, \, v \in W$, что
$$D (M) \leq |xy| + \varepsilon,\quad |xu| \leq |xW| + \varepsilon,\quad |yv| \leq |yW| + \varepsilon.$$
Тогда
$$D(W) \geq |uv| \geq |xy| - |xu| - |yv| \geq |xy| - |xW| - |yW| - 2 \varepsilon \geq$$
$$D(M) - 2 \alpha (M,W) - 3 \varepsilon.$$
Используя произвольность выбора $\varepsilon > 0$, получим
$$D(M) - D(W) \leq 2 \alpha (M,W).$$
Осталось использовать произвольность выбора $M, \, W \in B(X)$.
\vskip7pt
Найдем оценку изменения относительного чебышевского радиуса $R_W(M)$ при изменении множеств $M, \, W \in (B(X),\alpha)$.
\vskip7pt
{\bf Теорема 1} [s19], [s14]. {\it Пусть $X$ --- метрическое пространство, $M, \, W, \, A, \, B \in B(X)$. Тогда верны следующие утверждения.
$$|R_W (M) - R_B (A)|  \leq \max \{\beta (M,A) + \beta (B,W), \beta(A,M) + \beta (W,B)\}. \leqno (i)$$
$$R_W (M) \leq R(M) + |Z (M)W| \leqno (ii)$$
при $Z (M) \neq \emptyset$. В частном случае, когда $X$ --- банахово пространство, это неравенство известно {\rm \cite {Wisnicki}}.
$$R_W (M) \leq R_0 (M) + |Z_0 (M)W| \leqno (iii)$$
при $Z_0 (M) \neq \emptyset$.
$$|R_{\Sigma }(M) - R_{\Sigma }(W)| \leq \alpha (M,W). \leqno (iv)$$}
Из теоремы 1, определения метрики Хаусдорфа и функции $\theta$ получим
\vskip7pt
{\bf Следствие 1} [s19]. {\it Пусть $X$ --- метрическое пространство, 
\\
$M, \, W, \, A, \, B \in B(X)$. Тогда имеют место неравенства
$$|R_W (M) - R_B (A)|  \leq \alpha (M,A) + \alpha (W,B);$$
$$|R_0 (M) - R_0 (A)| \leq \beta (M,A) + \beta (A,M) \leq 2\alpha (M,A);$$
$$\theta (M,W) \leq \max \{R(M) + |Z(M)W|, R(W) + |Z(W)M|\},$$
при $Z (M) \neq \emptyset$, $Z (W) \neq \emptyset;$
$$\theta (M,W) \leq \max \{R_0 (M) + |Z_0 (M)W|, R_0 (W) + |Z_0 (W)M|\}$$
при $Z_0 (M) \neq \emptyset$, $Z_0 (W) \neq \emptyset$.}

{\bf Доказательство теоремы 1}.
\vskip7pt
$(i)\, \,$ Для любых $x \in W$, $y \in B$ справедливо неравенство
$$\beta (M,x) \leq \beta (M,y) + |yx|.$$
Тогда
$$R_W (M) \leq \beta (M,y) + |yW| \leq \beta (M,y) + \beta (B,W).$$
Следовательно,
$$R_W (M) \leq R_B (M) + \beta (B,W).$$
Кроме того, для любого $x \in B$ имеет место неравенство
$$\beta (M,x) \leq \beta (M,A) + \beta (A,x).$$
Тогда
$$R_B (M) \leq \beta (M,A) + R_B (A).$$
Учитывая эти неравенства, получим
$$R_W (M) - R_B (A) \leq R_W (M) - R_B (M) + R_B (M) - R_B (A) \leq$$
$$\beta (B,W) + \beta (M,A).$$
Из полученного неравенства следует требуемое неравенство.
\vskip7pt
$(ii)\, \,$ Для любых $y \in Z(M)$, $x \in W$ верно неравенство
$$\beta (M,x) \leq \beta (M,y) + |yx|.$$
Тогда из этого неравенства и определений $R(M)$, $Z(M)$ следует, что
$$R_W (M) \leq \beta (M,y) + |yW| = R(M) + |Z(M)W|.$$
$(iii)\, \,$ Для любых $y \in Z_0 (M)$, $x \in W$ верно неравенство
$$\beta (M,x) \leq \beta (M,y) + |yx|.$$
Тогда из этого неравенства и определений $R_0 (M)$, $Z_0 (M)$ следует, что
$$R_W (M) \leq \beta (M,y) + |yW| = R_0 (M) + |Z_0 (M)W|.$$
$(iv)\, \,$ Пусть $x \in M$ и  $\tau   > 0$. Тогда найдется точка $z \in W$ такая, что
$$|xz| \leq |xW| + \tau.$$
 Кроме того, получим неравенства
$$|xS| \leq |xz| + |zS| \leq |xW| + |zS| + \tau
 \leq$$
$$\sup\limits_{x \in M}{|xW|} + \sup\limits_{z \in W}{|zS|} + \tau   \leq
\alpha (M,W) + \sup\limits_{z \in W}{|zS|} + \tau $$
для $S \in \Sigma (X)$.
Отсюда следует, что
$$\inf\limits_{S \in \Sigma (X)}{\sup\limits_{x \in M}{|xS|}} \leq
\alpha (M,W) + \sup\limits_{z \in W}{|zS|} + \tau.$$
Переходя к точной нижней грани в правой части и используя
произвольность выбора $\tau  > 0$, получим
$$R_{\Sigma }(M) = \inf\limits_{S \in \Sigma (X)}{\sup\limits_{x \in M}{|xS| }} \leq
\alpha (M,W) + \inf\limits_{S \in \Sigma (X)}{\sup\limits_{z \in W}{|zS|}} = \alpha (M,W) +
R_{\Sigma }(W).$$
Неравенство
$$R_{\Sigma }(W) \leq \alpha (M,W) + R_{\Sigma }(M)$$
получается аналогично.
Таким образом, теорема 1 доказана.
\vskip7pt
{\bf Теорема 2} [s19]. {\it Пусть $X$ --- метрическое пространство и для каждого $n \in \mathbb{N}$
$$M_n \in K (X) \, \, \mbox{и} \, \, \lim \limits_{n \rightarrow \infty} \alpha (M_n,M) = 0.$$
Тогда найдется такая подпоследовательность
$(M_k)_{k \in \mathbb{N}}$ последовательности $(M_n)_{n \in \mathbb{N}}$, что верны равенства.
$$\lim \limits_{k \rightarrow \infty} \beta (Z_0 (M_k), Z_0 (M)) = 0. \leqno (i)$$
$$\lim \limits_{k \rightarrow \infty} \beta (H(M_k),H(M)) = 0. \leqno (ii)$$
$$\lim \limits_{k \rightarrow \infty} \beta (Q_0 (M_k),Q_0 (M)) = 0. \leqno (iii)$$}

{\bf Доказательство теоремы 2}.
\vskip7pt
Из свойства $5$ и компактности множеств $M, \, M_n \, \, (n \in \mathbb{N})$ следует, что для каждого $n \in \mathbb{N}$ найдется такой элемент
$$z_n \in Z_0 (M_n)\quad (z_n \in H(M_n),\, \, z_n \in Q_0 (M_n)),$$
что
$$|z_nZ_0 (M)| = \beta (Z_0 (M_n),Z_0 (M))\quad (|z_nH(M)| = \beta (H(M_n),H(M)),$$ 
$$|z_nQ_0 (M)| = \beta (Q_0 (M_n),Q_0 (M))).$$
$(i)\, \,$ Тогда для каждого $n \in \mathbb{N}$
$$\alpha(M_n,z_n) = \beta (M_n,z_n) = R_0 (M_n).$$
Используя условия теоремы 2, найдем такие подпоследовательность
$(z_k)_{k \in \mathbb{N}}$ последовательности $(z_n)_{n \in \mathbb{N}}$  и элемент $z \in M$, что
$$\lim \limits_{k \rightarrow \infty}|z_kz| = 0.$$
Тогда в силу следствия 1 и непрерывности метрики $\alpha$ получим
$$\beta (M,z) = \alpha (M,z) = R_0 (M).$$
Следовательно, $z \in Z_0 (M)$ и
$$\lim \limits_{k \rightarrow \infty} \beta (Z_0 (M_k),Z_0 (M)) = |z_kZ_0 (M)| = 0.$$

$(ii)\, \,$ Тогда для каждого $n \in \mathbb{N}$
$$\alpha(M_n,z_n) = \beta (M_n,z_n) = D(M_n).$$
Используя условия теоремы 2, найдем такие подпоследовательность
$(z_k)_{k \in \mathbb{N}}$ последовательности $(z_n)_{n \in \mathbb{N}}$ и элемент $z \in M$, что
$$\lim \limits_{k \rightarrow \infty} |z_kz| = 0.$$
Тогда в силу свойства 5 и непрерывности метрики $\alpha$ получим
$$\beta (M,z) = \alpha (M,z) = D(M).$$
Следовательно, $z \in H(M)$ и
$$\lim \limits_{k \rightarrow \infty} \beta (H(M_k),H(M)) = \lim \limits_{k \rightarrow \infty} |z_kH(M)| = 0.$$

$(iii)\, \,$ Тогда для каждого $n \in \mathbb{N}$
$$\beta (Z_0(M_n),z_n) = R_0 (M_n).$$
Используя условия теоремы 2 и доказанное утверждение $(i)$, найдем такие подпоследовательности
$(z_k)_{k \in \mathbb{N}}$, $(Z_0(M_k))_{k \in \mathbb{N}}$ последовательностей $(z_n)_{n \in \mathbb{N}}$, $(Z_0(M_n)))_{n \in \mathbb{N}}$ соответственно
и элемент $z \in M$, что
$$\lim \limits_{k \rightarrow \infty} |z_kz| = 0,\quad \lim \limits_{k \rightarrow \infty} \beta (Z_0 (M_k),Z_0 (M))= 0.$$
Но для каждого $k \in \mathbb{N}$ верно неравенство
$$|\beta (Z_0 (M_k),z_k) - \beta (Z_0 (M),z)| \leq \beta (Z_0 (M_k),Z_0 (M)) + |z_kz|.$$ 
Тогда в силу следствия 1 получим равенство
$$\beta (Z_0(M),z) = R_0 (M).$$
Следовательно, $z \in Q_0 (M)$ и
$$\lim \limits_{k \rightarrow \infty} \beta (Q_0 (M_k),Q_0 (M)) = \lim \limits_{k \rightarrow \infty} |z_kQ_0 (M)| = 0.$$
Таким образом, теорема 2 доказана.
\vskip7pt
{\bf Пример 1.}
\vskip7pt
Рассмотрим для каждого $n \in \{2,3,\ldots\}$ в евклидовой плоскости $\mathbb{R}^2$  $3$-сети
$$S = \{O(0;0),A(1/2;\sqrt{3}/2),B(-1/2;\sqrt{3}/2)\},$$ 
$$S_n = \{C_n(0;1/n),A(1/2;\sqrt{3}/2),B(-1/2;\sqrt{3}/2)\}.$$
 Тогда нетрудно проверить, что для каждого $n \in \{2,3,\ldots\}$ имеют место равенства:
$$\alpha(S_n,S) = 1/n,\quad Z_0 (S) = H(S) = Q_0(S) = S,$$ 
$$Z_0(S_n) = \{C_n\},\quad H(S_n) = Q_0(S_n) = \{A,B\},\quad \alpha (Z_0(S_n),Z_0(S)) = |BC_n|,$$ 
$$\alpha(H(S_n),H(S)) = \alpha(Q_0(S_n),Q_0(S)) = 1.$$
Таким образом, в теореме 2, в общем случае, нельзя заменить отклонение $\beta$ на метрику Хаусдорфа.

\section{Относительный чебышевский центр конечного множества геодезического пространства}

\vskip20pt

Рассмотрим метрическое пространство $(X,\rho)$ и примем следующие обозначения.
\vskip7pt
$B_{\alpha}(S,r)$ --- открытый шар с центром в точке $S \in (\Sigma^*_N (X),\alpha)$, радиуса $r > 0$.
\vskip7pt
Пусть $S \in \Sigma_N (X)$, $x \in S$. Если $S \neq \{x\}$, то $S(x) = S\backslash {\{x\}}$. Если $S = \{x\}$, то $S(x) = {x}$.
$$m(S) = \min \{|xy| : \, x, \, y \in S, \, x \neq y\}.$$
$$m_1(S) = \max \{|xS(x)| :\, x \in S\}.$$
$$h(S) = \{x \in S :\, |xS(x)| = m(S)\}.$$  
$$h_1(S) = \{x \in S :\, |xS(x)| = m_1 (S)\}.$$
$$d_0 (N) = \{S \in \Sigma^*_N (X) :\, D(S) = R_0 (S)\}$$ 
--- множество диаметральных $N$-сетей в $\Sigma^*_N (X)$.
\vskip7pt
Пусть $1 \leq k \leq N$.
$$Z_{k,N} = \{ S \in \Sigma^*_N (X) : \, \mbox{card\,}(Z_0 (S)) \leq k\}.$$
$$Z^*_{k,N} = \{ S \in \Sigma^*_N (X) : \, \mbox{card\,}(Z_0 (S)) = k\}.$$
$$md (N) = \{S \in \Sigma^*_N (X) :\, D(S) = m(S) \}.$$
$$mr_0 (N) = \{S \in \Sigma^*_N (X) :\, m(S) = R_0 (S) \}.$$
$$dm_1 (N) = \{S \in \Sigma^*_N (X) :\, D(S) = m_1 (S) \}.$$
$$mm_1 (N) = \{S \in \Sigma^*_N (X) :\, m(S) = m_1 (S) \}.$$
$$m_1r_0 (N) = \{S \in \Sigma^*_N (X) :\, m_1(S) = R_0 (S) \}.$$
В дальнейшем мы будем рассматривать геодезическое пространство $X$,
удовлетворяющее следующему условию $(CA)$.
\vskip7pt
$(CA)\quad$ Для любого сегмента $[x,y] \subset X$ и любой точки $z \in X$ верно неравенство
$$2|z\omega_{1/2} (x,y)| \leq |zx| + |zy|.$$
Причем, если точки $x$, $y$, $z$ не принадлежат одному сегменту, то неравенство строгое.
\vskip7pt
Из условия $(CA)$ нетрудно получить следующее свойство $(CB)$ в геодезическом пространстве $X$.
\vskip7pt
$(CB)\quad$ Для любого сегмента $[x,y] \subset X$, любой точки $z \in X$ и любого числа $\varepsilon \in (0,1)$ верно неравенство
$$|z\omega_{\varepsilon}(x,y)| \leq (1 - \varepsilon)|zx| + \varepsilon |zy|.$$
Причем, если точки $x$, $y$, $z$ не принадлежат одному сегменту, то неравенство строгое.
\vskip7pt
Приведем некоторые простые свойства и взаимосвязи рассмотренных понятий в метрическом пространстве $X$, доказательства первых двух из которых непосредственно следуют из их определений.
\vskip7pt
$1.\, \,$ Пусть $N > 2$. Тогда 
$$Z^*_{N-1,N} = \oslash,\quad Z_{N,N} = \Sigma^*_N (X),\quad Z^*_{N,N} = d_0 (N) = \Sigma^*_N (X) \backslash Z_{N-2,N}.$$

$2.\, \,$ Пусть $S \in \Sigma^*_N (X)$. Тогда следующие условия равносильны: 
$$(i)\, \, \, h(S) = S;\quad (ii)\, \, \, h_1 (S) = S;\quad (iii)\, \, \, m(S) = m_1(S).$$

$3.\, \,$ В геодезическом пространстве $X$, удовлетворяющем условию $(CA)$, для любых точек $x$, $y$ существует единственный сегмент $[x, y]$.
\vskip7pt
Действительно, в противном случае найдутся такие два сегмента $[x,y]$, $[x,y]'$ с концами в точках $x$, $y \in X$ и различными серединами $u = \omega_{1/2}(x,y)$ и $v = \omega'_{1/2}(x,y)$ соответственно. Тогда из условия $(CA)$
получим 
$$|x\omega_{1/2}(u,v)| + |y\omega_{1/2}(u,v)| < |xy|.$$ 
Получили противоречие.
\vskip7pt
$4.\, \,$ Пусть $X$ геодезическое пространство, удовлетворяющее условию $(CA)$, $M \in K(X)$, точки
$x$, $y \in Z_0 (M)$ различны и $\varepsilon \in (0,1)$. Тогда 
$$\omega_{\varepsilon}(x,y) \notin M.$$
\vskip7pt
Докажем это свойство.
Предположим противное, т.е. $\omega_{\varepsilon}(x,y) \in M$. Тогда в силу условий нашего утверждения и свойства $(CB)$ найдется такая точка $z \in M$, что  
$$\beta (M,\omega_{\varepsilon}(x,y)) = |z\omega_{\varepsilon}(x,y)| < \max \{|zx|,|zy|\} \leq R_0 (M).$$ 
Получили противоречие.
\vskip7pt
В следующих лемме 1, теореме 1 и следствии 1 установлены некоторые геометрические свойства рассмотренных классов фигур пространства $(\Sigma^*_N (X),\alpha)$ и их относительных чебышевских центров.
\vskip7pt
{\bf Лемма 1} [s19]. {\it Пусть $X$ метрическое пространство, 
\\
$S, \, \hat S \in (\Sigma^*_N (X),\alpha)$. Тогда 
$$mr_0 (N)= mm_1 (N)\cap m_1r_0 (N),\, \, \, md (N) = mm_1 (N)\cap m_1r_0 (N) \cap d_0 (N), \leqno (i)$$
$$dm_1 (N) = m_1r_0 (N) \cap d_0 (N)=$$
$$\{S \in d_0 (N) :   \,\mbox{существуют различные точки}\, \,$$ 
$$x_1,\ldots,x_N \in S \, (|x_1x_2| =\ldots =|x_1x_N| =D (S))\};$$
$$|m_1 (S) - m_1 (\hat S)| \leq 2 \theta (S,\hat S). \leqno (ii)$$}

{\bf Доказательство леммы 1}.
\vskip7pt
$(i)\, \,$ Докажем сначала, что для любого $S \in \Sigma^*_{N} (X)$ верно неравенство: 
$$m_1 (S) \leq R_0 (S).$$
Выберем такой элемент $x \in S$, что $m_1 (S) = |xS(x)|$. Тогда для каждого $y \in S(x)$
справедливы неравенства:
$$|xS(x)| \leq |xy| \leq \beta (S,y).$$ 
Кроме того, 
$$|xS(x)| \leq \beta (S,x).$$ 
Следовательно,
$$m_1 (S) \leq \min \{\beta (S,z) : z \in S \} = R_0 (S).$$ 
Используя полученное неравенство, а также
очевидные неравенства 
$$m(S) \leq m_1(S),\quad R_0 (S) \leq D(S)$$ 
для $S \in \Sigma_N (X)$, получим
$$mr_0 (N)= mm_1 (N)\cap m_1r_0 (N),\quad dm_1 (N) = m_1r_0 (N) \cap d_0 (N),$$
$$md (N) = mm_1 (N)\cap m_1r_0 (N) \cap d_0 (N).$$ 
Кроме того, из полученного представления множества $dm_1 (N)$ и его определения следует равенство
$$dm_1 (N) = \{S \in d_0(N) :   \,\mbox{существуют различные точки}\, \,$$ 
$$x_1,\ldots,x_N \in S \, (|x_1x_2| = \ldots = |x_1x_N| = D(S))\}.$$

$(ii)\, \,$ Заметим, что для любых $x \in S,\, \, y \in \hat S$ справедливы неравенства:
$$|xS(x)| - |y\hat S(y)| \leq |xy| + |yS(x)| \leq 2 \beta (S,y).$$ Следовательно, 
$$m_1 (S) - m_1 (\hat S) \leq 2 R_{\hat S} (S) \leq 2 \theta (S,\hat S).$$
Осталось использовать произвольность выбора $S, \, \hat S \in \Sigma^*_N (X)$. Лемма 1 доказана.
\vskip7pt
{\bf Теорема 1} [s19]. {\it Пусть $N > 2$, $X$ --- геодезическое пространство, удовлетворяющее условию $(CA)$. Тогда в пространстве $(\Sigma^*_N (X), \alpha)$  имеют место следующие соотношения.
\vskip7pt
$(i)$ Для всех $k \in \{1,\ldots,N\}$ 
$$\stackrel{\circ}{Z}_{k,N} = Z_{k,N}.$$

$$\overline{Z_{N-2,N}} = \Sigma^*_N (X); \leqno (ii)$$
$$\overline{Z_{1,N}} = \{S \in \Sigma^*_N (X) : D (Z_0 (S)) < R_0 (S)\}\cup  \leqno (iii)$$
$$\{S \in \Sigma^*_N (X) : \exists x, \, y \in S \, (|xy| = R_0 (S), D (\{x,y\},S\backslash \{x,y\})= R_0 (S))\}.$$
$$dm_1 (N) \cup d_{0,N-1} \subset \overline{Z_{1,N}}\cap d_0 (N), \leqno (iv)$$
где $$d_{0,N-1} = \{S \in d_0 (N) :   \,\mbox{существуют различные точки}\, \,$$ 
$$x_1,\ldots,x_{N-1} \in S \,
(|x_1x_2| =\ldots = |x_1x_{N-1}| = D (S))\}.$$ 
Причем, при $N \in \{3,4\}$ имеет место равенство.}
\vskip7pt
{\bf Замечание 1.} {\it Автор первоначально доказал утверждение $(i)$ теоремы 1 для случая $k = 1$.
На верность утверждения $(i)$ теоремы 1 для случая $k >1$ автору указал рецензент статьи} [s19].
\vskip7pt
{\bf Доказательство теоремы 1}.
\vskip7pt
$(i)\, \,$ В силу свойства 1 достаточно рассмотреть случай, когда $k \leq N-2$. Пусть
$S \in Z_{k,N}$ и 
$$\varepsilon = \frac {1}{8} \min \{R_0 (S\backslash Z_0 (S)) - R_0 (S), m(S)/2\}.$$ 
Рассмотрим произвольную
$N$-сеть $\hat S \in B_{\alpha} (S,\varepsilon)$. Тогда нетрудно проверить, что 
$$Z_0 (\hat S) = \hat S \cap B(Z_0(S),\varepsilon)$$ 
и $\hat S \in Z_{k,N}$. Следовательно, 
$$S \in \stackrel{\circ}{Z}_{k,N}$$ 
Таким образом,
$$\stackrel{\circ}{Z}_{k,N} = Z_{k,N}.$$

$(ii)\, \,$ Учитывая свойство 1, достаточно рассмотреть случай, когда $Z_0 (S) = S$ и
найти в произвольной окрестности $N$-сети $S$ $N$-сеть $\hat S \in Z_{N-2,N}$.
Выберем $\varepsilon = m(S)/4$. Пусть найдутся такие попарно различные элементы $x, \, y, \, z \in S$, что 
$$R_0 (S) = |xy| = |xz|.$$ 
Положим 
$$\hat y = \omega_{\varepsilon/2}(y,z),\quad \hat S = \{\hat y\}\cup S(y).$$ 
Тогда в силу свойства $(CB)$ получим  
$$x, \, z \notin Z_0(\hat S).$$ 
Следовательно,
$$\hat S \in Z_{N-2,N}\cap B_{\alpha} (S,\varepsilon).$$ 
Отметим, что для $N = 3$ наше утверждение доказано.
Рассмотрим случай, когда $N > 3$ и таких элементов нет. Тогда найдутся такие элементы $x, \, y, \, u, \, v  \in S$, что 
$$R_0 (S) = |xy| = |uv|$$ 
и $u, \, v \notin \{x,y\}$. Положим 
$$\hat x = \omega_{\varepsilon/2}(x,y),\quad \hat S = \{\hat x\}\cup S(x).$$ 
В силу свойства $(CB)$ получим  
$$u, \, v \notin Z_0(\hat S).$$ 
Таким образом, 
$$\hat S \in Z_{N-2,N}\cap B_{\alpha} (S,\varepsilon).$$

$(iii)\, \,$ Пусть 
$$S \in \{S \in \Sigma^*_N (X) : D (Z_0(S)) < R_0 (S)\}\cup $$
$$\{S \in \Sigma^*_N (X) : \exists x, \, y \in S \, (|xy| = R_0 (S),\, \, D (\{x,y\},S\backslash \{x,y\}) = R_0 (S))\}$$ 
и $S \notin Z_{1,N}$. Докажем, что $S \in \overline{Z_{1,N}}$. Выберем $\varepsilon = m(S)/4$ и рассмотрим два случая.
\vskip7pt
$1.\, \,$ Пусть
$$D (Z_0(S)) < R_0 (S),\quad u, \, v \in Z_0(S)$$ 
и $u \neq v$. Положим 
$$\hat u = \omega_{\varepsilon}(u,v),\quad \hat S = \{\hat u\}\cup S(u).$$ 
Тогда в силу свойства $(CB)$ получим 
$$\hat u \in Z_0(\hat S)$$ 
и $R_0 (\hat S) < R_0 (S)$. Если $z \in Z_0 (S)$, то найдется такой элемент $t \in S\backslash Z_0(S)$,
что 
$$|zt| = R_0 (S).$$ 
Если $z \in S\backslash Z_0 (S)$, то найдется такой элемент $t \in S(u)$,
что 
$$|zt| \geq R_0(S).$$ 
Следовательно, $\{\hat u\} = Z_0(\hat S)$.
\vskip7pt
$2.\, \,$ Пусть найдутся такие $x, \, y \in S$, что 
$$|xy| = R_0(S)\, \, \mbox{и}\, \, D(\{x,y\},S\backslash \{x,y\}) = R_0(S).$$ 
Для определенности считаем, что 
$$D(x,S\backslash \{x,y\}) = R_0 (S).$$ 
Если для каждого $z \in S\backslash \{x,y\}$ имеет место неравенство 
$$\beta(S\backslash \{y\},z) \geq R_0(S),$$ 
то выберем 
$$\hat y = \omega_{\varepsilon/2}(y,x),\quad \hat S = \{\hat y\}\cup S(y).$$ 
Тогда в силу свойства $(CB)$ получим 
$$\{\hat y\} = Z_0(\hat S).$$ 
Если найдется $z \in S\backslash \{x,y\}$ такой, что 
$$\beta(S\backslash \{y\},z) < R_0(S),$$ 
то $|zy| = R_0(S)$. Выберем 
$$\hat z = \omega_{\varepsilon/2}(z,x),\quad \hat S = \{\hat z\}\cup S(z).$$ 
Тогда в силу свойства $(CB)$ получим 
$$\{\hat z\} = Z_0 (\hat S).$$
Таким образом, в обоих случаях 
$$\hat S \in Z_{1,N}\cap B_{\alpha} (S,\varepsilon).$$
Пусть $S \in \overline{Z_{1,N}}$. Тогда либо 
$$D(Z_0(S)) < R_0 (S),$$ 
либо найдутся такие $x, \, y \in Z_0(S)$, что 
$$|xy| = R_0(S).$$ 
В первом случае $N$-сеть $S$ принадлежит множеству в правой части устанавливаемого равенства. Во втором случае предположим, что
$$D(\{x,y\},S\backslash \{x,y\}) < R_0(S)$$ 
и выберем
$$\delta = \frac {1}{8} \min \{R_0 (S) - D (\{x,y\},S\backslash \{x,y\}), m (S)/2 \}.$$ 
Тогда для любой $N$-сети
$\hat S \in B_{\alpha} (S,\delta)$ найдутся такие 
$$\{u\} = \hat S\cap B(x,\delta),\quad \{v\} = \hat S\cap B(y,\delta),$$ 
что $u \neq v$ и $u, \, v \in Z_0 (\hat S)$. Следовательно, 
$$\hat S \notin Z_{1,N}$$ 
и наше предположение приводит к противоречию с тем, что 
$$S \in \overline{Z_{1,N}}.$$
Таким образом, во втором случае 
$$D (\{x,y\},S\backslash \{x,y\}) = R_0 (S)$$ 
и $N$-сеть $S$ также принадлежит множеству в правой части устанавливаемого равенства. Следовательно, утверждение $(iii)$ доказано.
\vskip7pt
$(iv)\, \,$ Из равенств $(i)$ леммы 1 получим 
$$dm_1 (N) = \{S \in d_0 (N) :   \,\mbox{существуют различные точки}\, \,$$ 
$$x_1,\ldots,x_N \in S \, (|x_1x_2| =\ldots =|x_1x_N| =D (S))\}.$$ 
Следовательно, для произвольного
$S \in dm_1 (N)$ найдутся такие $x = x_1, \, y = x_2 \in S$, что
$$|xy| = D (S) = R_0 (S)\, \, \mbox{и}\, \, D (\{x,y\},S\backslash \{x,y\}) = D (S) = R_0 (S).$$ 
Тогда, учитывая доказанное утверждение $(iii)$,
получим 
$$dm_1 (N) \subset \overline{Z_{1,N}}\cap d_0 (N).$$ 
Пусть
$$S = \{x_1,\ldots,x_N\} \in d_{0,N-1}\, \, \mbox{и}\, \, |x_1x_2| =\ldots = |x_1x_{N-1}| = D (S).$$ Используя определения множеств $d_0 (N)$, $d_{0,N-1}$,
найдем такой элемент $x \in S(x_N)$, что $|xx_N| = D(S)$. Тогда
$$|xx_1| = D (S) = R_0 (S),\quad D (\{x,x_1\},S\backslash \{x,x_1\}) = D (S) = R_0 (S).$$ 
Учитывая доказанное утверждение $(iii)$,
получим 
$$d_{1,N-1} \subset \overline{Z_{1,N}}\cap d_0 (N).$$
Докажем обратное включение при $N \in \{3,4\}$. Используя доказанное утверждение $(iii)$, свойство 1 и определение множеств $d_0 (N)$, $d_{0,N-1}$, получим
$$\overline{Z_{1,N}}\cap d_0 (N) =$$ 
$$\{S \in d_0 (N) : \exists x, \, y \in S \, (|xy| = D(S),\, \, D(\{x,y\},S\backslash \{x,y\})= D(S))\} =$$
$$\{S \in d_0 (N) : \exists x, \, y, \, z \in S \, (|xy| = |xz| = D (S))\} = d_{0,N-1} \subset dm_1 (N) \cup d_{0,N-1}.$$
Таким образом, теорема 1 доказана.
\vskip7pt
Из определения множества $d_0 (N)$, теоремы 1 и свойства 1 получим
\vskip7pt
{\bf Следствие 1.} {\it Пусть $N > 2$, $X$ --- геодезическое пространство, удовлетворяющее условию $(CA)$. Тогда имеет место равенство
$$d_0 (N) = \partial Z_{N-2,N}.$$}

{\bf Пример 1.}
\vskip7pt
Пусть $4$-сеть $S$ состоит из вершин квадрата в евклидовой плоскости.
Tогда в силу утверждения $(iii)$ теоремы 1
$$S \in d_0 (4)\backslash \overline{Z_{1,4}}.$$

{\bf Пример 2.}
\vskip7pt
Рассмотрим $5$-сеть $S$, состоящую из вершин правильного пятиугольника в евклидовой плоскости. Используя утверждение $(iii)$ теоремы 1, нетрудно проверить, что
$$S \in \overline{Z_{1,5}}\cap d_0 (5).$$ 
С другой стороны, 
$$S \notin dm_1 (5) \cup d_{0,4}.$$ 
Таким образом, в общем случае, включение множеств в $(iv)$ теоремы 1 строгое уже при $N = 5$.
\vskip7pt
 Следующую теорему 2 вместе с доказательством автору сообщил рецензент при рецензировании статьи [s19], в которой можно найти ее доказательство.  
\vskip7pt
{\bf Теорема 2.} {\it Пусть $N > 2$, $X$ --- геодезическое пространство, удовлетворяющее условию $(CA)$,  
$$\Sigma^*_N (X)\backslash \overline{Z_{1,N}} \neq \oslash.$$ 
Тогда
$$\Sigma^*_N (X) = \overline{Z_{1,N}\cup \hat Z_{2,N}},$$ 
где 
$$\hat Z_{2,N} = \{S \in Z^*_{2,N} (X) : D(Z_0 (S),S\backslash Z_0 (S)) < R_0 (S), \, D(Z_0 (S)) = R_0 (S)\}$$ 
--- открытое множество в пространстве $(\Sigma^*_N (X), \alpha)$, непересекающееся с множеством $\overline{Z_{1,N}}$.}
\vskip7pt
Эта теорема 2 удачно дополняет теорему 1, поскольку из нее, в частности, следует, что при $N > 2$ множество 
$$Z_{1,N}\cup \hat Z_{2,N}$$ 
всюду плотно в пространстве $(\Sigma^*_N (X), \alpha)$.  

\section{Достаточные условия существования и единственности чебышевского центра непустого ограниченного множества геодезического пространства}

\vskip20pt

Достаточные условия существования и единственности чебышевского центра непустого ограниченного множества банахова пространства были получены Гаркави А. Л. \cite {Gark}, а достаточные условия принадлежности чебышевского центра замыканию выпуклой оболочки данного множества банахова пространства были получены Гаркави А. Л. \cite {Gark1} и Белобровым П. К. \cite {Belob1}. В случае геодезического пространства такие условия были установлены автором в [s9], [s21].
Перечислим отдельно условия, некоторые из которых будем налагать на метрическое пространство $(X,\rho)$.
\vskip7pt
$(A_4)\quad$  Для любых $x$, $y \in X$ существует единственная точка
\\
$\omega (x,y) \in X$ такая, что
$$|x\omega (x,y)| = |y\omega (x,y)| = \dfrac{1}{2}|xy|.$$

$(A_5)\quad$  Для всех точек $p,\, x,\, y$ пространства $X$, удовлетворяющего условию $(A_4)$, верно неравенство
$$|p\omega (x,y)| \leq \max \{|px|,|py|\}.$$

$(A_6)\quad$ Отображение $\omega : X\times X \to X$ равномерно непрерывно на каждом множестве $B \times B$, где В --- произвольный замкнутый шар пространства $X$, удовлетворяющего условию $(A_4)$.
\vskip7pt
$(A_7)\quad$  Для каждого $r > 0$ и для любых ограниченных последовательностей 
$$(p_n)_{n \in \mathbb{N}},\quad (x_n)_{n \in \mathbb{N}},\quad (y_n)_{n \in \mathbb{N}}$$ 
пространства $X$, удовлетворяющего условию $(A_4)$, таких, что для каждого $n \in \mathbb{N}$  
 $$|p_nx_n| \leq r, \quad |p_ny_n| \leq r, \quad \lim\limits_{n\to\infty}{|p_n\omega (x_n,y_n)| } = r,$$ 
верно равенство
$$\lim\limits_{n\to\infty}{|x_ny_n|} = 0.$$ 

$(A_8)\quad$  Для каждого $r > 0$, для каждой сходящейся к нулю последовательности неотрицательных вещественных чисел $(\varepsilon _n)_{n \in \mathbb{N}}$ и
для любых ограниченных последовательностей 
$$(p_n)_{n \in \mathbb{N}},\quad (x_n)_{n \in \mathbb{N}},\quad (y_n)_{n \in \mathbb{N}}$$
пространства $X$, удовлетворяющего условию $(A_4)$, таких, что для каждого $n \in \mathbb{N}$ 
$$|p_nx_n| \leq r + \varepsilon _n, \quad |p_ny_n| \leq r + \varepsilon _n, \quad
\lim\limits_{n\to\infty}{|p_n\omega (x_n,y_n)|} = r,$$ 
верно равенство 
$$\lim\limits_{n\to\infty}{|x_ny_n|} = 0.$$

Отметим, что условие $(A_4)$ является более сильным, чем условие  $(A_0)$, но более слабым, чем условие $(A_1)$. В геодезическом пространстве условие $(A_4)$ эквивалентно условию $(A_1)$. Полное метрическое пространство, удовлетворяющее условию $(A_4)$, является геодезическим пространством \cite {Efrem}. Очевидно, что в геодезическом пространстве, удовлетворяющем условию $(A_4)$, условие $(A_5)$ эквивалентно условию выпуклости произвольного замкнутого шара и слабее условия $(A_3)$. Примерами полных метрических пространств, в которых выполнены условия $(A_4-A_8)$, являются равномерно выпуклое банахово пространство и пространство Адамара \cite [c. 390] {Bur} (в частности, гильбертово пространство и пространство Лобачевского).
\vskip7pt
{\bf Лемма 1.} {\it Если геодезическое пространство $X$ удовлетворяет условиям $(A_4-A_7)$, то оно удовлетворяет и условию $(A_8)$.}
\vskip7pt
{\bf Доказательство леммы 1}.
\vskip7pt
Пусть заданы число $r > 0$ и последовательности
$$(\varepsilon _n)_{n \in \mathbb{N}},\quad  (p_n)_{n \in \mathbb{N}},\quad (x_n)_{n \in \mathbb{N}},\quad (y_n)_{n \in \mathbb{N}},$$
как указано в $(A_8)$.
Рассмотрим для каждого $n \in \mathbb{N}$ такие точки $u_n$, $v_n$ в сегментах $[p_n,x_n]$, $[p_n,y_n]$ соответственно, что
$$|u_nx_n| = \mu _n |p_nx_n|, \quad |v_ny_n| =  \mu_n |p_ny_n|,$$
где  $$\mu_n = \dfrac{\varepsilon_n}{r+\varepsilon_n}.$$
Тогда для каждого $n \in \mathbb{N}$ верны вспомогательные неравенства.
$$|u_nx_n| \leq \varepsilon _n, \quad |y_nv_n| \leq \varepsilon _n, \quad
|p_nu_n| \leq r, \quad |p_nv_n| \leq r,$$ 
$$|p_n\omega (x_n,y_n)| -  |\omega (x_n,y_n)\omega (u_n,v_n)| \leq  |p_n\omega (u_n,v_n)| \leq$$ 
$$|p_n\omega (x_n,y_n)| +  |\omega (x_n,y_n)\omega (u_n,v_n)|.$$
Используя эти неравенства и условие $(A_6)$, получим 
$$\lim\limits_{n\to\infty}{|p_n\omega (u_n,v_n)|} = r.$$
Следовательно, в силу условия $(A_7)$  
$$\lim\limits_{n\to\infty}{|u_nv_n|} = 0.$$
Из этого равенства и неравенств 
$$|x_ny_n|\leq |x_nu_n| + |u_nv_n| + |v_n y_n| \leq 2 \varepsilon _n + |u_nv_n|$$ для каждого $n \in \mathbb{N}$ получим 
$$\lim\limits_{n\to\infty}{|x_ny_n|} = 0.$$ 
Таким образом, лемма 1 доказана.
\vskip7pt
{\bf Лемма 2} ([s9], [s21]). {\it Для каждого непустого ограниченного множества метрического пространства, удовлетворяющего условиям $(A_4)$, $(A_5)$, $(A_7)$, существует не более одного чебышевского центра.}
\vskip7pt
{\bf Доказательство леммы 2}.
\vskip7pt
 Пусть $x$, $y$ чебышевские центры множества $M \in B(X)$. Тогда из определения чебышевского центра и условий $(A_4)$, $(A_5)$ следует, что верны неравенства
$$\sup\limits_{u\in M}{|\omega (x,y)u|} \leq \sup\limits_{u \in M}{\max \{|xu|,|yu|\}} \leq R (M).$$ 
Тогда $\omega (x,y)$ --- чебышевский центр множества $M$ и из определения чебышевского центра следует, что в множестве $M$ найдется такая последовательность $(p_n)_{n \in \mathbb{N}}$, что для каждого $n \in \mathbb{N}$
$$\lim\limits_{n\to\infty}{|p_n\omega (x,y)|} = R(M), \quad |p_nx| \leq R (M), \quad |p_ny| \leq R(M).$$
В силу условия $(A_7)$ получим $x=y$. 
Лемма 2 доказана.
\vskip7pt
{\bf Теорема 1} ([s9], [s21]). {\it Для каждого непустого ограниченного множества полного метрического пространства, удовлетворяющего условиям $(A_4)$, $(A_5)$, $(A_6)$, $(A_8)$, существует единственный чебышевский центр.}
\vskip7pt
{\bf Доказательство теоремы 1}.
\vskip7pt
Для множества $M \in B(X)$ рассмотрим семейство непустых замкнутых множеств
$$K_{\varepsilon} (M) = \{y \in X: \beta  (M,y) \leq R (M) + \varepsilon\},$$
где $\varepsilon > 0$. 
\\
Это семейство обладает следующими элементарными свойствами: 
\vskip7pt
$(i)\, \,$ если $\alpha \leq \varepsilon$, то 
$$K_{\alpha} (M) \subset K_{\varepsilon} (M) ;$$ 
\vskip7pt
$(ii)\, \,$ $\cap\{K_{\varepsilon} (M) : \varepsilon > 0\}$  есть множество всех чебышевских центров множества $M$. 
\vskip7pt
По лемме 2 существует не более одного чебышевского центра множества $M$.
Выберем сходящуюся к нулю последовательность положительных вещественых чисел $(\varepsilon_n)_{n \in \mathbb{N}}$ и для каждого $n \in \mathbb{N}$ такие точки $x_n$, $y_n \in K_{\varepsilon_n} (M)$, что
$$D (K_{\varepsilon _n} (M)) \leq |x_ny_n| + \varepsilon _n.$$
В силу условий $(A_4)$, $(A_5)$,  определения чебышевского радиуса и множества $K_{\varepsilon} (M)$ для каждого $n \in \mathbb{N}$ получим 
$$R (M) \leq \sup\limits_{z \in M}{|z\omega (x_n,y_n)|} \leq \sup\limits_{z\in M}{\max \{|zx_n|,|zy_n|\}} \leq R (M) + \varepsilon_n.$$ 
Следовательно, в множестве $M$ найдется такая последовательность $(z_n)_{n \in \mathbb{N}}$, что 
$$\lim\limits_{n\to\infty}{|z_n\omega (x_n,y_n)|} = R (M). $$
Теперь в силу условия $(A_8)$ получим 
$$\lim\limits_{n\to\infty}{|x_ny_n|} = 0.$$ 
Тогда по теореме Кантора \cite [c. 399]{Engel} множество 
$$\cap \{K_{\varepsilon} (M) : \varepsilon > 0\}$$
состоит из одной точки. Таким образом, теорема 1 доказана.
\vskip7pt
Из леммы 1 и теоремы 1 следует
\vskip7pt
{\bf Следствие 1} ([s9], [s21]). {\it Для каждого непустого ограниченного множества полного метрического пространства, удовлетворяющего условиям $(A_4-A_7)$, существует единственный чебышевский центр.}

\section{Обобщение некоторых теорем Б. Секефальви-Надь,
С. Б. Стечкина и Н. В. Ефимова}

\vskip20pt

Пусть $(X,\rho)$ метрическое пространство. Напомним следующие обозначения и определения из \cite{Vlas}, \cite{Brudn}, \cite{Singer}.
\vskip7pt
Пусть $M$ непустое подмножество пространства $X$.
$$E_M = \{ x \in X : P (x, M) \neq \emptyset \},\quad
T_M = \{ x \in X : P (x, M) - \mbox{одноточечно} \},$$
$$D_M (x) = \lim \limits_{\delta
\to 0+} D (P (x, M, \delta)),\quad T'_M = \{ x \in X : D_M (x) = 0 \}.$$

Множество $M$ из $X$ называется {\it множеством существования} ( {\it чебышевским множеством}), если  $E_M = X$  ($T_M = X$).
\vskip7pt
Чебышевское множество $M$ называется  {\it сильно чебышевским множеством}, если отображение
$$P_M : X \to M,\quad P_M (x) = P (x, M)$$
---  непрерывно \cite [определение 3.28]{Brudn}.
\vskip7pt
Последовательность $(y_n)_{n \in \mathbb{N}}$ множества $M$ называется {\it минимизирующей последовательностью} для точки $x \in X$, если
$$\lim \limits_{n \rightarrow \infty} |xy_n| = |xM|.$$

Множество $M \subset X$ называется {\it аппроксимативно компактным множеством}, если для каждого $x \in X$ каждая минимизирующая последовательность имеет подпоследовательность, сходящуюся к элементу из множества $M$.
\vskip7pt
Множество $M \subset X$ называется {\it телом} \cite{Vlas}, если
$$M \subset \overline{\stackrel{\circ}{M}}.$$

Множество $M \subset X$ называется {\it $a$-выпуклым} ($a > 0$), если для
каждого элемента $x \in X\backslash M$ найдется элемент $c \in X$ такой,
что $x \in B(c,a)$ и $M\bigcap B(c,a) = \emptyset$ \cite{Vlas}.
\vskip7pt
Замкнутый шар $B[x,a]$ называется {\it опорным шаром} к множеству $M$ в точке $y \in M$, если $|xM| = |xy| = a$ \cite{Vlas}.
\vskip7pt
Вместе с условиями $(A_4-A_6)$ на геодезическое пространство $X$ будем использовать также вместо условия $(A_7)$ более слабое условие $(A'_7)$.
\vskip7pt
$(A'_7)\quad$   Для каждого замкнутого шара $B[p,r]$ геодезического пространства $X$, удовлетворяющего условию $(A_4)$, и для любых
последовательностей 
$$(x_n)_{n \in \mathbb{N}},\quad (y_n)_{n \in \mathbb{N}}$$
из этого шара таких, что
$$\lim \limits_{n \rightarrow \infty} |p\omega (x_n,y_n)| = r,$$
верно равенство
$$\lim \limits_{n \rightarrow \infty} |x_ny_n| = 0.$$

Приведем несколько вспомогательных утверждений, доказательства первых
двух из них очевидны, а остальных известны.
\vskip7pt
$1.\, \,$ Каждый открытый шар геодезического пространства, через каждые две
различные точки которого можно провести единственную прямую, выпуклый тогда и только тогда, когда каждый замкнутый шар этого пространства выпуклый.
\vskip7pt
$2.\, \,$ В геодезическом пространстве, удовлетворяющем условиям $(A_4-A_6)$, $(A'_7)$, каждая сфера не содержит невырожденных отрезков.
\vskip7pt
$3.\, \,$ Каждое аппроксимативно компактное чебышевское множество метрического пространства является сильно чебышевским (\cite [следствие 2]{Singer}  или \cite [следствие 2.2]{Vlas}).
\vskip7pt
$4.\, \,$ Пусть $X$ геодезическое пространство, через каждые две
различные точки которого можно провести единственную прямую, $x \in X$, $M \in \Sigma (X)$. Если
$$y \in P (x, M),\quad z \in (x,y],$$
где  $(x,y] = [x, y] \backslash \{x\}$, то $P (z, M) = y$ \cite[предложение 0.3]{Vlas}.
\vskip7pt
$5.\, \,$ Пусть $M$ замкнутое подмножество полного геодезического пространства $X$, через каждые две различные точки которого можно провести единственную прямую. Тогда $M \subset T'_M \subset T_M$ и для каждого $x \in T'_M$ имеет место включение \cite[предложение 1.1]{Vlas}
$$[x,P (x, M)] \subset T'_M.$$
\vskip7pt
Следующая теорема 1 обобщает теорему Б. Секефальви - Надь \cite [теорема 3.35]{Brudn}, другие обобщения можно найти в \cite [предложения 2.5, 2.6]{Vlas}.
\vskip7pt
{\bf Теорема 1} [s7]. {\it
В полном геодезическом пространстве, удовлетворяющем условиям $(A_4-A_6)$, $(A'_7)$, каждое выпуклое замкнутое множество является аппроксимативно компактным, чебышевским {\rm (}а значит и сильно чебышевским{\rm )}.}
\vskip7pt
{\bf Замечание 1.} {\it
Условие $(A_6)$ в теореме 1 можно ослабить до следующего условия $(A'_6)$.
\vskip7pt
$(A'_6)\quad$ Если
$$\lim \limits_{n \rightarrow \infty} |x_ny_n| = 0$$
для ограниченных последовательностей $(x_n)_{n \in \mathbb{N}}$, $(y_n)_{n \in \mathbb{N}}$ пространства $X$, то
$$|\omega (x_n,x_m)\omega (y_n,y_m)| \to 0$$
при $n \to \infty,\quad m \to \infty$.}
\vskip7pt
{\bf Доказательство теоремы 1}.
\vskip7pt
Докажем аппроксимативную компактность
множества $M$. Пусть
$$p \in X\backslash M,\quad (x_n)_{n \in \mathbb{N}}$$
--- минимизирующая последовательность в множестве $M$ для точки $p$ и для каждого $n \in \mathbb{N}$
$$y_n = [p,x_n]\bigcap S(p,|pM|).$$
Тогда
$$\lim \limits_{n \rightarrow \infty} |x_ny_n| = \lim \limits_{n \rightarrow \infty} (|x_np| -
|py_n|) = \lim \limits_{n \rightarrow \infty} |x_np| - |pM| = 0.$$
Из условия $(A_6)$ (или $(A'_6)$) следует, что 
$$|\omega (x_n,x_m)\omega (y_n,y_m)| \to 0$$
при $n \to \infty,\quad m \to \infty$.
Кроме того,
$$\omega (y_n,y_m) \in B[p,|pM|],\quad \omega (x_n,x_m) \in M$$
при $n,\, m \in \mathbb{N}$, поскольку $M$, $B[p,|pM|]$ --- выпуклые множества.
Пусть
$$a_{nm} \in S(p,|pM|)\bigcap [\omega (x_n,x_m),\omega (y_n,y_m)]$$
при $n,\, m \in \mathbb{N}$. Тогда
$$|pM| - |\omega (y_n,y_m)a_{nm}| = |pa_{nm}| - |\omega (y_n,y_m)a_{nm}| \leq |p\omega (y_n,y_m) \leq |pM|$$
и значит,
$$|p\omega (y_n,y_m)| \to |pM|$$
при $n \to \infty$, $m \to \infty$.
Из условия $(A'_7)$ получим, что $|y_ny_m| \to 0$ при $n \to \infty$, $m \to \infty$. Но пространство $X$ --- полное, поэтому найдется точка $y \in
S(p,|pM|)$ такая, что
$$\lim \limits_{n \rightarrow \infty} y_n = y.$$
Из неравенства
$$|yx_n| \leq |x_ny_n| + |yy_n|$$
следует теперь, что
$$\lim \limits_{n \rightarrow \infty} x_n = y \in M\bigcap S(p,|pM|).$$
Значит, множество $M$ --- аппроксимативно компактно.

Пересечение $M\bigcap S(p,|pM|)$ содержит не более одной точки.
Действительно, в противном случае сфера содержала бы невырожденный
отрезок, поскольку множества $M$, $B[p,pM]$ выпуклые. А это
противоречит свойству 2. Следовательно, множество $M$ чебышевское и
теорема 1 доказана.
\vskip7pt
{\bf Лемма 1} [s7]. {\it
Пусть $X$ --- геодезическое пространство, удовлетворяющее условиям $(A_4-A_6)$, $(A'_7)$,   $x \in X$. Если
последовательности
$$(y_n)_{n \in \mathbb{N}},\quad (z_n)_{n \in \mathbb{N}},\quad (u_n)_{n \in \mathbb{N}}$$
удовлетворяют условиям:
$$|xy_n| = r,\quad \lim \limits_{n \rightarrow \infty} |xz_n| = r + h,\quad u_n \in [x,z_n],\quad |xu_n| = r,\quad \lim \limits_{n \rightarrow \infty} |y_nz_n| =h,$$
где $r>0$, $h>0$, $n \in \mathbb{N}$, то
$$\lim \limits_{n \rightarrow \infty} |y_nu_n| = 0.$$}

{\bf Доказательство леммы 1}.
\vskip7pt
Пусть для каждого $n \in \mathbb{N}$
$$c_n \in [z_n,y_n],\quad |c_nz_n| = |u_nz_n|.$$
Тогда
$$\lim \limits_{n \rightarrow \infty} |y_nc_n| = \lim \limits_{n \rightarrow \infty} (|y_nz_n| - |u_nz_n|) = 0$$
и по условию $(A_6)$
$$\lim \limits_{n \rightarrow \infty} |\omega (y_n,u_n)\omega (c_n,u_n)| = 0.$$
В силу условий леммы 1, $(A_5)$ и неравенства треугольника верны неравенства
$$|xz_n| \leq |x\omega (y_n,u_n)| + |\omega (y_n,u_n)\omega (c_n,u_n)| +$$
$$|\omega (c_n,u_n)z_n| \leq r + |\omega (y_n,u_n)\omega (c_n,u_n)| + h.$$
Следовательно,
$$\lim \limits_{n \rightarrow \infty} |x\omega (y_n,u_n)| = r.$$
Из условия $(A_7)$ следует теперь, что
$$\lim \limits_{n \rightarrow \infty} |u_ny_n| = 0.$$
Что и требовалось доказать.
\vskip7pt
Рассмотрим еще одно условие $(A_9)$ на геодезическое пространство $X$.
\vskip7pt
$(A_9)\quad$ Для любых вещественных чисел $r > 0$, $\lambda > 1$ и для каждой точки $p$ геодезического пространства $X$, через каждые две различные точки которого можно провести единственную прямую, отображение 
$$g : S(p, r) \rightarrow S(p, \lambda r),\quad g (x) = \omega_{\lambda} (p, x)$$
--- равномерно непрерывно.
\vskip7pt
Отметим, что совокупность условий $(A_6)$, $(A_9)$ в геодезическом пространстве $X$, через каждые две различные точки которого можно провести единственную прямую, слабее следующего условия $(A^*_6)$,  использованного автором в [s7].
\vskip7pt
$(A^*_6)\quad$ Для каждого $\lambda \in \mathbb{R}$ отображение
$$\omega_{\lambda} : X\times X \to X$$
равномерно непрерывно на каждом множестве вида $B \times B$, где В --- произвольный замкнутый шар геодезического пространства $X$, через каждые две различные точки которого можно провести единственную прямую.
\vskip7pt
Утверждение $(i)$ следующей леммы 2 является обобщением леммы 1.1 С. Б. Стечкина \cite{Vlas1}, утверждение $(ii)$ леммы 2 является обобщением утверждения в лемме 4  А. В. Маринова \cite{Marin1}.
\vskip7pt
{\bf Лемма 2} [s7], [s12]. {\it
Пусть в геодезическом пространстве $X$, через каждые две различные точки которого можно провести единственную прямую, выполнены условия $(A_5)$, $(A_6)$, $(A'_7)$, $(A_9)$,  $x \in X$ и даны такие вещественные числа $h$, $H$, что $0 < h <H$.
Тогда верны следующие утверждения.
$$D(B[z, H - h + \delta]\backslash B(x, H)) \to 0; \leqno (i)$$
$$D(B[z, H - h]\backslash B(x, H - \delta)) \to 0 \leqno (ii)$$
при $\delta \to 0+$ равномерно по всем $y,\, z \in X$, удовлетворяющим условиям:
$$z \in (x,y),\quad |xy| = H,\, |xz| = h,$$
где $(x,y) = [x, y] \backslash \{x, y\}.$}
\vskip7pt
{\bf Доказательство леммы 2}.
\vskip7pt
$(i)\, \,$ Пусть для каждого $\delta > 0$ в множестве
$$B[z, H - h + \delta]\backslash B(x,H)$$
произвольным образом выбран элемент $u$. Очевидно,
достаточно доказать, что
$$\lim \limits_{\delta \rightarrow 0+} |uy| = 0.$$
Обозначим
$$t = [x,u]\bigcap S(x, h),\quad v = [x,u]\bigcap S(x, H).$$
Тогда из условий леммы 2 и неравенства треугольника получим
$$H - h \leq |xu| - |xt| = |tu| \leq |zu| \leq H - h + \delta.$$
Следовательно, 
$$\lim \limits_{\delta \rightarrow 0+} |tu| = \lim \limits_{\delta \rightarrow 0+} |zu| = H - h.$$
В силу леммы 1 
$$\lim \limits_{\delta \rightarrow 0+} |tz| = 0.$$
Из условия $(A_9)$ следует теперь, что
$$\lim \limits_{\delta \rightarrow 0+} |vy| = 0.$$
Кроме того,  из условий леммы 2 и неравенства треугольника получим
$$|uy| \leq |vy| + |vu| \leq |vy| + \delta$$
Следовательно, 
$$\lim \limits_{\delta \rightarrow 0+} |uy| = 0.$$
$(ii)\, \,$ Пусть для каждого $0 < \delta < H - h$ в
множестве
$$B[z, H - h]\backslash B(x, H - \delta )$$
произвольным образом выбран элемент $u$. Очевидно, достаточно доказать, что $|uy| \to 0$ при $\delta \to 0+$.
Пусть $v$ точка пересечения прямой, проходящей через точки $x, \,  u$ со сферой
$$S(x,H),\quad t = [x,u]\bigcap S(x, h).$$
Из неравенств
$$H - \delta \leq |xu| \leq |xv| = H$$
следует, что
$$\lim \limits_{\delta \rightarrow 0+} |xu| = H, \quad \lim \limits_{\delta \rightarrow 0+} |uv| =\lim \limits_{\delta \rightarrow 0+} (H - |xu|) = 0.$$
 Тогда из неравенств
$$|zy| \leq |zv| \leq |zu| + |uv| \leq |zy| + |uv|$$
следует, что
$$\lim \limits_{\delta \rightarrow 0+} |zv| = H - h.$$
Но $|xt| = |xz| = h$. Поэтому по лемме 1, 
$$\lim \limits_{\delta \rightarrow 0+} |tz| = 0.$$
Из условия $(A_9)$ следует теперь, что 
$$\lim \limits_{\delta \rightarrow 0+} |vy| = 0.$$
Таким образом,
$$|uy| \leq |uv| +|vy| \to 0$$
при $\delta \to 0+$ и лемма 2 доказана.
\vskip7pt
Следующая теорема 2 обобщает теорему С. Б. Стечкина \cite [теорема 1.1]{Vlas}.
\vskip7pt
{\bf Теорема 2} [s7]. {\it
Пусть в полном геодезическом пространстве $X$, через каждые две различные точки которого можно провести единственную прямую, выполнены условия $(A_5)$, $(A_6)$, $(A'_7)$, $(A_9)$, $M$ --- непустое замкнутое множество в $X$. Тогда каждое из
множеств $T_M$, $T'_M$ является дополнением множества первой категории {\rm (}в частности, всюду плотно{\rm )} в пространстве $X$.}
\vskip7pt
{\bf Доказательство теоремы 2}.
\vskip7pt
Докажем, что множество
$$F_{\varepsilon} = \{ x \in X : D_M (x) \geq \varepsilon \}$$
замкнуто для каждого $\varepsilon > 0$.
Пусть
$$\lim \limits_{n \rightarrow \infty} x_n = x,\quad D_M (x_n) \geq \varepsilon$$
для каждого $n \in \mathbb{N}$. Покажем, что для каждого $n \in \mathbb{N}$
$$P (x_n, M, t) \subset P (x, M, \delta)$$
при
$$t < \delta /2,\quad |xx_n| < \delta /4.$$
Действительно, если для каждого $n \in \mathbb{N}$
$y \in P (x_n, M, t)$, то
$$|yx| \leq |yx_n| + |xx_n| \leq |x_nM| + t + |xx_n| \leq |xM| + t + 2 |xx_n|
\leq |xM| + \delta.$$
Следовательно, для каждого $n \in \mathbb{N}$ и для каждого $\delta > 0$
$$D (P (x, M, \delta)) \geq D(P (x_n, M, t)) \geq D_M (x_n) \geq \varepsilon.$$
Отсюда получаем
$$D_M (x) \geq \varepsilon, \quad x \in F_{\varepsilon}.$$
Докажем теперь, что множество $F_{\varepsilon}$ нигде не плотно в пространстве $X$. Если это не так, то замкнутое множество $F_{\varepsilon}$ содержит некоторый замкнутый шар $B[x,h]$ радиуса $h > 0$. Можно считать, что $h < |xM|$. Пусть точка $w \in M$ такая, что
$$|xw| \leq |xM| + \delta.$$
Положим
$$y = [x,w]\bigcap S(x,|xM|).$$
Пусть точка $z \in [x,y]$ такая, что $|xz| = h$. Тогда
$$|zw| = |xw| - |xz| \leq |xM| + \delta - h.$$
Кроме того, для каждого $u \in P (z, M, \delta)$
$$|zu| \leq |zM| + \delta \leq |zw| + \delta \leq |xM| - h + 2 \delta.$$
Значит,
$$P (z, M, \delta) \subset B[z,|xM| - h + 2 \delta]\backslash B(x,|xM|).$$
Из свойства 5 и утверждения $(i)$ леммы 2 следует, что
найдется такое $\delta > 0$, что
$$D_M (z) \leq D (P (z, M, \delta) \leq
D (B[z, |xM| - h + 2 \delta]\backslash B(x, |xM|)) < \varepsilon.$$
Получили противоречие с тем, что
$$z \in B[x,h] \subset F_{\varepsilon}.$$
Таким образом, для каждого $\varepsilon > 0$ множество $F_{\varepsilon}$
нигде не плотно в пространстве $X$. Кроме того, множество
$$\{x \in X : D_M (x) > 0 \}$$
является множеством первой категории, поскольку совпадает с объединением
множеств $F_{1/n}$, $(n \in \mathbb{N})$. Дополнением этого множества
является множество $T'_{M}$. Из свойства 5 следует, что $T'_M \subset
T_M$. Таким образом, теорема 2 доказана.
\vskip7pt
Следующая теорема обобщает теорему Н. В. Ефимова и С. Б. Стечкина \cite [теорема 1.2]{Vlas}.
\vskip7pt
{\bf Теорема 3} [s7]. {\it
Пусть в полном геодезическом пространстве $X$, через каждые две различные точки которого можно провести единственную прямую, выполнены условия $(A_5)$, $(A_6)$, $(A'_7)$, $(A_9)$, $M$ --- $a$-выпуклое тело этого пространства. Тогда множество точек границы $\partial M$, в которых существует единственный опорный шар радиуса $a$, является дополнением множества первой категории в $\partial M$.}
\vskip7pt
{\bf Доказательство теоремы 3.}
\vskip7pt
Доказательство в основном аналогично доказательству теоремы 1.2 в \cite{Vlas}, если учесть лемму 2.
\vskip7pt
$1.\, \,$ Для каждого $\delta > 0$ и для любого $y \in \partial M$ положим
$$C_{\delta} (y) = \{x \in X : B (x, a) \cap M = \emptyset, \, \, x \in B[y, a + \delta]\}.$$
Это множество непустое. Действительно, найдется такая точка
$z \in B(y, \delta)\backslash M$, для которой в силу определения $a$-выпуклого тела найдется открытый шар $B(x, a)$, удовлетворяющий условиям:
$$B (x, a) \cap M = \emptyset, \quad z \in  B(x, a).$$
Тогда 
$$|xy| \leq |xz| + |zy| < a + \delta,$$ 
 то есть $z \in C_{\delta} (y)$. Заметим, что при $0 < \delta < \delta'$
$$C_{\delta} (y) \subset C_{\delta'} (y). $$
Это простое свойство позволяет определить число 
$$d_y = \lim \limits_{\delta \rightarrow 0+} D (C_{\delta} (y)) = \inf \{D (C_{\delta} (y)) : \delta > 0\}$$
и для каждого $\varepsilon > 0$ множество
$$G_{\varepsilon} = \{y \in \partial M : d_y \geq \varepsilon\}.$$

$2.\, \,$ Покажем, что $G_{\varepsilon}$ замкнутое множество. Пусть последовательность $(y_n)_{n \in \mathbb{N}}$ из множества $G_{\varepsilon}$ сходится к точке $y \in X$, то есть для каждого $\delta > 0$
найдется такое $N \in  \mathbb{N}$, что для каждого $n > N$ 
$|y_ny| < \delta.$
Если $n > N$ и $x \in C_{\delta} (y_n)$, то 
$$B (x, a) \cap M = \emptyset, \quad x \in  B[y_n, a + \delta].$$
Кроме того,
$$|xy| \leq |xy_n| + |y_ny| \leq a + 2 \delta, \, \, \mbox{и} \, \, x \in C_{2 \delta} (y).$$
Следовательно, для любого $\delta > 0$ при $n > N$
$$C_{\delta} (y_n) \subset C_{2 \delta} (y).$$ 
Тогда для каждого $\delta > 0$ при $n > N$
$$D (C_{2 \delta} (y)) \geq D (C_{\delta} (y_n)) \geq d_{y_n} \geq
\varepsilon.$$
Следовательно, 
$$d_{y} \geq \varepsilon \, \, \, \mbox{и} \, \, \, y \in G_{\varepsilon}.$$

$3.\, \,$ Покажем, что для каждого $\varepsilon > 0$ множество $G_{\varepsilon}$ нигде не плотно в $\partial M$. Предположим, что это неверно. Тогда найдется такой замкнутый шар
$B[y, r]$, что
$$y \in \partial M, \quad B[y, r]\cap \partial M \subset G_{\varepsilon}.$$
В силу того, что $M$ --- тело, найдется такая точка
$x \in \stackrel{\circ}{M}$, что $|xy| < r /3$. Пусть
$(u_n)_{n \in \mathbb{N}}$ --- минимизирующая последовательность в
множестве $\partial M$ для точки $x$. Тогда найдется такое $N \in  \mathbb{N}$, что для каждого $n > N$
$$|xu_n| \leq 2 |xy|, \quad |yu_n| \leq |yx| + |xu_n| < r \, \, \mbox{и} \, \,
u_n \in  G_{\varepsilon}.$$
Пусть для каждого $n > N$
$$h = |x(\partial M)|,\quad z_n \in [x,u_n], \quad |xz_n| = h,$$
$y_n$ --- точка на луче с началом в точке $x$, содержащем точку $u_n$, такая, что $|xy_n| = a + h$ и $H = a + h$. Как установлено выше, для любого $\delta > 0$ и для каждого $n > N$
$$C_{\delta} (u_n) \neq \emptyset.$$
Следовательно, для каждого 
$n > N$  найдется $v_n \in C_{\delta} (u_n) $, что
$$a \leq |v_nM|, \quad \lim \limits_{n \rightarrow \infty} |v_nu_n| = a.$$
Кроме того, для каждого $n > N$
$$B (v_n, a) \cap M = \emptyset, \quad B(x, h) \subset M, \quad |xv_n| \geq a + h = H.$$
Следовательно,  для каждого $n > N$ 
$$v_n \neq B(x, |xy_n|).$$
Заметим, что для каждого $n > N$
$$a \leq |v_nz_n| \leq |v_nu_n| + |u_nz_n|, \quad \lim \limits_{n \rightarrow \infty} |u_nz_n| = 0.$$
Следовательно,
$$\lim \limits_{n \rightarrow \infty} |v_nz_n| = a.$$
Если сходящаяся к нулю последовательность положительных вещественных чисел $(\delta_n)_{n \in \mathbb{N}}$ удовлетворяет для каждого
$n > N$ условию
$$\delta_n > |v_nu_n| - a + |u_nz_n|,$$
то для каждого $n > N$ справедливо включение:
$$v_n \in C_{\delta_n} (u_n).$$
Таким образом, для каждого $n > N$
$$v_n,\, y_n \in B[z_n, a + \delta_n]\backslash B(x, H).$$
В силу леммы 2
$$\lim \limits_{n \rightarrow \infty} |v_ny_n| = 0.$$
Из включения  $u_n \in  G_{\varepsilon}$ для каждого $n > N$
следует, что найдется точка
$$q_n \in C_{\delta_n} (u_n),$$
для которой верно неравенство
$$|v_nq_n| \geq \varepsilon/3.$$
Снова, для каждого $n > N$
$$q_n \in B[z_n, a + \delta_n]\backslash B(x, H).$$
В силу леммы 2
$$\lim \limits_{n \rightarrow \infty} |q_ny_n| = 0, \quad 
\lim \limits_{n \rightarrow \infty} |v_nq_n| = 0.$$
Получили противоречие. Таким образом, для каждого $\varepsilon > 0$ множество $G_{\varepsilon}$ нигде не плотно в $\partial M$.  Тогда множество
$$\{y \in \partial M : d_y > 0\} = \bigcup \limits_{n = 1}^{\infty} G_{1/n}$$
--- первой категории в $\partial M$.
\vskip7pt
$4.\, \,$ Покажем, что в точке $y \in \partial M$ с  $d_y = 0$ существует единственный опорный шар радиуса $a$. Пусть $(\delta_n)_{n \in \mathbb{N}}$ сходящаяся к нулю строго убывающая последовательность вещественных чисел и для каждого $n \in \mathbb{N}$ имеет место включение:
$$ p_n \in C_{\delta_n} (y).$$
Последовательность $(p_n)_{n \in \mathbb{N}}$ --- фундаментальная, поскольку
$$\lim \limits_{n \rightarrow \infty} D (C_{\delta_n} (y)) = d_y = 0.$$
Но пространство $X$ полное, поэтому найдется такая точка $p \in X$, что
$$\lim \limits_{n \rightarrow \infty} p_n = p.$$
Для каждого $n \in \mathbb{N}$ имеет место неравенство: $|p_nM| \geq a.$
Следовательно, $|pM| \geq a.$
Кроме того, для каждого $n \in \mathbb{N}$ справедливы неравенства:
$$a \leq |pM| \leq |py| \leq |pp_n| + a + \delta_n.$$
Значит,
$$|pM| = |pa| = a\, \, \mbox{и} \, \, B[p, a]$$
--- искомый опорный шар. Предположим, что существует второй опорный шар $B[p', a]$ с тем же свойством. Тогда для каждого $n \in \mathbb{N}$ имеет место включение:
$$ p' \in C_{\delta_n} (y).$$
Откуда, $p = p'$.  Множество 
$$\{y \in \partial M : d_y = 0\}$$
является дополнением множества первой категории в $\partial M$.
Но это подмножество следующего множества
$$\{y \in \partial M : \mbox {в точке}\, \,  y \, \, \mbox {существует единственный опорный шар}\, \, B[p, a] \, \, \mbox {к}\, \,  M \}.$$
Следовательно, последнее множество также  является дополнением множества первой категории в $\partial M$. Теорема 3 доказана.

\section{Непрерывность и связность метрической \\
$\delta$-проекции в геодезическом пространстве}

\vskip20pt

Пусть $(X,\rho)$ --- геодезическое пространство, через каждые две
различные точки которого можно провести единственную прямую.

Напомним следующие обозначения и определения из \cite{Vlas}, \cite{Marin1}.
\vskip7pt
Пусть
$$x \in X,\quad M \in \Sigma (X),\quad t \geq 0,$$
$$\lambda (M) = \sup \{|AB| : \, A\cup B = M, \, A,\, B \neq \emptyset \},$$
$$\lambda_M (\delta ) = \sup \{\lambda (P (x, M, \delta)) : x \in X \}$$
--- {\it мера несвязности метрической $\delta$-проекции}.
$$\omega (t) = \sup \{\alpha (P (x, M, \delta), P (y, M, \delta)) : \, y \in X, \, |xy| \leq t, \,
P (y, M, \delta) \neq \emptyset \},$$
$$\hat {\omega}(t) = \sup \{\alpha (P (x, M, \delta), P (x, W, \delta))  : $$
$$ W \in \Sigma (X), \,
\alpha (M,W) \leq t, \, P (x, W, \delta) \neq \emptyset \},$$
$$\bar {\omega}(t) = \sup \{\alpha (P (x, M, \delta), P (x, M, \varepsilon)) : \, \varepsilon \geq 0, \,
|\varepsilon  - \delta | \leq t, \, P (x, M, \varepsilon) \neq \emptyset \},$$
$$\Omega (t) = \sup \{\alpha (P (x, M, \delta), P (y, W, \varepsilon)) : \,
\varepsilon \geq 0, \,
y \in X, \, W \in \Sigma (X), $$ 
$$|xy| \leq t, \, |\varepsilon  - \delta | \leq t, \, \alpha (M,W) \leq t, \,
P(y, W, \varepsilon) \neq \emptyset \}$$
--- четыре {\it модуля непрерывности оператора} $P$ в точке  $(x,M,\delta)$ по $x$, по $M$, по $\delta $ и по
 совокупности этих переменных, причем $P (x, M) \neq \emptyset$ при $\delta = 0$.
$$u_M (\delta ) = \sup \{\omega (0+) : x \in X \},\quad \hat {u}_M(\delta ) =
\sup \{\hat {\omega }(0+) : x \in X \},$$
$$\bar {u}_M (\delta ) = \sup \{\bar {\omega}(0+) : x \in X \},\quad
U_M (\delta ) = \sup \{\Omega (0+) : x \in X \},$$
где
$$\omega (0+) = \lim \limits_{t \rightarrow 0+}{\omega (t)},$$ --- {\it верхние грани разрывов
оператора} $P$ на множестве
$$\{(x,M,\delta ) : x \in X \}.$$

Оператор
$$P_M : X \rightarrow 2^X,\quad P_M (x) = P (x, M)$$
метрического проектирования на множество $M$ называется {\it непрерывным в слабом смысле} \cite{Vlas1}, если
$$\lim \limits_{n \rightarrow \infty} x_n = x$$
влечет $$\lim \limits_{n \rightarrow \infty} |P_M (x_n)P_M (x)| = 0.$$

Множество $M \subset X$ называется {\it $P$-связным}, если для каждого $x \in X$
множество $P_M (x)$ непустое и связное.
\vskip7pt
Множество $M \subset X$ называется {\it $B$-связным} ({\it $\stackrel{\circ}{B}$-связным}), если пересечение множества $M$ с каждым замкнутым (открытым) шаром является связным.
\vskip7pt
Приведем одно вспомогательное утверждение.
\vskip7pt
{\bf Предложение 1} \cite [теорема 1]{Marin1}. {\it
Пусть $M$ --- подмножество геодезического пространства $X$.
Тогда
$$\lambda_M (\delta + 0) \leq u_M (\delta )$$
при
$$0 < \delta < \lambda_M (\delta + 0),$$
а если $M$ --- множество существования, то справедливы неравенства
$$\lambda _M (\delta ) \leq u_M (\delta )$$
при $0 \leq  \delta < \lambda _M (\delta );$
$$\max \{\lambda_M (0+),\lambda_M (0)\} \leq  u_M (0).$$}

Следующая теорема 1 является обобщением теоремы 1 Л. П. Власова \cite{Vlas1} и леммы 2 А. В. Маринова \cite{Marin1}.
\vskip7pt
{\bf Теорема 1} [s12]. {\it
Пусть в полном геодезическом пространстве $X$, через каждые две
различные точки которого можно провести единственную прямую, каждый замкнутый шар выпуклый, $M \in \Sigma (X)$ и кроме того, выполняется одно из следующих условий.
\vskip7pt
$(i)\, \,$ Множество $M$ --- $P$-связно и оператор $P_M$ непрерывен в слабом смысле.
\vskip7pt
$(ii)\, \,$ Для каждого $x \in X$
$$\lim \limits_{\delta \to 0+} \lambda (P (x, M, \delta)) = 0.$$
Тогда множество $M$ --- $\stackrel{\circ}{B}$-связно.}
\vskip7pt
{\bf Доказательство теоремы 1}.
\vskip7pt
Доказательство пункта $(i)$ аналогично
доказательству теоремы 1 Л. П. Власова \cite{Vlas1}. Отличие состоит в содержании понятия деления отрезка в данном отношении, в своевременном упоминании условия выпуклости каждого замкнутого шара и обозначениях.
 Пусть, напротив, найдется такой открытый шар $B_0 = B(x,R_0)$, что множество $M\cap B_0$ --- несвязно. Значит,
$$M\cap B_0 = A\cup C,$$
где $A, \, C$ непустые, замкнутые в $B_0$ и непересекающиеся множества. Выберем $R_1$ так, чтобы
$$\max \{|xA|, |xC|\} < R_1 < R_0.$$
Тогда замкнутый шар $B_1 = B[x,R_1]$ обладает тем
свойством, что множество $B_1\cap M$ распадается  на непустые, замкнутые и
непересекающиеся множества
$$A\cap B_1,\quad C\cap B_1.$$
Покажем, что найдется замкнутый шар
$$B = B[y,r]\subset B_0$$
такой, что
$$|(B\cap A)(B\cap C)| > 0.$$
Допустим, что это не так. Пусть $R_i > 0$ при $i \geq 2$
такие, что
$$\sum \limits_{i=2}^{\infty }{R_i} < R_0 - R_1.$$
Построим последовательность замкнутых шаров $(B_n)_{n \in \mathbb{N}}$
индуктивно. Шар $B_1$ уже построен. Пусть
$$B_n \cap A \neq \emptyset, \quad  B_n\cap C \neq \emptyset$$
и радиус $B_n$ равен $R_n$. Тогда, по предположению
$$|(B_n\cap A)(B_n\cap C)| = 0.$$
Значит, найдутся точки
$$a_n \in B_n\cap A,\quad c_n \in B_n\cap C$$
такие, что
$$|a_nc_n| < 2R_{n+1}.$$
В качестве $B_{n+1}$ возьмем шар $B[q_{n+1},R_{n+1}]$, где
$$q_{n+1} = \omega _{1/2}(a_n,c_n).$$
Из условия выпуклости каждого замкнутого шара и неравенства треугольника следует,
что для каждого $n \in \mathbb{N}$
$$|q_nq_{n+p}| \leq |q_nq_{n+1}| + \ldots + |q_{n+p-1}q_{n+p}| \leq
\sum \limits_{k=n}^{n+p-1}{R_k} \rightarrow 0$$
при $p\rightarrow \infty, \, n\rightarrow \infty $.
Значит, $(q_n)_{n \in \mathbb{N}}$ --- фундаментальная последовательность и найдется такая точка
$q \in X$, что
$$\lim \limits_{n \rightarrow \infty} q_n= q.$$
Кроме того, $x = q_1$ и
$$|xq| = |q_1q| \leq |q_1q_n| + |q_nq| \leq \sum \limits_{k=1}^{n-1}{|q_kq_{k+1}|} + |q_nq|
 \leq \sum \limits_{k=1}^{\infty }{R_k} + |q_nq|.$$
Но
$$\sum \limits_{k=1}^{\infty }{R_k} < R_0, \quad \lim \limits_{n \rightarrow \infty} |q_nq| = 0.$$
Следовательно, $q \in B_0$. Кроме того,
$$\lim \limits_{n \rightarrow \infty} R_n = 0.$$
Значит,
$$\lim \limits_{n \rightarrow \infty} a_n= q, \quad  \lim \limits_{n \rightarrow \infty} c_n = q$$
и $q \in A\cap C\cap B_0.$
Получили противоречие. Таким образом, вышеуказанный шар существует. Тогда
$$B\cap M = A'\cup C',$$
где
$$A' = A\cap B, \quad C' = C\cap B$$
--- непустые замкнутые множества такие, что $|A'C'| > 0$.
\vskip7pt
$(i)\, \,$ Следовательно,
$$P_M (y) \subset A'\cup C'.$$
Пусть, для определенности,  $P_M (y) \subset A'$. Тогда
$P_M (y) \not \subset C'$, поскольку $|A'C'| > 0$ и множество $P_M (y)$ связно. Пусть $c \in C'$.
Тогда
$$P_M (\omega_{\lambda }(y,c)) \subset B$$
при $0 \leq \lambda \leq 1$, поскольку
$$|\omega _{\lambda }(y,c)M| \leq |\omega _{\lambda }(y,c)c|  \leq R - |y\omega _{\lambda }(y,c)| =
|\omega _{\lambda }(y,c)(X\backslash B)|.$$
Значит,
$$P_M (\omega_{\lambda }(y,c)) \subset A'\cup C'.$$
Выберем
$$\lambda_0 = \sup \{\lambda : 0 \leq \lambda \leq 1, \,
P_M (\omega_{\lambda }(y,c)) \subset A' \}.$$
Пусть последовательность $(\lambda _n)_{n \in \mathbb{N}}$ такая, что
$$\lim \limits_{n \rightarrow \infty} \lambda _n = \lambda _0$$
и для каждого $n \in \mathbb{N}$ справедливо неравенство $\lambda _n > \lambda _0$.
Тогда
$$P_M (\omega_{\lambda_n }(y,c)) \subset C',\quad \lim \limits_{n \rightarrow \infty} \omega_{\lambda_n }(y,c)= \omega_{\lambda_0 }(y,c),$$
$$\lim \limits_{n \rightarrow \infty} |P_M (\omega_{\lambda_n }(y,c))P_{M}(\omega_{\lambda_0 }(y,c))| = 0.$$
Следовательно,
$$|A'P_M (\omega_{\lambda_{0} }(y,c)| = 0.$$
Но возможны только два случая:
$$P_M (\omega_{\lambda_{0} }(y,c)) \subset A'$$
или
$$P_{M}(\omega_{\lambda_{0} }(y,c))
 \subset C'.$$
Следовательно, $|A'C'| = 0$. Получили противоречие.
\vskip7pt
$(ii)\, \,$ Если множество $M$ не является $\stackrel{\circ}{B}$-связным, то, как было доказано выше, найдутся $x \in X$ и
$\delta > 0$  такие, что
$$P (x, M, \delta) = A'\cup C',\quad |A'C'| > 0,\quad \max \{|xA'|, |xC'|\} < |xM| + \delta.$$
Нетрудно заметить, что в случае, когда
$|xA'| = |xC'|$, верно неравенство
$$|A'C'| \leq \lambda (P (x, M, \delta))$$
при $\delta > 0$. Но это приводит к противоречию с условием $(ii)$. Для завершения
доказательства осталось привести случай, когда $|xA'| < |xC'|$ к предыдущему случаю.
Пусть
$$0 < 2t < |xM| + \delta - |xC'|,\quad y \in P (x, C', t).$$
Найдется такая точка $v \in [x,y]$,
что $|vA'| = |vC'|$, поскольку непрерывная функция $f(z) = |zA'| - |zC'|$ меняет знак на сегменте $[x,y]$.
Тогда
$$|xv| = |xy| - |vy| \leq |xC'| + t - |vC'|$$
и $B[v,|vC'| + t] \subset B[x,|xM| + \delta ]$. Следовательно,
$$P (v, M, t) \subset A'\cup C',\quad |A'C'| \leq \lambda (P (v, M, t)).$$
Вместе с условием $(ii)$ это приводит к противоречию. Теорема 1 доказана.
\vskip7pt
Следующая теорема 2 является обобщением теоремы 4.2 Л. П. Власова \cite{Vlas}.
\vskip7pt
{\bf Теорема 2} [s12]. {\it
Пусть в полном геодезическом пространстве $X$, через каждые две различные точки которого можно провести единственную прямую, выполнены условия $(A_5)$, $(A_6)$, $(A'_7)$, $(A_9)$, $M$ --- $P$-связное множество. Тогда  $M$ --- $B$-связное множество.}
\vskip7pt
{\bf Доказательство теоремы 2}
\vskip7pt
Доказательство нетрудно получить из доказательства теоремы 4.2
Л. П. Власова \cite{Vlas} заменой ссылок на лемму 1.1 и теорему 1.1 С. Б. Стечкина \cite{Vlas} на их обобщения: лемму 2.4.2 и теорему 2.4.2 соответственно, а также заменой
обозначений.
\vskip7pt
Следующая лемма 1 является обобщением леммы 3 А. В. Маринова \cite{Marin1}.
\vskip7pt
{\bf Лемма 1} [s12]. {\it
Пусть в полном геодезическом пространстве $X$, через каждые две различные точки которого можно провести единственную прямую, выполнены условия $(A_5)$, $(A_6)$, $(A'_7)$, $(A_9)$, $M$ ---  замкнутое множество. Тогда следующие утверждения эквивалентны.
\vskip7pt
$(i)\, \,$ Множество $M$ --- $\stackrel{\circ}{B}$-связно.
\vskip7pt
$(ii)\, \,$ $P (x, M, \delta)$ --- связно для каждого $x \in X$ и каждого $\delta > 0$.
\vskip7pt
$(iii)\, \,$ Для каждого $\delta > 0$
$$\lambda_M(\delta ) = 0.$$
\vskip7pt
$(iv)\, \,$ Для каждого $x \in X$
$$\lim \limits_{\delta \to 0+}{\lambda (P (x, M, \delta))} = 0.$$}

{\bf Доказательство леммы 1}.
\vskip7pt
Импликации $(ii) \Rightarrow (iii) \Rightarrow (iv)\,$ очевидны.
\vskip7pt
Импликация $(iv) \Rightarrow (i)\,$ доказана в пункте $(ii)$ теоремы 1.
\vskip7pt
Импликация $(i) \Rightarrow (ii)\,$ отличается от доказательства в лемме 3 А. В. Маринова \cite{Marin1} только тем, что ссылка на доказательство теоремы 4.2  Л. П. Власова должна быть заменена ссылкой на видоизмененное доказательство согласно указанию в доказательстве теоремы 2.
\vskip7pt
Следующая лемма 2 является обобщением неравенства $(3)$ А. В. Маринова \cite{Marin1}.
\vskip7pt
{\bf Лемма 2} [s12]. {\it
Пусть $X$ геодезическое пространство, через каждые две
различные точки которого можно провести единственную прямую. Тогда для каждой точки
$$(x,M,\delta) \subset X\times \Sigma (X)\times \mathbb{R}_{+},$$
где $\delta   > 0$, выполняется неравенство
$$\bar {\omega}(t) < \hat {\omega}(t) + t$$
при $0 < t < \delta $.}
\vskip7pt
{\bf Доказательство леммы 2}.
\vskip7pt
Поставим в соответствие каждой точке $z \in P (x, M, t)$ точку
$$v(z) = S(x, |xM| +t)\cap \{\omega _l (x, z) : l \geq 0 \}.$$
Положим
$$W = \{v(z) : z \in P (x, M, t) \} \cup (M\backslash P (x, M, t)).$$
Тогда
$$\alpha (M, W) = \beta (P (x, M, t + \delta),P (x, W, \delta)) = t.$$
Следовательно,
$$\beta (P (x, M, t + \delta), P (x, M, \delta)) \leq \beta (P (x, M, t + \delta), P (x, W, \delta)) +$$
$$\beta (P (x, W, \delta), P (x, M, \delta)) \leq t + \hat {\omega}(t).$$
Положим теперь
$$W_1 = \{\omega_l (x,z) : z \in P (x, M, \delta) \backslash P (x, M, \delta - t), \,  l |xz| = |xz| + t \},$$
$$K = [M\backslash (P (x, M, \delta) \backslash P (x, M, \delta - t))]\cup W_1.$$
Тогда
$$W_1\cap P (x, M, \delta) = \emptyset,\quad P (x, K, \delta) = P (x, M, \delta - t),\quad
\alpha (M,K) \leq t.$$
Отсюда получим
$$\beta (P (x, M, \delta), P (x, M, \delta - t)) \leq \beta (P (x, M, \delta), P (x, K, \delta)) +$$
$$\beta (P (x, K, \delta), P (x, M, \delta - t)) =\beta (P (x, M, \delta), P (x, K, \delta)) \leq \hat {\omega}(t).$$
Лемма 2 доказана.
\vskip7pt
{\bf Замечание 1.} {\it
Из леммы $2$ и неравенства $(2)$ А. В. Маринова {\rm \cite{Marin1}:}
$$\Omega (t) \leq \bar {\omega}(t) + t$$
для $(x,M,\delta) \subset X\times \Sigma (X)\times \mathbb{R}_{+}$, где $\delta   > 0$,
в геодезическом пространстве $X$, через каждые две
различные точки которого можно провести единственную прямую,
следуют соотношения, аналогичные соотношениям $(4)$ в {\rm \cite{Marin1}}
$$u_M (\delta ) \leq \hat {u}_M (\delta ) = \bar {u}_M (\delta )= U_M (\delta )$$
для каждого $\delta > 0.$}
\vskip7pt
Следующая теорема 3 является обобщением теоремы 2  А. В. Маринова \cite{Marin1}.
\vskip7pt
{\bf Теорема 3} [s12]. {\it
Пусть в геодезическом пространстве $X$, через каждые две различные точки которого можно провести единственную прямую, выполнены условия $(A_5)$, $(A_6)$, $(A'_7)$, $(A_9)$, $M \in \Sigma (X), \, \delta > 0$.
Тогда
$$U_M (\delta ) \leq \lambda (P (x, M, \delta - 0)).$$}

{\bf Доказательство теоремы 3}
\vskip7pt
Доказательство нетрудно получить из доказательства теоремы 2
А. В. Маринова \cite{Marin1} заменой ссылки на лемму 4 в \cite{Marin1} ссылкой на более общую лемму 2.4.2 данной работы и заменой обозначений.
\vskip7pt
Следующая теорема 4 является обобщением теоремы 4  А. В. Маринова \cite{Marin1}.
\vskip7pt
{\bf Теорема 4} [s12]. {\it
Пусть в геодезическом пространстве $X$, через каждые две различные точки которого можно провести единственную прямую, выполнены условия $(A_5)$, $(A_6)$, $(A'_7)$, $(A_9)$, $M \in \Sigma (X)$.
Тогда нижеследующие условия $(a - e)$ эквивалентны. Если же, кроме того, $M$ замкнутое множество и  $X$ полное пространство, то условия $(a - h)$
эквивалентны.
\vskip7pt
$(a)\, \,$ $U_M (\delta ) = 0$ для каждого $\delta > 0$.
\vskip7pt
$(b)\, \,$ $\lim \limits_{\overline{\delta \to 0+}}{U_M (\delta )} = 0$.
\vskip7pt
$(c)\, \,$ $u_M (\delta ) = 0$ для каждого $\delta > 0$.
\vskip7pt
$(d)\, \,$ $\lim \limits_{\overline{\delta \to 0+}}{u_M (\delta )} = 0$.
\vskip7pt
$(e)\, \,$ $\lambda_M (\delta) = 0$ для каждого $\delta > 0$.
\vskip7pt
$(f)\, \,$  $\lim \limits_{\delta \to 0+}{\lambda (P (x, M, \delta))} = 0$ для каждого $x \in X$.
\vskip7pt
$(g)\, \,$ $M$ --- $\stackrel{\circ}{B}$-связное множество.
\vskip7pt
$(h)\, \,$ $P (x, M, \delta)$ --- связное множество для каждого $\delta > 0$ и для каждого $x \in X$.}
\vskip7pt
{\bf Доказательство теоремы 4}.
\vskip7pt
Из предложения 1 и теоремы 3 следует
справедливость импликаций $(d) \Rightarrow (e) \Rightarrow (a).$ Но из первых четырех
условий первое самое сильное, а последнее самое слабое. Значит, условия
$(a - e)$ эквивалентны. Эквивалентность условий $(e - h)$ при дополнительном
предположении полноты пространства $X$ и замкнутости множества $M$ следует теперь из теоремы 3.
\vskip7pt
Из теорем 2, 3 получаем
\vskip7pt
{\bf Следствие 1} [s12]. {\it
Пусть в полном геодезическом пространстве $X$, через каждые две различные точки которого можно провести единственную прямую, выполнены условия $(A_5)$, $(A_6)$, $(A'_7)$, $(A_9)$, $M$ --- чебышевское множество. Тогда оператор $P$ метрического $\delta $-проектиро\-вания на
множество $M$ непрерывен по совокупности переменных, то есть для каждого $\delta > 0$
$$U_M(\delta ) = 0.$$}

Из предложения 1  и теоремы 3 получаем
\vskip7pt
{\bf Следствие 2} [s12]. {\it
Пусть в геодезическом пространстве $X$, через каждые две различные точки которого можно провести единственную прямую, выполнены условия $(A_5)$, $(A_6)$, $(A'_7)$, $(A_9)$, $M$ --- множество существования и оператор $P_M$ метрического проектирования на
множество $M$ непрерывен по Хаусдорфу, то есть $u_M (0) = 0$. Тогда для каждого $\delta > 0$
$$U_M (\delta ) = 0.$$}

\section{Непрерывность метрической $\delta$-проекции на выпуклое множество в геодезическом пространстве}
\vskip20pt

Пусть метрическое пространство $(X,\rho )$ с выделенным семейством $S$ сегментов, называемых хордами,  является $U$-множеством, удовлетворяющем следующему условию $(C)$.
\vskip7pt
$(C)\quad$ Для всех $ p, \, x, \, y \in X$ справедливо неравенство
$$2 |p\omega_{1/2}(x, y)| \leq |px| + |py|,$$
где  $\omega_{1/2} (x,y)$ --- середина хорды $[x,y]$ с концами $x, \, y$.
\vskip7pt
Известно (см. \cite [с. 63]{Bus2}, \cite [с. 304]{Bus}), что в прямых хордовых пространствах неположительной кривизны это условие выполнено. Аналогичные локальные условия налагались на метрическое пространство в \cite{Alex}.

Приведем два элементарных свойства в $U$-множестве $X$, удовлетворяющем условию $(C)$.
\vskip7pt
$1.\, \,$ Для каждого $\lambda \in [0,1]$ и для всех точек $p$, $x$, $y$ из $U$-множества $X$, удовлетворяющего условию $(C)$, справедливо неравенство
$$|p\omega_{\lambda} (x,y)| \leq (1 - \lambda )|px| + \lambda |py|, \eqno (1)$$
где точка $\omega_{\lambda} (x,y) \in [x,y]$ такая, что
$$|x\omega_{\lambda}(x,y)| = \lambda |xy|.$$

$2.\, \,$ Каждый замкнутый (открытый) шар $U$-множества $X$, удовлетворяющего условию $(C)$, является выпуклым множеством.
\vskip7pt
В дальнейшем используем следующее обозначение
$$\beta (P (x, M, \delta \pm 0),F) = \lim \limits_{t \to 0+}{\beta (P (x, M, \delta \pm t), F)},$$
где $F, \, M \in \Sigma (X)$ \cite{Marin}, причем бесконечные
значения для отклонения и расстояния по Хаусдорфу между непустыми множествами из $X$ не исключаются.
Нам также понадобится следующая лемма 1 из \cite{Marin}.
\vskip7pt
{\bf Лемма 1} \cite{Marin}. {\it
Пусть $(X,\rho )$ метрическое пространство,
$$x, \, y \in X,\quad F, \, M, \, W \in \Sigma (X),\quad \delta \geq 0.$$
Тогда
$$\beta (P (y, W, \delta) ,F) \leq \beta (W,M) + \beta (P (x, M, \delta_1 + 0), F),$$
$$\beta (F, P (y, W, \delta)) \leq \beta (M, W) + \beta (F, P (x, M, \delta_2 - 0)),$$
где
$$\delta_1 = \delta + |xy| + |yW| + \beta (W, M) - |xM|,$$
$$\delta_2 =\delta + |yW| - |xy| - |xM| - \beta (M, W).$$
Причем предполагается, что $\delta_2 > 0$.}
\vskip7pt
Следующая лемма обобщает леммы 6, 7, 8 из \cite{Marin2}.
\vskip7pt
{\bf Лемма 2} [s13]. {\it
Пусть $M$ выпуклое подмножество $U$-множества $X$, удовлетворяющего условию $(C)$,
$$x \in X,\quad 0 < t < \varepsilon < \delta,\quad 0 < \varepsilon ' < \delta ',\quad
\varepsilon ' \leq \varepsilon,\quad \delta ' \leq \delta,$$
$$(t \geq 0 \, \, \, \mbox {при}\, \, \, P_M (x) \neq \emptyset),\quad F \subset  P (x, M, \delta),
\quad G  \subset  P (x, M, t).$$
Тогда имеют место следующие неравенства.
$$\dfrac{\beta (F, P (x, M, \varepsilon ))}{\delta - \varepsilon } \leq \dfrac{\beta (F, P (x, M, t))}{\delta - t}; \leqno (a)$$
$$\dfrac{\beta (P (x, M, \delta ),G)}{\delta - t} \leq \dfrac{\beta ( P (x, M, \varepsilon ),G)}{\varepsilon - t}; \leqno (b)$$
$$\dfrac{\beta (P (x, M, \delta ),P (x, M, \varepsilon ))}{\delta - \varepsilon } \leq \dfrac{\beta (P(x, M, \delta'), P (x, M, \varepsilon'))}{\delta' - \varepsilon' }. \leqno (c)$$}

{\bf Доказательство леммы 2}.
\vskip7pt
$(a)\, \,$ Пусть
$$\beta (F,P (x, M, \varepsilon )) > 0$$
и число $s > 0$ достаточно мало. Выберем точку $z \in F \backslash P (x, M, \varepsilon )$ так, что
$$|zP (x, M, \varepsilon )| \geq \beta (F,P (x, M, \varepsilon )) - s. \eqno (3)$$
Затем выберем точку $p \in P (x, M, t)$ так, что
$$|zp| \leq |zP (x, M, t)| + s. \eqno (4)$$
Введем обозначения
$$a = [x,z]\cap S(x,|xM| + t),\, \, \, b = [x,z]\cap S(x,|xM| + \varepsilon ),\, \, \, c = [p,z]\cap S(x,|xM| + \varepsilon ).$$
Из $(3)$ следует, что
$$|zc| \geq \beta (F,P (x, M, \varepsilon )) - s. \eqno (5)$$
А из $(1)$ следует неравенство
$$|xc| \leq |px| \dfrac{|cz|}{|pz|} + |xz| \dfrac{|cp|}{|pz|}.$$
Но $|xc| = |xb| = |xz| - |bz|$, $|ax| \geq |xp|$, $|pc| = |pz| - |cz|$. Поэтому
$$|xz| - |bz| \leq |ax| \dfrac{|cz|}{|pz|} + |xz| \dfrac{|pz| - |cz|}{|pz|}.$$
Откуда следует, что
$$|cz| \leq |pz| \dfrac{|bz|}{|az|}.$$
Кроме того, $\delta - t > \delta - \varepsilon $ и величина
$$\tau = \delta - \varepsilon - |bz| = \delta - t - |az|$$
неотрицательная. Следовательно,
$$|cz| \leq |pz| \dfrac{|bz|}{|az|} \leq |pz|\frac{\delta - \varepsilon - \tau }{\delta - t - \tau } \leq |pz|\dfrac{\delta - \varepsilon}{\delta - t}.$$
Из этих неравенств, а также из $(4)$, $(5)$ следует
$$\beta (F,P (x, M, \varepsilon )) -s \leq cz \leq (|zP (x, M, t)| + s) \dfrac{\delta - \varepsilon}{\delta - t} \leq $$
$$(\beta (F,P (x, M, t)) + s)\dfrac{\delta - \varepsilon}{\delta - t}$$
Переходя к пределу при $s \to 0+$, получим неравенство $(a)$.
\vskip7pt
$(b)\, \,$ Пусть $\beta (P (x, M, \delta ),G) > 0$ и число $s > 0$ достаточно мало. Выберем точку
$$z \in P (x, M, \delta ) \backslash G$$
так, что
$$|zG| \geq \beta (P (x, M, \delta ),G) - s. \eqno (6)$$
Затем выберем точку $p \in G$ так, что
$$|zp| \leq |zG| + s. \eqno (7)$$
Тогда
$$\beta (P (x, M, \delta ),G) -s \leq |zp| \leq \beta (P (x, M, \delta ),G) + s. \eqno (8)$$
Если для каждого достаточно малого $s$ точка $z$ принадлежит $P (x, M, \varepsilon )$, то
$$\beta (P (x, M, \delta ),G) = \beta (P (x, M, \varepsilon ),G)$$
и имеет место неравенство $(b)$. Предположим, что
$$z \in P (x, M, \delta ) \backslash P (x, M, \varepsilon )$$
и введем обозначения
$$b = [x,z]\cap S(x,|xM| + \varepsilon ),\quad c = [p,z]\cap S(x,|xM| + \varepsilon ). \eqno (9)$$
Из $(1)$ следует, что
 $$|xc| \leq |px| \dfrac{|cz|}{|pz|} + |xz| \dfrac{|cp|}{|pz|} =  |px| \dfrac{|pz| - |pc|}{|pz|} + |xz| \dfrac{|cp|}{|pz|} = |px| + |pc| \dfrac{|xz| - |px|}{|pz|}.$$
Но $|xc| = |xb|$. Следовательно,
$$|pz| \leq |pc| \dfrac{|xz| - |xp|}{|xb| - |px|} \leq |pc| \dfrac{\delta - t }{\varepsilon - t}. \eqno (10)$$
Кроме того, из $(7)$ и $(9)$ следует, что
$$|cp| = |pz| - |cz| \leq |zG| + s - |cz| \leq |cG| + s
\leq \beta (P (x, M, \varepsilon ),G) + s. \eqno (11)$$
Тогда из $(8)$, $(10)$, $(11)$ получим
$$\beta (P (x, M, \delta ),G) -s \leq |zp| \leq |pc| \dfrac{\delta - t }{\varepsilon - t}
 \leq (\beta (P (x, M, \delta ),G) + s) \dfrac{\delta - t }{\varepsilon - t}.$$
Переходя к пределу при $s \to 0+$, получим неравенство $(b)$.
\vskip7pt
$(c)\, \,$ Достаточно применить $(b)$ при
$$G = P (x, M, \varepsilon' ),\quad t = \varepsilon ',\quad \varepsilon = \delta ',$$
а затем $(a)$ при $F =  P (x, M, \delta )$. Лемма 2 доказана.
\vskip7pt
{\bf Следствие 1} [s13]. {\it
Пусть $M$ выпуклое подмножество $U$-множества $X$, удовлетворяющего условию $(C)$, $x \in X$. Тогда
имеют место следующие неравенства.
$$\dfrac{\beta (P (x, M, \delta ),P (x, M, 0+))}{\delta } \leq \dfrac{\beta (P (x, M, \delta' ),P (x, M,  0+))}{\delta' }, \leqno (a)$$
где $\delta > \delta' > 0$;
$$\beta (P (x, M, \delta ),P (x, M, \varepsilon )) \leq  \dfrac{\beta (P (x, M, \delta ),P (x, M,  0+))}{\delta }(\delta - \varepsilon), \leqno (b)$$
где $\delta \geq \varepsilon > 0$;
$$\alpha (P (x, M, \delta ),P (x, M, \varepsilon )) \leq
  \dfrac{\beta (P (x, M, \mu), P (x, M, 0+)) }{\mu }|\delta - \varepsilon| \leq \leqno (c)$$
$$\dfrac{\beta (P (x, M, \delta ),P (x, M, 0+)) }{\delta }|\delta - \varepsilon|,$$
где $\delta > 0$, $\mu = \max \{\delta,\varepsilon\}$;
$$\beta (P (x, M, \delta ),P (x, M, \varepsilon )) \leq (\dfrac{2 |xM|}{\delta } + 1)(\delta - \varepsilon ), \leqno (d)$$
где $\delta \geq \varepsilon \geq 0$. При $P_M (x) \neq \emptyset $ выражение $P (x, M,  0+)$ в этих
неравенствах можно заменить на $P_M (x)$.}
\vskip7pt
{\bf Доказательство следствия 1}.
\vskip7pt
$(a)\, \,$ Достаточно в неравенстве $(c)$ леммы 2 перейти к пределам при
$\varepsilon' \to 0$, $\varepsilon \to 0$.
\vskip7pt
$(b)\, \,$ В том же неравенстве нужно положить $\delta ' = \delta $ и перейти к пределу при $\varepsilon ' \to 0$.
\vskip7pt
$(c)\, \,$ Следует из $(b)$.
\vskip7pt
$(d)\, \,$ Следует из $(b)$ и простого неравенства
$$\beta (P (x, M, \delta ),P (x, M, 0+)) \leq 2 |xM| + \delta $$
\cite [(3.10)]{Marin2}. Следствие 1 доказано.
\vskip7pt
Следующая  лемма 3 является обобщением леммы 5 из \cite{Marin2} и имеет аналогичное доказательство, если использовать утверждение $(a)$ леммы 2.
\vskip7pt
{\bf Лемма 3} [s13]. {\it
Пусть $M$ выпуклое подмножество $U$-множества $X$, удовлетворяющего условию $(C)$,
$$x \in X,\quad \delta >0,\quad F \in \Sigma (X).$$
Тогда имеют место следующие четыре равенства.
$$\beta (P (x, M, \delta \pm 0),F) = \beta (P (x, M, \delta ),F),$$
$$\beta (F,P (x, M, \delta \pm 0)) = \beta (F,P (x, M, \delta )).$$}

Следующая лемма аналогична лемме 10 из \cite{Marin2}.
\vskip7pt
{\bf Лемма 4} [s13]. {\it
Пусть $M$ выпуклое подмножество $U$-множества $X$, удовлетворяющего условию $(C)$,
$$x, \, y \in X,\quad \delta \geq 0,\quad \varepsilon \geq 0,\quad W \in  \Sigma (X),$$
$$\mu =\max \{\varepsilon,\delta\},\quad \lambda = \mu + 2 |xy| + 2\alpha (M,W).$$ Тогда
\vskip7pt
$(a)\, \,$ при $\lambda > 0$ верно неравенство
$$\beta (P (y, W, \delta ), P (x, M, \varepsilon )) \leq \beta (W,M) +$$
$$\dfrac{\beta (P (x, M, \lambda ),P (x, M, 0+))}{\lambda }(|\varepsilon - \delta | + 2 |xy| + 2\alpha (M,W));$$

$(b)\, \,$ при
$$|\varepsilon - \delta | + 2 |xy| + 2 \alpha (M,W) < \mu $$
верно неравенство
$$\alpha (P (y, W, \delta ),P (x, M, \varepsilon )) \leq \alpha (W,M) +$$
$$ \dfrac{\beta (P (x, M, \mu ), P (x, M, 0+))}{\mu }(|\varepsilon - \delta | + 2 |xy| + 2\alpha (M,W)).$$}

{\bf Доказательство леммы 4}.
\vskip7pt
$(a)\, \,$ Из простых оценок \cite [(1.5)]{Marin} получим
$$\delta_1 = \delta + |xy| + |yW| + \beta (W,M) - |xM| \leq \delta + 2 |xy| + 2 \alpha (M,W) \leq \lambda.$$
Тогда из лемм 1, 3 и утверждения $(b)$ следствия 1 получим.
$$\beta (P (y, W, \delta ),P (x, M, \varepsilon )) \leq \beta (W,M) +
 \beta (P (x, M, \lambda ),P (x, M, \varepsilon )) \leq$$
$$\beta (W,M) +
 \dfrac{\beta (P (x, M, \lambda ),P (x, M, 0+))}{\lambda }(\lambda - \varepsilon )  \leq \beta (W,M) +$$
$$\dfrac{\beta (P (x, M, \lambda ),P (x, M, 0+))}{\lambda }(|\varepsilon - \delta | + 2 |xy| + 2 \alpha (M,W)).$$

$(b)\, \,$ Из простых оценок (\cite{Marin}, (1.5)) получим
$$|\mu - \delta_2| = |\mu - \delta - |yW| + |xy| + |xM| + \beta (M,W)| \leq  |\mu - \delta | + 2 |xy| + 2\alpha (M,W) \leq $$
 $$|\varepsilon - \delta | + 2 |xy| + 2 \alpha (M,W) < \mu.$$
Следовательно, $\delta_2  > 0$ и можно применить
леммы 1, 3 и утверждение $(b)$ следствия 1.
$$\beta (P (x, M, \varepsilon ),P (y, W, \delta )) \leq \beta (M,W) +
 \beta (P (x, M, \varepsilon ),P (x, M, \delta_2 - 0)) \leq$$
$$\beta (M,W) +
 \beta (P (x, M, \mu ),P (x, M, \delta_2)) \leq \beta (M,W) +$$
$$ \dfrac{\beta (P (x, M, \mu ),P (x, M, 0+))}{\mu } (|\varepsilon - \delta | + 2 |xy| + 2 \alpha (M,W)).$$
Лемма 4 доказана.
\vskip7pt
Следующая теорема обобщает теорему 4 из \cite{Marin2}.
\vskip7pt
{\bf Теорема 1} [s13]. {\it
Пусть $M$ --- выпуклое подмножество $U$-множества $X$, удовлетворяющего условию $(C)$,
$$x, \, y \in X,\quad W \in \Sigma (X),\quad \mu =\max \{\varepsilon,\delta\} > 0.$$
Тогда верно неравенство
$$\alpha (P (y, W, \delta ),P (x, M, \varepsilon )) \leq \alpha (W,M) +$$
$$\left(\dfrac{2 |xM|}{\mu } + 2\right)(|\varepsilon - \delta | + 2 |xy| + 2\alpha (M,W)).$$}

{\bf Доказательство теоремы 1}.
\vskip7pt
Из утверждения $(a)$ леммы 4 и утверждения $(d)$ следствия 1 получим
$$\beta (P (y, W, \delta ),P (x, M, \varepsilon )) \leq \beta (W,M) +
 \dfrac{\beta (P (x, M, \lambda ),P (x, M, 0+))}{\lambda } (|\varepsilon - \delta | +$$
$$2 |xy| + 2 \alpha (M,W)) \leq \alpha (W,M) +
\left(\dfrac{2 |xM|}{\mu } + 1\right) (|\varepsilon - \delta | + 2 |xy| + 2 \alpha (M,W)).$$
А при
$$|\varepsilon - \delta | + 2 |xy| + 2 \alpha (M,W) < \mu $$
из утверждения $(a)$ леммы 4 и утверждения $(d)$ следствия 1 получим
$$\alpha (P (y, W, \delta ),P (x, M, \varepsilon )) \leq \alpha (W,M) +$$
$$\left(\dfrac{2 |xM|}{\mu } + 1\right) (|\varepsilon - \delta | + 2 |xy| + 2 \alpha (M,W)).$$
Пусть теперь
$$|\varepsilon - \delta | + 2 |xy| + 2 \alpha (M,W) \geq  \mu.$$
Тогда
$$\beta (P (x, M, \varepsilon ),P (y, W, \delta )) \leq \beta (P (x, M, \varepsilon ),x) + |xy| + |yP (y, W, \delta )| \leq$$
$$|xM| + \varepsilon + |xy| + |yW| \leq 2 |xM| + ||xM| - |yW|| + \varepsilon + |xy| \leq$$
$$2 |xM| + 2 |xy| + \alpha (M,W) +  \varepsilon \leq
2 |xM| + \mu + |\varepsilon - \delta | + 2 |xy| + 2 \alpha (M,W)  \leq$$
$$\left(\dfrac{2 |xM|}{\mu } + 2\right) (|\varepsilon - \delta | + 2 |xy| + 2 \alpha (M,W)).$$
Из полученных неравенств следует теорема 1.
\vskip7pt
Следующая теорема аналогична теореме 5 из \cite{Marin2}.
\vskip7pt
{\bf Теорема 2} [s13]. {\it
Пусть $M$, $W$ --- выпуклые подмножества $U$-множества $X$, удовлетворяющего условию $(C)$,
$$x, \, y \in X,\quad \lambda =\max \{\varepsilon,\delta\} + 2 |xy| + 2\alpha (M,W).$$
Тогда при $\lambda > 0$ выполняется неравенство
$$\alpha (P (y, W, \delta ),P (x, M, \varepsilon )) \leq \alpha (W,M) +$$
$$\left(\dfrac{2 \min \{|xM|,|yW|\}}{\lambda } + 1\right)(|\varepsilon - \delta | + 2 |xy| + 2\alpha (M,W)) \eqno(2).$$}

{\bf Доказательство теоремы 2}.
\vskip7pt
Обозначим правую часть в $(2)$ через $A$. Рассмотрим три случая для оценки
выражения
$$\beta (P (y, W, \delta ),P (x, M, \varepsilon )).$$
$$\delta_1 = \delta + |xy| + |yW| + \beta (N,M) - |xM| = 0. \leqno 1.$$
Тогда из неравенства
$$|xM| - |yW| \leq |xy| + \beta (W,M)$$
следует, что $\delta = 0$. Из леммы 1 получим
$$\beta (P (y, W, \delta ),P (x, M, \varepsilon )) \leq \beta (W,M) + \beta (P (x, M, 0+),P (x, M, \varepsilon )) \leq$$
$$\beta (W,M) + 2 |xM| \leq \beta (W,M) + 2 \min \{|xM|,|yW|\} + 2 ||xM| - |yW|| \leq$$
$$\beta (W,M) + 2 \min \{|xM|, |yW|\} + 2 |xy| + 2 \alpha (M,W) \leq A.$$
$$0 < \delta_1 \leq \varepsilon. \leqno 2.$$
Тогда из лемм 1, 3 получим
$$\beta (P (y, W, \delta ),P (x, M, \varepsilon )) \leq \beta (W,M) \leq A.$$
$$\delta_1 > \varepsilon. \leqno 3.$$
Сначала получим простые неравенства
$$\delta_1 < \lambda,\quad 2 |xM| + \delta_1 = \delta + |xy| + 2 \min \{|xM|,|yW|\} + ||xM| - |yW|| +$$ 
$$\beta (W,M) \leq \mu + 2 |xy| +  2 \min \{|xM|, |yW|\} + 2 \alpha (M,W) \leq$$ 
$$\lambda + 2 \min \{|xM|,|yW|\}, \quad\dfrac{\delta_{1} - \varepsilon }{\delta_{1}} \leq \dfrac{\lambda - \varepsilon }{\lambda }.$$
Используя эти неравенства, а также лемму 1 и утверждение $(d)$ следствия 1, получим
$$\beta (P (y, W, \delta ),P (x, M, \varepsilon )) \leq \beta (W,M) + \beta (P (x, M, \delta_1),P (x, M, \varepsilon ))  \leq$$
$$\beta (W,M) + (2 |xM| + \delta_1) \frac{\delta_1 - \varepsilon }{\delta_1} \leq$$
$$\beta (W,M) + (2\min \{|xM|, |yW|\} + \lambda) \dfrac{\lambda - \varepsilon }{\lambda } \leq A.$$
Осталось использовать симметричность правой части неравенства. Таким образом, теорема 2 доказана.

\section{Наилучшие $N$-сети ограниченных замкнутых выпуклых множеств в геодезическом пространстве}
\vskip20pt

Рассмотрим метрическое пространство $(X,\rho)$ и примем следующие обозначения и определения.
\vskip7pt
{\it Радиусом покрытия множества} $M \in B(X)$ $N$-сетью $S \in \Sigma_N (X)$  называется число
$\beta (M,S)$ \cite{Gark}. 
\vskip7pt
При $\Sigma = \Sigma_N (X)$  используем обозначение 
$$R_N (M) = \inf\{\beta (M, S) : S \in \Sigma_N (X)\}$$ 
для {\it наименьшего радиуса покрытия} множества $M \in B(X)$ $N$-сетями. 
\vskip7pt
$N$-сеть $S^{*}_N$ называется {\it наилучшей  $N$-сетью} множества $M \in B(X)$ \cite{Gark}, если 
$$\beta (M,S^{*}_N) = R_N (M).$$
В частности, наилучшая $1$-сеть множества $M \in B(X)$ есть чебышевский центр этого множества, а число $R (M) = R_1 (M)$ --- чебышевский радиус этого множества \cite{Gark}. 
\vskip7pt
$co (M)$ --- выпуклая оболочка множества $M$ (т.е. пересечение всех выпуклых множеств, содержащих множество $M$) \cite{Gark}.

 Кроме того, будем использовать  следующее условие  $(A_{10})$ в геодезическом пространстве $X$. 
\vskip7pt
$(A_{10})\quad$  Для каждого $x \in X$ и для любой точки $p$ из произвольного сегмента $[u,v] \subset X$ верно неравенство 
$$|pP_{[u,v]}(x)| \leq  |px|.$$

Примером полного метрического пространства, в котором выполнены условия $(A_4-A_{10})$,  является пространство Адамара \cite [c. 390]{Bur} (в частности, гильбертово пространство и пространство Лобачевского).
\vskip7pt
{\bf Лемма 1} [s14]. {\it Пусть $X$ --- полное геодезическое пространство, удовлетворяющее условиям $(A_4-A_7)$,  $(A_{10})$, $W \in \Sigma (X)$, $M$ --- непустое замкнутое ограниченное выпуклое множество в $X$. Тогда имеют место неравенства.
\vskip7pt
$(i)\, \,$ Для каждого $x \in M$
$$|xP_M (W)| \leq |xW|.$$ 

$$\beta (M,P_M (W)) \leq \beta (M,W). \leqno (ii)$$}

{\bf Доказательство леммы 1}. 
\vskip7pt
$(i)\, \,$  Из теоремы 2.4.1 и условий леммы 1 следует, что оператор $P_M$ является однозначным и непрерывным в пространстве $X$.
Выберем произвольно $x \in M$, $y \in W$. Тогда найдется единственная точка $y_1 = P_M (y)$. В силу выпуклости множества $M$ имеет место включение $[x,y_1] \subset M$. Тогда, используя условие $(A_{10})$, получим, что 
$$|xP_M (y)| = |xP_{[x,y_1]} (y)| \leq |xy|.$$ 
Следовательно, 
$$|xP_M (W)| \leq |xy|\, \, \mbox {и} \, \, |xP_M (W)| \leq |xW|.$$
\vskip7pt
$(ii)\, \,$ Неравенство $(ii)$ доказывается переходом к точной верхней грани по всем $x \in M$ сначала в правой части неравенства $(i)$, а затем в левой части полученного неравенства.
Лемма 1 доказана.
\vskip7pt
{\bf Следствие 1} [s14]. {\it Пусть $X$ --- полное геодезическое пространство, удовлетворяющее условиям $(A_4-A_7)$,  $(A_{10})$, $Z$ --- наилучшее аппроксимирующее множество из семейства $\Sigma $ 
для непустого замкнутого ограниченного выпуклого множества $M \subset X$ и $P_M (Z) \in \Sigma$. Тогда $P_M (Z)$ --- наилучшее аппроксимирующее множество для множества $M$.}
\vskip7pt
Следующая лемма 2 обобщает теорему 2 из \cite{Gark1}, теоремы 1, 2 из \cite{Belob1}.
\vskip7pt
{\bf Лемма 2} ([s14], [s21]). {\it Пусть $X$ --- полное геодезическое пространство, удовлетворяющее условиям $(A_4-A_7)$,  $(A_{10})$, $M \in B(X)$. Тогда чебышевский центр множества $M$ принадлежит замыканию выпуклой оболочки множества $M$, совпадает с чебышевскими центрами замыкания выпуклой оболочки множества $M$ и выпуклой оболочки множества $M$.}
\vskip7pt
{\bf Доказательство леммы 2}. 
\vskip7pt 
Существование и единственность чебышевского центра $Z (M)$ непустого ограниченного множества $M$ следует из следствия 2.6.1, а из условий $(A_4)$, $(A_5)$ следует, что замкнутые шары пространства $X$ выпуклые. Тогда 
$$M \subset co (M) \subset \overline{co (M)} \subset B[Z (M), R(M)]$$
Следовательно,
$$R (M) = R (co (M)) = R (\overline{co (M)}).$$
Из следствия 1 и единственности чебышевского центра следует, что 
$$Z (M) \subset \overline{co (M)}, \, Z (M) = Z (co (M)) = Z (\overline{co (M)}).$$
Лемма 2 доказана.
\vskip7pt 
{\bf Лемма 3} [s14]. {\it Пусть $X$ --- полное геодезическое пространство, удовлетворяющее условиям $(A_4-A_7)$,  $(A_{10})$,  $M$ --- непустое ограниченное замкнутое выпуклое множество пространства $X$, 
$(W_n)_{n \in \mathbb{N}}$ --- последовательность непустых ограниченных замкнутых множеств пространства $X$, сходящаяся в метрике Хаусдорфа к некоторому компакту $W \subset X$. Тогда 
$$\lim\limits_{n\to\infty}{\alpha (P_M (W_n), P_M (W))} = 0.$$}

{\bf Доказательство леммы 3}. 
\vskip7pt 
Из теоремы 2.4.1 и условий леммы 3 следует, что 
$M$ является аппроксимативно компактным, сильно чебышевским множеством. Кроме того, в силу  леммы 1 из \cite{Lisk} оператор $P_M$ является $\beta$-непрерывным на компакте $W \subset X$, поэтому $$\lim\limits_{n\to\infty}{\beta (P_M (W_n), P_M (W))} = 0.$$ 
Докажем, что 
$$\lim\limits_{n\to\infty}{\beta (P_M (W),P_M (W_n))} = 0$$ 
методом от противного. Пусть это утверждение неверно, тогда найдутся константа $c > 0$ и подпоследовательность 
$(W_m)_{m \in \mathbb{N}}$ последовательности 
$(W_n)_{n \in \mathbb{N}}$ такие, что для каждого $m \in \mathbb{N}$
$$\beta (P_M (W),P_M (W_m)) > c.$$ 
Следовательно, найдется такая последовательность $(z_m)_{m \in \mathbb{N}}$  в компакте $W$, что для каждого $m \in \mathbb{N}$
$$|P_M (z_m) P_M (W_m)| > c.$$
Кроме того, найдется подпоследовательность $(z_k)_{k \in \mathbb{N}}$ последовательности  $(z_m)_{m \in \mathbb{N}}$, сходящаяся к некоторой 
точке $z \in W$, так как $W$ --- компакт. Но оператор $P_M$ является непрерывным, поэтому найдется такое число $k_0 \in \mathbb{N}$, что для каждого $k > k_0$
$$|P_M (z) P_M (W_k)| > c.$$
Из условия
$$\lim\limits_{k\to\infty}{\alpha (W_k, W)} = 0$$ 
следует, что найдется сходящаяся к точке $z$ последовательность $(y_k)_{k \in \mathbb{N}}$ такая, что для каждого $k > k_0$ $y_k \in W_k$.
Но тогда 
$$\lim\limits_{k\to\infty}{|P_M (z)P_M (W_k)|} \leq \lim\limits_{k\to\infty}{|P_M (z)P_M (y_k)|} = 0.$$ 
Получили противоречие. Таким образом, лемма 3 доказана.
\vskip7pt
{\bf Лемма 4} [s14]. {\it Пусть геодезическое пространство $X$ удовлетворяет условиям  $(A_3)$, $(A_4)$, $B = B[p,r]$. Тогда отображение $P_B$ является липшицевым отображением с липшицевой константой $Lip (P) \leq 2$. Если, кроме того, выполняется условия $(A_6)$, $(A_7)$,  $(A_{10})$, то $Lip (P) = 1$.}
\vskip7pt
{\bf Доказательство леммы 4}. 
\vskip7pt
 Из условий $(A_3)$, $(A_4)$ нетрудно получить неравенство:
$$|\omega_{\lambda}(p,x)\omega_{\lambda}(p,y)| \leq \lambda |xy|$$ 
для любых
$$0 \leq \lambda \leq 1;\quad x,\, y,\, p \in X.$$
Очевидно, что для любых $x,\,  y \in B$
$$|P_B (x)P_B (y)| = |xy|.$$
Докажем, что для любых $x,\,  y \in X$
$$|P_B (x)P_B (y)| \leq 2|xy|.$$
Пусть для определенности 
$$|py| > r, \quad 0 < |px|  \leq |py|.$$
Если $|px| \leq r$, то 
$$|P_B (x)P_B (y)| \leq |xy| + |yP_B (y)| \leq 2|xy|.$$
Если $|px| > r$, то 
$$|P_B (x)P_B (y)| = |\omega_{\lambda} (p,x)\omega_{\lambda}(p,z)|,$$
где точка $z \in [p,y]$ такая, что 
$$|pz| = |px| \quad \mbox {и} \, \, \lambda =\frac{r}{|px|} < 1.$$ 
Тогда 
$$|\omega_{\lambda}(p,x)\omega_{\lambda}(p,z)| \leq \lambda|xz| \leq |xy| + |yz| \leq 2|xy|.$$
Докажем второе утверждение. Пусть для определенности 
$$|py| > r, \quad 0 < |px|  \leq |py|.$$
Если $|px| \leq r$, то необходимое неравенство следует из
неравенства $(i)$ леммы 1. Если $|px| > r$, то
$$|\omega_{\lambda}(p,x)\omega_{\lambda}(p,y)| \leq \lambda|xy| < |xy|,$$
где 
$$\lambda =\dfrac{r}{|px|} < 1.$$ 
Но $\omega_{\lambda}(p,x) = P_B (x)$ и
$$|p\omega_{\lambda}(p,y)| \geq r.$$ 
Следовательно, в силу неравенства $(i)$ леммы 1 получим
$$|P_B (x)P_B (y)| \leq |\omega_{\lambda}(p,x)\omega_{\lambda}(p,y)|\leq \lambda|xy| < |xy|.$$
Лемма 4 доказана.
\vskip7pt
{\bf Теорема 1} [s14]. {\it Пусть $(M_n)_{n \in \mathbb{N}}$ --- последовательность непустых ограниченных замкнутых множеств пространства $(X,\rho)$. Тогда верны следующие утверждения.
\vskip7pt
$(i)\, \,$ Если 
$$M \in \Sigma,\quad \lim\limits_{n\to\infty} \alpha (M_n,M) = 0,\quad
(K_n)_{n \in \mathbb{N}} \subset \Sigma $$ 
--- последовательность наилучших аппроксимирующих множеств
для последовательности $(M_n)_{n \in \mathbb{N}}$ такая, что для каждого натурального $n\, \,$  $K_n \subset M_n$, то 
$$\lim\limits_{n\to\infty} \alpha (K_n, M)= 0.$$ 

$(ii)\, \,$ Если $M$ --- компакт в $X$, 
$$\lim\limits_{n\to\infty} \alpha (M_n, M) = 0,\quad (S^n)_{n \in \mathbb{N}}$$ 
--- последовательность наилучших $N$-сетей для
последовательности $(M_n)_{n \in \mathbb{N}}$ такая, что для каждого натурального $n$  
$S^n \subset M_n$, то найдется подпоследовательность последовательности
$(S^n )_{n \in \mathbb{N}}$, сходящаяся в метрике Хаусдорфа к некоторой наилучшей $N$-сети множества $M$.}
\vskip7pt
Непосредственое следствие 2 из следствия 1 и теоремы 1 обобщает теорему 3 из \cite{Belob}.
\vskip7pt
{\bf Следствие 2}  [s14]. {\it Пусть $X$ --- полное геодезическое пространство, удовлетворяющее условиям $(A_4-A_7)$,  $(A_{10})$, $(M_n)_{n \in \mathbb{N}}$ --- последовательность непустых ограниченных замкнутых выпуклых множеств пространства $X$, сходящаяся в метрике Хаусдорфа к некоторому компакту $M \subset X$. 
Тогда верны следующие утверждения.
\vskip7pt
$(i)\, \,$ Если $(K_n)_{n \in \mathbb{N}}$ --- последовательность наилучших аппроксимирующих компактов для
последовательности $(M_n)_{n \in \mathbb{N}}$, то найдется подпоследовательность  $(M_m)_{m \in \mathbb{N}}$ последовательности $(M_n)_{n \in \mathbb{N}}$,  для которой найдется последовательность $(\hat K_m)_{m \in \mathbb{N}}$ наилучших аппроксимирующих компактов, сходящаяся в метрике Хаусдорфа к множеству $M$. 
\vskip7pt
$(ii)\, \,$ Если $(S^n )_{n \in \mathbb{N}}$ --- последовательность наилучших $N$-сетей для последовательности $(M_n)_{n \in \mathbb{N}}$, то найдется подпоследовательность 
\\
$(M_m)_{m \in \mathbb{N}}$ последовательности $(M_n)_{n \in \mathbb{N}}$, для которой найдется последовательность 
$(\hat S^m)_{m \in \mathbb{N}}$ наилучших $N$-сетей, сходящаяся в метрике Хаусдорфа к некоторой наилучшей $N$-сети множества $M$.}
\vskip7pt
{\bf Доказательство теоремы 1}. 
\vskip7pt
$(i)\, \,$  Из условий теоремы 1 и утверждения $(iv)$ теоремы 2.1.1 нетрудно получить: 
$$\beta (K_n,M) \leq \beta (M_n,M) \leq \alpha (M_n,M),$$
$$\beta (M,K_n) \leq \beta (M,M_n) + \beta (M_n,K_n) = \beta (M,M_n) + R_{\Sigma }(M_n) \leq 2\alpha (M,M_n)$$ 
для каждого $n \in \mathbb{N}$. Следовательно,
$$\lim\limits_{n\to\infty}\alpha (K_n,M)= 0.$$

$(ii)\, \,$ Из условий теоремы 1 следует, что найдется такая подпоследовательность 
$$(S^m= \{y_1^m,\ldots y_N^m\})_{m \in \mathbb{N}}$$ 
последовательности $(S^n)_{n \in \mathbb{N}}$, что существуют пределы
$$y_i^{*} = \lim \limits_{m \to \infty}{y_i^m} \quad (1 \leq i \leq N).$$ 
Пусть 
$$S^{*}= \{y_1^{*},\ldots y_N^{*}\}.$$ 
Тогда для каждого $m \in \mathbb{N}$ справедливы неравенства:
$$\beta (M,S^{*}) \leq \beta (M,M_m) + \beta (M_m,S^m) + \beta (S^m,S^{*})
 \leq$$ 
$$\alpha (M,M_m) + R_N (M_m) + \beta (S^m,S^{*}).$$ 
Из условия теоремы, утверждения $(iv)$ теоремы 2.1.1 и определения $N$-сети $S^{*}$ следует, что правая часть полученного неравенства стремится к $R_N (M)$ при $(m \to \infty)$. Следовательно, $S_N^{*}$ является наилучшей $N$-сетью для множества $M$.
Теорема 1 доказана.

\section{Наилучшая $N$-сеть и наилучшее сечение ограниченного множества в бесконечномерном пространстве Лобачевского}

\vskip20pt

Рассмотрим в вещественном гильбертовом пространстве $V$ открытый шар $B(0,1)$ единичного радиуса, с центром в нуле и зададим в множестве $X = B(0,1)$ метрику $\rho$ модели Бельтрами -- Клейна бесконечномерного пространства Лобачевского \cite [с. 48]{Nut}:
$$|xy|= k \mathop{\rm Arch} \left(\frac{1- (x,y)}{((1 -x^2)(1 - y^2))^{ 1/2}}\right), \eqno (1) $$
где (x,y) --- скалярное произведение векторов $x, \, y $ из $X$, $k > 0$ ---
константа. Напомним следующие определения.
\vskip7pt
Наилучшим $N$-мерным сечением множества $M \in\Sigma (X)$ называется $N$-мерная плоскость $H^{*}$ пространства $X$, для которой
$$\beta (M, H^{*})
= \inf\limits_{H \subset X}{\beta (M, H)} < \infty, \eqno (2)$$
где точная нижняя грань берется по всевозможным $N$-плоскостям пространства $X$ \cite {Gark}.
\vskip7pt
Пусть $X$ --- метрическое пространство, в котором любое множество $M \in B (X)$ имеет единственный чебышевский центр $Z (M)$. Чебышевский центр $Z (M)$ множества $M \in B (X)$ называется {\it сильно устойчивым} \cite {Belob}, если для каждой последовательности непустых ограниченных множеств $(M_n)_{n \in \mathbb{N}}$, сходящейся в метрике Хаусдорфа к множеству $M$
$$\lim\limits_{n\to\infty}|Z (M_n) Z (M)| =0.$$
\vskip7pt
{\bf Теорема 1} [s8]. {\it
Для каждого непустого ограниченного множества $M$ пространства Лобачевского $X$ и каждого натурального $N$ существует наилучшая $N$-сеть.}
\vskip7pt
{\bf Доказательство теоремы 1}.
\vskip7pt
 Пусть $ (\varepsilon_l)_{l \in \mathbb{N}}$ --- последовательность положительных вещественных чисел,
сходящаяся к нулю. Тогда для
каждого $l \in \mathbb{N}$ найдется такая $N$-сеть 
$$S^l = \{y^l_1,\,\ldots,\, y^l_N\},$$ 
что
$$\beta (M,S^l) \leq R_N (M) + \varepsilon_l.\eqno (3)$$
Можно считать, что для каждого $l \in \mathbb{N}$   $N$-сеть $S^l$ принадлежит некоторому ограниченному замкнутому множеству  $P$ в пространстве $X$. Действительно, пусть $P$ --- замыкание объединения множества всех тех шаров пространства  $X$, центры которых принадлежат ограниченному множеству
$M$ и радиусы которых равны:
$$1 + R_N (M) + \max \{\varepsilon_l :  l  \in \mathbb{N}\}.$$
Для каждого $l \in \mathbb{N}$  $N$-сеть $S^l$  можно изменить (при сохранении прежних обозначений ) следующим образом:
 элементы $N$-сети $S^l$, принадлежащие множеству  $X\backslash P$
( если такие найдутся ), заменяем на произвольные элементы из множества $P$. Неравенство $(3)$ при этом сохранится.
  Нетрудно понять, что множество $P$ ограничено и
замкнуто в гильбертовом пространстве $V$ и, следовательно, слабо компактно в нем. Пусть 
$$u^{*} = \{y^{*}_1,\,\ldots,\, y^{*}_N\}$$
--- $N$-сеть, составленная из $N$ предельных точек последовательностей
$$(y^l_1)_{l \in \mathbb{N}},\,\ldots,\, (y^l_N)_{l \in \mathbb{N}}$$ 
в слабой топологии на множестве $P$.
Покажем, что $N$-сеть $u^{*}$ наилучшая. Допустим, что это не
выполняется. Тогда найдутся такие $x \in M$ и $c > 0$, что для каждого $n \in \mathbb{N}$
$$|xy^{*}_n| > R_N (M) + c. \eqno (4) $$
В силу $(3)$ имеем:
$$|xy^l_{n_l}|  \leq R_N (M) + \varepsilon_l, \eqno (5) $$
где $n_l = n_l (x)$ принимает одно из значений $1, 2,\, \ldots, \,N$. В
дальнейшем, вместо двойного индекса $n_l$ (или $n_t$) всюду пишем
индекс $n$. Пусть $(y^t_n)_{t \in \mathbb{N}}$ --- подпоследовательность последовательности
$(y^l_n)_{l \in \mathbb{N}}$, слабо сходящаяся к $y^{*}_n$. Тогда существуют сходящиеся к нулю последовательности положительных вещественных чисел
$ (\delta_t)_{t \in \mathbb{N}}$ и $(\tau_t)_{t \in \mathbb{N}}$, для которых, начиная с некоторого $t \in \mathbb{N}$, выполняются неравенства:
$$|(y^t_n - y^{*}_n,x)| < \tau_t, \quad |(y^t_n, y^{*}_n)
 - (y^{*}_n)^2| < \delta_t ||y^{*}_n||,$$ 
$$(y^{*}_n)^2 < (y^t_n, y^{*}_n) + \delta_t ||y^{*}_n||
 \leq ||y^t_n|| ||y^{*}_n|| + \delta_t ||y^{*}_n||, \quad ||y^{*}_n|| \leq ||y^t_n|| + \delta_t  < 1. \eqno (6)$$
Используя эти неравенства, а также $(1)$, $(4)$ и $(5)$, получим:
$$\ch \left(\dfrac{R_N (M) + c}{k}\right) < \dfrac{1- (x,y^{*}_n)}
{((1 -x^2)(1 - (y^{*}_n)^2))^{ 1/2}}
\leq \dfrac{1- (x,y^t_n) + \tau_t}
{[(1 -x^2)(1 - (y^{*}_n)^2)]^{ 1/2}}
\leq$$
 $$\ch \left(\dfrac{R_N (M) + \varepsilon_t}{k} \right) (1 - (y^t_n)^2)^{1/2}(1 - (|y^t_n| +
\delta_t)^2)^{- 1/2} + \tau_t ((1 -x^2)(1 - (y^{*}_n)^2))^{-
1/2}. \eqno (7)$$
Очевидно, что правая часть этого неравенства стремится к
$$\ch \left(\dfrac{R_N (M)}{k}\right)$$ 
при $t \to \infty$. Но 
$$\ch \left(\dfrac{R_N (M) + c}{k}\right) > \ch \left(\dfrac{R_N (M)}{k}\right)$$ 
при $c > 0, \quad R_N (M) > 0, \quad k > 0$. Получили
противоречие. Таким образом, теорема 1 доказана.
\vskip7pt
{\bf Теорема 2} [s8]. {\it
Для каждого множества $M$ пространства Лобачевского $X$, для которого нижняя грань из $(2)$ конечна {\rm{(}}в частности, для множества
$M \in B (X))$,  существует наилучшee $N$-мерное сечение для каждого натурального $N$.}
\vskip7pt
{\bf Доказательство теоремы 2}. 
\vskip7pt
Очевидно, что каждый элемент $y$
произвольной $N$-плоскости $H_N$ пространства $X$ можно представить
в виде:
$$y = y_0 + a_1 y_1 + \ldots +a_N y_N,$$
 где $y_0$ --- вектор,
ортогональный $N$-плоскости $H_N$, определяющий точку из $H_N$
(пишем: $y_0 \in H_N$), 
$$\{y_1,\, \ldots,\, y_N\}$$ 
--- семейство из $N$ ортонормированных векторов пространства $V$, ортогональных вектору
$y_0$ и 
$$|a_m| < 1,\quad \mbox {где} \, \, m \in \{1, \ldots,N\}.$$ 
Последнее
неравенство следует из того, что $y^2 < 1$.
 Пусть для множества $M$ пространства $X$
$$r = \inf\limits_{H \subset X}{\beta (M,
H)} < \infty,$$
где точная нижняя грань берется по всевозможным $N$-плоскостям пространства $X$.
Рассмотрим последовательность $N$-мерных плоскостей 
$(H^n)_{n \in \mathbb{N}}$ такую, что
$$\beta (M,
H^n)  \leq r + \varepsilon_n,$$
где $(\varepsilon_n)_{n \in \mathbb{N}}$ --- последовательность положительных вещественных чисел,
сходящаяся к нулю. Предположим, что 
$$y^n_0,\,y^n_1,\,\ldots,\,y^n_N$$  
--- элементы, введенные ранее для $N$-мерных
плоскостей $H^n$ $(n \in \mathbb{N})$. Пусть $x$ произвольная
точка из $M$ и
$$y^n_x = y^n_0 + a^n_{x,1} y^n_1 + \ldots
+a^n_{x,N} y^n_N$$ 
--- ближайший к точке $x$ в метрике $\rho$ элемент
$N$-мерной плоскости $H^n$ $(n \in \mathbb{N})$. Наделим замкутый шар $B[0,1]$ с
центром в нуле, единичного радиуса индуцированной слабой топологией
гильбертова пространства $V$, а отрезок  $I_x = [-1,1]$ наделим индуцированной стандартной
топологией прямой. Тогда множество
$$Q = B^{N + 1}[0,1] \times \prod \{I^N_x : x \in M \}$$ является компактным множеством в топологии произведения. Пусть $(u_n)_{n \in \mathbb{N}}$ --- последовательность точек из
множества $Q$, проекциями которых являются точки:
$$y^n_0,y^n_1,
\ldots ,y^n_N \in B[0,1]; \quad a^n_{x,1}, \ldots ,a^n_{x,N} \in
I_x,$$ 
где $x \in M$, $n \in \mathbb{N}$. И пусть точка $u \in Q$ с проекциями:
$$y_0,y_1,
\ldots ,y_N \in B[0,1]; \quad a_{x,1}, \ldots ,a_{x,N} \in I_x,$$ 
где $x
\in M$, является предельной для последовательности $(u_n)_{n \in \mathbb{N}}$.
Докажем, что $N$-мерная плоскость $H^{*}$, элементы которой имеют
вид:
$$y = y_0 + a_1 y_1 + \ldots +a_N y_N,$$
 является наилучшим
$N$-мерным сечением для $M$. Пусть $x$ --- произвольная точка из $M$.
Тогда 
$$|xy_x| \leq r, \eqno (8)$$
 где $y_x = y_0 + a_{x,1} y_1 +
\ldots +a_{x,N} y_N$ --- ближайший к точке $x$ элемент $N$-мерной
плоскости $H^{*}$ в метрике $\rho$. Действительно, если это
неравенство не выполняется, то найдется такое число $c > 0$, что
$$|xy_x| > r + c.$$
 Очевидно, что для последовательности
$(y^n_x)_{n \in \mathbb{N}}$ найдется подпоследовательность
$$(y^t_x =
 y^t_0 + a^t_{x,1} y^t_1 + \ldots +a^t_{x,N} y^t_N\})_{t \in \mathbb{N}},$$
слабо сходящаяся к точке $y_x$. Рассуждая так же, как и в теореме 1,
получим неравенства $(6)$ и $(7)$, где вместо точек $y^t_n$,
$y^{*}_n$ нужно подставить точки $y^t_x$ и $y^{*}_x$
соответственно. Следовательно, получаем противоречие. Таким образом,
неравенство $(8)$ верно и $N$-мерное сечение $H^{*}$ действительно наилучшее. Теорема 2 доказана.
\vskip7pt
{\bf Теорема 3} [s8]. {\it
Чебышевский центр произвольного множества $M\in B (X)$ пространства Лобачевского $X$ сильно устойчив.}
\vskip7pt
{\bf Доказательство теоремы 3}.
\vskip7pt
Поскольку при различных коэффициентах
$k$ в формуле $(1)$ получаются подобные метрики, то в теореме достаточно рассмотреть случай, когда $k = 1$.
 Пусть теорема неверна. Тогда найдутся непустые
ограниченные множества $M$, $(M_n)_{n \in \mathbb{N}}$ пространства $X$
с чебышевскими центрами $y^{*}$, $(y^{*}_n)_{n \in \mathbb{N}}$ и чебышевскими радиусами
 $r$, $(r_n)_{n \in \mathbb{N}}$
соответственно, для которых выполняются следующие условия.
$$\lim\limits_{n\to\infty}\alpha(M_n,M) = 0. \leqno (i)$$
$(ii)\, \,$ Найдется подпоследовательность 
$(y^{*}_l)_{l \in \mathbb{N}}$ последовательности
$(y^{*}_n)_{n \in \mathbb{N}}$ и число $A > 0$, для которых, начиная с некоторого номера
$l_0$, выполняется следующее неравенство.
$$|y^{*}y^{*}_l| \geq 2A . \eqno (9)$$
В силу утверждения $(iv)$ теоремы 2.1.1 имеет место неравенство: 
$$|r_n - r| \leq \alpha(M_n,M) \eqno (10)$$
Пусть $x$ --- произвольная точка в множестве $M$. Выберем число
$\varepsilon > 0$ и номер $l > l_0$ такие, что:
$$\ch (r + 2 \alpha(M_l,M) + \varepsilon) < (2 \ch (A) - 1) \ch (r). \eqno (11)$$
Кроме того, $|xy^{*}| \leq r$ и найдется точка $x_l \in M_l$ такая, что 
$$|xx_l| \leq \alpha(M_l,M) + \varepsilon.$$
Тогда, учитывая $(10)$, получим:
$$|xy^{*}_l| \leq |xx_l| + |x_ly^{*}_l| \leq
\alpha(M_l,M) + \varepsilon + r_l \leq r + 2 \alpha(M_l,M) + \varepsilon. \eqno (12)$$
Пусть $w_l$ середина отрезка с концами $y^{*}$ и $y^{*}_l$ в метрике $\rho$.  Используем формулу
$$\ch |zw| = \dfrac{\ch |zu| + \ch |zv|}{2 \ch \left(\dfrac{|uv|}{2}\right)}$$
для нахождения длины медианы $|zw|$  в треугольнике с вершинами $z$, $u$ и $v$ в пространстве Лобачевского $X$ (эту формулу нетрудно получить из теоремы косинусов \cite [с. 60]{Shir}). Учитывая выражение для длины медианы в
треугольнике с вершинами $x$, $y^{*}$, $y^{*}_l$ и неравенства (9), (11), (12), получим:
$$\ch|xw_l| = \dfrac{\ch|xy^{*}| +
\ch |xy^{*}_l|}{2 \ch \left(\dfrac{|y^{*}y^{*}_l|}{2}\right)} \leq
\dfrac{\ch (r) + \ch (r + 2 \alpha(M_l,M) + \varepsilon)}{2 \ch (A)} < \ch (r).$$
Таким образом, $|xw_l| < r$ для каждой точки $x \in M$. Получили
противоречие с  тем, что $r$ --- чебышевский радиус множества $M$.
Следовательно, теорема 3 доказана.

\section{Наилучшее приближение выпуклого компакта геодезического пространства шаром}
\vskip20pt

В метрическом пространстве  $(X,\rho)$ используем следующие обозначения.
\vskip7pt
$BC [X]$ --- множество всех замкнутых шаров в пространстве $X$.
$$\psi (M,x) = \dfrac{1}{2}(\beta (M,x) + |xM| - |x(X\backslash M)|),$$
$$r (M,x) = \dfrac{1}{2} (\beta (M,x) - |xM| + |x(X\backslash M)|),$$
где $x \in X$, $M \in B(X)$, $M \neq X$.
$$\chi (M) = \{x \in X : \alpha (B[x, r], M) = \alpha (M,BC [X]) \quad
\mbox {для некоторого}\, \, r \geq 0 \}$$ 
\cite {Dudov}, где
$$\alpha (M, BC[X])= \inf \{\alpha (M,B[y,t]) : y \in X, t \geq 0 \}.$$

Напомним, что множество $M$ в пространстве $X$, удовлетворяющем условию $(A_1)$, называется строго выпуклым, если для любых различных
$x,\, y \in \overline{M}$ $[x, y]\backslash \{x,y\} \subset \stackrel{\circ}{M}$ \cite [с. 153]{Bus}.
\vskip7pt
Нашими основными задачами в данном параграфе являются следующие задачи.
\vskip7pt
Получить оценку сверху для расстояния Хаусдорфа от непустого
ограниченного множества до множества всех замкнутых шаров специального геодезического пространства $X$ неположительной
кривизны. Доказать, что множество всех центров $\chi (M)$
замкнутых шаров, наилучшим образом приближающих в метрике
Хаусдорфа выпуклый компакт $M \subset X$, непустое и принадлежит
множеству $M$. Исследовать некоторые другие свойства множества
$\chi (M)$. Таким образом, необходимо обобщить некоторые результаты С. И. Дудова и И. В. Златорунской \cite {Dudov}, \cite {Dudov1} на случай специального геодезического пространства неположительной кривизны по Буземану. Полученные основные результаты
содержатся в теоремах 1, 2, 3.

 Сначала исследуем некоторые простые свойства вышеперечисленных понятий.
\vskip7pt
{\bf Лемма 1} [s17]. {\it
\vskip7pt
$(i)\, \,$ Для любых $M,\, \, W \in B(X)$ имеют место неравенства:
$$|\beta(M,BC[X])-\beta(W,BC[X])| \leq \alpha (M,W),$$
$$|\alpha (M,BC[X]) - \alpha (W,BC[X])| \leq \alpha (M,W),$$
где $\beta (M,BC[X]) = \inf \{\beta (M,V) : V \in BC[X] \}.$
\vskip7pt
$(ii)\, \,$ Для любых 
$$M \in B(X),\quad M \neq X;\quad x,\, y \in X$$
имеют место неравенства:
$$|\psi(M,x) - \psi(M,y)| \leq \dfrac{3|xy|}{2},$$
$$|r(M,x) - r(M,y)| \leq \dfrac{3|xy|}{2}.$$

$(iii)\, \,$ Для любых 
$$M \in B(X),\quad M \neq X,\quad x \in M$$
имеет место неравенство:
$$\psi(M,x)\leq \dfrac{D(M)}{2}.$$
$(iv)\, \,$ Для любых $M \in B(X)$, $x \in X\backslash M$ имеет место неравенство:
$$r(M,x) \leq \dfrac{D(M)}{2}.$$
$(v)\, \,$ Для любых 
$$M \in B(X),\quad M \neq X,\quad x \in X$$
имеет место равенство:
$$\psi (M,x) = \inf \{\max \{r + |xM| - |x(X\backslash M)|, \beta (M, x) - r\} :
r \geq 0 \},$$
где инфимум достигается при $r = r(M,x)$.}
\vskip7pt
{\bf Доказательство леммы 1}.
\vskip7pt
$(i)\, \,$ Из неравенства треугольника для функции $\beta$ следует
$$\beta (M, BC[X]) = \inf \{\beta (M,V) : V \in BC[X] \} \leq$$ 
$$\inf \{\beta (M,W) + \beta (W,V) : V \in BC[X] \} \leq  \alpha (M,W) +
\beta (W,BC[X]).$$ 
Неравенство 
$$\beta (W,BC[X]) \leq \alpha (M,W) +
\beta (M,BC[X])$$ 
доказывается аналогично.
\vskip7pt
$(ii)\, \,$ Из неравенства треугольника для функции $\beta$ и неравенства
из \cite [с. 219]{Kurat} следует
$$|\psi (M,x) - \psi (M,y)| \leq \dfrac{1}{2} (|\beta (M, x) - \beta (M,y)|+ ||xM| - |yM|| +$$
$$||y(X\backslash M)| - |x(X\backslash M)||) \leq \dfrac{3 |xy|}{2}.$$ Для второго неравенства доказательство аналогично.
\vskip7pt
$(iii)\, \,$ Свойство следует из неравенств:
$$\psi (M,x) = \dfrac{(\beta (M, x) - |x(X\backslash M)|)}{2} \leq \dfrac{\beta (M, x) }{2} \leq \dfrac{D(M)}{2}.$$

$(iv)\, \,$ Пусть $x \in X\backslash M$. Тогда для каждого
$\varepsilon > 0$ найдутся такие
$z, \, y \in M$, что
$$\beta (M, x) \leq |zx| + \varepsilon,\quad |xy| \leq |xM| + \varepsilon.$$ Утверждение следует теперь из неравенств
$$r(M,x) = \dfrac{\beta (M, x) - |xM|}{2} \leq \dfrac{|zx| - |xy|}{2} \leq \dfrac{|zy|}{2} \leq \dfrac{D(M)}{2}.$$

$(v)\, \,$ Неравенство 
$$\psi (M,x) \geq \inf \{\max \{r + |xM| - |x(X\backslash M)|, \beta (M, x) - r\} : r \geq 0 \}$$ 
получается, если положить $r = r(M,x)$. Если 
$$r + |xM| - |x(X\backslash M)| \geq \beta (M, x) -r,$$ 
то $r \geq  r(M,x)$. Следовательно,
$$\max \{r + |xM| - |x(X\backslash M)|, \beta (M, x) - r\} = r + |xM| -
|x(X\backslash M)| \geq \psi (M,x).$$
Если $$r + |xM| - |x(X\backslash M)| \leq  \beta (M, x) - r,$$
то  $r \leq  r(M,x)$. Следовательно,
$$\max \{r + |xM| - |x(X\backslash M)|, \beta (M, x) - r\} = \beta (M, x) - r  \geq \psi (M,x).$$
Таким образом, получено неравенство:  
$$\psi (M,x) \leq \inf \{\max
\{r + |xM| - |x(X\backslash M)|, \beta (M, x) - r\} : r \geq 0 \}.$$ 
Лемма 1 доказана.
\vskip7pt
{\bf Лемма 2} [s17]. {\it Пусть $X$ геодезическое пространство, через каждые две различные точки которого можно провести единственную прямую. Тогда имеют место следующие утверждения.
\vskip7pt
$(i)\, \,$ Для любых $M \in B(X)$, $x \in X$ справедливо неравенство: 
 $$r(M,x) \leq \dfrac{D (M)}{2}.$$

$(ii)\, \,$ Пространство $(BC[X],\alpha)$ изометрично пространству
$$(X\times  \mathbb{R}_{+},\, \rho + d),$$ 
где $d$ --- стандартная метрика на $ \mathbb{R}_{+}.$
\vskip7pt
$(iii)\, \,$ Для любых 
$$x,\, y \in X;\quad r_1,\, r_2 \in  \mathbb{R}_{+},
\quad r = \dfrac{r_1+ r_2}{2}$$
верны равенства:
$$\alpha (B[x,r_1],B[\omega (x,y),r]) =
\alpha (B[y,r_2],B[\omega (x,y),r]) =\dfrac{1}{2} \alpha (B[x,r_1],B[y,r_2]).$$

$(iv)\, \,$ Если, кроме того, пространство $(X,\rho)$ удовлетворяет
условию $(A_3)$, то каждый замкнутый шар пространства $X$ является
строго выпуклым множеством.}
\vskip7pt
{\bf Доказательство леммы 2}.
\vskip7pt
$(i)\, \,$ В силу утверждения $(iv)$ леммы 1 достаточно рассмотреть случай,
когда $x \in M$. Если $M = \{x\}$, то доказательство очевидно.
Если $M \neq \{x\}$, то для каждого $\varepsilon > 0$ найдется
такая точка $u \in M\backslash \{x\}$, что 
$$|ux| \geq \beta (M, x) - \varepsilon.$$ 
Пусть $\Lambda [u,x)$ --- полупрямая с началом в точке
$u$, содержащая точку $x$, и точка 
$$y \in \overline{M}\cap \Lambda [u,x)$$
такая, что 
$$|yu| = D (\overline{M}\cap \Lambda [u,x)).$$ 
Тогда
$$r(M,x) = \dfrac{1}{2} (\beta (M, x) + |x(X\backslash M)|) \leq \dfrac{1}{2} (|ux| + |xy| - \varepsilon) =$$ 
$$\dfrac{1}{2} (|uy| - \varepsilon) \leq \dfrac{1}{2} |D (M)- \varepsilon|.$$
Таким образом, второе утверждение верно в силу произвольности
выбора $\varepsilon > 0$.
\vskip7pt
$(ii)\, \,$ Определим отображение 
$$f : (BC[X],\alpha) \rightarrow (X\times  \mathbb{R}_{+}, \, \rho + d)$$ 
формулой $f(B[x,r]) = (x; r)$.
Тогда для любых $x,\, y \in X$, $r,\, R \in  \mathbb{R}_{+}$
$$(\rho + d)(f(B[x,r]),f(B[y,R])) = |xy| + |r - R| =
\alpha (B[x,r],B[y,R]).$$
Следовательно, $f$ --- изометрия.
\vskip7pt
$(iii)\, \,$ Докажем одно из равенств непосредственным вычислением:
$$\alpha (B[x,r_{1}],B[\omega (x,y),r]) = |x\omega (x,y)| +
|r_1 - r|=$$ 
$$\dfrac{1}{2} (|xy| + |r_1 - r_2|) = \dfrac{1}{2} \alpha (B[x,r_1],B[y,r_2]).$$
\vskip7pt
$(iv)\, \,$ Доказательство выпуклости шара почти очевидно. Действительно, пусть $B[p,r]$ --- произвольный замкнутый шар в пространстве
$X$ и $x,\, y \in B[p,r]$. Рассмотрим функцию $f = |pz(t)|$, где
$z = z(t)$ --- параметризация сегмента $[x,y]$ длиной дуги. Из
неравенств
$$||pz(t)| - |pz(t_1)|| \leq |z(t)z(t_1)| = |t - t_1|$$
для любых $t,\, t_1 \in [0,|xy|]$ следует, что функция $f$
непрерывна. Кроме того, из условия $(A_3)$ и неравенства
треугольника для любых $t, \, t_1 \in [0,|xy|]$ получим
$$f \left(\dfrac{t+t_1}{2}\right) = |p\omega(z(t),z(t_1))| \leq |p\omega(p,z(t))| +
|\omega(p,z(t))\omega(z(t),z(t_1))|\leq$$ 
$$\dfrac{1}{2} (|pz(t)| +  |pz(t_1)|) = \dfrac{1}{2} (f(t)+f(t_1)).$$ 
Следовательно, функция $f$ выпуклая и 
$[x,y] \subset B[p,r]$. Таким образом, каждый замкнутый шар пространства
$X$ является выпуклым множеством (это свойство нетрудно также
получить, используя доказательства теорем (36.4 - 36.6) из \cite [c. 304-307]{Bus}). Докажем строгую выпуклость замкнутого шара $B[p,r]$
методом от противного. Пусть этот шар не является строго выпуклым
множеством. Тогда найдутся такие попарно различные точки $x$, $y$,
$z$, что 
$$|px| = |py| = |pz| = r,\quad z \in [x,y].$$ 
Выберем на прямой
$\Lambda (p,z)$, содержащей точки $p$, $z$, точку $u$, удовлетворяющую
условиям: $u \neq z$ и $p \in [z,u]$. Тогда из условия леммы 2 и
неравенства треугольника следует, что 
$$|uz| = |up| + r > \max \{|ux|,|uy|\}.$$
Следовательно, замкнутый шар $B[u,\max \{|ux|,|uy|\}]$ невыпуклый.
Получили противоречие. Лемма 2 доказана.
\vskip7pt
{\bf Лемма 3} [s17]. {\it Пусть геодезическое пространство $X$
удовлетворяет условиям $(A_1)$, $(A_3)$. Тогда имеют место следующие
утверждения.
\vskip7pt
$(i)\, \,$ Для любых $x, y \in X$, $M \in B(X)$ 
$$2\beta (M,\omega (x,y)) \leq \beta (M, x) + \beta (M,y).$$

$(ii)\, \,$ Если 
$$\beta (M, x) = \beta (M,y) + |xy|,$$ 
то 
$$2\beta (M,\omega (x,y)) = \beta (M, x) + \beta (M,y)$$
$(iii)\, \,$ Если, кроме того, $M$ --- выпукло, то 
$$2 |\omega (x,y)M| \leq |xM| + |yM|.$$}

{\bf Доказательство леммы 3}.
\vskip7pt
$(i)\, \,$ Для любого $\varepsilon > 0$ найдется такое $z \in M$, что
$$|z\omega (x,y)| \geq \beta (M,\omega (x,y)) - \varepsilon .$$ 
Кроме того, из
условия $(A_3)$ следует, что
$$2 |z\omega (x,y)| \leq |xz| + |yz| \leq \beta (M, x) + \beta (M,y).$$ 
Тогда  
$$2 \beta (M,\omega (x,y)) - 2\varepsilon \leq \beta (M, x) + \beta (M,y).$$
Осталось использовать произвольность $\varepsilon > 0$.
\vskip7pt
$(ii)\, \,$ Используем рассуждение от противного. Пусть свойство неверно.
Тогда в силу доказанного неравенства $(i)$ 
$$2\beta (M,\omega (x,y)) < \beta (M, x) + \beta (M,y).$$ 
Из условия леммы 3 следует, что 
$$\beta (M, x) + \beta (M,y) = 2 \beta (M, x) - |xy|.$$ 
Следовательно,
$$|xy| < 2\beta (M, x) - 2\beta (M,\omega (x,y)) \leq 2|x\omega (x,y)| = |xy|.$$ Получили противоречие. Таким образом, равенство верно. 
\vskip7pt
$(iii)\, \,$ Для любого $\varepsilon > 0$ найдутся такие $u,\, v \in M$,
что 
$$|xu| \leq |xM| + \varepsilon,\quad |yv| \leq |yM| + \varepsilon.$$ 
Тогда из условия $(A_3)$ следует, что
$$2|\omega (x,y)\omega (u,v)| \leq 2|\omega (x,y)\omega (u,y)| +
2|\omega (u,y)\omega (u,v)| \leq$$ 
$$|xu| + |yv| \leq |xM| + |yM| + 2 \varepsilon.$$ 
Но множество $M$ --- выпуклое, поэтому $\omega (u,v) \in M$ и 
$$2|\omega (x,y)M| \leq 2|\omega (x,y)\omega (u,v)|
\leq |xM| + |yM| + 2\varepsilon.$$
Осталось использовать произвольность $\varepsilon > 0$. Лемма 3
доказана.
\vskip7pt
Следующая лемма 4 обобщает лемму 4.1 из \cite {Dudov}.
\vskip7pt
{\bf Лемма 4} [s17]. {\it Пусть $X$ геодезическое пространство, через каждые две различные точки которого можно провести единственную прямую, 
$$M \in B [X],\quad \psi (M,x^{*}) = \inf \{\psi (M,x) : x \in X\}.$$ 
Тогда верны следующие утверждения.
\vskip7pt
$(i)\, \,$ Если  
$$x^{*} \in X\backslash M, \quad y \in P_M (x^{*}),$$ 
то найдется такая точка $z \in X$, что 
$$S(x^{*},\beta (M,x^{*}))\cap S(y,\beta (M,y)) = \{z\}.$$ 
Если, кроме того, существуют  такие точки $u,\, v \in X$, что 
$$\beta (M,x^{*}) = |ux^{*}|,\quad \beta (M,y) = |vy|,$$ 
то $x^{*} \in [y,z]$ и $u = v = z$.
\vskip7pt
$(ii)\, \,$  Если множество $M$ --- выпуклое, $y \in P_M (x^{*})$ и
существуют такие точки $u, \, v \in X$, что 
$$\beta (M,x^{*}) = |ux^{*}|,\quad \beta (M,y) = |vy|,$$ 
то $x^{*} = y$.}
\vskip7pt
{\bf Доказательство леммы 4}.
\vskip7pt
$(i)\, \,$ Из условия леммы 4 и неравенства треугольника для функции
$\beta$ получим:
$$\beta (M,y) \leq \beta (M,x^{*}) + |x^{*}y|,$$ 
$$2 \psi (M,x^{*}) = \beta (M,x^{*}) + |x^{*}M| \leq
\beta (M,y) + |yM| - |y(X\backslash M)| = \beta (M,y).$$ 
Следовательно,
$$\beta (M,y) = \beta (M,x^{*}) + |x^{*}y|\, \, \mbox { и} \, \,
S(x^{*},\beta (M,x^{*}))\cap S(y,\beta (M,y)) \neq \emptyset.$$ 
Предположим, что
$$u,\, v \in S(x^{*},\beta (M,x^{*})) \cap S(y,\beta (M,y)),\quad u \neq v$$ 
 и точка $u$ принадлежит пересечению прямой $\Lambda (x^{*},y)$, проходящей через точки  $x^{*}$, $y$, со сферой $S(y,\beta (M,y))$.
Тогда точки $y,\, v $ являются концами двух различных сегментов
$[y,v]$ и $[y,x^{*}]\cup [x^{*},v]$. Получили противоречие с
условием леммы 4. Равенства $u = v = z$ следуют теперь из доказанной
единственности точки $z$. Кроме того, 
$$|zy| = |vy| = \beta (M,y) = \beta (M,x^{*}) +
|x^{*}y| = |ux^{*}| + |x^{*}y| = |zx^{*}| + |x^{*}y|.$$ 
Следовательно, $x^{*} \in [y,z]$.
\vskip7pt
$(ii)\, \,$ Если предположить, что $x^{*} \in X\backslash M$, то из
предыдущего пункта получим, что $x^{*} \in [y,z]$. Но множество
$M$ --- выпуклое, поэтому  $[y,z] \subset M$. Получили
противоречие. Лемма 4 доказана.
\vskip7pt
В следующей теореме 1 обобщаются леммы 3.1, 3.3 из \cite {Dudov}.
\vskip7pt
{\bf Теорема 1} [s17]. {\it Пусть геодезическое пространство $X$ удовлетворяет условию $(A_1)$. Тогда для любых 
$$r \geq 0,\quad x \in X,\quad M \in B(X),\quad M \neq X$$ 
верны следующие утверждения.
$$\beta (M,B[x,r]) = \max \{0,\beta (M, x) -r\}. \leqno (i)$$
$$\beta (B[x,r],M) \leq \max \{0,r + |xM| - |x(X\backslash M)|\}. \leqno (ii)$$
$$\alpha (B[x,r],M) \leq
\max \{r + |xM| - |x(X\backslash M)|,\beta (M, x) - r\}. \leqno (iii)$$
$$\alpha (M,BC[X]) \leq \inf \{\psi (M,x) : x \in X\}. \leqno (iv)$$
$(v)\, \,$ Если 
$$\alpha (B[x,r],M) = \alpha (M,BC[X]),$$ 
то найдется такое $r_0 \geq 0$ , что 
$$r(M,x) \leq r_0 \leq r,$$ 
$$\alpha (B[x,r_0],M) = \beta (M, x) - r_0 = \beta (B[x,r_0],M) =
\alpha (M,BC[X]).$$
$(vi)\, \,$ Множество $\chi (M)$ замкнуто и ограничено. Если, кроме
того, пространство $X$ --- собственное, то это множество компактно.
\vskip7pt
Если $X$ геодезическое пространство, через каждые две различные точки которого можно провести единственную прямую, удовлетворяет условию $(A_3)$ и $M$ --- выпуклое ограниченное множество существования  в $X$, то
для любых $r \geq 0$, $x \in X\backslash M$ верны следующие утверждения.
$$\beta (B[x,r],M) = r + |xM|. \leqno (vii)$$
$$\alpha (B[x,r],M) = \max \{r + |xM|,\beta (M, x) - r\}. \leqno (viii)$$
$$\psi (M,x) = \inf \{\alpha (B[x,r],M) : r \geq 0 \}, \leqno (i\mathrm x)$$
$$\inf \{\alpha (B[y,r],M) : y \in X\backslash M, \, r \geq 0 \} =
\inf \{\psi (M,y) : y \in X\backslash M \}.$$
$(\mathrm x)\, \,$ Если 
$$\tau \geq 0,\quad \alpha (B[x,r],M) \leq \alpha (M,BC[X]) + \tau,$$ 
то
$$\alpha (B[P_M (x),r(M,P_M (x))],M)  \leq \alpha (M,BC[X]) + \tau,$$
где $r(M,P_M (x)) = \dfrac{1}{2} \beta (M,P_M (x))$.}
\vskip7pt
{\bf Замечание 1.}

{\it $1.\, \,$ Нетрудно заметить, что если множество $M$ является замкнутым
шаром, то правая часть в неравенстве $(iv)$ равна нулю и это
неравенство становится равенством.
\vskip7pt
$2.\, \,$ Из леммы 3.4, теоремы 1.1 и утверждения в начале $\S 4$ из {\rm \cite {Dudov}} следует, что для выпуклого компакта $M$ из пространства
$\mathbb{R}^n$ с произвольной нормой в $(iv)$ имеет место равенство.}
\vskip7pt
{\bf Доказательство теоремы 1}.
\vskip7pt
$(i)\, \,$  Из неравенства треугольника для функции $\beta$ получим:
$$\beta (M,B[x,r]) \geq \beta (M, x) - \beta (B[x,r],x) \geq \beta (M, x) - r.$$
С другой стороны, для любого $\varepsilon > 0$ найдется такое
$y \in M$, что 
$$\beta (M,B[x,r]) - \varepsilon \leq |yB[x,r]| =
\max \{0,|xy| - r\} \leq \max \{0,\beta (M, x) - r\}.$$
Осталось использовать произвольность $\varepsilon > 0$.
\vskip7pt
$(ii)\, \,$ Для любых 
$$r \geq 0,\quad x \in X\backslash M,\quad M \in B(X)$$ 
из неравенства треугольника для функции $\beta$ получим: 
$$\beta (B[x,r],M) \leq  \beta (B[x,r],x) +|xM| \leq r + |xM|.$$ 
Пусть теперь
$$r \geq 0,\quad x \in M,\quad M \in B(X)$$ 
произвольны. Тогда для любого
$\varepsilon > 0$ найдется такое $z \in B[x,r]$, что 
$$|zM| \geq \beta (B[x,r],M) - \varepsilon.$$
Кроме того, найдется такое $u \in [z,x]$, что
$$|zu| = |z(\overline{M}\cap [z,x])|.$$ 
Тогда
$$\beta (B[x,r],M) - \varepsilon \leq |zM| \leq |zu| = |zx| - |xu| \leq
\max \{0,r - |x(X\backslash M)|\}.$$ 
Осталось использовать произвольность $\varepsilon > 0$.
\vskip7pt
$(iii)\, \,$ Очевидно, что 
$$\max \{r + |xM| - |x(X\backslash M)|,\beta (M, x) - r\} \leq$$
$$\max \{\max \{0,r + |xM| - |x(X\backslash M)|\},\max \{0,\beta (M, x) -r\}\}.$$
Предположим, что неравенство строгое. Тогда 
$$r + |xM| - |x(X\backslash M)| < 0\, \, \mbox {и} \, \, \beta (M, x) - r < 0.$$ Следовательно,
$$x \in M\, \, \mbox {и} \, \,\beta (M, x) < |x(X\backslash M)|.$$ 
Пусть $z \in X\backslash M$.
Тогда найдется такое $u \in [z,x]$, что 
$$|zu| = |z(\overline{M}\cap [z,x])|.$$ 
Следовательно, 
$$\beta (M, x) < |xu|.$$ 
Получили противоречие. Таким образом, равенство установлено. Осталось использовать соотношения из утверждений $(i)$, $(ii)$.
\vskip7pt
$(iv)\, \,$ Неравенство следует из неравенства $(iii)$ и утверждения $(v)$
леммы 1.
\vskip7pt
$(v)\, \,$ Заметим сначала, что функция 
$$g : \mathbb{R}_{+} \rightarrow \mathbb{R}_{+}, \quad g (R) = \beta (B[x,R], M)$$ 
--- липшицева в силу неравенства:  
$$|\beta (B[x,r],M) - \beta (B[x,R],M)| \leq |r - R|.$$ 
Предположим, что 
$$\beta (M, x) - r > \beta (B[x,r],M).$$ 
Тогда из утверждения $(i)$, липшицевости функции $g$ и
равенства 
$$\alpha (B[x,r],M) = \alpha (M,BC[X])$$ 
следует, что найдется такое $r_0 > r$, что 
$$\beta (M, x) - r_0 = \beta (B[x,r_0],M) < \beta (M, x) - r = \alpha (M,BC[X]).$$ Получили противоречие. Предположим
теперь, что 
$$\beta (M, x) - r < \beta (B[x,r],M).$$ 
Тогда найдется такое $r_0 < r$, что 
$$\beta (M, x) - r_0 = \beta (B[x,r_0],M) \leq  \beta (B[x,r],M).$$ Кроме того, из неравенств $(ii)$, $(iii)$ получим 
$$\beta (M, x) - r_0 = \beta (B[x,r_0],M) = \beta (B[x,r],M) = \alpha (M,BC[X])
=$$ 
$$\alpha (B[x,r_0],M) \leq \max \{0, r_0 + |xM| - |x(X\backslash M)|\}$$
и $r(M,x) \leq r_0$.
\vskip7pt
$(vi)\, \,$ Ограниченность множества $\chi (M)$ следует из предыдущего
свойства. Пусть последовательность $(x_n)_{n \in \mathbb{N}}$ из множества $\chi (M)$ сходится к точке $x \in X$. Тогда в силу
предыдущего свойства найдется такая ограниченная
последовательность неотрицательных вещественных чисел $(r_n)_{n \in \mathbb{N}}$, что 
$$\alpha (B[x_n,r_n],M) = \beta (M,x_n) - r_n = \beta (B[x_n,r_n],M) = \alpha (M,BC[X]).$$ 
Из этой последовательности
выберем сходящуюся к числу $r \in  \mathbb{R}_{+}$
\\
подпоследовательность $(r_m)_{m \in \mathbb{N}}$. Тогда
$$|\alpha (B[x_m,r_m],M) - \alpha (B[x,r],M)| \leq
|x_mx| + |r_m - r| \rightarrow 0$$ 
при $n \rightarrow \infty$. В силу непрерывности псевдометрики $\alpha$ получим 
$$\alpha (B[x,r],M) =  \alpha (M,BC[X])$$ 
и $x \in \chi (M)$. Второе утверждение следует из того, что в собственном метрическом пространстве каждое замкнутое и ограниченное множество компактно \cite [c. 18]{Bus}.
\vskip7pt
$(vii)\, \,$ Если пространство $X$ изометрично множеству вещественных
чисел со стандартной метрикой, то утверждение $(vii)$, очевидно,
верное. Рассмотрим случай, когда такой изометрии нет. Из условий и
утверждения $(iv)$ леммы 2 следует, что множество $P_M (x)$ ---
одноточечное. Пусть $u$ --- точка из пересечения луча
$\Lambda [P_M (x),x)$ с началом в точке $P_M (x)$, содержащего точку $x$, со
сферой $S(x,r)$ такая, что $x \in [P_M (x),u]$. Если 
$P_M (u) = P_M (x)$, то 
$$r + |xM| = |uP_M (x)| \leq \beta (B[x,r],M).$$ 
Осталось использовать неравенство $(ii)$. Докажем, что предположение о
справедливости неравенства
$$P_M (u) \neq P_M (x)$$ 
приводит к противоречию. Из этого неравенства и
утверждения $(iv)$ леммы 2 следует, что 
$$|xP_M (u)| > |xP_M (x)|.$$ Пусть
$z = z(t)$ --- параметризация луча 
$$\Lambda = \Lambda [P_M (u),P_M (x))$$ 
с началом в точке $P_M (u)$, содержащего точку $P_M (x)$, длиной
дуги. Из доказательства утверждения $(iv)$ леммы 2 следует
выпуклость функции $f = |xz(t)|$. Тогда найдется такая точка
$v = P_{\Lambda} (x)$, что либо 
$$v \in \Lambda\backslash [P_M (u),P_M (x)]\, \, \mbox {и}\, \, |xv| < |xP_M (x)|,$$ 
либо $v = P_M (x)$. Рассмотрим первый случай. Тогда
$$|uv| \leq |xv| + |xu| < |xP_M (x)| + |xu| = |uP_M (x)|.$$ 
Кроме того, в силу выпуклости функции $g = |uz(t)|$ и неравенства $(36.8)$ из \cite [с. 307]{Bus} имеют место неравенства
$$|uP_M (x)|\cdot |vP_M (u)| \leq |uP_M (u)|\cdot |vP_M (x)| + |uv|\cdot
|P_M (x)P_M (u)| <$$ 
$$|uP_M (x)|\cdot |vP_M (u)|.$$ 
Получили противоречие. Во втором
случае также получаем противоречие, так как из теоремы 1.3.6
следует, что существует единственный перпендикуляр к полупрямой
$P$, содержащий точку $u$ (см. также доказательство утверждения
(20.10) из \cite [c. 158]{Bus}).
\vskip7pt
$(viii)\, \,$ Это следствие утверждений $(i)$ и $(vii)$.
\vskip7pt
$(i\mathrm x)\, \,$ Это следствие утверждений $(v)$ леммы 1 и $(viii)$.
\vskip7pt
$(\mathrm x)\, \,$ Воспользуемся равенством $(viii)$. Тогда
$$\alpha (B[x,r],M) = \max \{r + |xM|,\beta (M, x) - r\} \leq \alpha (M,BC[X]) +
\tau.$$
Рассмотрим два случая.
\vskip7pt
$1.\, \,$ Если 
$$r + |xM| \geq \beta (M, x) - r,$$ 
то 
$$r \geq \dfrac{1}{2} (\beta (M, x) - |xM|)$$
и
$$\alpha (M,BC[X]) + \tau \geq r + |xM| \geq \dfrac{1}{2}  (\beta (M, x) + |xM|) =
\dfrac{1}{2} (\beta (M, x) + |xP_M (x)|) \geq$$ 
$$\dfrac{1}{2} \beta (M,P_M (x)) = \alpha (B[P_M (x),r(M,P_M (x))],M).$$

$2.\, \,$ Если 
$$r + |xM| \leq \beta (M, x) - r,$$ 
то 
$$r \leq \dfrac{1}{2} (\beta (M, x) - |xM|)$$
и
$$\alpha (M,BC[X]) + \tau \geq \beta (M, x) - r \geq \dfrac{1}{2} (\beta (M, x) + |xM|) =$$
$$\dfrac{1}{2} (\beta (M, x) + |xP_M (x)|) \geq \dfrac{1}{2} \beta (M,P_M (x)) = \alpha (B[P_M (x),r(M,P_M (x))],M).$$
Теорема 1 доказана.
\vskip7pt
Следующие лемма 5 и теорема 2 являются обобщением теоремы 4.3 из
\cite {Dudov}.
\vskip7pt
{\bf Лемма 5} [s17]. {\it Пусть $X$ геодезическое пространство, через каждые две различные точки которого можно провести единственную прямую, удовлетворяет условию $(A_3)$, $M$ --- выпуклое ограниченное множество
существования  в $X$, $x \in \chi (M)$. Тогда 
$$P_M (x) \in \chi (M).$$
Если, кроме того, существуют такие точки $u, \, v \in X$, что
$$\beta (M, x) = |ux|,\quad \beta(M,P_M (x)) = |vP_M (x)|,$$ 
то $x \in M$.}
\vskip7pt
{\bf Доказательство леммы 5}.
\vskip7pt 
В случае, когда $x \in M$ доказательство леммы
очевидно. Предположим, что 
$$x \in \chi (M)\backslash M.$$ Тогда из
утверждений $(v)$, $(vii)$ теоремы 1 следует, что найдется такое
$r_0 \geq 0$, что 
$$r(M,x) \leq r_0 \leq r,\quad \alpha (B[x,r_0],M) = \beta (M, x) - r_0 =$$ 
$$\beta (B[x,r_0],M) = \alpha (M,BC[X]) = r_0 + |xM|.$$ 
Следовательно, 
$$r_0 = \dfrac{1}{2}(\beta (M, x) - |xM|) = r(M,x)$$ 
и в силу утверждения $(iv)$ теоремы 1 имеют место равенства
$$\alpha (M,BC[X]) = r(M,x) + |xM| = \psi (M,x) \leq
\psi (M,P_M (x)) = \dfrac{1}{2} \beta(M,P_M (x)).$$ 
Но
$$2\psi (M,P_M (x)) = \beta(M,P_M (x)) \leq \beta (M, x) + |xP_M (x)| = 2\alpha (M,BC[X]),$$
поэтому
$$P_M (x) \in \chi (M).$$ 
Если, кроме того, существуют такие
точки $u,\, v \in X$, что 
$$\beta (M, x) = |ux|,\quad \beta(M,P_M (x)) = |vP_M (x)|,$$ 
то из леммы 4 получим $x = P_M (x)$. Получили противоречие с нашим
предположением. Следовательно, $x \in M$. Лемма 5 доказана.
\vskip7pt 
{\bf Теорема 2} [s17]. {\it Пусть $X$ геодезическое пространство, через каждые две различные точки которого можно провести единственную прямую, удовлетворяет условию $(A_3)$, $M$ --- выпуклый компакт в $X$. Тогда
$$\chi (M) \neq \emptyset\, \, \mbox {и} \, \, \chi (M) \subset M.$$}

{\bf Доказательство теоремы 2}.
\vskip7pt 
Очевидно, что в пространстве $X$ найдется такая
последовательность замкнутых шаров $(B[x_n,r_n])_{n \in \mathbb{N}}$, что $$\alpha (B[x_n,r_n],M) \leq \alpha (M,BC[X]) + 1/n.$$ 
В силу утверждения $(\mathrm x)$ теоремы 1 можно считать, что для каждого $n \in \mathbb{N}$
$x_n \in M.$ 
Тогда последовательность
$(r_n)_{n \in \mathbb{N}}$ ограничена и найдутся
подпоследовательности 
$$(x_m)_{m \in \mathbb{N}},\quad (r_m)_{m \in \mathbb{N}}$$
 последовательностей 
$$(x_n)_{n \in \mathbb{N}},\quad (r_n)_{n \in \mathbb{N}}$$
соответственно, точка $x \in M$ и число $r \in  \mathbb{R}_{+}$ такие,
что
$$\lim \limits_{m \rightarrow \infty} |x_mx| = 0,\quad \lim \limits_{m \rightarrow \infty} r_m = r.$$ 
Если теперь учесть непрерывность метрики $\alpha$ и для каждого 
$m \in \mathbb{N}$ равенство 
$$\alpha (B[x_m,r_m],B[x,r]) = |x_mx| + |r_m - r|,$$ 
то получим 
$$\alpha (B[x,r],M) = \alpha (M,BC[X])\, \, \mbox {и} \, \, x \in \chi (M).$$
Второе утверждение следует из леммы 5. Теорема 2 доказана.
\vskip7pt
Из теорем 1, 2 получим следствие 1.
\vskip7pt
{\bf Следствие 1} [s17]. {\it Пусть $X$ геодезическое пространство, через каждые две различные точки которого можно провести единственную прямую, удовлетворяет условию $(A_3)$, $x,\, y \in X$. Тогда
$$\chi ([x,y]) = \{\omega (x,y)\}\, \, \mbox {и} \, \, r(M,\omega (x,y)) = \dfrac{1}{4}|xy|$$ 
--- радиус шара наилучшего приближения сегмента $[x,y]$ в метрике
$\alpha$.}
\vskip7pt
Следующая теорема 3 является аналогом теоремы 4.6 из \cite {Dudov}.
\vskip7pt
{\bf Теорема 3} [s17]. {\it Пусть $X$ геодезическое пространство, через каждые две различные точки которого можно провести единственную прямую, удовлетворяет условию $(A_3)$, $(M_n)_{n \in \mathbb{N}}$ ---
последовательность ограниченных выпуклых множеств существования,
сходящаяся в метрике Хаусдорфа к компакту $M \subset X$, и для каждого $n \in \mathbb{N}$
$$\chi (M_n) \neq \emptyset.$$ 
Тогда
$$\lim \limits_{n \rightarrow \infty} \beta (\chi (M_n)\cap M_n,\chi (M)) = 0.$$}

{\bf Доказательство теоремы 3}.
\vskip7pt 
Заметим прежде всего, что в силу условий теоремы 3
и леммы 5 для каждого $n \in \mathbb{N}$ 
$$\chi (M_n)\cap M_n \neq \emptyset.$$
Используем метод доказательства от противного. Пусть теорема 3 неверна. Тогда найдутся константа $c > 0$ и
подпоследовательность $(\chi (M_m)\cap M_m)_{m \in \mathbb{N}}$
 последовательности  $ (\chi (M_n)\cap M_n)_{n \in \mathbb{N}}$ такие, что $$\beta (\chi (M_m)\cap M_m,\chi (M)) > c.$$ 
Следовательно, для каждого 
$m \in \mathbb{N}$ найдется точка 
$$y_m \in \chi (M_m)\cap M_m$$ 
такая, что $|y_m\chi (M)| > c$. Из условий теоремы 3 следует, что найдутся
подпоследовательность $(y_l)_{l \in \mathbb{N}}$
последовательности $(y_m)_{m \in \mathbb{N}}$ и точка $y \in M$
такие, что 
$$\lim \limits_{l \rightarrow \infty}|y_ly| = 0.$$ 
Кроме того, в силу утверждения $(v)$ теоремы 1 для каждого 
$l \in \mathbb{N}$ найдется такое число $r_l \in  \mathbb{R}_{+}$, что
$$\alpha (B[y_l,r_l],M_l) = \beta(M_l,y_l) - r_l= \beta
(B[y_l,r_l],M_l) = \alpha (B[y_l,r_l],BC[M]).$$
Следовательно, последовательность $(r_l)_{l \in \mathbb{N}}$ ограничена. Выделим подпоследовательность $(r_k)_{k \in \mathbb{N}}$ последовательности $(r_l)_{l \in \mathbb{N}}$, сходящуюся к
числу $r \in  \mathbb{R}_{+}$. Тогда 
$$\lim \limits_{k \rightarrow \infty} \alpha (B[y_k,r_k],M_k)
= \alpha (B[y,r],M)\, \, \mbox {и} \, \, \alpha (B[y,r],M) = \alpha (B[y,r],BC[M]).$$ Таким образом,
$$\lim \limits_{k \rightarrow \infty}|y_k\chi (M)| = 0.$$ 
Получили противоречие. Теорема 3 доказана.

\section{Касательное пространство по Буземану в точке геодезического пространства}

\vskip20pt

Рассмотрим метрическое пространство $(X,\rho)$ с выделенным семейством $S$ сегментов, называемых хордами, удовлетворяющее условиям  $(C_1)$, $(C_4)$, и для определения касательного пространства по
Буземану в точке $p \in X$ введем на множествах $B(p, r_p)$, $B(p, r_p) \times \mathbb{R}_{+}$ специальные псевдометрики.
\vskip7pt
{\bf Лемма 1} [s18]. {\it Пусть пространство $(X,\rho)$ удовлетворяет
условиям $(C_1)$, $(C_4)$, $p \in X$. Тогда верны следующие утверждения.
\vskip7pt

$(i) \, \,$ Функция
$$\overline{m}_p : B(p, r_p)\times B(p, r_p) \rightarrow \mathbb{R}_{+}, \quad \overline{m}_p (x, y) =
\overline{\lim \limits}_{\lambda \to 0+}{\dfrac{|\omega_{\lambda} (p,
x)\omega_{\lambda}(p, y)|}{\lambda}}$$ 
является псевдометрикой, обладающей следующими свойствами: для любых 
$$\lambda \in [0, 1]; \quad x, \, y \in  B(p, r_p)$$ 
$$\overline{m}_p (p, x) = |px|,\quad\overline{m}_p (\omega_{\lambda}(p, x),\omega_{\lambda} (p, y)) = \lambda \overline{m}_p( x, y).$$

$(ii) \, \,$ Функция
$$m_p : (B(p, r_p) \times \mathbb{R}_{+}) \times (B(p, r_p) \times
\mathbb{R}_{+}) \rightarrow \mathbb{R}_{+},$$
$$m_p ((x;\lambda),(y;\mu)) =
\tau\overline{m}_p (\omega_{\lambda/\tau}(p,x),\omega_{\mu/\tau}(p,
y)),$$
где $\tau \in  (max \{\lambda,\mu\},+\infty)$, является псевдометрикой.}
\vskip7pt
{\bf Доказательство леммы 1.} 
\vskip7pt
Сводится к непосредственной проверке.
\vskip7pt
Введем на множестве $B(p, r_p) \times \mathbb{R}_{+}$ следующее отношение эквивалентности.
$$(x;\lambda)\sim (y;\mu),\, \, \mbox {если} \, \, m_p ((x;\lambda),(y;\mu)) = 0.$$ Тогда псевдометрика $m_p$ индуцирует на фактор-пространстве   
$$X_p = (B(p, r_p) \times \mathbb{R}_{+})/\!\!\sim$$ 
метрику, которую обозначим тем же символом: $m_p$.  Сначала исследуем свойства псевдометрики $\overline{m}_p$ и метрики $m_p$ в пространстве 
неположительной кривизны по Буземану.
\vskip7pt
Пространство $(X,\rho)$, удовлетворяющее условиям $(C_1)$, $(C_4)$, назовем
{\it пространством неположительной кривизны по Буземану} (\cite [с. 304]{Bus},  \cite [c. 63]{Bus2}), если для любых $x, \,y,\, z \in  B(p, r_p)$ имеет место неравенство 
$$2 |\omega_{1/2}(z,x)\omega_{1/2}(z,y)| \leq |xy|.$$

{\bf Теорема 1} [s18]. {\it Пусть $(X,\rho)$ пространство неположительной
кривизны по Буземану. Тогда верны следующие утверждения.
\vskip7pt
$(i)\, \, $ Для любых $x, \, y \in B(p,r_p)$
$$\overline{m}_p (x, y) = \lim \limits_{\nu \to 0+}
{\dfrac{|\omega_{\nu} (p, x)\omega_{\nu}(p,y)|}{\nu}}.$$

$(ii)\, \, $ Для любых $x, \, y \in B(p,r_p)$
$$\overline{m}_p (x, y) \leq |xy|.$$

$(iii)\, \, $ $m_p$ есть внутренняя метрика на $X_p$.
\vskip7pt
$(iv)$ Если, кроме того, пространство $X$ локально компактно, то для любого
$p \in X \, \,$ $(X_p,m_p)$ --- собственное геодезическое пространство.}
\vskip7pt
{\bf Доказательство теоремы 1.} 
\vskip7pt
$(i-ii)\, \,$ Пусть геодезическое пространство $X$ удовлетворяет условиям $(C_1)$, $(C_4)$ и для любых $x,\,  y,\,  z \in B(p, r_p)$ имеет место неравенство: 
$$2 |\omega_{1/2} (z, x)\omega_{1/2} (z, y)|\leq |xy|.$$ 
Из этих условий и неравенства нетрудно получить
следующее неравенство:
$$|\omega_{\lambda}(z,x)\omega_{\lambda}(z,y)| \leq \lambda |xy| \eqno (1)$$ 
для любых $0 \leq \lambda \leq 1$, $x,\,  y,\,  z \in B(p,r_p)$. Следовательно, для любых 
$x,\,  y \in  B(p, r_p)$ справедливо неравенство: 
$$\overline{m}_p ( x, y) \leq |xy|.$$ 
Кроме того, из доказательств теорем
(36.4 - 36.6) в  \cite [с. 304-307]{Bus} следует, что для любых $x, \, y,\,  z
\in  B(p,r_p)$ непрерывная функция 
$$f : [0,1] \rightarrow \mathbb{R},\quad f (\lambda) = |\omega_{\lambda}(z,x)\omega_{\lambda}(z,y)|$$
выпуклая. Следовательно, для любых $x, \, y \in  B(p,r_p)$ существует предел 
$$\lim\limits_{\nu \to 0+}{\dfrac{|\omega_{\nu}(p, x)\omega_{\nu}(p,y)|}{\nu}}.$$

$(iii)\, \,$ Выберем произвольно $[x_1;\lambda],\,  [y;\mu] \in X_p$,
где для определенности $\lambda \leq \mu$. Если $\mu = 0$, то эти классы совпадают и утверждение очевидно. Пусть $\mu > 0$. Тогда
$$[x_1;\lambda] = [\omega_{\lambda/\mu}(p,x_1);\mu].$$ 
Обозначим $x = \omega_{\lambda/\mu} (p,x_1)$. В силу леммы 1.1.1 достаточно
доказать, что для каждого $\varepsilon >0$ найдется такая точка
$[z;\tau] \in X_p$, что 
$$2 \max \{m_p ([x;\mu],[z;\tau]), m_p([y;\mu],[z;\tau ])\} < m_p ([x;\mu],[y;\mu]) + \varepsilon.$$
Очевидно, найдется такое $\nu_0 \in (0,1]$, что для каждого $\nu \in (0,\nu_0)$
$$z_{\nu} = \omega_{1/2} (\omega_{\nu} (p,x),\omega_{\nu} (p,y)) \in B(p,r_p).$$
В силу неравенства $(ii)$ получим 
$$\max \{\overline{m}_p (\omega_{\nu} (p,x),z_{\nu}),\overline{m}_p (\omega_{\nu} (p,y),z_{\nu})\} \leq$$ 
$$\max \{|\omega_{\nu} (p,x)z_{\nu}|,|\omega_{\nu} (p,y)z_{\nu}|\} =
\dfrac{1}{2} |\omega_{\nu} (p,x)\omega_{\nu} (p,y)|. \eqno  (2)$$
Отсюда и из леммы 1 получим
$$\nu\overline{m}_p (x,y) = \overline{m}_p (\omega_{\nu} (p,x),\omega_{\nu} (p,y)) \leq$$
$$\overline{m}_p (\omega_{\nu} (p,x), z_{\nu}) + \overline{m}_p (\omega_{\nu} (p,y),z_\nu) \leq |\omega_{\nu} (p,x)\omega_{\nu} (p,y)|.$$ 
Из этих неравенств и из утверждения $(i)$ получим
$$\lim\limits_{\nu \to 0+}{\dfrac{\overline{m}_p (\omega_{\nu} (p,x),z_{\nu}) +
\overline{m}_p (\omega_{\nu} (p,y),z_{\nu})}{\nu}} = \overline{m}_p (x,y).$$
Если учтем $(2)$, то получим 
$$\lim\limits_{\nu \to 0+}{\dfrac{\overline{m}_p (\omega_{\nu} (p,x),z_{\nu})}{\nu}}
= \lim\limits_{\nu \to 0+}{\dfrac{\overline{m}_p (\omega_{\nu} (p,y),z_{\nu})}{\nu}} =
 \dfrac{1}{2} \overline{m}_p (x,y).  \eqno (3)$$
Следовательно, для каждого
$\varepsilon > 0$ найдется такое $\nu \in (0,\min\{\delta,1\})$,
что
$$\dfrac{\overline{m}_p(\omega_{\nu} (p,x),z_{\nu})}{\nu} <
\dfrac{1}{2} \overline{m}_p (x,y) + \varepsilon,$$ 
$$\dfrac{\overline{m}_p (\omega_{\nu} (p,y),z_{\nu})}{\nu} <
\dfrac{1}{2} \overline{m}_p(x,y) + \varepsilon.$$
Тогда для каждого $\varepsilon > 0$ найдется такое $\nu \in (0,\min \{\delta,1\})$,
что 
$$2 \max \{m_p ([x;\mu],[z_{\nu};\mu/\nu]), m_p ([y;\mu],[z_{\nu};\mu/\nu])\} =$$
$$\dfrac{2\mu}{\nu}
\max\{\overline{m}_p (\omega_{\nu} (p,x),z_{\nu}),\overline{m}_p (\omega_\nu (p,y),z_{\nu})\} <$$ 
$$\mu \overline{m}_p (x,y)+ 2\mu\varepsilon = m_p ((x;\mu],[y;\mu]) +
2\mu\varepsilon.$$

$(iv)\, \,$ Пусть $X$ локально компактное
пространство неположительной кривизны по Буземану и
$([x_n;\lambda_n])_{n \in \mathbb{N}}$ произвольная ограниченная
последовательность пространства $(X_p,m_p)$. Заметим, что
для $\tau \geq \max \{1,\lambda_n\}$ имеют место равенства
$$m_p ([x_n;\lambda_n],[p;1]) = \tau \overline m_p (\omega_{\lambda_n/\tau}(p,x_n),p)=  \lambda_n |px_n|,$$
поэтому последовательность $(\lambda_n)_{n \in \mathbb{N}}$ --- ограничена. 
Тогда найдутся такие вещественные числа $r \in (0,r_p]$, $\nu > 0$, что замкнутый шар $B[p,r]$ компактен и для каждого $n \in \mathbb{N}$ 
$$\omega_{\lambda_n/\nu} (p,x_n) \in B[p,r].$$
Выделим подпоследовательность
$(\omega_{\lambda_k/\nu} (p,x_k))_{k \in \mathbb{N}}$ последовательности
$(\omega_{\lambda_n/\nu}(p,x_n))_{n \in \mathbb{N}}$, сходящуюся в исходной метрике $\rho$ к точке $y \in B[p,r]$. 
Отсюда, из леммы 1 и неравенства $(1)$ следует, что
$$m_p ([x_k;\lambda_k],[y;\nu]) =
\nu \overline{m}_p (\omega_{\lambda_k/\nu} (p,x_k),y)  \leq
\nu |\omega_{\lambda_k/\nu} (p,x_k)y| \rightarrow 0$$
при $k \rightarrow \infty.$ Таким образом, последовательность $([x_n;\lambda_n])_{n \in \mathbb{N}}$ обладает
сходящейся подпоследовательностью и пространство
$(X_p,m_p)$ собственное. 
Из следствия обобщенной альтернативы Хопфа-Ринова \cite [с. 297]{Konf} и того, что пространство $(X_p,m_p)$ собственное следует, что пространство $(X_p,m_p)$ --- геодезическое. Покажем, что в данном случае можно дать и непосредственное доказательство этого факта, избежав ссылки на обобщенную альтернативу Хопфа-Ринова. 

Выберем произвольно $[x_1;\lambda], \,  [y;\mu] \in X_p$,
где для определенности $\lambda \leq \mu$. Если $\mu = 0$, то эти классы совпадают и утверждение очевидно. Пусть $\mu > 0$. Тогда
$$[x_1;\lambda] = [x;\mu],\, \, \mbox {где} \, \, x = \omega_{\lambda/\mu} (p,x_1)$$
и найдется такое число $N \in \mathbb{N}$, что для каждого натурального $n > N$
$$z_n = \omega_{1/2} (\omega_{1/n} (p,x),\omega_{1/n} (p,y)) \in B[p,r_p].$$
Кроме того, в силу леммы 1 для каждого 
$$n > \max \{N,\dfrac{1}{\mu}\}$$
$$m_p ([p;1],[z_n;n\mu]) = n \mu \overline{m}_p (p,z_n) = n\mu |pz_n| \leq$$ 
$$n\mu (|p\omega_{1/n} (p,x)| + |\omega_{1/n} (p,x) z_n|) \leq
\mu (|px| + \dfrac{n}{2} |\omega_{1/n} (p,x)\omega_{1/n} (p,y)|) \leq $$ 
$$\mu (|px| + \dfrac{1}{2}(|px| + |py|)) \leq 2 \mu r_p.$$ 
Таким образом, последовательность $([z_n;n\mu])_{n \in \mathbb{N}}$ ограничена в собственном пространстве $(X_p,m_p)$. Следовательно, у нее найдется подпоследовательность
$$([z_{n_k}; n_k \mu])_{k \in \mathbb{N}},$$ 
сходящаяся относительно
метрики $m_p$ к некоторой точке $[z;\tau]  \in X_p$. Повторяя
рассуждения в доказательстве утверждения $(iii)$,  получим
равенства $(3)$, где нужно подставить $\nu = 1/n_k$.  Тогда
$$\dfrac{1}{2} m_p ([x;\mu],[y;\mu]) = \dfrac{\mu}{2} \overline{m}_p (x,y) =$$
$$\mu \lim\limits_{k \to \infty}{(n_k \overline{m}_p (\omega _{1/n_k} (p,x), z_{n_k})) } =
\lim\limits_{k \to \infty}{m_p ([x;\mu],[z_{n_k};n_k \mu])}
= m_p ([x;\mu],[z;\tau]), $$
$$\dfrac{1}{2} m_p ([x;\mu],[y;\mu]) = \dfrac{\mu}{2} \overline{m}_p (x,y) =$$
$$\mu \lim\limits_{k \to \infty}{(n_k\overline{m}_p (\omega _{1/n_k} (p,y), z_{n_k})) } =
\lim\limits_{k \to \infty}{m_p ([y;\mu],[z_{n_k};n_k \mu])}
= m_p ([y;\mu],[z;\tau]). $$
Но известно \cite {Efrem}, что если для любых двух точек
$u, \, v$ из полного метрического  пространства  $(Y,\rho)$  найдется точка
$w \in Y$ такая, что $|uw| = |vw| = \dfrac{1}{2}|uv|$, то это пространство
является геодезическим пространством. Теорема 1 доказана.
\vskip7pt
Напомним следующие определения.
\vskip7pt
Пространство $(X,\rho)$, удовлетворяющее условиям $(C_1)$, $(C_4)$ назовем
{\it дифференцируемым пространством в точке $p \in X$  по Буземану} (\cite [с. 293] {Bus}, \cite [с. 21]{Bus1}), если для каждого $\varepsilon > 0$ найдется такое $\delta \in (0, r_p)$, что для любого $\lambda \in [0, 1]$, для всех $x,\,  y \in B(p,\delta)$ имеет место неравенство:
$$||\omega_{\lambda}(p, x)\omega_{\lambda}(p,y)| - \lambda |xy|| \leq
\varepsilon \lambda |xy|.$$

Отображение $f$ из метрического пространства $(X,\rho)$ в
метрическое пространство $(Y,d)$ называется {\it билипшицевым
отображением}, если найдется такая константа
$\lambda \geq 1$, что для любых $x,\, y \in X$ справедливы неравенства: 
$$\dfrac{1}{\lambda}|xy| \leq d(f(x),f(y)) \leq \lambda|xy|.$$

Метрические пространства $(X,\rho)$, $(Y,d)$ называются {\it липшицево
эквивалентными}, если существует сюръективное билипшицево
отображение $f : (X,\rho) \rightarrow (Y,d)$ \cite [с. 269]{Bonk}.
\vskip7pt
Исследуем теперь свойства псевдометрики $\overline{m}_p$ и метрики $m_p$ в пространстве, дифференцируемом в точке $p \in X$ по Буземану.
\vskip7pt
{\bf Теорема 2} [s18]. {\it Пусть $(X,\rho)$ пространство, дифференцируемое в
точке $p \in X$ по \mbox{Буземану}. Тогда верны следующие утверждения.
\vskip7pt
$(i)\, \, $ Для всех $x, \, y \in B(p,r_p)$
$$\overline{m}_p (x, y) = \lim \limits_{\nu \to 0+}
{\dfrac{|\omega_{\nu}(p, x)\omega_{\nu} (p,y)|}{\nu}}.$$

$(ii)\, \, $ $\overline{m}_p$ есть метрика на открытом шаре $B(p, r_p)$
и для каждого $\varepsilon > 0$ найдется такое $\delta \in (0,
r_p)$, что для всех $x, \, y \in B(p,\delta )$ справедливо неравенство:
$$|\overline{m}_p (x, y) - |xy|| \leq \varepsilon |xy|.$$

$(iii)\, \, $ $m_p$ есть внутренняя метрика на $X_p$.
\vskip7pt
$(iv)\, \, $ Если, кроме того, $X$ является локально полным пространством, то
$(X_p,m_p)$ --- полное пространство.
\vskip7pt
$(v)\, \, $ Пусть $\delta$ то же, что и в пункте $(ii)$ при $0 < \varepsilon < 1$.
И пусть найдется такое число $t(p) \in (0,\delta)$, что для
каждой точки 
\\
$x \in B(p,t(p)) \backslash \{p\}$ существует точка $\hat x \in S(p,t(p))$
такая, что $x \in [p,\hat x]$, где $[p,\hat x] \in S$. Тогда для каждого
$R>0$ пространства $(B[p,t(p)],\rho)$, $(B[[p;1],R],m_p)$ липшицево
эквивалентны.}
\vskip7pt
{\bf Замечание 1.} {\it Нетрудно проверить, что если  $X$
дифференцируемое в точке $p$
$G$-пространство Буземана, то пространство $(X_p,m_p)$  изометрично
касательному пространству в точке $p$ по Буземану (при наличии в
пространстве $(X_p,m_p)$ линейной структуры эта изометрия
является линейной) {\rm \cite [c. 22]{Bus1}}. Поэтому метрическое пространство
$(X_p,m_p)$ естественно называть касательным пространством по
Буземану в точке $p \in X$.}
\vskip7pt
{\bf Доказательство теоремы 2.} 
\vskip7pt
$(i-ii)\, \,$ Пусть $X$ пространство, дифференцируемое в точке
$p \in  X$ по Буземану. Тогда для каждого $\varepsilon > 0$
найдется такое
$\delta(\varepsilon) \in (0, r_p)$, что для любых 
$$x,\, y \in B(p,\delta(\varepsilon)),\quad 0 < \lambda < \mu < 1$$ 
имеет место неравенство
$$||\omega_{\lambda}(p,x)\omega_{\lambda}(p,y)| -
\dfrac{\lambda}{\mu} |\omega_{\mu}(p,x)\omega_{\mu}(p,y)| |\leq \varepsilon
\dfrac{\lambda}{\mu} |\omega_{\mu}(p,x)\omega_{\mu}(p,y)|.$$
Тогда для любых 
$$x, \, y \in B(p,\delta(1/2)),\quad 0 < \mu < \mu_0 <
\dfrac{\delta(1/2)}{r_p}$$ 
справедливы неравенства
$$\dfrac{|\omega_{\mu_0} (p,x)\omega_{\mu_0} (p,y)|}{2\mu_0}\leq
\dfrac{|\omega_{\mu} (p,x)\omega_{\mu}(p,y)|}{\mu} \leq
\dfrac{3 |\omega_{\mu_0}(p,x)\omega_{\mu_0}(p,y)|}{2\mu_0}.$$
И для любых 
$$0 < \varepsilon < \dfrac{1}{2};\quad x, \, y \in B(p,\delta(\varepsilon));\quad
\mu < \min \{\dfrac{\delta(\varepsilon)}{r_p},\mu_0\}$$ 
имеют место неравенства
$$|\dfrac{|\omega_{\lambda}(p,x)\omega_{\lambda}(p,y)|}{\lambda }  -
\dfrac{|\omega_{\mu}(p,x)\omega_{\mu}(p,y)|}{\mu} |\leq 
\dfrac{\varepsilon |\omega_{\mu}(p,x)\omega_{\mu}(p,y)|}{\mu} \leq$$
$$\dfrac{3 \varepsilon |\omega_{\mu_0}(p,x)\omega_{\mu_0}(p,y)|}{2\mu_0}.$$
Из этих неравенств следует существование предела
$$\overline{m}_p ( x, y) = \lim\limits_{\nu \to 0+}{\dfrac{|\omega_{\nu} (p, x)\omega_{\nu} (p,
y)|}{\nu}}$$ 
и то, что $\overline{m}_p ( x, y) \neq 0$ при $x \neq y$. Если обе части неравенства из определения дифференцируемости в
точке $p \in X$  по Буземану поделить на $\lambda > 0$ и перейти к
пределу при $\lambda \to 0+$, то для любых
$x, \, y \in  B(p,\delta)$ получим неравенство
$$|\overline{m}_p ( x, y) - |xy|| \leq \varepsilon |xy|.$$

$(iii)\, \,$ Выберем произвольно $[x_1;\lambda],\, [y;\mu] \in X_p$ и для
определенности считаем, что $\lambda \leq \mu$. Тогда
$$[x_1;\lambda] = [\omega_{\lambda/\mu}(p,x_1);\mu].$$ 
Обозначим через $x = \omega_{\lambda/\mu} (p,x_1)$. В силу леммы 1.1.1 достаточно
доказать, что для каждого $\varepsilon >0$ найдется такая точка
$[z;\tau ] \in X_p$, что 
$$2 \max \{m_p ([x;\mu],[z;\tau]), m_p([y;\mu],[z;\tau ])\} < m_p([x;\mu],[y;\mu]) + \varepsilon.$$
Очевидно, найдется такое $\nu_0 \in (0,1]$, что для каждого $\nu \in (0,\nu_0]$
$$z_\nu =
\omega_{1/2} (\omega_{\nu} (p,x),\omega_{\nu} (p,y)) \in B(p,r(p)).$$
Из утверждения $(ii)$ получим, что для каждого
$\varepsilon > 0$ существует такое $\delta > 0$, что найдется
$\nu_0 \in (0,1]$ такое, что для любого $\nu \in (0,\nu_0]$ 
$$\omega_{\nu} (p,x), \, \omega_{\nu} (p,y), \,
\omega_{1/2} (\omega_{\nu} (p,x),\omega_{\nu} (p,y)) \in B(p,\delta)$$ 
и
$$\overline{m}_p (\omega_{\nu}(p,x),\omega_{\nu}(p,y)) \leq (1 +
\varepsilon) |\omega_{\nu} (p,x)\omega_{\nu}(p,y)|.$$
Тогда для любых $\varepsilon > 0$, $\nu  \in (0,\nu_0]$ справедливы неравенства:
$$\max \{\overline{m}_p (\omega_{\nu} (p,x), z_{\nu}),
\overline{m}_p (\omega_{\nu} (p,y), z_{\nu})\} \leq$$
$$(1 + \varepsilon) \max \{|\omega_{\nu} (p,x)z_{\nu}|,|\omega_{\nu} (p,y)z_{\nu}|\} = \dfrac{(1 + \varepsilon)}{2} |\omega_{\nu}(p,x)\omega_{\nu} (p,y)|. \eqno (4)$$
Отсюда и из леммы 1 получим для любых $\varepsilon
> 0$, $\nu  \in (0,\nu_0]$ справедливы неравенства:
$$\nu \overline{m}_p (x,y) = \overline{m}_p (\omega_{\nu} (p,x),\omega_{\nu} (p,y)) \leq$$ 
$$\overline{m}_p (\omega_{\nu} (p,x), z_{\nu}) +
\overline{m}_p(\omega_{\nu} (p,y),z_{\nu}) \leq (1 +
\varepsilon)|\omega_{\nu} (p,x)\omega_{\nu} (p,y)|.$$ 
Если учтем $(4)$, то для любого $\varepsilon > 0$ получим
$$\overline{m}_p (x,y) \leq
\overline{\lim\limits_{\nu \to 0+}}{\dfrac{\overline{m}_p (\omega_{\nu} (p,x),z_{\nu}) +
\overline{m}_p (\omega_{\nu} (p,y),z_{\nu})}{\nu}} \leq (1 + \varepsilon)
\overline{m}_p(x,y),$$
$$\overline{\lim\limits_{\nu \to
0+}}{\dfrac{\max \{\overline{m}_p (\omega_{\nu} (p,x),z_{\nu}),
\overline{m}_p (\omega_{\nu} (p,y), z_{\nu})\}}
{\nu}} \leq \dfrac{(1 + \varepsilon)}{2}\overline{m}_p (x,y).$$ 
Аналогичные неравенства верны и для нижних
пределов. В силу произвольности $\varepsilon > 0$
$$\lim\limits_{\nu \to 0+}{\dfrac{\overline{m}_p (\omega_{\nu}(p,x),z_{\nu})}{\nu}} =
\lim\limits_{\nu \to 0+}{\dfrac{\overline{m}_p (\omega_{\nu} (p,y), z_{\nu})}{\nu}} = \dfrac{1}{2}\overline{m}_p(x,y).$$ 
Оставшаяся часть доказательства
совпадает с частью доказательства утверждения $(iii)$ теоремы 1
после равенств $(3)$.
\vskip7pt
$(iv)\, \,$  Пусть $X$ локально полное дифференцируемое в точке
$p \in X$ по Буземану пространство и  $(B[p,r],\rho)$, где $0 < r <
\delta$, полное подпространство пространства $(X,\rho)$. Из неравенств
$$(1-\varepsilon)|xy| \leq \overline{m}_p ( x, y) \leq (1+\varepsilon)|xy|$$  
для любых 
$$x,\, y \in  B[p,r],\quad 0< \varepsilon <1,$$ 
полученных в пункте $(ii)$, следует липшицева эквивалентность пространств
$$(B[p,r],\overline{m}_p),\quad (B[p,r],\rho)$$ 
и полнота пространства $(B[p,r],\overline{m}_p)$.  Пусть $([x_n;\lambda_n])_{n \in \mathbb{N}}$ фундаментальная последовательность в пространстве $(X_p,m_p)$. Тогда отсюда, из ограниченности фундаментальной последовательности
$([x_n;\lambda_n])_{n \in \mathbb{N}}$ и леммы 1 следует, что найдется число $\tau > 0$ такое, что 
$$\omega_{\lambda_n/\tau}(p,x_n) \in
(B[p,r],\overline{m}_p)\, \, \mbox {и} \, \, (\omega_{\lambda_n/\tau} (p,x_n))_{n \in \mathbb{N}}$$
является фундаментальной последовательностью. Пусть эта
последовательность сходится к точке $z \in B[p,r]$. Тогда
$$m_p ([x_n;\lambda_n],[z;\tau]) = \tau \overline
m_p (\omega_{\lambda_n/\tau}(p,x_n),z) \rightarrow 0$$ 
при $n \rightarrow \infty$. Таким образом, метрическое пространство
$(X_p,m_p)$ полное.
\vskip7pt
$(v)\, \,$ Докажем, что искомую липшицеву эквивалентность определяет
следующая композиция
$$h\circ f\circ id : (B[p,t(p)],\rho) \rightarrow (B[[p;1],R],m_p)$$
отображений
$$id : (B[p,t(p)],\rho) \rightarrow (B[p,t(p)],\overline m_p);$$ 
$$f : (B[p,t(p)],\overline m_p) \rightarrow (B[[p;1],t(p)],m_p),\quad f(x)=[x;1];$$
$$h : (B[[p;1],t(p)],m_p) \rightarrow (B[[p;1],R],m_p),\quad h ([x;\lambda])=[x;\dfrac{R \lambda}{t(p)}].$$
Рассмотрим каждое из этих отображений. Тождественное отображение
$$id : (B[p,t(p)],\rho) \rightarrow (B[p,t(p)],\overline m_p)$$ 
определяет, как было отмечено в начале доказательства пункта $(iv)$, липшицеву эквивалентность соответствующих пространств. Докажем, что $f$ есть изометрия. В силу равенств
$$m_p (f(x),f(y)) = m_p ([x;1],[y;1]) = \overline m_p (x,y)$$ 
для всех $x$, $y \in B[p,t(p)]$, отображение $f$ есть изометрическое вложение.
Осталось доказать, что оно сюръективно. Выберем произвольно
$[x;\mu] \in B[[p;1],t(p)]$. Если $\mu \leq 1$, то 
$$f (\omega_{\mu}[p,x]) = [\omega_{\mu}(p,x);1] = [x;\mu].$$ 
Если $\mu > 1$ и $\lambda = \dfrac{\mu |px|}{t(p)}$, то 
$$f(\omega_{\lambda} (p,\hat x)) = [\omega_{\lambda} (p,\hat x);1] =
[x;\mu].$$ 
Следовательно, отображение $f$ сюръективно. Отображение
$h$ есть подобие (определение см. в \cite [с. 286]{Bus}). Действительно, отображение $h$, очевидно, обратимо и
для любых $[x;\lambda]$, $[y;\mu] \in B[[p;1],t(p)]$ имеют место равенства:
$$m_p(h([x;\lambda]),h([y;\mu])) = m_p ([x;\dfrac{R\lambda}{t(p)}],[y;\dfrac{R\mu}{t(p)}]) =
\dfrac{R}{t(p)} m_p ([x;\lambda], [y;\mu]).$$  
Следовательно, отображение
$h\circ f\circ id$ определяет липшицеву эквивалентность пространств
$(B[p,t(p)],\rho)$, $(B[[p;1],R],m_p)$. Теорема 2 доказана.
\vskip7pt
{\bf Пример.} 
\vskip7pt
Рассмотрим модель Бельтрами-Клейна геометрии Лобачевского в шаре $B(0,r)$ вещественного гильбертова пространства. В силу теоремы 1.4.4 
$$\overline{m}_0 (x, y) = ||y \dfrac{\rho (0,y)}{||y||} - x \dfrac{\rho (0,x)}{||x||}|| = \leqno (ii)$$
$$\sqrt{\rho^2 (0,x) + \rho^2 (0,y) - 2 \rho (0,x) \rho (0,y) \cos \alpha},$$
где $x,\, y \in B(0,r)\backslash \{0\}$, $\alpha$ --- величина угла между векторами $x,\, y$.

%% file: ch33.tex
\chapter{Специальные отображения метрических пространств}
\vskip20pt

В параграфе 3.1 рассматривается метрическое пространство слабо ограниченных отображений метрических пространств с метрикой Куратовского $\delta$. Доказано, что пространство всех слабо ограниченных гомеоморфизмов $(HB(X),\delta)$, каждый из которых равномерно непрерывен на произвольном замкнутом шаре с центром в фиксированной точке метрического пространства $X$ вместе со своим обратным гомеоморфизмом, является паратопологической группой, непрерывно действующей на пространстве $X$. Установлено, что $(HB(X),\delta)$ является топологической группой при связности произвольного замкнутого шара с центром в фиксированной точке метрического пространства $X$.

В параграфе 3.2 рассматриваются геодезические отображения специальных геодезических пространств. Исследованы некоторые геометрические свойства таких отображений. Теорема Банаха об обратном операторе и принцип равностепенной непрерывности для  $F$-пространств \cite [c. 99, с. 104]{Sadov},  обобщены на случай специального геодезического отображения специальных геодезических пространств. 

В параграфе 3.3 метрика Буземана $\delta_p$ распространяется на множество $\Phi (X, Y)$ всех непрерывных отображений 
$$f : (X,\rho) \rightarrow (Y, d),$$ 
удовлетворяющих следующему условию: для любых $x, \, y \in X$
$$d (f (x), f (y)) \leq B_f e^{|xy|},$$
где $B_f$ --- неотрицательная константа. Доказано, что пространство 
\\
$(H_B (X,Y,\alpha),\delta_p)$ всех отображений из метрического пространства $X$ в метрическое пространство $Y$, удовлетворяющих равномерному условию Гельдера с фиксированными показателем $\alpha \in (0,1]$ и коэффициентом $B \in \mathbb{R}_{+}$, является полным (собственным) метрическим пространством, если $Y$ --- полное метрическое пространство ($X$, $Y$ --- собственные метрические пространства). Установлено, что если $X$ --- собственное метрическое пространство, то топология пространства
$(H_B (X,Y,\alpha),\delta_p)$ совпадает как с топологией поточечной сходимости, так и с компактно-открытой топологией.

В параграфе 3.4 рассматривается метрическое пространство всех подобий $(Sim (X, Y),\delta_p)$ из метрического пространства $X$ на метрическое пространство $Y$ с метрикой Буземана $\delta_p$. Доказано, что: 

-- если группа всех подобий действует транзитивно на полном метрическом пространстве, то и группа изометрий действует на нем транзитивно;

-- если $X$, $Y$ --- полные (собственные) метрические пространства, то пространство $(Sim (X,Y)\cup Const (X,Y),\delta_p)$ --- полное (собственное), где $Const (X,Y)$ --- множество всех постоянных отображений из $X$ в $Y$; 

-- если $X$ --- собственное метрическое пространство, то топология пространства $(Sim (X,Y)\cup Const (X,Y),\delta_p)$ совпадает как с топологией поточечной сходимости, так и с компактно-открытой топологией;
 
-- $(Sim (X),\delta_p)$ --- топологическая группа, действующая непрерывно на пространстве $X$;

-- группы подобий $Sim (X)$ и изометрий $Iso (X)$ с метрикой Куратовского $\delta$ являются топологическими группами, непрерывно действующими на пространстве $X$. 

Найдено замыкание группы подобий полного метрического пространства в пространстве $(\Phi (X,X),\delta_p)$.

В параграфе 3.5 в специальном метрическом пространстве введены два аналога слабой сходимости последовательности в вещественном гильбертовом пространстве и исследованы их геометрические свойства.

Основные результаты по перечисленным темам опубликованы в статьях
автора [s1, s2, s5, s15].

\section{Метрическое пространство слабо ограниченных отображений метрических пространств}
\vskip20pt
Пусть $(r_i)_{i \in \mathbb{N}}$ --- неограниченная строго возрастающая последовательность положительных вещественных чисел, $(B_i = B[p,r_i])_{i \in \mathbb{N}}$ --- последовательность замкнутых шаров с общим центром в точке $p$ метрического пространства $(X,\rho)$. Тогда
$$X =\bigcup \limits_{i = 1}^{\infty} B_i.$$
Рассмотрим множество $BW(X,Y)$ всех отображений из метрического пространства
$(X,\rho)$ в метрическое пространство $(Y,d)$, отображающих ограниченные множества в ограниченные множества. Элементы этого множества назовем {\it слабо ограниченными отображениями} из пространства
$(X,\rho)$ в пространство $(Y,d)$. Зададим инъекцию
$$\theta :  BW(X,Y) \rightarrow \prod \limits_{i = 1}^{\infty}  BW(B_i,Y), \quad
\theta (f) = (f_1, f_2, \ldots ),$$
где $f_i = f |_{B_i}$ для каждого $i \in \mathbb{N}\,$.
Очевидно, что  $\theta$ --- биекция множества $BW(X,Y)$ на предел обратной последовательности \cite [с. 160]{Engel}:
$$\varprojlim \{BW(B_i,Y), \pi^i_j\},\quad \pi^i_j : BW(B_j,Y) \rightarrow BW(B_i,Y),\, \, \pi^i_j (f_i) = f_j,$$
где $j \leq i$. Для каждого $i \in \mathbb{N}$ метризуем множество $BW(B_i,Y)$ метрикой
$\delta_i$:
$$\delta_i (f_i, \varphi_i) = \sup \{ d (f(x), \varphi (x)) : x \in B_i\},$$
где $f_i, \, \varphi_i \in BW(B_i,Y)$. Определим метрику на множестве $BW(X,Y)$ следующим образом
$$\delta (f, \varphi) = \sum \limits_{i = 1}^{\infty} 2^{-i}\dfrac{\delta_i (f_i, \varphi_i)}{1 + \delta_i (f_i, \varphi_i)}.$$
Если метрическое пространство $X$ собственное, то для каждого $i \in \mathbb{N}$ шар $B_i$ компактен и множество $CBW(X,Y)$ всех непрерывных отображений из пространства $(X,\rho)$ в пространство $(Y,d)$, отображающих ограниченные множества в ограниченные множества, совпадает с множеством $C(X,Y)$ всех непрерывных отображений из $(X,\rho)$ в $(Y,d)$. Поэтому для собственного пространства метрика $\delta$     совпадает с метрикой Куратовского на пространстве $C(X,Y)$ \cite [$\S$44, п. 7]{Kurat}. Из результатов Куратовского следует
\vskip7pt
{\bf Теорема 1} [s1]. {\it Пусть $(X,\rho)$ --- собственное метрическое пространство. Тогда верны следующие утверждения.
\vskip7pt
$(i\, \,)$ Метрика $\delta$ определяет на множестве $CBW(X,Y)$ компактно-открытую топологию $\tau_{\delta}$.
\vskip7pt
$(ii\, \,)$ Определенное выше отображение $\theta$ является гомеоморфизмом пространства $C(X,Y) = CBW(X,Y)$ с компактно-открытой топологией на предел обратной последовательности
$$\varprojlim \{BW(B_i,Y), \pi^i_j\}$$
с топологией, определяемой метрикой $\delta$.
\vskip7pt
$(iii\, \,)$ Если $(Y,d)$ --- полное {\rm (}сепарабельное{\rm )} метрическое пространство, то 
$(CBW(X,Y), \delta)$ --- полное {\rm (}сепарабельное{\rm )} метрическое пространство.}
\vskip7pt
Отметим некоторые элементарные свойства метрики $\delta$.
\vskip7pt
$1.\, \,$ Пусть для каждого $n \in \mathbb{N}$ 
$$f, \, f_n \in BW(X,Y).$$ 
Равенство
$$\lim \limits_{n \rightarrow \infty} \delta (f_n, f) = 0,$$ 
имеет место тогда и только тогда, когда для каждого $i \in \mathbb{N}$
$$\lim \limits_{n \rightarrow \infty} \delta_i (f_{n,i}, f_i) = 0.$$

$2.\, \,$ Если для каждого $n \in \mathbb{N}$
$$f, \, f_n \in BW(X,Y) \, \, \mbox {и} \, \, \lim \limits_{n \rightarrow \infty} \delta (f_n, f) = 0,$$ 
то для каждого $x \in X$ 
$$\lim \limits_{n \rightarrow \infty} d(f_n (x), f (x)) = 0,$$
то есть топология $\tau_{\delta}$, определяемая метрикой $\delta$ на множестве $BW(X,Y)$, сильнее топологии поточечной сходимости.
\vskip7pt
$3.\, \,$ Подпространство отображений, удовлетворяющих равномерному условию Гельдера с фиксированным показателем $\alpha \in (0,1]$ и коэффициентом $B \geq 0$,
$$H_B (X,Y,\alpha) = \{f : X \rightarrow Y :\, 
\forall x, \, y \in X ( d(f (x), f (y)) \leq
B |xy|^{\alpha})\}$$
замкнуто в пространстве $(CBW(X,Y), \delta)$.
\vskip7pt
$4.\, \,$ Топология $\tau_{\delta}$ на множестве $BW(X,Y)$ совместно непрерывная (\cite {Kelly}, c. 290), то есть отображение вычисления
$$\Omega :  BW(X,Y)\times X \rightarrow Y, \quad  \Omega (f, x) = f (x)$$ непрерывно.
\vskip7pt
Докажем свойство 4. Пусть для всех $n, \, m \in \mathbb{N}$
$$x, \, x_n \in X; \quad f, \, f_m \in BW(X,Y).$$  
Кроме того, 
$$\lim \limits_{n \rightarrow \infty} |x_nx| = 0,\quad \lim \limits_{m \rightarrow \infty} \delta (f_m, f) = 0.$$
Тогда для всех $n \in \mathbb{N}$ найдется такое $i \in \mathbb{N}$, что
$x, \, x_n \in B_i$ и по свойству 1 имеем
$$\lim \limits_{m \rightarrow \infty} \delta_i (f_{m,i}, f_i) = 0.$$
Следовательно,
$$d(f_{m,i} (x_n), f_i (x_n)) \rightarrow 0$$
при $n \rightarrow  \infty, \quad m \rightarrow \infty$.
Тогда
$$d(f_m (x_n), f (x)) \leq d(f_{m,i} (x_n), f_i (x_n)) + d (f_i (x_n), f_i (x))\rightarrow 0$$
при $n \rightarrow  \infty, \quad m \rightarrow \infty$. Таким образом, свойство 4 доказано.
\vskip7pt
Рассмотрим множество $HB (X)$ всех тех гомеоморфизмов из множества $BW(X,X)$, каждый из которых при любом $i \in \mathbb{N}$ равномерно непрерывен на замкнутом шаре $B_i \subset X$ вместе со своим обратным гомеоморфизмом. Напомним, что {\it паратопологической группой} называется
множество, наделенное структурой группы и топологией, в которой отображение произведения непрерывно \cite [c. 20]{Burb1}. 
\vskip7pt
{\bf Теорема 2} [s1]. {\it 
\vskip3pt
$(i)\, \,$ Пространство $(HB (X), \delta)$, наделенное операцией композиции отображений, является паратопологической группой, непрерывно действующей на пространстве $X$.
\vskip7pt
$(ii)\, \,$ Если для каждого $i \in \mathbb{N}$ шар $B_i \subset X$ связен, то $(HB (X), \delta)$ --- топологическая группа, непрерывно действующая на пространстве $X$.}
\vskip7pt
{\bf Доказательство теоремы 2}.
\vskip7pt
$(i)\, \,$ Пусть для  каждого $n \in \mathbb{N}$
$$f, \, \varphi, \, f_n, \, \varphi_n \in HB (X).$$  
Кроме того,
$$\lim \limits_{n \rightarrow \infty} \delta (f_n, f)  = 0,\quad \lim \limits_{n \rightarrow \infty} \delta (\varphi_n, \varphi) = 0.$$ 
В силу свойства 1 достаточно доказать, что для каждого $i \in \mathbb{N}$ 
$$\lim \limits_{n \rightarrow \infty} \delta_i (f_n\circ \varphi_n)_i, (f\circ \varphi)_i) = 0.$$
Фиксируем  $i \in \mathbb{N}$, тогда для любых 
$x, \, y \in B_i$ получим
$$D (\varphi_n (B_i)) = \sup \{|\varphi_n (x)\varphi_n (y)| : \, x, \, y \in B_i\} \leq$$
$$\sup \{|\varphi_n (x)\varphi (x)| + |\varphi_n (y)\varphi (y)| + |\varphi (x)\varphi (y)|: \, x, \, y \in B_i\} \leq$$ 
$$2 \delta_i (\varphi_{n,i}, \varphi_i) + D (\varphi (B_i)).$$ 
Следовательно, найдется такое $s \in \mathbb{N}$, что для каждого $n \in \mathbb{N}$
$$\varphi_n (B_i) \subset B_s, \quad  \varphi (B_i)\subset B_s.$$
Кроме того,
$$\delta_i ((f_n\circ\varphi_n)_i, (f\circ\varphi_n)_i) \leq \sup \{|f_n (y)f (y)| : \, y \in \varphi_n (B_i)\} \leq$$ 
$$\sup \{|f_n (y)f (y)| : \, y \in B_s\} = \delta_s (f_{n,s}, f_s) \rightarrow 0$$
при $n \rightarrow \infty$. 
Предположим, что для любого числа
$A > 0$ существует такое $n_0 \in \mathbb{N}$, что для каждого $n > n_0$  
$$\delta_i ((f\circ\varphi_n)_i, (f\circ\varphi)_i) > A > 0.$$
Тогда существует такая последовательность $(x_n)_{n \in \mathbb{N}}$, что
для каждого $n > n_0$
$$|(f\circ\varphi_n (x_n))(f\circ\varphi (x_n))| > A > 0.$$
Но
$$|\varphi_n (x_n)\varphi (x_n)| \leq \sup \{|\varphi_n (x)\varphi (x)| : \, x \in B_i\} = \delta_i (\varphi_{n,i}, \varphi_i) \rightarrow 0$$
при $n \rightarrow \infty$
и отображение $f \in HB (X)$ равномерно непрерывно на шаре $B_s$.
Значит, 
$$0 < A < |f (\varphi_n (x_n))f(\varphi (x_n))|\rightarrow 0$$
при $n \rightarrow \infty$. Получили противоречие. Следовательно, существует подпоследовательность $(\varphi_k)_{k \in \mathbb{N}}$ последовательности $(\varphi_n)_{n \in \mathbb{N}}$ такая, что 
$$\delta_i ((f\circ\varphi_k)_i, (f\circ\varphi)_i)\rightarrow 0$$
при $k \rightarrow \infty$.
Утверждение $(i)$ следует теперь из свойства 1 и неравенств:
$$\delta_i ((f_k\circ\varphi_k)_i, (f\circ\varphi)_i) \leq \delta_i ((f_k\circ\varphi_k)_i, (f\circ\varphi_k)_i) + \delta_i ((f\circ\varphi_k)_i, (f\circ\varphi)_i) \leq$$ 
$$\delta_s (f_{k,s}, f_s) + \delta_i ((f\circ\varphi_k)_i, (f\circ\varphi)_i).$$

$(ii)\, \,$ Из первой части теоремы следует, что достаточно доказать следующее утверждение:
если для  каждого $n \in \mathbb{N}$ 
$$ f_n \in HB (X) \, \, \mbox {и}\, \, \lim \limits_{n \rightarrow \infty} \delta (f_n, id)  = 0,$$ 
то
$$\lim \limits_{n \rightarrow \infty} \delta (f^{-1}_n, id)  = 0.$$
Фиксируем $i \in \mathbb{N}$, тогда найдется такое $n_0 \in \mathbb{N}$, что для каждого $n > n_0$ 
$$\delta_{i+1} (f_{n,i+1}, id_{i+1}) =\sup \{|f_n (x)x| : \, x \in B_{i+1}\} <$$
$$\varepsilon = \min \{r_i, r_{i+1} - r_i\} .$$
Покажем, что для каждого $n > n_0$ 
$$f^{-1}_n (B_i) \subset B_{i+1}.$$
Предположим противное, то есть найдется такое $k > n_0$, для которого
$$f^{-1}_k (B_i) \nsubseteq B_{i+1}.$$ 
Ясно, что
$$f^{-1}_k (B_i) \cap B_{i+1} \neq \emptyset,$$
поскольку $f_n (p) \in B_i$. Кроме того, множество $f^{-1}_k (B_i)$ связно, так как шар $B_i$ связен. Следовательно, найдется точка
$$z \in f^{-1}_k (B_i) \cap (\partial B_{i+1}).$$
Получаем противоречие, поскольку
$$|pz| \leq |pf_k (z)| + |f_k (z)z|< r_i + \varepsilon \leq r_{i+1}$$
и  $z \in \stackrel{\circ}{B}_{i+1}$. Таким образом, для каждого $n > n_0$  
$$f^{-1}_n (B_i) \subset B_{i+1}.$$
Кроме того, для каждого $n > n_0$ 
$$\delta_i (f^{-1}_{n,i}, id_i) =\sup \{|f^{-1}_n (x) x| : \, x \in B_i\} \leq$$
$$\sup \{|yf_n (y)| : \, y \in f^{-1} (B_i)\} \leq \sup \{|yf_n (y)| : \, y \in B_{i+1}\} =  \delta_{i+1} (f_{n,i+1}, id_{i+1}).$$
Следовательно,
$$\lim \limits_{n \rightarrow \infty} \delta (f^{-1}_n, id)  = 0.$$
Теорема 2 доказана.
\vskip7pt
{\bf Следствие 1} [s1]. {\it 
Если для каждого $i \in \mathbb{N}$ шар $B_i \subset X$ связен, то группа липшицевых гомеоморфизмов
$$LH (X) = \{f : X \rightarrow X : \, f(X) = X,$$ 
$$\forall x,\, y \in X (A_f |xy| \leq|f(x)f(y)|\leq B_f |xy|, \, \, A_f, \, B_f \in \mathbb{R}_{+}, \, \, A_f > 0)\}$$
является топологической группой, непрерывно действующей на пространстве $X$.}
\vskip7pt
{\bf Замечание 1.} {\it В статье {\rm \cite {Dieud}} Ж. Дьедонне использовал общий подход для введения равномерных структур на группах гомеоморфизмов хаусдорфовых равномерных пространств. 
Он, в частности, доказал, что введенная им равномерная структура на группе 
$HB (X)$ определяет, в случае связности для каждого $i \in \mathbb{N}$ шара $B_i$ структуру топологической группы, непрерывно действующей на пространстве $X$. Можно показать, что топология, определяемая равномерной структурой Дьедонне на группе $HB (X)$ слабее топологии, определяемой на этой  группе метрикой $\delta$.}

\section{Геодезические отображения специальных геодезических пространств}
\vskip20pt

Пусть метрические пространства $(X,\rho )$, $(Y,d)$ с выделенными семействами сегментов $S$ и $\tilde S$ соответственно, являются $U$-множествами. Непрерывное отображение 
$$f : (X,\rho ) \rightarrow (Y,d)$$
назовем {\it геодезическим отображением}, если оно отображает сегменты семейства $S$ в сегменты семейства $\tilde S$, то есть $f (S) \subset \tilde S$.
\vskip7pt
Будем использовать следующие обозначения. 
\vskip7pt
Если выполнено условие $(PG)$, то $K = \mathbb{R}$, в ином случае $K = [0,1]$.
\vskip7pt
$G(X,Y)$ --- множество всех геодезических отображений из пространства $(X,\rho)$  в  пространство $(Y,d)$. 
\vskip7pt
$G_{\omega} (X,Y)$ --- подмножество множества $G(X,Y)$, содержащее только такие геодезические отображения $f \in G(X,Y)$, которые удовлетворяют условию 
$$f (\omega_{\lambda} (x, y)) = \tilde \omega_{\lambda} (f (x),f (y))$$
для всех $x, \, y \in X$ и для каждого $\lambda \in K$. 
\vskip7pt
Приведем некоторые элементарные свойства геодезических отображений.
\vskip7pt
$1.\, \,$ Образ выпуклого множества при геодезическом отображении является выпуклым множеством.
\vskip7pt
$2.\, \,$ Изометрическое отображение 
$$f : (X, \rho) \rightarrow (Y,d)$$ 
является геодезическим отображением, если $f (S) \subset \tilde S$. 
\vskip7pt
$3.\, \,$ Если множество $\tilde S$ замкнуто относительно топологии поточечной сходимости в множестве всех сегментов пространства $(Y,d)$, последовательность геодезических отображений поточечно сходится к некоторому непрерывному отображению, то это отображение геодезическое.
\vskip7pt
$4.\, \,$ Если пространство $(Y,d)$ удовлетворяет условию $(C_6)$, последовательность отображений из множества $G_{\omega} (X,Y)$ поточечно сходится к некоторому непрерывному отображению, то это отображение принадлежит множеству $G_{\omega} (X,Y)$.
\vskip7pt
$5.\, \,$ Если гомеоморфизм $f$ принадлежит множеству $G(X,Y)$ $(G_{\omega} (X,Y))$, то обратный гомеоморфизм $f^{-1}$ принадлежит множеству 
\\
$G(Y,X)$ $(G_{\omega} (Y,X))$.
\vskip7pt
{\bf Теорема 1} [s5]. {\it Если $f \in G_{\omega} (X,Y)$, то $f$ отображает ограниченные множества в ограниченные множества, то есть $G_{\omega} (X,Y) \subset CBW (X,Y)$.}
\vskip7pt
{\bf Доказательство теоремы 1}.
\vskip7pt
Допустим противное. Пусть существует такая ограниченная последовательность $(x_n)_{n \in \mathbb{N}}$ в пространстве $X$, что
$$\lim \limits_{n \rightarrow \infty} d(f (x_n), f(p)) = \infty,$$ 
где $p$ --- фиксированная точка из $X$. Тогда, начиная с некоторого номера 
$$n \in \mathbb{N}, \quad \lambda_n =  \dfrac{1}{d(f (x_n), f(p))} \leq 1\, \, \mbox {и} \, \, \lim \limits_{n \rightarrow \infty} \lambda_n = 0.$$
Кроме того, 
$$\lim \limits_{n \rightarrow \infty} \omega_{\lambda_n} (p, x_n) = 0,$$
поскольку
$$\lim \limits_{n \rightarrow \infty} |p\omega_{\lambda_n} (p, x_n)| = \lim \limits_{n \rightarrow \infty} \lambda_n |px_n| = 0.$$
Получили противоречие с условием непрерывности отображения $f$, поскольку
$$d(f (p), f (\omega_{\lambda_n} (p, x_n))) = d(f (p), \tilde \omega_{\lambda_n} (f(p), f (x_n))) = \lambda_n d (f(p), f (x_n)) = 1.$$
Теорема 1 доказана.
\vskip7pt
{\bf Теорема 2} [s5]. {\it Пусть пространства $X$, $Y$ удовлетворяют условию $(PG)$, отображение $f : X \rightarrow Y$ отображает ограниченные множества в ограниченные множества и для всех $x,\, y \in X$; $\lambda \in \mathbb{R}$ 
имеет место равенство:
$$ f (\omega_{\lambda} (x, y)) = \tilde \omega_{\lambda} (f (x),f (y)).$$
Тогда $f \in G_{\omega} (X,Y)$.}
\vskip7pt
{\bf Доказательство теоремы 2}.
\vskip7pt
Докажем непрерывность отображения $f$. Пусть для любого $n \in \mathbb{N}$
$$x,\, x_n \in X\quad \mbox {и} \, \,  \lim \limits_{n \rightarrow \infty} x_n = x.$$
Выберем такую последовательность действительных чисел $(k_n)_{n \in \mathbb{N}}$, что
$$\lim \limits_{n \rightarrow \infty} k_n = \infty \quad \mbox {и} \, \,  \lim \limits_{n \rightarrow \infty} (k_n |x_nx|) = 0.$$
Тогда 
 $$\lim \limits_{n \rightarrow \infty} |x\omega_{k_n} (x,x_n)| = \lim \limits_{n \rightarrow \infty} (k_n |xx_n|) = 0.$$ 
Кроме того, последовательность $(f(\omega_{k_n} (x,x_n)))_{n \in \mathbb{N}}$  ограничена, так как отображение $f$ отображает ограниченные множества в ограниченные множества. Следовательно,
$$\lim \limits_{n \rightarrow \infty} d (f(x_n), f(x)) = \lim \limits_{n \rightarrow \infty} (k^{-1}_n d (f(\omega_{k_n} (x,x_n)), f(x))) = 0.$$
Таким образом, отображение $f$ непрерывно. Теорема 2 доказана.
\vskip7pt
{\bf Теорема 3} [s5]. {\it Пусть полные метрические пространства $X$, $Y$ удовлетворяют условиям $(PG)$, $(C_6)$ и все замкнутые шары в этих пространствах выпуклые. Если сюръекция $f$ принадлежит множеству $G_{\omega} (X,Y)$ , то $f$ --- открытое отображение.}
\vskip7pt
{\bf Доказательство теоремы 3}.
\vskip7pt
Из представления  
$$X = \bigcup \limits_{n = 1}^{\infty} \omega_n (p, B(p, \varepsilon)),$$
где $p \in X, \, \, \varepsilon > 0$, имеем
$$Y = \bigcup \limits_{n = 1}^{\infty} \tilde \omega_n (f (p), f (B(p, \varepsilon))) = \bigcup \limits_{n = 1}^{\infty} \overline {\tilde \omega_n (f (p), f (B(p, \varepsilon)))},$$
где черта означает замыкание множества. По известной теореме о категориях  \cite [c. 44]{Sadov} существует такой номер $n \in \mathbb{N}$, что множество
$$\overline {\tilde \omega_n (f (p), f (B(p, \varepsilon)))}$$
содержит непустое открытое множество. Но из условий $(PG)$, $(C_6)$ следует, что для любого $z \in Y$ и для каждого $\lambda \in \mathbb{R}\backslash \{0\}$ отображение 
$$\tilde \omega_{\lambda} (z, \cdot) : (Y,d) \rightarrow (Y,d)$$
является гомеоморфизмом. Поэтому множество
$$\overline {f (B(p, \varepsilon))}$$
содержит открытое множество $V \neq \emptyset$. Кроме того,
$$\overline {f (B(p, \varepsilon))} = \overline {f (\omega_{1/2} (B(p, \varepsilon),\omega_{-1} (p, B(p, \varepsilon))))} =$$
$$\overline {\tilde \omega_{1/2} (f (B(p, \varepsilon)),\tilde \omega_{-1} (f (p), f (B(p, \varepsilon))))} \supset \tilde \omega_{1/2} (V, \omega_{-1} (f (p), V)).$$
Последнее множество содержит точку $f (p)$ и открыто, поскольку в силу условий  $(PG)$, $(C_6)$ для каждого $\lambda \in \mathbb{R}$  открыто отображение 
$$\tilde \omega_{\lambda} : (Y,d)\times (Y,d) \rightarrow (Y,d).$$
Выберем числа
$$\varepsilon_k = \varepsilon 2^{-k}, \quad k \in \{0, 1, 2, \ldots\}.$$ 
Из приведенных выше рассуждений следует, что найдутся такие вещественные числа
$$\delta_k > 0, \quad  k \in \{0, 1, 2, \ldots\}, \, \mbox {что} \, \, \lim \limits_{k \rightarrow \infty} \delta_k = 0$$ 
и  для всех $x \in X$, $k \in \{0, 1, 2, \ldots\}$
$$\overline {f (B(p, \varepsilon_k))} \supset B(f(x), \delta_k).$$
Пусть $x_0 \in X$, $y \in B(f(x_0),\delta_0)$. Тогда найдется такая точка 
$$x_1 \in B(x_0,\varepsilon_0),\, \mbox {что} \, \, d (f (x_1), y) < \delta_1.$$
Следовательно,
$$y \in B(f (x_1),\delta_1).$$
Повторяя это рассуждение, построим такую последовательность $(x_k)_{k \in \mathbb{N}}$, что для каждого $k \in \mathbb{N}$ 
$$x_{k+1} \in B(x_k, \varepsilon_k), \quad y \in B(f(x_k),\delta_k).$$
Эта последовательность фундаментальная, поскольку
для всех $k \in \{0, 1, 2, \ldots\}, \, l \in \mathbb{N}$ 
$$|x_kx_{k+l}| \leq |x_kx_{k+1}| + \cdots + |x_{k+l-1}x_{k+l}| <$$
$$\varepsilon_k + \cdots + \varepsilon_{k+l-1} = \varepsilon 2^{1-k} (1 - 2^{-l}).$$
Но пространство $X$ полное, поэтому последовательность $(x_k)_{k \in \mathbb{N}}$ сходится к некоторому пределу $x \in X$.
Кроме того,  
$$y = \lim \limits_{k \rightarrow \infty} f(x_k) = f(x),$$
поскольку отображение $f$ непрерывно.
Таким образом, 
$$B(f(x_0),\delta_0) \subset f (B(x_0, 3 \varepsilon))$$
так как точка $y$ из $B(f(x_0),\delta_0)$  произвольная и 
$$|x_0x| = \lim \limits_{k \rightarrow \infty} |x_0x_k| \leq 2 \varepsilon < 3 \varepsilon.$$
Теперь нетрудно заметить, что отображение $f$ открытое. Теорема 3 доказана.
\vskip7pt
Следствием свойства 5 и теоремы 3 является
\vskip7pt
{\bf Теорема 4} [s5]. {\it Пусть полные метрические пространства $X$, $Y$ удовлетворяют условиям $(PG)$, $(C_6)$ и все замкнутые шары в этих пространствах выпуклые. Если отображение $f$ из множества $G_{\omega} (X,Y)$ биективно, то обратное отображение $f^{-1}$ принадлежит множеству $G_{\omega} (Y,X)$.}
\vskip7pt
{\bf Теорема 5} [s5]. {\it Пусть полные метрические пространства $X$, $Y$ удовлетворяют условиям $(PG)$, $(C_6)$  и $(f_\alpha)_{\alpha \in A}$ --- семейство отображений из множества $G_{\omega} (X,Y)$. Если для каждого $x \in X$ множество 
$$\{f_\alpha (x) : \, \alpha \in A\}$$ 
ограничено, то
$$\lim\limits_{x\to p} f_\alpha (x) = f_\alpha (p)$$
равномерно относительно $\alpha \in A$.}
\vskip7pt
{\bf Доказательство теоремы 5}.
\vskip7pt
Пусть заданы $\varepsilon > 0$, $p \in X$.
Очевидно, что для каждого $n \in \mathbb{N}$ множество
$$X_n = \{x \in X ; \, d (f_{\alpha} (x), f_{\alpha} (p)) \leq n \varepsilon, \, \alpha \in A\}$$
замкнуто в пространстве $X$ и
$$X = \bigcup \limits_{n = 1}^{\infty} X_n.$$
По теореме о категориях \cite [c. 44]{Sadov} существует такое натуральное число $n_0$, что для некоторого $x_0 \in X_{n_0}$ и некоторого $\delta > 0$
$$B(x_0,\delta) \subset X_{n_0}.$$
В силу условий $(PG)$, $(C_6)$ для любого $z \in X$ и для каждого $\lambda \in \mathbb{R}\backslash \{0\}$ отображение 
$$\tilde \omega_{\lambda} (z, \cdot) : X \rightarrow X$$
является гомеоморфизмом. Поэтому $x$ сходится к $x_0$ тогда и только тогда, когда
$$y = \omega_{1/n_0} (p, \omega_{1/2} (\omega_{-1}(p, x_0), x))$$
сходится к $p$. Но если $|xx_0| < \delta$, то
$$d (f_{\alpha} (y), f_{\alpha} (p)) = \dfrac{1}{n_0} d (f_{\alpha} (\omega_{1/2} (\omega_{-1}(p, x_0), x)), f_{\alpha} (p)) \leq$$
$$\dfrac{1}{n_0} \{ d (\tilde \omega_{1/2} (f_{\alpha} (\omega_{-1}(p, x_0)),f_{\alpha} (x)), f_{\alpha} (\omega_{-1}(p, x_0))) + d (f_{\alpha} (\omega_{-1}(p, x_0)),f_{\alpha} (p))\} =$$
$$\dfrac{1}{n_0} \{\dfrac{1}{2} d (f_{\alpha} (\omega_{-1}(p, x_0)),f_{\alpha} (x)) + d (f_{\alpha} (\omega_{-1}(p, x_0)),f_{\alpha} (p))\}
\leq$$
$$\dfrac{1}{2 n_0} \{ d (f_{\alpha} (p), f_{\alpha} (x)) + 3 d (f_{\alpha} (\omega_{-1}(p, x_0)),f_{\alpha} (p))\} =$$
$$\dfrac{1}{2 n_0} \{ d (f_{\alpha} (p), f_{\alpha} (x)) + 3 d (f_{\alpha} (p), f_{\alpha} (x_0)) \} \leq 2 \varepsilon.$$
Следовательно, если $y$ сходится к $p$, то $f_{\alpha} (y)$ сходится к $f_{\alpha} (p)$. Теорема 5 доказана.
\vskip7pt
Пусть $(r_i)_{i \in \mathbb{N}}$ --- неограниченная строго возрастающая последовательность положительных вещественных чисел и для каждого $i \in \mathbb{N}$
$B_i = B [p, r_i]$--- замкнутый шар с центром в точке $p \in X$, радиуса $r_i$. Из теоремы 1 следует, что 
$$G_{\omega} (X,Y) \subset (CBW (X,Y),\delta)$$ 
и значит на множестве $G_{\omega} (X,Y)$ можно задать индуцированную метрику $\delta$
$$\delta (f, \varphi) = \sum \limits_{i = 1}^{\infty} 2^{-i}\dfrac{\sup \{ d (f(x), \varphi (x)) : x \in B_i\}}{1 + \sup \{ d (f(x), \varphi (x)) : x \in B_i\}},$$
для всех $f, \, \varphi \in G_{\omega} (X,Y)$.
\vskip7pt
Обозначим через $H_{\omega} (X)$ множество всех тех гомеоморфизмов из множества $G_{\omega} (X,X)$, каждый из которых при любом $i \in \mathbb{N}$ равномерно непрерывен на замкнутом шаре $B_i \subset X$ вместе со своим обратным гомеоморфизмом. 
\vskip7pt
Из теорем 1, 3.1.2 следует
\vskip7pt
{\bf Теорема 6} [s5]. {\it Метрическое пространство $(H_{\omega} (X),\delta)$, наделенное операцией композиции отображений, является топологической группой, непрерывно действующей на пространстве $(X,\rho)$.}

\section{Полнота и собственность некоторых пространств отображений с метрикой Буземана}
\vskip20pt

Рассмотрим множество $\Phi (X, Y)$ всех непрерывных отображений 
$$f : (X,\rho) \rightarrow (Y, d),$$ 
удовлетворяющих следующему условию: для любых $x, \, y \in X$ имеет место неравенство 
$$d (f (x), f (y)) \leq B_f e^{|xy|},$$
где $B_f$ --- неотрицательная константа.
Очевидно, что множество $\Phi (X, Y)$ содержит подмножество $Const (X,Y)$ всех постоянных отображений из $X$ в $Y$. Кроме того, в нем содержится множество $H_B (X,Y,\alpha)$ всех отображений, удовлетворяющих равномерному условию Гельдера с фиксированными показателем $\alpha \in (0,1]$ и коэффициентом $B \in \mathbb{R}_{+}$, то есть
$$H_B (X,Y,\alpha) = \{f : (X,\rho) \rightarrow (Y, d) :\, \forall x, \, y \in X (
 d (f (x), f (y)) \leq B |xy|^{\alpha})\},$$  
где $B \in \mathbb{R}_{+}$, $\alpha \in (0,1]$ --- фиксированные константы.
Фиксируем точку $p \in X$ и зададим метрику $\delta_p$ на множестве $\Phi (X, Y)$, введенную Буземаном для группы изометрий \cite [с. 30]{Bus}:
$$\delta_p (f, \varphi) = \sup \{ d (f (x), \varphi (x)) e^{-|px|} : \, x \in X\},$$
для всех $f, \, \varphi \in \Phi (X, Y).$
Приведем в лемме 1 некоторые элементарные свойства перечисленных множеств и метрики. 
\vskip7pt
{\bf Лемма 1} [s2]. {\it 
\vskip3pt
$(i)\, \,$ Для любых $p, \, q \in X$ $\delta_p$ --- метрика и для любых $f, \, \varphi \in \Phi (X, Y)$ справедливы неравенства
$$ e^{-|pq|}\delta_p (f, \varphi) \leq \delta_q (f, \varphi)  \leq e^{|pq|}\delta_p (f, \varphi),$$ 
то есть метрика $\delta_p$ липшицево эквивалентна метрике $\delta_q$.
\vskip7pt
$(ii)\, \,$ Для любого $p \in X$ метрическая топология пространства 
\\
$(\Phi (X, Y), \delta_p)$ является совместно непрерывной и больше компактно-открытой топологии.
\vskip7pt
$(iii)\, \,$ Для любого $p \in X$ и для любого $W \in \Sigma (Y)$ отображение вычисления в точке $x \in X$
$$E_x : (Const (X,W), \delta_p) \rightarrow (W, d), \quad E_x (f) = f (x)$$
является изометрией и 
$$\overline {Const (X,W)} = Const (X,\overline {W}),$$
где в левой части замыкание выполнено в пространстве $(\Phi (X, Y), \delta_p)$.
\vskip7pt
$(iv)\, \,$ Для любого $p \in X$ множество $H_B (X,Y,\alpha)$ замкнуто в
пространстве $(\Phi (X, Y), \delta_p)$.}
\vskip7pt
{\bf Доказательство леммы 1}.
\vskip7pt
$(i)\, \,$ Для каждого $p\in X$ и для любых 
$f, \, \varphi \in \Phi (X, Y)$ конечность расстояния $\delta_p (f, \varphi)$
следует из неравенств:
$$\delta_p (f, \varphi) \leq \sup \{ (d (f (x), f(p)) + d (f(p), \varphi (p)) + d(\varphi (p), \varphi (x))) e^{-|px|} : \, x \in X\} \leq$$
$$B_f + B_{\varphi} + d (f(p), \varphi (p)).$$
Определение метрики нетрудно проверить непосредственно. Для любых $p, \, q \in X$ липшицева эквивалентность метрик  $\delta_p$,  $\delta_q$ устанавливается c помощью неравенства:
$$\delta_p (f, \varphi) \leq \sup \{d (f (x), \varphi (x))  e^{-|qx| +|pq|} : \, x \in X\} = e^{|pq|}\delta_q (f, \varphi),$$
а также последующих замен $p$ на $q$ и $q$ на $p$. 
\vskip7pt
$(ii)\, \,$ Пусть 
$$x, \, x_0 \in X,\quad f, \, f_0 \in \Phi (X, Y).$$
Из неравенств
$$d (f (x), f_0(x_0)) \leq  d (f (x), f_0(x)) + d (f_0 (x), f_0(x_0)) \leq$$ $$e^{|xx_0|}\delta_{x_0} (f, f_0) + d (f_0 (x), f_0(x_0))$$
следует, что отображение вычисления
$$\Omega :  \Phi (X, Y)\times X \rightarrow Y, \quad  \Omega (f, x) = f (x)$$ непрерывно. Поэтому для любого $p \in X$ метрическая топология пространства $(\Phi (X, Y), \delta_p)$ является совместно непрерывной и больше 
\\
компактно-открытой топологии \cite [с. 290-292]{Kelly}.
\vskip7pt
$(iii)\, \,$ Пусть $x \in X$, $f, \, \varphi \in Const (X, W)$ и для любого $z \in X$ 
$$f (z) = y_f, \quad \varphi (z) = y_{\varphi}.$$
Первое утверждение следует из равенств:
$$d (E_x (f), E_x (\varphi)) = d (f (x), \varphi (x)) = \sup \{d (y_f, y_{\varphi})  e^{-|pz|} : \, z \in X\} = \delta_p (f, \varphi) .$$
Второе утверждение нетрудно проверить, используя утверждение $(ii)$. 
\vskip7pt
$(iv)\, \,$ Это следствие утверждения $(ii)$ и определения множества 
\\
$H_B (X,Y,\alpha)$. Лемма 1 доказана.
\vskip7pt
{\bf Теорема 1} [s2]. {\it Если  $(Y,d)$ --- полное метрическое пространство, $p\in X$, то 
$$(H_B (X,Y,\alpha),\delta_p)$$ 
--- полное метрическое пространство.}
\vskip7pt
{\bf Доказательство теоремы 1}.
\vskip7pt
Пусть $p \in X$ и последовательность $(f_n)_{n \in \mathbb{N}}$ пространства $H(X,Y,\alpha)$ фундаментальная, то есть для каждого $\varepsilon > 0$ найдется такое $n_0 \in  \mathbb{N}$, что для всех $n, \, m > n_0$ 
$$\delta_p (f_n, f_m) < \varepsilon. \eqno (1)$$
Тогда для всех $n, \, m > n_0$ и для любого $x \in X$
$$d (f_n (x), f_m (x)) e^{-|px|} < \varepsilon.$$
Поэтому последовательность $(f_n (x))_{n \in \mathbb{N}}$ пространства $Y$
фундаментальная при фиксированном $x \in X$. В силу полноты пространства $Y$ определено отображение
$$f : X \rightarrow Y, \quad f (x) = \lim \limits_{n \rightarrow \infty} f_n (x).$$
Тогда $f \in H(X,Y,\alpha)$, поскольку для всех $x, \, y \in X$
$$d (f (x), f (y)) = \lim \limits_{n \rightarrow \infty} d (f_n (x), f_n (y)) \leq
B |xy|^{\alpha}.$$
В неравенстве $(1)$ $\varepsilon$ не зависит от $x$, поэтому, переходя в его левой части к пределу при $m \rightarrow \infty$ и вычисляя точную верхнюю грань, получим для любого $n > n_0$
$$\delta_p (f_n, f) \leq \varepsilon.$$
Следовательно, 
$$\lim \limits_{n \rightarrow \infty} \delta_p (f_n, f) =0$$
и $(H_B (X,Y,\alpha),\delta_p)$ --- полное метрическое пространство. Теорема 1 доказана.
\vskip7pt
{\bf Теорема 2} [s2]. {\it Если $(X,\rho)$ --- собственное метрическое пространство, то топология пространства
$(H_B (X,Y,\alpha),\delta_p)$ совпадает как с топологией поточечной сходимости, так и с компактно-открытой топологией.}
\vskip7pt
{\bf Доказательство теоремы 2}.
\vskip7pt
Пусть $p \in X$. В силу утверждения $(ii)$ леммы 1 достаточно доказать следующее утверждение: если для любого $f \in H_B (X,Y,\alpha$ и любой последовательности $(f_n)_{n \in \mathbb{N}}$ пространства $H(X,Y,\alpha)$
$$\lim \limits_{n \rightarrow \infty} d (f_n (x), f (x)) = 0$$
для каждого $x \in X$, то 
$$\lim \limits_{n \rightarrow \infty} \delta_p (f_n, f) =0.$$
Докажем это утверждение методом от противного. Пусть это утверждение неверно. Тогда найдутся подпоследовательность $(f_l)_{l \in \mathbb{N}}$
последовательности $(f_n)_{n \in \mathbb{N}}$, последовательность $(x_l)_{l \in \mathbb{N}}$ пространства $X$ и константа $A > 0$ такие, что
$$d (f_l (x_l), f (x_l)) e^{-|px_l|} > A. \eqno (2)$$

Рассмотрим два случая.
\vskip7pt
$1. \, \,$ Допустим сначала, что последовательность $(x_l)_{l \in \mathbb{N}}$ неограничена. Тогда в силу собственности пространства $X$ найдется такая подпоследовательность $(x_t)_{t \in \mathbb{N}}$ последовательности $(x_l)_{l \in \mathbb{N}}$, что 
$$\lim \limits_{t \rightarrow \infty} |px_t| = \infty.$$
Кроме того,
$$d (f_t (x_t), f (x_t)) e^{-|px_t|} \leq \{d (f_t (x_t), f_t (p)) + d (f_t (p), f (p)) + $$
$$d (f (p), f (x_t))\} e^{-|px_t|} \leq d (f_t (p), f (p)) + 2 B |x_tp|^{\alpha}  e^{-|px_t|} \rightarrow 0$$
при $t \rightarrow \infty$. Учитывая неравенство $(2)$, получили противоречие. 
\vskip7pt
$2. \, \,$ Допустим теперь, что последовательность $(x_l)_{l \in \mathbb{N}}$ ограничена. Тогда в силу собственности пространства $X$ найдется такая подпоследовательность $(x_k)_{k \in \mathbb{N}}$ последовательности $(x_l)_{l \in \mathbb{N}}$ и точка $x \in X$, что 
$$\lim \limits_{k \rightarrow \infty} |px_k| = 0.$$
Кроме того,
$$d (f_k (x_k), f (x_k)) e^{-|px_k|} \leq d (f_k (x_k), f_k (x)) + d (f_k (x), f (x)) + d (f (x), f (x_k)) \leq$$
$$d (f_k (x), f (x)) + 2 B |x_kx|^{\alpha} \rightarrow 0$$
при $k \rightarrow \infty$. Учитывая неравенство $(2)$, снова получили противоречие. Теорема 2 доказана.
\vskip7pt
Напомним, что замкнутое непрерывное отображение метрических пространств, прообразы точек которого компактны, называется {\it совершенным отображением} \cite [с. 277]{Engel}.
\vskip7pt
{\bf Теорема 3} [s2]. {\it Если $(X,\rho)$, $(Y,d)$ --- собственные метрические пространства, то $(H_B (X,Y,\alpha),\delta_p)$ --- собственное метрическое пространство, а отображение вычисления в точке $p \in X$
$$E_p : (H_B (X,Y,\alpha),\delta_p) \rightarrow (Y, d), \quad E_p (f) = f (p)$$
сюръективно и совершенно. Если, кроме того, $(Y,d)$ --- компактное пространство, то $(H_B (X,Y,\alpha),\delta_p)$ --- компактное пространство.}
\vskip7pt
{\bf Доказательство теоремы 3}.
\vskip7pt
$1. \, \,$ Пусть $p \in X$ и последовательность $(f_n)_{n \in \mathbb{N}}$ пространства 
\\
$(H_B (X,Y,\alpha),\delta_p)$ ограничена. Тогда для любого $x \in X$ ограничена последовательность $(f_n (x))_{n \in \mathbb{N}}$ пространства $Y$. 
Для любого $\varepsilon > 0$ выберем $\delta = \varepsilon$ при $B = 0$ и
$$\delta = \left(\dfrac{\varepsilon}{B}\right)^{1/\alpha}$$
при $B > 0$ так, что для любого $n \in \mathbb{N}$
$$d (f_n (x), f_n (y)) \leq B |xy|^{\alpha} < \varepsilon,$$
как только $|xy| <  \delta$ для $x, \, y \in X$. Следовательно, последовательность  $(f_n)_{n \in \mathbb{N}}$ равностепенно равномерно непрерывна в пространстве $X$ (определение см. в \cite [с. 20]{Bus}). В силу условий теоремы 3 и теоремы 2.17 \cite [с. 20]{Bus} найдется 
подпоследовательность  $(f_l)_{l \in \mathbb{N}}$ последовательности  $(f_n)_{n \in \mathbb{N}}$, поточечно сходящаяся к равномерно непрерывному 
отображению 
$$f : X \rightarrow Y.$$
Кроме того, $f \in H_B (X,Y,\alpha)$, поскольку для всех $x, \, y \in X$
$$d (f (x), f (y)) = \lim \limits_{l \rightarrow \infty} d (f_l (x), f_l (y)) \leq
B |xy|^{\alpha}.$$ 
Из теоремы 2 следует теперь, что 
$$\lim \limits_{l \rightarrow \infty} \delta_p (f_l, f) =0.$$
и пространство $(H_B (X,Y,\alpha),\delta_p)$ собственное. 
\vskip7pt
$2. \, \,$
Ясно, что отображение $E_p$ сюръективно и непрерывно. Поэтому для каждого
$\varphi \in H_B (X,Y,\alpha)$ множество $E^{-1}_p (\varphi (p))$ замкнуто. Из доказательства утверждения $(i)$ леммы 1 следует, что для каждого $\theta \in E^{-1}_p (\varphi  (p))$
$$\delta_p (\varphi , \theta) \leq 2 B.$$
Следовательно, множество $E^{-1}_p (\varphi  (p))$ компактно, как замкнутое и ограниченное множество в собственном пространстве $(H_B (X,Y,\alpha),\delta_p)$. Пусть $y \in Y$ и $(\varphi _n)_{n \in \mathbb{N}}$ такая последовательность в замкнутом множестве $Z$ пространства $(H_B (X,Y,\alpha),\delta_p)$, что
$$\lim \limits_{n \rightarrow \infty} d (\varphi _n (p), y) = 0.$$
Последовательность $(\varphi _n)_{n \in \mathbb{N}}$ ограниченная, поскольку последовательность $(\varphi_n (p))_{n \in \mathbb{N}}$ сходится и
для фиксированного $\psi \in H_B (X,Y,\alpha)$ для каждого $n \in \mathbb{N}$ 
$$\delta_p (\varphi_n, \psi) \leq 2 B + d (\varphi_n (p), \psi (p)).$$
В силу собственности пространства $(H_B (X,Y,\alpha),\delta_p)$ найдутся
подпоследовательность $(\varphi _l)_{l \in \mathbb{N}}$ последовательности $(\varphi _n)_{n \in \mathbb{N}}$ и отображение $\varphi  \in H_B (X,Y,\alpha)$ такие, что
$$\lim \limits_{l \rightarrow \infty} \delta_p (\varphi_l, \varphi ) = 0.$$
Кроме того, 
$$d (\varphi (p), y) \leq d (\varphi (p), \varphi_l (p)) + d (\varphi_l(p), y) \rightarrow 0$$
при $l \rightarrow \infty$. Следовательно,
$$E_p (\varphi ) = \varphi (p) = y \, \, \mbox {и} \, \, \varphi \in Z,$$
поскольку множество $Z$ замкнуто.
Таким образом, $E_p$ --- замкнутое отображение, прообразы точек из $Y$ у которого компактны, то есть это отображение совершенное. 
\vskip7pt
$3. \, \,$
Если пространство $(Y,d)$ компактно, то $D (M) < + \infty$ и любая последовательность $(f_n)_{n \in \mathbb{N}}$ пространства $(H_B (X,Y,\alpha),\delta_p)$ ограничена. Используя доказательство первой части, получим, что пространство $(H_B (X,Y,\alpha),\delta_p)$ компактно. Теорема 3 доказана.

\section{Пространство всех подобий метрического пространства с метрикой Буземана}
\vskip20pt

Напомним некоторые определения.
\vskip7pt
Сюръекция $f : (X,\rho) \rightarrow (Y,d)$ называется {\it подобием с коэффициентом подобия
$\sigma_f > 0$}, если для любых $x, \, y \in X$ 
$$d(f(x),f(y)) = \sigma_f |xy|.$$

$Sim (X, Y)$ --- множество всех подобий из метрического пространства $(X,\rho)$ на метрическое пространство  $(Y,d)$. $Sim (X) = Sim (X, X)$.
\vskip7pt
Подобие $f$ называется {\it собственным подобием {\rm(}изометрией{\rm)}}, если $\sigma_f \neq 1$ $(\sigma_f = 1)$. 
\vskip7pt
$Iso (X, Y)$ --- множество всех изометрий из метрического пространства $(X,\rho)$ на метрическое пространство  $(Y,d)$. $Iso (X) = Iso (X, X)$.
\vskip7pt
Неподвижная точка собственного подобия называется {\it центром подобия}, а множество всех центров подобий пространства $X$ обознается $Cen (X)$. 
\vskip7pt
В следующей лемме 1 установлены некоторые простые геометрические свойства группы подобий метрического пространства.
\vskip7pt
{\bf Лемма 1} [s2]. {\it 
\vskip3pt
$(i)\, \,$ Отображение 
$$\sigma : (Sim (X),\circ) \rightarrow \mathbb{R}^{\ast}_{+}, \quad
\sigma (f) = \sigma_f$$ 
является непрерывным гомоморфизмом в группу положительных вещественных чисел $\mathbb{R}^{\ast}_{+}$ со стандартной операцией умножения, ядро которого есть группа изометрий.
\vskip7pt
$(ii)\, \,$ Для любых 
$$p \in X; \quad f,\, \varphi \in Sim (X,Y);\quad \theta \in Sim (Y),\quad \psi \in Iso (X)$$
имеют место равенства 
$$\delta_p (f\circ \psi,\varphi\circ \psi) = \delta_{\psi (p)} (f,\varphi),\quad
\delta_p (\theta\circ f,\theta\circ \varphi) = \sigma_{\theta} \delta_p (f,\varphi).$$
\vskip7pt
$(iii)\, \,$  Множества $Cen (X)$,  $X\backslash Cen (X)$ инвариантны относительно действия группы $Sim(X)$ на пространстве $X$.
\vskip7pt
$(iv)\, \,$ Если группа $Sim (X)$ действует транзитивно на полном метрическом пространстве $X$, то и группа $Iso (X)$ действует транзитивно на пространстве $X$.}
\vskip7pt
{\bf Доказательство леммы 1}.
\vskip7pt
$(i)\, \,$ Пусть $x$, $y \in X$ различны, $f$, $\varphi \in Sim (X,Y)$ произвольные. Тогда
$$d(f\circ \varphi(x), f\circ \varphi(y)) = \sigma_f \sigma_{\varphi} |xy|.$$ Следовательно,
$$ \sigma (f\circ \varphi) = \sigma(f) \sigma (\varphi).$$ 
Последнее утверждение очевидно. Докажем непрерывность $\sigma$. Пусть
$$\lim \limits_{n \rightarrow \infty} \delta_p (f_n,f) = 0,$$ 
где $f$, $f_n \in Sim (X,Y)$ для любого $n \in \mathbb{N}$. 
Тогда для любого $x \in X$
$$\lim \limits_{n \rightarrow \infty} d(f_n(x),f(x)) = 0.$$
Следовательно, для любых различных $x, \, y \in X$ 
$$\lim \limits_{n \rightarrow \infty} \sigma (f_n) |xy| = \lim \limits_{n \rightarrow \infty} d(f_n(x),f_n(y))$$ 
$$= d(f(x),f(y)) = \sigma (f) |xy|\, \, \mbox {и} \, \, \lim \limits_{n \rightarrow \infty}\sigma (f_n) = \sigma (f).$$
\vskip7pt
$(ii)\, \,$ Нетрудно проверить непосредственно.
\vskip7pt
$(iii)$ Пусть $x \in Cen (X)$, то есть найдется собственное подобие $f$ такое, что $f(x) = x$. Тогда для каждого $\psi \in  Sim (X)$
$$\sigma (\psi\circ f\circ \psi^{-1}) =\sigma (f).$$ 
Следовательно, $\psi\circ f\circ \psi^{-1}$ собственное подобие. Кроме того,
$$\psi\circ f\circ \psi^{-1}(\psi(x)) = \psi(x).$$ 
Таким образом, $\psi(x) \in Cen (X)$.
\vskip7pt
$(iv)$ Пусть $x$, $y \in X$. Тогда найдется такое подобие $f \in  Sim (X)$, что $f(x) = y$. Рассмотрим нетривиальный случай, когда подобие $f$ собственное. Тогда в силу полноты пространства найдется такая точка $z \in X$, что $f(z) = z$. Кроме того, найдется такое подобие $g$, что
$g(x) =z$. Тогда 
$$\sigma (f\circ g^{-1}\circ f^{-1}\circ g) = 1\, \, \mbox {и} \, \, f\circ g^{-1}\circ f^{-1}\circ g(x) = y.$$ 
Лемма 1 доказана.
\vskip7pt
{\bf Теорема 1} [s2]. {\it Если $(X,\rho)$, $(Y,d)$ --- полные метрические пространства, то пространства
$$(Sim (X,Y)\cup Const (X,Y),\delta_p),\quad (Iso (X,Y)\cup Const (X,Y),\delta_p)$$ 
--- полные.}
\vskip7pt
{\bf Доказательство теоремы 1}.
\vskip7pt
Пусть последовательность $(f_n)_{n \in \mathbb{N}}$ пространства 
\\
$(Sim (X,Y)\cup Const (X,Y),\delta_p)$ фундаментальная. Рассуждая так же, как и в доказательстве теоремы 3.3.1, получим отображение 
$f \in \Phi (X, Y)$ такое, что 
$$\lim \limits_{n \rightarrow \infty} \delta_p (f_n,f) = 0$$ 
и для любых $x, \, y \in X$
$$d(f(x), f(y)) = \lim \limits_{n \rightarrow \infty} d (f_n (x), f_n (y)) = \lim \limits_{n \rightarrow \infty} \sigma_{f_n} |xy| = \sigma_f |xy|.$$
Возможны следующие два случая.
\vskip7pt
$1. \, \,$ Если $\sigma_f = 0$, то $f \in Const (X,Y)$.
\vskip7pt
$2. \, \, $ Если $\sigma_f > 0$, то найдется такое $N \in  \mathbb{N}$, что $\sigma_{f_n}> 0$ для каждого $n > N$ и 
$f_n \in Sim (X,Y)$. Пусть $z \in Y$, $n >N$, тогда 
$$\lim \limits_{n \rightarrow \infty} |f^{-1}_n (z)p| = \lim \limits_{n \rightarrow \infty} (\sigma^{-1}_{f_n} d (z, f_n(p)))= \dfrac{1}{\sigma_f} d(z, f(p).$$
Следовательно, последовательность
$$(x_n = f^{-1}_n (z))_{n \in \mathbb{N}}$$
ограничена. Кроме того, для каждого
$\varepsilon > 0$ найдется такое 
$N_1 \in \mathbb{N}$, что для каждого $n > N_1$
$$d (f_n(x_n), f(x_n)) e^{-|px_n|}
= d (z_n, f(x_n)) e^{-|px_n|}< \varepsilon.$$
Поэтому
$$\lim \limits_{n \rightarrow \infty} f (x_n) = z.$$
Последовательность $(x_n)_{n \in \mathbb{N}}$ фундаментальная, так как для любых $n, \, m \in \mathbb{N}$
$$|x_nx_m|= \dfrac{1}{\sigma_f} d (f (x_n),f (x_m)).$$
Следовательно, эта последовательность сходится к некоторой точке $y \in X$ и
$$f (y) = f(\lim \limits_{n \rightarrow \infty} x_n) = \lim \limits_{n \rightarrow \infty} f(x_n) = z.$$
Таким образом, отображение $f$ сюръективно и $f \in Sim (X,Y).$
Заметим, что множество $Iso (X,Y)$ замкнуто в пространстве
$(Sim (X,Y)\cup Const (X,Y),\delta_p)$ и поэтому полное метрическое пространство. Теорема 1 доказана.
\vskip7pt 
{\bf Теорема 2} [s2]. {\it Если метрическое пространство $(X,\rho)$ --- собственное, то топология пространства
$(Sim (X,Y)\cup Const (X,Y),\delta_p)$ совпадает как с топологией поточечной сходимости, так и с компактно-открытой топологией.}
\vskip7pt
{\bf Замечание 1.} {\it Можно показать, что топология метрического пространства $(Sim (X), \delta)$ совпадает с топологией, определяемой на группе подобий $Sim (X)$ для каждой точки $p$ метрикой Буземана 
$\delta_p$ {\rm [s2]}.}
\vskip7pt
{\bf Доказательство теоремы 2}.
\vskip7pt
Пусть для всех $n \in \mathbb{N},\, \, x \in X$
$$\varphi_n, \, \varphi \in Sim (X,Y)\cup Const (X,Y)\, \, \mbox {и}\, \,
\lim \limits_{n \rightarrow \infty} d (\varphi_n (x),\varphi (x)) = 0.$$
Тогда для любых $x, \, y \in X$
$$\lim \limits_{n \rightarrow \infty} \sigma_{\varphi_n} |xy| = \lim \limits_{n \rightarrow \infty} d(\varphi_n (x), \varphi_n (y)) = d (\varphi (x), \varphi (y)) = \sigma_{\varphi} |xy|.$$
Следовательно, найдется такая константа $B \geq 0$, что
для каждого $n \in \mathbb{N}$
$$\sigma_{\varphi_n} \leq B \, \, \mbox {и}\, \, \sigma_{\varphi}\leq B.$$
Дальнейшие рассуждения аналогичны рассуждениям, приведенным в доказательстве теоремы 3.3.2. Теорема 2 доказана.
\vskip7pt
{\bf Теорема 3} [s2]. {\it Если метрические пространства $(X,\rho)$, $(Y,d)$ --- собственные, то пространство
$(Sim (X,Y)\cup Const (X,Y),\delta_p)$ --- собственное. Если кроме того, пространство $(X,\rho)$ --- компактно, то пространство $(Sim (X,Y)\cup Const (X,Y),\delta_p)$ --- компактно.}
\vskip7pt
{\bf Доказательство теоремы 3}.
\vskip7pt
Пусть последовательность $(f_n)_{n \in \mathbb{N}}$ пространства
\\
$(Sim (X,Y)\cup Const (X,Y),\delta_p)$ ограничена. Тогда для каждого $x \in X$ последовательность $(f_n(x))_{n \in \mathbb{N}}$ ограничена.
Кроме того, для всех 
$$n \in \mathbb{N}; \quad q \in Y; \quad x,\, y \in X$$
$$\sigma_{f_n} |xy| = d (f_n(x), f(_n(y)) \leq d (f_n(x),q) + d(q,f_n(y)).$$
Поэтому последовательность $(\sigma_{f_n})_{n \in \mathbb{N}}$ ограничена.  Дальнейшие рассуждения аналогичны рассуждениям, приведенным в доказательстве теоремы 3.3.3, если учесть теорему 2 и доказательство теоремы 1. Теорема 3 доказана.
\vskip7pt
{\bf Теорема 4} [s1], [s2]. {\it 
\vskip3pt
$(i)\, \,$ $(Sim (X),\delta_p)$ --- топологическая группа, действующая непрерывно на пространстве $X$.
\vskip7pt
$(ii)\, \,$ Группы подобий $Sim (X)$ 
 и изометрий $Iso (X)$ с метрикой  $\delta$ являются топологическими группами, непрерывно действующими на пространстве $X$.}
\vskip7pt
{\bf Замечание 1.} {\it Доказательство того факта, что $(Iso (X),\delta_p)$ --- топологическая группа, действующая непрерывно на пространстве $X$, известно} (см. \cite {Beres2}). 
\vskip7pt
{\bf Доказательство теоремы 4}.
\vskip7pt
$(i)\, \,$ Пусть $p \in X$, для каждого $n \in \mathbb{N}$ 
$$f, \, g, \, f_n, \, g_n \in Sim (X)\, \, \mbox {и}\, \,  \lim \limits_{n \rightarrow \infty} \delta_p (f_n,f) = 0, \quad \lim \limits_{n \rightarrow \infty} \delta_p (g_n,g) =0.$$
Рассмотрим неравенство
$$\delta_p (f_n\circ g_n,f\circ g) \leq \delta_p (f_n\circ g_n,f_n\circ g) + \delta_p (f_n\circ g,f\circ g)=$$ 
$$\sigma_{f_n} \delta_p (g_n,g) + \sup \{|f_n(y)f (y)| e^{-\sigma_{g^{-1}} |g(p)y|} : \, y \in X \}.$$
По утверждению $(i)$ леммы 1 первое слагаемое в правой части стремится к нулю при $n \rightarrow \infty$. Докажем методом от противного, что и второе слагаемое стремится к нулю. Пусть, напротив, найдется подпоследовательность $(f_l)_{l \in \mathbb{N}}$ последовательности $(f_n)_{n\in \mathbb{N}}$, последовательность $(x_l)_{l \in \mathbb{N}}$ пространства $X$ и константа $A > 0$ такие, что для каждого $l \in \mathbb{N}$
$$|f_l(x_l)f(x_l)| e^{-\sigma_{g^{-1}} |g(p)x_l|} > A.$$   
Тогда для каждого $l \in \mathbb{N}$
$$0 < A < |f_l(x_l)f(x_l)| e^{-\sigma_{g^{-1}} |g(p)x_l|} \leq
\{|f_l(x_l)f_l(g(p))| +$$ 
$$|f_l(g(p))f(g(p))| + |f(g(p))f(x_l)|\} e^{-\sigma_{g^{-1}} |g(p)x_l|} =$$ 
$$\{(\sigma_{f_l} + \sigma_f) |x_lg(p)|+ |f_l(g(p))f(g(p))|\}e^{-\sigma_{g^{-1}} |g(p)x_l|}. \eqno (1)$$
$1.\, \,$ Пусть последовательность $(x_l)_{l \in \mathbb{N}}$ неограничена. Значит, найдется такая подпоследовательность $(x_t)_{t \in \mathbb{N}}$ последовательности $(x_l)_{l \in \mathbb{N}}$, что
$$\lim \limits_{t \rightarrow \infty} |x_tg(p)| =\infty.$$
Тогда правая часть неравенства $(1)$ после замены $x_l$ на $x_t$ стремится к нулю при $t \rightarrow \infty$. Получаем противоречие.
\vskip7pt
$2.\, \,$ Пусть последовательность $(x_l)_{l \in \mathbb{N}}$ ограничена, то есть найдется такая константа $c > 0$, что  
$|px_l| < c$ для каждого $l \in \mathbb{N}.$ 
Тогда для каждого 
$\varepsilon > 0$ найдется такое $N \in \mathbb{N}$, что 
$$\delta_p (f_l,f) < \varepsilon, \quad 0 < A < |f_l(x_l)f(x_l)| e^{-\sigma_{g^{-1}} |g(p)x_l|} \leq$$
$$|f_l(x_l)f(x_l)| < \varepsilon e^{|px_l|} < \varepsilon e^c$$
при $l > N$.
Снова получили противоречие. Следовательно, 
$$\lim \limits_{n \rightarrow \infty} \delta_p (f_n\circ g_n,f\circ g) = 0.$$
Кроме того, 
$$\lim \limits_{n \rightarrow \infty} \delta_p (f^{-1}_n,f^{-1}) = \lim \limits_{n \rightarrow \infty} (\sigma_{f^{-1}_n} \delta_p (id,f_n \circ f^{-1})) = \lim \limits_{n \rightarrow \infty} (\dfrac{1}{\sigma_{f_n}} \delta_p (f \circ f^{-1},f_n \circ f^{-1})) = 0.$$
Непрерывность действия группы $Sim (X)$ на пространстве $X$ следует из неравенства:
$$|f(x)f_0(x_0)| \leq |f(x)f_0(x)| + |f_0 (x)f_0(x_0)|\leq \delta_{x_0} (f, f_0) e^{|xx_0|} + \sigma_{f_0} |xx_0|$$
для всех
$$x, \, x_0 \in X\quad f, \, f_0 \in Sim (X).$$ 
$(ii)\, \,$ Согласно теореме 3.1.2 необходимо доказать лишь следующее утверждение: 
если для каждого $n \in \mathbb{N}$
$$f_n \in Sim (X) \, \, \mbox \, \, \lim \limits_{n \rightarrow \infty} \delta (f_n, id)  = 0,$$
то  
$$\lim \limits_{n \rightarrow \infty} \delta (f^{-1}_n, id)  = 0.$$
В силу свойства 3.1.2 получим для любого $x \in X$ 
$$\lim \limits_{n \rightarrow \infty} |f_n (x)x|  = 0.$$
Следовательно, для любых $x, \, y \in X$ 
$$\lim \limits_{n \rightarrow \infty} \sigma_{f_n} |xy| =\lim \limits_{n \rightarrow \infty}|f_n (x)f_n (y)|  = |xy|,$$
$$\lim \limits_{n \rightarrow \infty} \sigma_{f_n} = 1, \quad
\lim \limits_{n \rightarrow \infty} \sigma_{f^{-1}_n} = \lim \limits_{n \rightarrow \infty} \dfrac{1}{\sigma_{f_n}} = 1.$$
Кроме того, для каждого $i \in \mathbb{N}$ 
$$\lim \limits_{n \rightarrow \infty}\delta_i (f^{-1}_{n,i}, id_i) =\lim \limits_{n \rightarrow \infty}\sup \{|f^{-1}_n (x)x| : \, x \in B_i\} =$$ $$\lim \limits_{n \rightarrow \infty} (\sigma_{f^{-1}_n} \sup \{|xf_n (x)| : \, x \in B_i\}) =$$
$$\lim \limits_{n \rightarrow \infty} \sigma_{f^{-1}_n} \delta_i (f_{n,i}, id_i) = 0.$$
Осталось применить свойство 3.1.1. Теорема 4 доказана.
\vskip7pt
{\bf Теорема 5} [s2]. {\it Пусть пространство $(X,\rho)$ --- полное. Тогда верны следующие утверждения.
\vskip7pt
$(i)\, \,$ Если $Cen (X) = X$, то
$$\overline{Sim (X)} = Sim (X)\cup Const (X,X),$$
где замыкание выполняется в пространстве $(\Phi (X,X),\delta_p)$.
\vskip7pt
$(ii)\, \,$ Если $Cen (X) \neq \emptyset, \, X$, то 
$$Cen (X) \subset \partial (Cen (X))\, \, \mbox {и}\, \, \overline{Sim (X)} = Sim (X)\cup Const (X,\partial (Cen (X))).$$
$(iii)\, \,$ Если $Cen (X) = \emptyset$, то $\overline{Sim (X)} = Sim (X) = Iso (X)$.}
\vskip7pt
{\bf Доказательство теоремы 5}.
\vskip7pt
$(i)\, \,$ Пусть $p \in X$, $\psi \in Const (X,X)$, то есть найдется такое $x_0 \in X = Cen (X)$, что $\psi (x) = x_0$ для любого $x \in X$. Пусть $f$ такое собственное подобие с коэффициентом $\sigma_f < 1$, что $f(x_0) = x_0$. Тогда для всех $k \in \mathbb{N},\quad x \in X$
$$|xf^k(x)| \leq |xx_0| + |x_0f^k(x)| = (1 + \sigma^k_f) |xx_0|.$$
Последовательность $(f^n)_{n \in \mathbb{N}}$ --- фундаментальная, поскольку
для любых $n, \, k \in \mathbb{N}$ 
$$\delta_p (f^n,f^{n+k}) = \sigma^n_f  \delta_p (id,f^k) \leq \sigma^n_f  \sup \{|xf^k (x)| e^{-|px|} :\, x \in X\} \leq$$
$$\sigma^n_f  \sup \{(1 + \sigma^k_f) |x_0x| e^{-|px|} :\, x \in X\} \leq
\sigma^n_f  (1 + \sigma^k_f) (1 + |px_0|).$$
По теореме 1 найдется такое $\varphi \in Sim (X)\cup Const (X,X)$, что
$$\lim \limits_{n \rightarrow \infty} \delta_p (f^n, \varphi) = 0.$$
Но для любого $x \in X$ 
$$\lim \limits_{n \rightarrow \infty}|f^n(x)\psi(x)| =
\lim \limits_{n \rightarrow \infty}|f^n(x)x_0| = 0.$$
Поэтому для любого $x \in X$
$$\varphi(x) = \lim \limits_{n \rightarrow \infty} f^n (x) = \psi(x) \, \, \mbox {и} \, \, \varphi = \psi.$$

$(ii)\, \,$ Пусть $\psi \in Const (X,X)\cap \overline{Sim (X)}$, то есть найдутся такое $x_0 \in X$ и последовательность $(f_n)_{n \in \mathbb{N}}$ пространства $Sim (X)$, что для любого $x \in X$
$$\psi (x) = x_0\, \, \mbox {и} \, \, \lim \limits_{n \rightarrow \infty} \delta_p (f^n, \psi) = 0.$$
Следовательно, для любого $x \in X$
$$\lim \limits_{n \rightarrow \infty} |f^n (x)\psi (x)| = \lim \limits_{n \rightarrow \infty} |f^n (x)x_0| = 0.$$
Если $x \in Cen (X)$, то по утверждению $(iii)$ леммы 1 $f^n (x) \in Cen (X)$ для любого $n \in \mathbb{N}$.
Следовательно, $x_0 \in \overline{Cen (X)}.$

Если $x \in X\backslash Cen (X)$, то по утверждению $(iii)$ леммы 1 $ f^n (x) \in X\backslash Cen (X)$ для любого $n \in \mathbb{N}$.
Следовательно, $x_0 \in \overline{(X\backslash Cen (X))}.$
Поэтому 
$$x_0 \in \partial (Cen (X))\, \, \mbox {и} \, \, \psi \in Const (X,\partial (Cen (X))).$$
Заметим также, что 
$$Cen (X) \subset \partial (Cen (X)).$$ 
Действительно, если $y \in Cen (X)$, то найдется такое собственное подобие $f$, что для каждого $x \in X$ 
$$\lim \limits_{n \rightarrow \infty} |f^n (x)y| = 0.$$
Тогда $y \in \partial (Cen (X))$, поскольку по утверждению $(iii)$ леммы 1:
если $x \in \partial (Cen (X))$, то для любого $n \in \mathbb{N}$ 
$$f^n (x) \in \partial (Cen (X)).$$
Пусть 
$$\psi \in Const (X,\partial (Cen (X))) \subset Const (X,\overline{Cen (X)}) =
\overline{Const (X,Cen (X))}.$$
Докажем, что 
$$\overline{Const (X,Cen (X))} \subset \overline{Sim (X)}.$$
Действительно, если $\varphi \in Const (X,Cen (X))$, то найдется такое $x_0 \in X$, что для каждого $x \in X$
$$\varphi (x) = x_0 \in Cen (X).$$
Следовательно, найдется такое собственное подобие $f$, что для каждого $x \in X$
$$\lim \limits_{k \rightarrow \infty} |f^k (x)\varphi (x)| = \lim \limits_{k \rightarrow \infty} |f^k (x)x_0| = 0.$$
Повторяя доказательство пункта $(i)$, получим
$$\lim \limits_{n \rightarrow \infty} \delta_p (f^n, \varphi) = 0.$$
Следовательно,
$$\varphi \in \overline{Sim (X)}\, \, \mbox {и} \, \, Const (X,\partial (Cen (X))) \subset  \overline{Sim (X)}.$$

$(iii)\, \,$ Это следствие теоремы 1 (см. также \cite [с. 32]{Bus}). Теорема 5 доказана.
\vskip7pt
{\bf Пример.} 
\vskip7pt
Пусть $X$ --- замкнутое полупространство вещественного гильбертова пространства $l_2$ с граничной гиперплоскостью $X_0$, содержащей нулевой вектор. Тогда нетрудно проверить, что
$$Cen (X) = \partial (Cen (X)) = X_0, \quad \overline{Sim (X)} = Sim (X)\cup Const (X,X_0).$$

\section{Два аналога слабой сходимости в специальном метрическом пространстве}
\vskip20pt

Пусть в метрическом пространстве $(X,\rho)$ выделено семейство $S_p$
сегментов с общим началом в точке $p \in X$, удовлетворяющее следующему условию $(A_{11})$.
\vskip7pt
$(A_{11})\quad$ Для каждой точки $x \in X$ найдется единственный сегмент $[p,x]\in S_p$ с концами $p$ и $x$.
\vskip7pt
Нам понадобятся также некоторые дополнительные условия на пространство $X$. 
\vskip7pt
$(A_{12})\quad$ Пусть пространство $X$ удовлетворяет условию $(A_{11})$. Существует такое семейство ориентированных прямых $\hat S_p$ пространства, проходящих через точку $p$, что каждый невырожденный ориентированный сегмент $[p,x]$ из семейства $S_p$ содержится в единственной ориентированной прямой из $\hat S_p$, ориентация которой определяется этим сегментом. 
\vskip7pt
Пусть $\Delta OAB$ --- треугольник в евклидовой плоскости с длинами
сторон 
$$||OA|| = |px|,\quad ||OB|| =|py|,\quad ||AB||=|xy|,$$
где $p$, $x$, $y \in X$.
\vskip7pt
$\gamma (OA,OB)$ --- величина угла при вершине $O$ в треугольнике  $\Delta OAB$.
$$\bar \gamma (x,y) = \varlimsup \limits_ {||OA|| \rightarrow 0, \, ||OB|| \rightarrow 0} {\gamma (OA,OB)}$$ 
--- верхний угол между сегментами $[p,x]$, $[p,y]\in S_p$ пространства $X$ удовлетворяющего условию $(A_{11})$ \cite {Beres}.
\vskip7pt
Пусть для любого $y \in X$
$$\psi _y : X \rightarrow \mathbb{R},\quad \psi _y (x) = |px||py| \cos \bar \gamma (x,y)$$
при $x\in X\backslash \{p\}$ и $\psi _y (x) = 0$ при $x = p$ или $y = p$.
$$\Psi_p (X) = \{\psi_y : X \rightarrow \mathbb{R} :\, y \in X\}.$$

$(A_{13})\quad$ Пространство $X$ удовлетворяет условию $(A_1)$ и все замкнутые шары пространства $X$ строго выпуклые.
\vskip7pt
$(A_{14})\quad$ Пространство $X$ удовлетворяет условию $(A_{11})$ и тождественное отображение подпространства $(X\backslash \{p\}, \rho)$ на псевдометрическое пространство $(X\backslash \{p\}, \bar \gamma)$ является непрерывным.
\vskip7pt
Отметим, что пространство $(X\backslash \{p\}, \bar \gamma)$ является псевдометрическим в силу известных свойств верхнего угла между сегментами с общим концом и условие $(A_{14})$ выполняется для полных пространств неположительной кривизны по А. Д. Александрову, являющихся областями $R_K$ при $K \leq 0$ \cite {Beres}. 

Пусть выполняются условия $(A_{12}-A_{13})$ и $L \in \hat S_p$, тогда для любой точки $x \in X$ cуществует единственная точка $P_L (x)$ \cite [c. 22]{Bus} и определено отображение 
$$\varphi_L : X \rightarrow \mathbb{R}, \quad \varphi_L (x) = \varepsilon_L (x) |pP_L (x)|,$$  
где $\varepsilon_L (x)$ равняется $0$, $1$ или $-1$, если соответственно 
$P_L (x) =p$,  сегмент $[p,P_L (x)]$ и ориентированная прямая $L$ сонаправлены, сегмент $[p,P_L (x)]$ и ориентированная прямая $L$ противоположно направлены.
Введем следующие обозначения:
$$\Phi_p (X) = \{\varphi_L : X \rightarrow \mathbb{R} : \, L \in
\hat S_p \};$$
$$F_p (X) = \{ f : X \rightarrow \mathbb{R} : \, f (p) = 0, \,
||f|| = \sup \{\dfrac{|f (x)|}{|px|} : x \neq p \} < +\infty \};$$
$$FC_p (X) = \{ f  \in F_p (X) : \, f \, \mbox{--- непрерывно} \}.$$

$\hat \Psi_p (X)$ $(\hat \Phi_p (X))$ --- замкнутая линейная оболочка множества $\Psi_p (X)$ $(\Phi_p (X))$ в линейном нормированном пространстве $(F_p (X), ||.||)$.
\vskip7pt
Отметим, что $FC_p (X)$ является замкнутым линейным подпространством в
банаховом пространстве $(F_p (X), ||.||)$. Это непосредственное следствие
теоремы 1 из \cite [c. 181]{Burb} и следствия 2 теоремы 2 из \cite [c. 184]{Burb}.
\vskip7pt
{\bf Лемма 1} [s15]. {\it 
\vskip3pt
$(i)\, \,$ Пусть пространство $X$ удовлетворяет условию $(A_{14})$. Тогда $\Psi_p (X) \subset FC_p(X)$ и
$\psi(\omega_{\lambda}(p,x)) = \lambda\psi(x)$ для любых $\psi \in \hat \Psi_p(X)$,
$\lambda \in \mathbb{R}_{+}$, $x \in X$. 
\vskip7pt
$(ii)\, \,$ Пусть пространство $X$ удовлетворяет условиям $(A_{12})$, $(A_{13})$. Тогда $\Phi_p (X) \subset FC_p (X)$. 
\vskip7pt
$(iii)\, \,$ Пусть пространство $X$ удовлетворяет условию $(A_{11})$ и для любых $x,\, y \in X\backslash \{p\}$  из условия $\bar \gamma (x, y) = 0$ следует, что один из сегментов $[p,x]$, $[p,y] \in S_p$ принадлежит другому сегменту. Тогда для любых различных
$x, \, y \in X$ найдется такая точка $z \in X\backslash \{p\}$, что $\psi _z (x) \neq \psi _z (y)$. 
\vskip7pt
$(iv)\, \,$ Пусть пространство $X$ удовлетворяет условиям $(A_{12})$, $(A_{13})$ и для всех $x \in X, \, \, L \in \hat S_p$ 
$$|pP_L(x)| \leq  |px|.$$ 
Тогда для любых различных
$x, \, y \in X$ найдется $L \in \hat S_p$ такое, что
$\varphi _L(x) \neq \varphi _L(y)$.}
\vskip7pt
{\bf Доказательство леммы 1}. 
\vskip7pt
$(i)\, \,$ Для любого $y \in X$ непрерывность функции $\psi_y$ следует из непрерывности метрики $\rho$, косинуса и условия $(A_14)$. Кроме того,
$||\psi_y || = 0$ при $y = p$ и
$$||\psi_y || =\sup \{|py| \cos \bar \gamma (x,y) : \, x \neq p\} = |py| < + \infty$$ 
при $y\neq p$.
Второе утверждение нетрудно получить из определений 
$$\hat \Psi_p(X),\quad \omega_{\lambda} (p,x), \quad \psi_y $$ 
для $y \in X$.
\vskip7pt
$(ii)\, \,$ Утверждение следует из равенств
$$|\varphi _L(x) - \varphi _L(y)| = |\varepsilon_L(x) |pP_L (x)| -
\varepsilon_L(y) |pP_L (y)||= |P_L (x)P_L (y)|$$ 
и непрерывной зависимости $P_L (x) \in L$ от $x \in X$ \cite [c. 22]{Bus}.
\vskip7pt
$(iii)\, \,$ Докажем методом от противного. Пусть найдутся различные $x,\, y
\in X$ такие, что $\psi_z (x) = \psi_z (y)$ для каждого $z \in X\backslash \{p\}.$ Если $x = p$ ($y = p$), то при $z =y$ ($z = x$) найдем
$$0 = \psi_z (p) =  \psi _z (z) = |pz|^2.$$
Получили противоречие.
Если $x \neq  p$, $y \neq p$, то из предположения и равенств
$$\psi _x (x) =|px|^2 = \psi _x (y) =  |px| |py| \cos \bar \gamma (x,y),$$
$$\psi _y (y) = |py|^2 = \psi _y (x) =  |px| |py| \cos \bar \gamma (x,y),$$
следует $x=y$. Получили противоречие.
\vskip7pt
$(iv)\, \,$ Докажем методом от противного. Пусть найдутся различные $x, \, y
\in X$ такие, что для каждого $L \in\hat S_p$
$\varphi _L (x) = \varphi _L (y)$. Если $x = p$ ($y = p$), то
при $L = L_y$ ($L = L_x$), где ориентация прямой
$L_y$, содержащей точку $y$, определяется сегментом $[p,y] \in S_p$ (аналогично для $L_x$), найдем 
$$0 = \varphi _{L_y} (p) = \varphi _{L_y} (y) = |py|\quad (0 = \varphi _{L_x} (p) = \varphi _{L_x}(x) = |px|).$$ 
Получили противоречие. Если $x \neq  p$, $y \neq p$, то
$$\varphi _{L_x}(x) = |px| = \varphi _{L_x}(y) = \varepsilon_{L_x}(y) |pP_{L_x}(y)| \leq |py|,$$
$$\varphi _{L_y}(y) = |py| =
\varphi _{L_y}(x) = \varepsilon_{L_y}(x) |pP_{L_y}(x)|\leq |px|.$$
Следовательно,
$$P_{L_x}(y) =x,\quad P_{L_y}(x)=y,\quad |px| = |py|.$$ 
Пусть $\lambda \in (0,1)$. Тогда из полученных равенств, неравенства
$$|pP_{L_x} (\omega_{\lambda}(x,y))| \leq  |p\omega_{\lambda}(x,y)|$$
и строгой  выпуклости шара $B[p,|px|]$ получим
$$|pP_{L_x} (\omega_{\lambda}(x,y))| \leq  |p\omega_{\lambda}(x,y)| < |px|.$$ 
Но 
$$P_{L_x}(\omega_{\lambda}(x,y)) = x, \, \, P_{L_{\omega_{\lambda} (x,y)}} (x) = \omega_{\lambda}(x,y)\quad \, \, \mbox {и}\, \, |px| = |p\omega_{\lambda}(x,y)|.$$ 
Получили противоречие. Лемма 1 доказана.
\vskip7pt
Последовательность $(x_n)_{n \in \mathbb{N}}$ пространства $X$, удовлетворяющего условию $(A_{14})$,
назовем {\it $\psi$-сходящейся} к точке $x \in X$,
если для каждого  $\psi \in \hat \Psi_p (X)$ 
$$\lim \limits_{n\rightarrow \infty}{\psi (x_n)} = \psi(x).$$

Последовательность $(x_n)_{n \in \mathbb{N}}$ пространства $X$, удовлетворяющего условиям $(A_{12})$, $(A_{13})$,
назовем {\it  $\varphi$-сходящейся} к точке $x \in X$,
если для каждого  $\varphi \in \hat \Phi_p(X)$
$$\lim \limits_{n\rightarrow \infty}{\varphi (x_n)} =  \varphi (x).$$
Следующие теоремы 1, 2 аналогичны соответствующим теоремам 1, 2 из
\cite [с. 225-227]{Kolm} о свойствах слабо сходящихся последовательностей
в нормированном пространстве.
\vskip7pt
{\bf Теорема 1} [s15]. {\it 
\vskip3pt
$(i)\, \,$ Если последовательность $(x_n)_{n \in \mathbb{N}}$ пространства $X$, удовлетворяющего условию $(A_{14})$,
$\psi$-сходится к точке $x \in X$, то она ограничена.
\vskip7pt
$(ii)\, \,$ Если последовательность $(x_n)_{n \in \mathbb{N}}$ пространства $X$, удовлетворяющего условиям $(A_{12})$, $(A_{13})$,
$\varphi$-сходится к точке $x \in X$, то она ограничена.}
\vskip7pt
{\bf Доказательство теоремы 1}. 
\vskip7pt
Рассмотрим в пространстве $\hat \Psi_p(X)$
$(\hat \Phi_p(X))$ для всех $k, \, n \in \mathbb{N}$ замкнутые множества
$$A_{kn} = \{f : \,|f (x_n)| \leq k \},\quad A_k =\bigcap \limits_{n=1}^{\infty} A_{kn}.$$
Для каждого $f \in \hat \Psi_p(X)$ $(f \in \hat \Phi_p(X))$
последовательность $(f (x_n))_{n \in \mathbb{N}}$ ограничена,
так как последовательность
$(x_n)_{n \in \mathbb{N}}$ $\psi$-сходится ($\varphi$-сходится). Тогда
$$\hat \Psi_p (X)= \bigcup \limits_{k=1}^{\infty}A_k \quad (\hat \Phi_p(X) = \bigcup \limits_{k=1}^{\infty}A_{k}).$$ 
По теореме Бэра \cite [c.83]{Kolm} найдется замкнутое множество $A_{k_0}$, плотное в некотором замкнутом шаре
$$B[f_0,r] \subset \hat \Psi_p(X)\quad (B[f_0,r] \subset \hat \Phi_p(X))$$ 
и, следовательно, содержащее этот шар. Рассмотрим отображение 
$$\pi_1 : X \rightarrow \hat \Psi_p^{*}(X), \quad \pi_1 (y) (f) = f (y)$$
$$(\pi_2 : X \rightarrow \hat \Phi_p^{*}(X),\quad \pi_2 (y) (f) = f (y)).$$ 
Тогда последовательность 
$$(\pi_1 (x_n))_{n \in \mathbb{N}}\quad  ((\pi_2 (x_n))_{n \in \mathbb{N}})$$
ограничена на шаре $B[f_0,r]$ и, следовательно, на каждом шаре
пространства $\hat \Psi_p(X)$ $(\hat \Phi_p(X))$. Кроме того, $||\psi _y||= |py|$ для каждого $y \in X$.
Следовательно,
$$||\pi_1 (x_n)|| = \sup\{\dfrac{|\psi (x_n)|}{||\psi||} : \psi \neq 0\} = |px_n|$$
и последовательность $(x_n)_{n \in \mathbb{N}}$ ограничена.
Отметим, что для каждой ориентированной прямой $L \in \hat S_p$ имеют место неравенства
$$ 1 \leq ||\varphi_L|| \leq 2,$$ 
поскольку первое неравенство очевидно, а второе следует из определения нормы и неравенств 
$$|pP_L (x)| \leq |px| + |xP_L (x)| \leq 2|px|$$
для каждого $x \in X$. 
Следовательно, 
$$||\pi_2 (x_n)|| = \sup\{\dfrac{|\varphi (x_n)|}{||\varphi||} : \varphi \neq 0\} \geq \dfrac{|px_n|}{2}$$ 
и последовательность $(x_n)_{n \in \mathbb{N}}$ ограничена.
Теорема 1 доказана.
\vskip7pt
{\bf Теорема 2} [s15]. {\it Пусть последовательность  $(x_n)_{n \in \mathbb{N}}$ ограничена в пространстве $X$, $x \in X$. Тогда верны следующие утверждения.
\vskip7pt
$(i)\, \, $ Если пространство $X$ удовлетворяет условию $(A_{14})$ и для каждого $y \in X$ 
$$\lim \limits_{n\rightarrow \infty}{\psi_y (x_n)} =
\psi_y (x),$$
то последовательность $(x_n)_{n \in \mathbb{N}}$ $\psi$-сходится к точке $x \in X$.
\vskip7pt
$(ii)\, \, $ Если пространство $X$ удовлетворяет условиям $(A_{12})$, $(A_{13})$  и для каждого $L \in \hat S_p$ 
$$\lim \limits_{n\rightarrow \infty}{\varphi_L (x_n)} = \varphi_L (x),$$
то последовательность $(x_n)_{n \in \mathbb{N}}$ $\varphi$-сходится к точке $x \in X$.}
\vskip7pt
Из утверждений $(i)$, $(ii)$ леммы 1 и теоремы 2 получим следствия 1, 2.
\vskip7pt
{\bf Следствие 1} [s15]. {\it  Если последовательность $(x_n)_{n \in \mathbb{N}}$ пространства $X$, удовлетворяющего условию $(A_{14})$ {\rm {(}}условиям $(A_{12})$, $(A_{13}))$, сходится к точке $x \in X$, то она $\psi$-сходится
$(\varphi$-сходится{\rm {)}} к точке $x \in X$.}
\vskip7pt
{\bf Следствие 2} [s15]. {\it Если пространство $X$, удовлетворяет условию $(A_{14})$ {\rm {(}} условиям $(A_{12})$, $(A_{13}))$, то
$\hat \Psi_p(X) \subset FC_p(X)$ $(\hat \Phi_p(X) \subset FC_p(X))$.}
\vskip7pt
{\bf Доказательство теоремы 2}.  
\vskip7pt
Пусть $f$ из линейной оболочки множества
$$\Phi_p(X) \subset \hat \Phi_p(X)\quad (\Psi_p(X) \subset \hat \Psi_p(X)).$$ Тогда в силу условия теоремы 2
$$\lim \limits_{n\rightarrow \infty }{f (x_n)} = f(x).$$ 
Для произвольного элемента
$$f \in \hat \Phi_p(X)\quad (f \in \hat \Psi_p(X))$$ 
найдется последовательность
$(f_n)_{n \in \mathbb{N}}$ из линейной оболочки множества
$$\Phi_p(X) \subset \hat \Phi_p(X)\quad (\Psi_p(X) \subset \hat \Psi_p(X))$$ такая, что 
$$f = \lim \limits_{n\rightarrow \infty}{f_n}.$$
Выберем такую константу $M  > 0$, чтобы для каждого $n \in \mathbb{N}$  
$$|px_n| \leq M,\quad|px| \leq M.$$
Для каждого $\varepsilon > 0$ найдется такое число
$K \in \mathbb{N}$, что для каждого $k \geq K$
$$ ||f - f_k|| < \dfrac{\varepsilon}{6M}.$$ 
Фиксируем  $k \geq K$. Тогда найдется такое натуральное число $N$, что для каждого $n \geq N$
$$|f_k (x_n) - f_k (x)| < \dfrac{\varepsilon}{2}.$$ 
Следовательно,  для каждого $n \geq N$ 
$$|f (x_n) - f (x)| \leq |f (x_n) - f_k (x_n)| +|f_k (x_n) - f_k (x)| + |f_k (x) - f(x)| \leq$$
$$||f - f_k||(|px| + |px_n|) + \dfrac{\varepsilon}{2} < \varepsilon.$$
Теорема 2 доказана. 
\vskip7pt
Пусть пространство $X$ удовлетворяет условию $(A_{14})$. Последовательность $(\psi_n)_{n \in \mathbb{N}}$ пространства $\hat \Psi_p(X)$
назовем {\it слабо сходящейся к
$\psi \in \hat \Psi_p (X)$}, если для каждого $x \in X$ 
$$\lim \limits_{n\rightarrow \infty}{\psi_n (x)} = \psi (x).$$

Пусть пространство $X$ удовлетворяет условиям $(A_{12})$, $(A_{13}))$. Последовательность $(\varphi_n)_{n \in \mathbb{N}}$ пространства $\hat \Phi_p(X)$ назовем {\it слабо сходящейся к $\varphi \in \hat \Phi_p(X)$}, если для каждого $x \in X$ 
$$\lim \limits_{n\rightarrow \infty}{\varphi_n (x)} =  \varphi (x).$$
Следующая лемма 2 аналогична соответствующей теореме 1* из \cite [с. 230]{Kolm}.
\vskip7pt
{\bf Лемма 2} [s15]. {\it Если $X$ --- полное метрическое пространство, удовлетворяющее условию $(A_{14})$ {\rm {(}} условиям $(A_{12})$, $(A_{13}))$, и последовательность 
$(\psi_n)_{n \in \mathbb{N}}$ $((\varphi_n)_{n \in \mathbb{N}})$ пространства $\hat \Psi_p(X)$  $(\hat \Phi_p(X))$
слабо сходится к $\psi \in \hat \Psi_p(X)$ $(\varphi \in \hat \Phi_p(X))$, то эта последовательность ограничена на некотором шаре пространства $X$.}
\vskip7pt
{\bf Доказательство леммы 2}. 
\vskip7pt
Для всех $n, \, k \in \mathbb{N}$ рассмотрим замкнутые множества
$$A_{kn} = \{x \in X : \,|\psi_n (x)| \leq k \}\quad (A_{kn} = \{x \in X : |\varphi_n(x)| \leq k \}),$$
$$A_k =\bigcap \limits_{n=1}^{\infty}A_{kn}.$$ 
Для каждого $x \in X$ последовательность 
$$(\psi_n (x))_{n \in \mathbb{N}}\quad ((\varphi _n (x))_{n \in \mathbb{N}})$$ ограничена, так как последовательность
$$(\psi_n)_{n \in \mathbb{N}}\quad ((\varphi _n)_{n \in \mathbb{N}})$$ пространства $\hat \Psi_p(X)$  $(\hat \Phi_p(X))$ слабо сходится. Тогда 
$$X = \bigcup \limits_{k=1}^{\infty}A_k.$$ 
По теореме Бэра \cite [c.83]{Kolm} найдется замкнутое множество $A_{k_0}$ плотное в некотором замкнутом шаре
$B[x_0,r]$ пространства $X$ и, следовательно, содержащее этот шар. Значит,
последовательность 
$$(\psi_n)_{n \in \mathbb{N}}\quad ((\varphi _n)_{n \in \mathbb{N}})$$ ограничена на этом шаре. Лемма 2 доказана.
\vskip7pt
{\bf Пример 1.} 
\vskip7pt
Нетрудно проверить, что в вещественном гильбертовом
пространстве $X$ последовательность $(x_n)_{n \in \mathbb{N}} \subset X$
$\psi$-сходится ($\varphi$-сходится) к точке $x \in X$ тогда и
только тогда, когда она слабо сходится к точке $x \in X$.
\vskip7pt
{\bf Пример 2.}
\vskip7pt
 Пусть $X$ --- открытый шар единичного радиуса с
центром в нулевом векторе вещественного гильбертова пространства
$V$ с метрикой
$$\rho : X \times X\rightarrow \mathbb{R}_{+},\quad |xy| =
k Arch \left(\dfrac{1 - (x,y)}{((1 - x^2)(1 - y^2))^{1/2}}\right),$$
где $(x,y)$ --- скалярное произведение векторов $x, \, y \in X$, $k
>0$ константа. Это бесконечномерный вариант модели
Бельтрами -- Клейна геометрии Лобачевского \cite [с. 48]{Nut}. Пусть $p = 0 \in X$. Нетрудно
проверить, что в пространстве Лобачевского $(X,\rho)$
последовательность $(x_n)_{n \in \mathbb{N}}$ пространства $X$
$\psi$-сходится ($\varphi$-сходится) к точке $0 \in X$ ($x \in X$)
тогда и только тогда, когда она слабо сходится к этой точке в пространстве $V$. Кроме того, последовательность $(x_n)_{n \in \mathbb{N}}$ пространства $X$ $\psi$-сходится к точке $x \in X\backslash\{0\}$
тогда и только тогда, когда последовательность 
$$\left(\dfrac{\rho (0,x_n) x_n}{||x_n||}\right)_{n \in \mathbb{N}}$$ 
слабо сходится к точке 
$$\dfrac{\rho (0,x) x}{||x||}$$
в пространстве $V$. Отметим также, что условия $(A_{11}-A_{14})$ верны в пространстве Лобачевского $(X,\rho)$.
\vskip7pt
Пусть $\varepsilon > 0$. Введем следующие множества
$$W_{x_1,\ldots,x_N;\, \varepsilon} =
\{\varphi \in \hat \Phi_p(X) :  \dfrac{|\varphi (x_k)|}{|px_k|} < \varepsilon, \,
k = 1,\ldots,N \},$$ 
где точки $x_1,\ldots,x_N$ отличны от точки $p$.
$$V_{x_1,\ldots,x_N;\, \varepsilon} = \{\psi \in \hat
\Psi_p(X) : |\psi (x_k)| < \varepsilon, \, k = 1,\ldots,N \}.$$
Если выбрать множества 
$$V_{x_1,\ldots,x_N;\, \varepsilon}\quad (W_{x_1,\ldots,x_N;\, \varepsilon})$$
в качестве сиcтемы окрестностей нуля в пространстве $\hat \Psi_p(X)$ $(\hat \Phi_p(X))$, то они определят топологию $\tau_{\psi}$ $(\tau_{\varphi})$ в этом пространстве. Это аналоги $*$-слабой топологии в сопряженном
банаховом пространстве $(\cite {Kolm}, c. 230)$.
Следующая теорема 3 аналогична теореме 4 из \cite [с. 232-234]{Kolm}.
\vskip7pt
{\bf Теорема 3} [s15]. {\it Пусть $(X,\rho)$ --- сепарабельное метрическое пространство и
$B[0,1]$ --- замкнутый шар с центром в нуле, радиуса $1$ в пространстве
$\hat \Psi_p(X)$ $(\hat \Phi_p(X))$. Тогда топология в $B[0,1]$, индуцированная топологией $\tau_{\psi}$ $(\tau_{\varphi})$, сильнее
топологии, определенной с помощью метрики
$$d : B[0,1] \times B[0,1] \rightarrow \mathbb{R}_{+},$$
$$d (f,g) = \sum \limits_{n=1}^{\infty} 2^{- n}|(f - g) (x_n)|\quad
(d (f,g) = \sum \limits_{n=1}^{\infty} 2^{- n} \dfrac{|(f - g) (x_n)|}{|px_n|}),$$
где 
$$\{x_n\}_{n \in \mathbb{N}}\quad ( \{x_n\}_{n \in \mathbb{N}}, \, \,  x_n \neq p)$$ 
--- фиксированное счетное всюду плотное множество в замкнутом шаре
$B[p,1] \subset X$  {\rm {(}}пространстве $X)$.}
\vskip7pt
{\bf Доказательство теоремы 3}. 
\vskip7pt
Нетрудно проверить, что $d$ является метрикой, удовлетворяющей условию $$d(f+h,g+h) = d(f,g) \quad \mbox {для всех}\, \, f,\,  g,\,  f + h,\, g + h \in B[0,1]).$$ 
Следовательно, достаточно
проверить, что любой открытый шар
$$Q_{\varepsilon} = \{ f : d (f,0) < \varepsilon \}$$ 
в пространстве $\hat \Psi_p(X)$ $(\hat \Phi_p(X))$ содержит пересечение замкнутого шара $B[0,1]$ с некоторой окрестностью нуля в
топологии $\tau _{\psi }$ $(\tau_{\varphi})$. Выберем такое натуральное число $N$, что $2^{-N} < \varepsilon/2$ и рассмотрим
следующую окрестность нуля в топологии $\tau _{\psi}$ $(\tau_{\varphi})$
$$V = V_{x_1,\ldots,x_N;\, \varepsilon/2} = \{f : |f (x_k)| < \varepsilon/2, \, k = 1,\ldots,N \}$$
$$(W = W_{x_1,\ldots,x_N;\, \varepsilon/2} = \{f : \dfrac{|f(x_k)|}{|px_k|}
< \varepsilon/2, \, k = 1,\ldots,N \}).$$
Пусть 
$$f \in B[0,1]\cap V \quad (f \in B[0,1]\cap W).$$ 
Тогда 
$$d (f,0) = \sum \limits_{n=1}^{N} 2^{- n}|f (x_n)| +
\sum \limits_{n=N+1}^{\infty } 2^{- n}|f(x_n)| \leq \dfrac{\varepsilon}{2} \sum \limits_{n=1}^{N} 2^{- n} +
\sum \limits_{n=N+1}^{\infty } 2^{- n} < \varepsilon$$
$$(d (f,0) = \sum \limits_{n=1}^{N} 2^{- n}\dfrac{|f (x_n)|}{|px_n|} +
\sum \limits_{n=N+1}^{\infty } 2^{- n}\dfrac{|f (x_n)|}{|px_n|} \leq \dfrac{\varepsilon}{2}
\sum \limits_{n=1}^{N} 2^{- n} +
\sum \limits_{n=N+1}^{\infty } 2^{- n} < \varepsilon).$$
Следовательно,
$$B[0,1]\cap V \subset Q_{\varepsilon}\quad (B[0,1]\cap W \subset Q_{\varepsilon}).$$
Теорема 3 доказана.
\vskip7pt
Теорема 4 является следующим шагом в получении достаточных условий на пространство $X$ для аналога теоремы Банаха-Алаоглу \cite [c. 80]{Rudin}. 
\vskip7pt
{\bf Теорема 4} [s15]. {\it Пусть $U_p$ --- окрестность точки $p \in X$
и в множестве 
$$P = \prod \limits_{x \in X} D_x,$$ 
где
$$D_x = \{\alpha \in \mathbb{R} : |\alpha| \leq \lambda (x) \}\, \, \mbox {и}\, \, \lambda(x) \in \mathbb{R}_{+}$$ 
такое, что
$$x \in \omega_{\lambda (x)} (p,U_p),$$ 
задана топология произведения $\tau$. Тогда топология $\tau_{\psi}$ на подмножестве
$$K= \{\psi \in \hat \Psi_p(X) : \, \forall x \in U_p \, (|\psi (x)| \leq 1 )\} \subset P$$
совпадает с топологией, индуцированной топологией произведения
$\tau$. Кроме того, замыкание множества $K$ в пространстве
$(P,\tau)$ принадлежит множеству 
$$\{ f \in F_p(X) : \, \forall  x \in X, \, \forall \beta \in \mathbb{R}_{+}
\, (f (\omega_{\beta } (p,x)) = \beta f(x)); \, \forall x \in U_p
\, (|f(x)| \leq 1)\}.$$}

{\bf Замечание.} {\it При выполнении условий теоремы Банаха-Алаоглу {\rm \cite [c. 80]{Rudin}} множество $K$ замкнуто в пространстве $(P,\tau)$. Задача о
нахождении структуры подмножества $\overline{K}\backslash K
\subset (P,\tau)$ в условиях теоремы 4 остается пока нерешенной.}
\vskip7pt
{\bf Доказательство теоремы 4}. 
\vskip7pt
Из утверждения $(i)$ леммы 1 следует, что для всех $x \in X,\, \, \psi \in K$ 
$$|\psi (x)| \leq \lambda (x).$$ 
По теореме Тихонова \cite [с. 413]{Rudin} пространство $(P,\tau)$ является компактным. Кроме того,
$$P =  \{f : X \rightarrow  \mathbb{R} : \,\forall x \in X \, (|f (x)| \leq \lambda (x)) \}.$$
Фиксируем $\psi_0 \in K$ и выберем произвольно 
$$n \in \mathbb{N},\quad \delta >0, \quad x_i \in X \quad \mbox {при} \, \, 1\leq i \leq n.$$
Множества   
$$W_1 = \{\psi \in \hat \Psi_p(X) : \,|\psi (x_i) - \psi_0 (x_i)| < \delta, \, 1\leq i \leq n \}$$
$$(W_2 = \{f \in P : \, |f  (x_i) - \psi_0 (x_i)| < \delta, \, 1\leq i \leq n \})$$ образуют локальную базу топологии  $\tau_{\psi}$ $(\tau)$  пространства
$\hat \Psi_p(X)$ $(P)$ в точке $\psi_0$.
Тогда 
$$W_1\cap K = W_2\cap K,\, \, \mbox {поскольку} \, \, K \subset \hat \Psi_p (X)\cap P.$$
Таким образом, в множестве $K$ топология $\tau_{\psi}$ совпадает с топологией $\tau$. Пусть $f_0$ принадлежит  $\tau$-замыканию множества $K$. Выберем произвольно 
$$x \in X,\quad \beta \in \mathbb{R}_{+}\, \, \mbox {и}\, \, \varepsilon > 0.$$ Множество всех $f \in P$, для которых 
$$|(f - f_0)(x)| < \varepsilon,\quad |(f - f_0) (\omega_{\beta } (p,x))| < \varepsilon,$$
является $\tau$-окрестностью точки $f_0$, поэтому множество $K$ содержит
хотя бы одну такую точку $\psi$ из этой окрестности. В силу утверждения $(i)$ леммы 1 получим 
$$|f_0 (\omega_{\beta } (p,x)) - \beta f_0 (x)| =
|(f_0 - \psi) (\omega_{\beta } (p,x)) + \beta(\psi - f_0)(x)| < (1 + \beta) \varepsilon.$$
Поскольку $\varepsilon$ произвольно, то
$$f_0 (\omega_{\beta } (p,x)) = \beta f_0 (x) \quad \mbox {для}\, \, x \in X,\, \, \beta \in \mathbb{R}_{+}.$$
Если $x \in U_p$ и $\varepsilon > 0$, то аналогичное рассуждение показывает, что в множестве $K$ найдется функция $\psi$, для которой $$|\psi (x) - f_0 (x)| < \varepsilon.$$ 
Но $|\psi (x)| \leq 1$, следовательно $|f_0 (x)| \leq 1$.
Теорема 4 доказана.

%% file: liter1.tex
\addtocontents{toc}{Список публикаций автора по теме диссертации и литература \hfill {246}\par}
\section*{\bf Список публикаций автора по теме диссертации}

\vskip20pt

\noindent 
[s1] Сосов Е.~Н. О метрическом пространстве слабо ограниченных
          отображений метрических пространств / Е.~Н. Сосов // Изв. вузов. Математика.
          - 1993. - № 9. - C. 61-64.

\noindent
[s2] Сосов Е.~Н. О конечной компактности и полноте метрических
          пространств с метрикой Буземана / Е.~Н. Сосов // Изв. вузов. Математика. -
          1993. - № 11. - С. 62-68.

\noindent
[s3] Сосов Е.~Н. Об одном одуле в геометрии Гильберта / Е.~Н. Сосов // Изв.
          вузов. Математика. - 1995. - № 5. - C. 78-82.

\noindent
[s4] Сосов Е.~Н. О выпуклых множествах в специальном метрическом
          пространстве / Е. Н. Сосов // Изв. вузов. Математика. - 1996. - № 3. - С. 77-79.

\noindent
[s5] Сосов Е.~Н. О геодезических отображениях специальных
          метрических пространств / Е.~Н. Сосов // Изв. вузов. Математика. - 1997. -
          № 8. - C. 46-49.

\noindent
[s6] Сосов Е.~Н. О выпуклых множествах в обобщенном хордовом
          пространстве / Е.~Н. Сосов // Изв. вузов. Математика. - 1998. - № 7. - C. 47-52.

\noindent
[s7] Сосов Е.~Н. Об аппроксимативных свойствах множеств в
          специальном метрическом пространстве / Е.~Н. Сосов // Изв. вузов. Математика.
          - 1999. - № 6. - C. 81-84.

\noindent
[s8] Сосов Е.~Н. О наилучшей сети, наилучшем сечении и
          чебышевском центре ограниченного множества в бесконечномерном
          пространстве Лобачевского / Е.~Н. Сосов // Изв. вузов. Математика. - 1999. - № 9.
          - C. 78-83.

\noindent
[s9] Sosov E.~N. On existence and uniqueness of Chebyshev center
          of a bounded set in a special geodesic space / E.~N. Sosov // Lobachevskii
          Journal of Mathematics. - 2000. - Vol. 7. - С. 43-46.

\noindent
[s10] Сосов Е.~Н. О существовании и единственности чебышевского
          центра ограниченного множества в специальном геодезическом
          пространстве / Е.~Н. Сосов // Труды Матем. центра имени Н.И. Лобачевского.
          - Т. 5. Актуальные проблемы матем. и механики. Материалы
          Международной научной конференции (Казань, 1-3 октября 2000 г.).
          - Казань: " УНИПРЕСС ". - 2000. - С. 198-199.

\noindent
[s11] Sosov E.~N. On Hausdorff intrinsic metric / E.~N. Sosov // Lobachevskii J.
          of Math. - 2001. - V. 8. - P. 185-189.

\noindent
[s12] Сосов Е.~Н. О непрерывности и связности метрической
          $\delta $-проекции в равномерно выпуклом геодезическом
          пространстве / Е. Н. Сосов // Изв. вузов. Математика. - 2001. - № 3. - С. 55-59.

\noindent
[s13] Сосов Е.~Н. О непрерывности метрической $\delta$-проекции
          на выпуклое множество в специальном метрическом пространстве / Е.~Н. Сосов //
          Изв. вузов. Математика. - 2002. - № 1. - C. 71-75.

\noindent
[s14] Сосов Е.~Н. О наилучших $N$-сетях ограниченных замкнутых
          выпуклых множеств в специальном метрическом пространстве / Е.~Н. Сосов // Изв.
          вузов. Математика. - 2003. - № 9. - С. 42-45.

\noindent
[s15] Сосов Е.~Н. Об аналогах слабой сходимости в специальном
          метрическом пространстве / Е.~Н. Сосов // Изв. вузов. Математика. - 2004.
          - № 5. - С. 69-73.

\noindent
[s16] Сосов Е.~Н. О метрическом пространстве всех $2$-сетей пространства
          неположительной кривизны / Е.~Н. Сосов // Изв. вузов. Математика. - 2004.  - № 10.
          - С. 57-60.

\noindent
[s17] Сосов Е.~Н. Наилучшее приближение в метрике Хаусдорфа выпуклого
          компакта шаром / Е.~Н. Сосов // Матем. заметки. - 2004. - Т. 76. - Вып. 2. - С. 226-236.

\noindent
[s18] Сосов Е.~Н. Касательное пространство по Буземану / Е.~Н. Сосов // Изв. вузов. Математика. - 2005. - № 6. - С. 71-75.

\noindent
[s19] Сосов Е.~Н. Относительный чебышевский центр конечного множества геодезического пространства / Е.~Н. Сосов // Изв. вузов. Математика. - 2008. - № 4. - С. 66-72.

\noindent
[s20] Сосов Е.~Н. Пространство всех $N$-сетей и симметризованная степень порядка $N$ метрического пространства / Е.~Н. Сосов // Уч. зап. Казанск. ун-та. - 2009. - Т. 123. - Кн.10. - С. 21-30.

\noindent
[s21] Сосов Е.~Н. Достаточные условия существования и единственности чебышевского
          центра непустого ограниченного множества геодезического пространства / Е.~Н. Сосов // Изв. вузов. Математика. - 2010. - № 6. - С. 47-51.

\noindent
[s22] Ivanshin P.~N. Local Lipschits property for the Chebyshev center mapping over $N$-nets / P.~N. Ivanshin, E.~N. Sosov // Matematicki Vesnik. - 2008. - № 60. - P. 9-22. 

\noindent
[s23] Ivanshin P.~N. The mapping of compact into the set
          of its Chebyhev centres is lipschitz in the space $l^n_{\infty}$ / P.~N. Ivanshin, E.~N. Sosov // www.arXiv.org : 0609229v4 [math.MG].